\documentclass[12pt]{report}
\usepackage{amssymb}
\usepackage{graphicx}
\usepackage{natbib}
\bibpunct{(}{)}{;}{a}{,}{,}
\usepackage{hyperref}

\newcommand{\commentout}[1]{}

\newcommand{\ba}{\begin{array}}
        \newcommand{\ea}{\end{array}}
\newcommand{\bc}{\begin{center}}
        \newcommand{\ec}{\end{center}}
\newcommand{\bdm}{\begin{displaymath}}
        \newcommand{\edm}{\end{displaymath}}
\newcommand{\bds} {\begin{description}}
        \newcommand{\eds} {\end{description}}%17Apr01
\newcommand{\ben}{\begin{enumerate}}
        \newcommand{\een}{\end{enumerate}}
\newcommand{\beq}{\begin{equation}}
        \newcommand{\eeq}{\end{equation}}
\newcommand{\bfg} {\begin{figure}[t]}
        \newcommand{\efg} {\end{figure}}%Nov 5,99
\newcommand{\bi} {\begin {itemize}}
        \newcommand{\ei} {\end {itemize}}
\newcommand{\bqn}{\begin{eqnarray}}
        \newcommand{\eqn}{\end{eqnarray}}
\newcommand{\bqs}{\begin{eqnarray*}}
        \newcommand{\eqs}{\end{eqnarray*}}
\newcommand{\bsl} {\begin{slide}[8.8in,6.7in]}
        \newcommand{\esl} {\end{slide}}
\newcommand{\bss} {\begin{slide*}[9.3in,6.7in]}
        \newcommand{\ess} {\end{slide*}}
\newcommand{\btb} {\begin {table}[t]}
        \newcommand{\etb} {\end {table}}%Nov 10,99
\newcommand{\bvb}{\begin{verbatim}}
        \newcommand{\evb}{\end{verbatim}}

\newcommand{\m}{\mbox}

\newcommand{\mat}[1]{{{\left[ \ba #1 \ea \right]}}}
\newcommand{\cas}[1]{{{\left \{ \ba #1 \ea \right. }}}
\newcommand {\nm} [1] {{\parallel #1\parallel }}
\newcommand{\reffig}[1] {{{Figure \ref {#1}}}}
\newcommand{\refe}[1] {{{Equation \ref {#1}}}}%Nov 5
\newcommand{\reft}[1] {{{Table \ref {#1}}}}
\newcommand{\refet}[2] {{{Equations \ref {#1} and \ref {#2}}}}%Nov 10'99
\def\l      {{\lambda}}

\def\pmb#1{\setbox0=\hbox{$#1$}%
   \kern-.025em\copy0\kern-\wd0
   \kern.05em\copy0\kern-\wd0
   \kern-.025em\raise.0433em\box0 }

\def\eop{{\hfill $\blacksquare$}}%17Apr01
\def\r{{\rho}}
\newtheorem{theorem}{Theorem}[section]%17Apr01
\def\dx     {{\Delta x}}
\def\dt     {{\Delta t}}
\def\dxt    {{\frac{\dx}{\dt}}}
\def\dtx    {{\frac{\dt}{\dx}}}

\def\r {\rho}
\def\vs {V}
\def\ra {{\frac {\r}a}}
\def\l {{\lambda}}

\def\oc{{{OC}}}
\def\uc{{{UC}}}

\def\l {{\lambda}} \def\la {{\lambda}}
\def\vecr{{\textbf{R}}}
\def\vecl{{\textbf{l}}}
\begin{document}
\pagenumbering{roman}
\bc
{  \bf Kinematic Wave Models of Network Vehicular Traffic}
\vskip 0.1in
\centerline{By}
%\vskip 0.1in
\centerline{   Wenlong Jin}
\vskip 0.1in
\centerline{  B.S. (University of Science and Technology of China, Anhui, China) 1998}
\centerline{  M.A. (University of California, Davis) 2000}
\vskip 0.1in
\centerline{   DISSERTATION}
\vskip 0.1in
\centerline{  Submitted in partial satisfaction of the requirements for the
degree of}
%\vskip 0.05in
\centerline{   DOCTOR OF PHILOSOPHY}
\vskip 0.1in
\centerline{  in}
%\vskip 0.1in
\centerline{  APPLIED MATHEMATICS}
\vskip 0.1in
\centerline{  in the}
%\vskip 0.1in
\centerline{   OFFICE OF GRADUATE STUDIES}
%\vskip 0.1in
\centerline{  of the}
%\vskip 0.1in
\centerline{   UNIVERSITY OF CALIFORNIA}
%\vskip 0.05in
\centerline{   DAVIS}
\end{center}
\vskip 0.1in
{  Approved:}
\begin{center}
\centerline{\underbar{Dr. H. Michael Zhang}}
\vskip 0.1in
\centerline{\underbar{Dr. Elbridge Gerry Puckett}}
\vskip 0.1in
\centerline{\underbar{Dr. Zhaojun Bai}}
%\vskip 0.05in
\centerline{  Committee in Charge}
\vskip 0.1in
\centerline{  September 2003}
\end{center}
%%%%%%%%%%%%%%%%%%%%%%%%%%%%%%%%%%%%%%%%%%%%%%%%%%%%%%%%%%%%%%%%%%
\newpage
%Wenlong Jin: 2003
\begin{flushright}
\begin{tabular}{l}
Wenlong Jin\\
September 2003\\
Applied Mathematics\\
\end{tabular}
\end{flushright}
\bc {\bf Kinematic Wave Models of Network Vehicular Traffic} \ec
\bc\underline{\bf Abstract}\ec
The kinematic wave theory, originally proposed by \citep{lighthill1955lwr,richards1956lwr}, has been a good candidate for studying vehicular traffic. In this dissertation, we study kinematic wave models of network traffic, which are expected to be theoretically rigorous, numerically reliable, and computationally efficient.

For inhomogeneous links, we reformulate the Lighthill-Whitham-Richards model into a nonlinear resonant system. In addition to shock and rarefaction waves, standing (transition) waves appear in the ten basic wave solutions. The solutions are consistent with those by the supply-demand method \citep{daganzo1995ctm,lebacque1996godunov}.

For merging traffic, we examine existing supply-demand models and, particularly, distribution schemes. Further, we propose a new distribution scheme, which captures key merging characteristics and leads to a model that is computationally efficient and easy to calibrate.

For diverging traffic, we propose an instantaneous kinematic wave model, consisting of nonlinear resonant systems. After studying the seven basic wave solutions, we show that this model is equivalent to a supply-demand model with modified definitions of traffic demands.

For traffic with mixed-type vehicles, we show the existence of contact waves. Using simulations by the developed Godunov method, we demonstrate that First-In-First-Out (FIFO) principle is observed in this model.

For network traffic flow, we propose a multi-commodity kinematic wave (MCKW) model, in which we combine kinematic wave models of different network components and a commodity-based kinematic wave theory. We also propose an implementation of the MCKW simulation and carefully design the data structure for network topology, traffic characteristics, and simulation algorithms. The solutions are consistent with FIFO principle in the order of a time interval.

For a road network with a single origin-destination (O/D) pair and two routes, we first demonstrate the formation of an equilibrium state and find multiple equilibrium status for different route distributions. We then show the formation of periodic oscillations and discuss their structure and properties.

Finally, we summarize our work and discuss future research directions. 

\newpage
%Wenlong Jin: 2003
\chapter*{Acknowledgements}
I owe many thanks to my wife, Ling. Her selfless support and encouragement is indispensable to the completion of this dissertation. In particular, she brings me the best gift for ever, our baby Laurel, who makes my life more joyful than ever. I also want to thank my family in China, including my parents, sister, and brothers. They have been a constant source of encouragement and support for me.

I'm very grateful to my advisor, Dr. Michael Zhang, for his financial support and academic guidance. Four years ago, he introduced me to the wonderful land of transportation studies. Since then, I've worked on several projects on traffic flow models and ramp metering methods. I always enjoy discussing research questions with him and have been inspired by his advices all the time. During these years when I work with him, he gave me many helpful suggestions on research and career development.

I'm also grateful to Dr. Elbridge Gerry Puckett, my academic advisor and a committee member of my dissertation. He spent a great deal of time in answering all kinds of questions about study and research. His careful revisions to my master thesis and this dissertation have helped me a lot on efficient technical writings. His comments on my research and career development will continue to influence my research in the future. He also introduced me to his colleagues, including Dr. Randall J. LeVeque and Dr. Phillip Colella, to whom I also owe my thanks for their interests in my research.

I'm also grateful to Dr. Zhaojun Bai for serving on the committee of this dissertation and his supportive and encouraging comments on my career development.

I'd like to thank Dr. John Hong for introducing the concept of resonant nonlinear waves, which form a foundation of the kinematic wave theories of inhomogeneous links (Chapter 2) and diverges (Chapter 4). I'd also like to thank Dr. Blake Temple, who generously offered many suggestions on my studies of vacuum problems and merging and diverging models. His research results in resonant nonlinear waves is a key reference in my studies.

Thanks go to Dr. Albert Fannjiang, for his comments and suggestions on my study, research and career, and Dr. Carlos Daganzo, for his suggestions on my research directions.

Thanks go to my fellow classmate, Scott Beaver, for his friendliness and help.

I also owe my thanks to many other colleagues who are not mentioned above but have also contributed to this dissertation in different aspects.

Finally, I offer my sincere thanks the University of California Transportation Center for their financial support through a dissertation grant. 
%%%%%%%%%%%%%%%%%%%%%%%%%%%%%%%%%%%%%%%%%%%%%%%%%%%%%%%%%%%%%%%%%%%%%%%%%%%%%%%%%%%%%%%%%%%%%%%%%%%%%%%%%%%%%%%%%%%%%%%%%%%%%%%%%%%%%%%%%%%%%%%%%%%%%%%%%%

\newpage
\tableofcontents
\newpage
\listoftables
\newpage
\listoffigures

%%%%%%%%%%%%%%%%%%%%%%%%%%%%%%%%%%%%%%%%%%%%%%%%%%%%%%%%%%%%%%%%%%

\newpage
\pagenumbering{arabic}
\pagestyle{myheadings}
\markright{  \rm \normalsize CHAPTER 1. \hspace{0.5cm}
 INTRODUCTION}
\chapter{Introduction}
%\thispagestyle{myheadings}
%Wenlong Jin: 2003
%Wenlong Jin: 2003
\section{Background}
\subsection{Traffic congestion}
As the backbone of the intermodal transportation network in the United States, road networks - consisting of highways, arterial roads, surface streets, and other kinds of roadways - connect air, transit, rail, and port facilities and terminals. In particular, highways carry 90 percent of passenger travel and 72 percent value of freight. Therefore, the performance of its road networks largely defines the mobility of the nation and affects economic and social activities in the United States. However, traffic conditions on road networks in many metropolitan areas are becoming increasingly congested: in 68 major American urban areas, the percentage of un-congested periods was two thirds of the whole peak period in 1982, while the percentage drops to one third in 1997 \citep{schrank1999mobility}.

In a road network, traffic congestion can be recurrent and non-recurrent. Recurrent traffic congestion is generally caused by limited physical infrastructure, increasing travel demand, rush hours, and toll booths. Non-recurrent congestion is associated with accidents, work zones, and weather. Traffic congestion increases travel delay and fuel consumption and adversely affects safety, mobility, productivity, the human and natural environment. As a result, increasing delay has seriously damaged the speed and reliability of road networks, which are vital to emerging industries, like warehousing and logistics.

It has been a major goal of transportation scientists and engineers to alleviate traffic congestion. Many strategies have been pursued to achieve this goal. One method is to expand existing road facilities by adding new roads or lanes, rebuilding key network components, and enhancing the physical condition of roadways. Another method is to apply traffic management and operation technologies, including automatic highways, travel demand management, freeway management, incident management, emergency response management, weather response management, value pricing, arterial signal control, on-ramp metering, traveler information,  changeable message signs, and so on.

All these approaches are intended to either increase the capacity of road networks or reduce their load. Among them, the method of expanding existing road networks has seen limited use due to huge construction costs and the difficulty of addressing public and environmental concerns. Therefore, the idea of using existing infrastructure more efficiently by building Intelligent Transportation Systems (ITS) is a more preferable option.

\subsection{The role of traffic models}
Whichever of the aforementioned strategies are taken to improve the mobility of a road network, their success relies on a better understanding of the properties of congestion and the overall performance of the road network. Widely available sensing, information and communication technology has had a major influence in detecting traffic conditions in a timely manner. However, collecting meaningful traffic data is not an easy job. Moreover, many management and control decisions are grounded on an estimate of traffic conditions in the near or long-term future. Thus, it is essential to accurately estimate traffic conditions in a road network during certain time periods and understand the evolution pattern of traffic conditions, i.e., the traffic dynamics. This calls for the development of traffic models.

A traffic model is a function which relates the movement of a vehicle to driver's behavior, vehicle type, network characteristics, weather conditions, traffic signals, guidance information, and interaction with other vehicles. The movement of a vehicle can be represented by its position at any time, i.e., its trajectory, from which its speed and acceleration rate can be obtained. Given all vehicles' trajectories, one can measure the performance of a road network, e.g., travel time, level of service, congestion level, etc. As shown in \reffig{traffic_model}, traffic models can considered as a theoretical substitution of a real traffic system.

\bfg
\bc
\includegraphics[height=12cm]{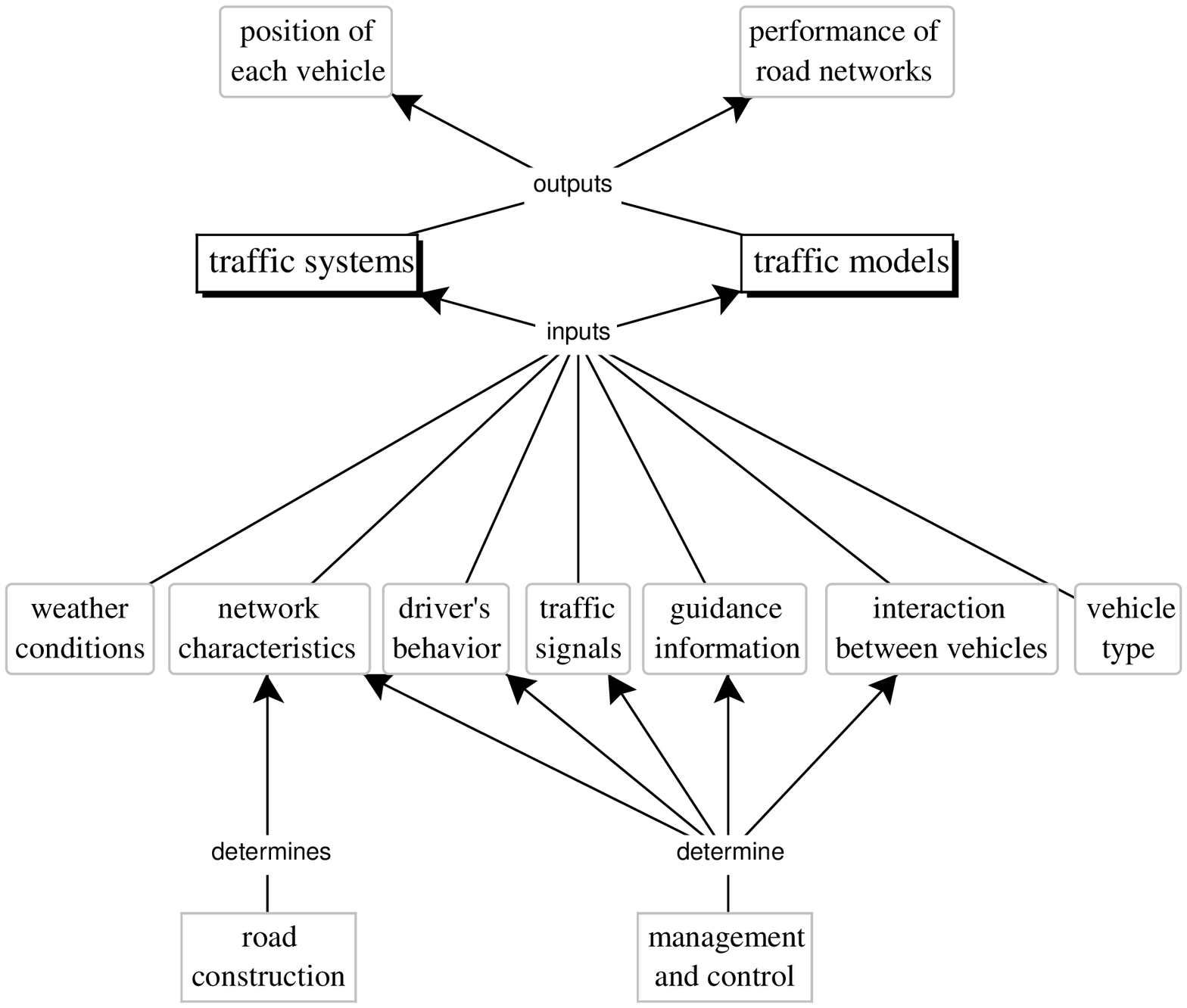}\caption{The definition and role of traffic models in traffic studies}\label{traffic_model}
\ec
\efg

Besides the constraints in \reffig{traffic_model}, traffic models are also subject to both space and time limitations, i.e., initial traffic conditions and boundary conditions. Moreover, all traffic models are subject to calibration and validation. In calibration, parameters in a traffic model are adjusted so that it acts as closely to the real traffic system as possible. In validation, the output of a traffic model is compared with observation of real traffic.

The role of traffic models in transportation engineering is two-fold. First, they provide better understanding of traffic dynamics, in particular the formation and propagation of traffic congestion. Hence, traffic researchers can use traffic models to identify possible bottlenecks. Second, they can serve as a simulation platform, on which different strategies for improving mobility, as shown in \reffig{traffic_model}, can be developed and evaluated. For example, in the plan for expanding a road network, traffic models can be used to simulate proposed expanded networks and help chose the most cost-effective strategy. For another example, they are also helpful in determining the best location of tolling booths, which are designed to divert traffic away from busy roads by charging fees. Finally, traffic models can be used to evaluate previously implemented strategies.

\subsection{Traffic models and simulation packages}
In a vehicular traffic system, a number of trips - defined by their origin/destination, departure time/arrival time and travel route - interact on the road network and generate various dynamics and phenomena. To study these traffic phenomena and the corresponding applications, many traffic models and simulation packages have been proposed in literature.

From a microscopic point of view, a stream of vehicular traffic is the combination of the movements of all vehicles. In a microscopic model of traffic, the movement of each vehicle and interactions between the vehicle pairs are studied. Three types of microscopic models have been suggested. One approach is the GM family of car-following models developed in the 1960's (e.g., \citep{gazis1961follow}). In these models, the movement of a vehicle is described by an ordinary differential equation in time. In another approach, known as coupled-map lattice models, a vehicle's movement is defined by equations that are discrete in time \citep{chowdhury2000statistical}. Yet another approach,  cellular automata (CA) models are based on the framework of statistical physics, in which not only independent variables (space and time) but also descriptive variables (speed and acceleration) are discrete \citep{chowdhury2000statistical}.

From a macroscopic point of view, one approach to modeling traffic is to treat the traffic as it were a gas of interacting particles, in which each particle represents a vehicle. The resulting model is called a kinetic model of vehicular traffic. Another approach is to consider traffic flow as a compressible fluid (continuum). Such models are known as continuum models, or fluid-dynamical models. In continuum models, the basic characteristics are flow rate $q$, traffic density $\rho$, and travel speed $v$, which are all functions in time and space. The first-order continuum models are generally called kinematic wave models, in which traffic dynamics are regarded as a combination of ``kinematic" waves in these quantities.

Based on the aforementioned theoretical models of vehicular traffic, many simulation packages have been developed; for example, the PASSER series \citep{passer1991}, TRANSYT \citep{courage1991transyt}, INTEGRATION \citep{aerde1995integration}, DYNASMART \citep{jayakrishnan1994dynasmart}, NETSIM \citep{fhwa1998netsim}, FRESIM \citep{smith1995fresim}, FREFLO \citep{payne1979freflo}, MITSIM \citep{mitsim1999}, TRANSIMS \citep{transims1992}, and PARAMICS \citep{cameron1996paramics}. In some simulation packages, however, the underlying traffic flow theories lack mathematical rigor, particularly for traffic dynamics at highway junctions. Furthermore, numerical methods used in some simulations have not been well justified and may cause numerical instability. Another common drawback with many simulation packages is their enormous computational cost. Therefore, many of these simulation packages have serious limitations in applications. That said, continuing efforts to develop theoretically rigor, numerically sound, and computationally effective models are still necessary, in particular for applications to dynamic traffic assignment and other advanced traffic engineering strategies. 
%Wenlong Jin: 2003
\section{Continuum models}
\subsection{Kinematic wave models}
In the kinematic wave models, traffic is viewed as a continuous media and characterized by traffic density ($\r$), travel speed ($v$), and flow-rate ($q$). The movements of vehicles on a road network are considered as combinations of kinematic waves in either of these three quantities.

The different types of kinematic waves are associated with different components of a traffic system. For example, on a homogeneous link with uniform conditions in vehicles, drivers, weather, etc., one can observe two basic waves: decelerating shock waves, generally seen when lights turn red, and accelerating rarefaction waves, generally seen when lights turn green. In a traffic system with inhomogeneous links, merges, diverges, or different vehicle types, more complicated waves can be observed, such as standing transitional waves, contact waves, and periodic waves.

In kinematic wave models, a hyperbolic conservation law is derived from traffic conservation:
\bqn
\r_t+q_x=0.\label{traffic_conservation}
\eqn
 A fundamental assumption in kinematic wave models is that the flow-rate $q$ is a function of the traffic density $\r$; i.e., $q=Q(\r)$, which is called the fundamental diagram. Hence, $v=V(\r)\equiv Q(\r)/\r$. Generally, flow-rate is a concave function in density and retains its maximum, the capacity, at the critical density $\r_c$. When traffic density is higher than the critical density, it is in the over-critical region and in the under-critical region, otherwise. The fundamental diagram, or the speed-density relation, varies with link characteristics, vehicle types, and so on.

\citet{lighthill1955lwr} and \citet{richards1956lwr} first proposed and analyzed kinematic waves on homogeneous links. The corresponding model is known as the LWR model and is described by a first-order, nonlinear PDE:
\bqn
\r_t+Q(\r)_x&=&0. \label{org:lwr}
\eqn
This model is an hyperbolic conservation law, whose Riemann problems \footnote{In a Riemann problem, the initial traffic conditions are described by a Heaviside function or step function.} can be analyzed using well developed tools \citep [e.g.,][]{lax1972shock, smoller1983shock}. Analysis shows that solutions to \refe{org:lwr}, where $f(x,\r)=\r \vs (\r)$, have wave properties analogous to those of water flow in channels \citep{lighthill1955flood}. In other words, the Riemann problem consists of either shock or rarefaction waves.

Numerically, the LWR model can be solved with a first-order Godunov method \citep{godunov1959}, in which a link is partitioned into a number of cells, a time duration into a number of time steps, and traffic conditions in each cell at a time step are uniform. In a Godunov method, traffic conditions are updated according to the conservation equation, \refe{org:lwr}; i.e., during a time interval, the increasing number of vehicles in a cell are the difference between the in-flow through its upstream boundary and the out-flow through its downstream boundary. Traditionally \citep{leveque2002fvm}, flows through cell boundaries are computed from wave solutions of Riemann problems.

As an alternative of the Godunov method, the celebrated supply-demand method was first proposed in \citep{daganzo1995ctm,lebacque1996godunov}. In this intuitive, engineering method, the flow through a boundary equals the minimum of the traffic demand of its upstream cell and the traffic supply of the downstream cell. Here the traffic demand \citep{lebacque1996godunov}, called sending flow in \citep{daganzo1995ctm},  of a cell is defined as its flow-rate when traffic is under-critical or its capacity when over-critical, and its traffic supply \citep{lebacque1996godunov}, receiving flow in \citep{daganzo1995ctm}, is the capacity when traffic is under-critical or the flow-rate when under-critical. Since one does not have to understand wave solutions of Riemann problems, this method has been widely applied in traffic studies.

Many attempts have been made to model more complicated traffic systems under the framework of the LWR model. For example, \citet{lighthill1955lwr} discussed the kinematic wave theory of traffic dynamics on an inhomogeneous link with lane-drops or curvatures. Lately, \citet{lebacque1996godunov} gave a detailed analysis of the kinematic wave solutions of such a traffic system and summarized the solutions into the supply-demand method. In \citep{daganzo1995ctm}, the concepts of sending flow (traffic supply) and receiving flow (traffic demand) was intuitively extended for simulating traffic on an inhomogeneous link.

Another attempt is to model traffic dynamics on a highway network, where the focus is on highway junctions, including merges, diverges, and other intersections. \citet{holden1995unidirection} studied kinematic waves initiated at highway junctions, by assuming the existence of an optimization problem at each junction and excluding route choice behavior. Without route choice, this model sees limited applications in reality. In the network traffic flow model by \citet{kunhe1992continuum}, on-ramps and off-ramps are considered as sources and sinks respectively. Although including the influence of ramps on mainline freeways, this models omits the other side: the influence of mainline freeways on ramps. Thus it also has limitations without giving a full picture of traffic dynamics. To capture overall traffic phenomena in a traffic system, \citep{daganzo1995ctm,lebacque1996godunov} proposed some discrete models, in which traffic demands and supplies for cells around an intersection are defined the same as those for cells inside a link, and route choice behavior and certain optimization rules for flows are incorporated in order to determine unique flows through the intersection. In literature, there have been little progress in analyzing the kinematic waves initiated at a highway junction considering route choice behavior, mainly due to the difficulties in formulating them into a system of continuous partial differential equations as in the LWR model.

In the framework kinematic wave theories, other extensions have been proposed for addressing different concerns in traffic systems. In \citep{daganzo1997special,daganzo1997it}, traffic on special lanes is investigated. In \citep{daganzo2002behavior}, driver behavior are incorporated. \citet{wong2002multiclass} discuss differentiated vehicle types. Other interesting studies can be found in \citep{vaughan1984od,newell1993sim,newell1999exit,leveque2001night,lebacque2003acceleration}

Compared to other traffic models, kinematic wave models have the following appealing features. First, they have inherit compliance with many applications in large-scale road networks, in which aggregate quantities such as traffic counts, flows, and space-mean travel speed are more important than characteristics of individual vehicles. Second, kinematic wave models can generally be written into a system of hyperbolic partial differential equations, or hyperbolic conservation laws. Thus, one can better understand the formation and structure of a traffic phenomenon on a road network through theoretical analysis of these equations. Finally, there exist many sound, efficient numerical methods for solving hyperbolic conservation laws, and one can carry out efficient and trustful simulations of large-scale road networks.

\subsection{Higher-order models}
Among many efforts to extend the LWR theory to capture instabilities in practical traffic flow, one direction leads to higher-order, or nonequilibrium, models.
\citet{payne1971PW} and \citet{whitham1974PW} introduced a momentum equation to capture the change in travel speed in addition to the traffic conservation equation, \refe{traffic_conservation}:
\bqn
v_t+vv_x+\frac{c_0^2}\r \r_x &=&\frac {V (\r)-v}{\tau}.\label{org:pw2}
\eqn
Here the constant $c_0$ is the traffic sound speed, the source term $\frac {\vs (\r)-v}{\tau}$ is called a relaxation term, and $\tau$ is the relaxation time. With the momentum equation,
the PW model attempts to model driver behavior by accounting for drivers' anticipation and inertia. One can show that the LWR model is an asymptotic approximation of the PW model \citep{schochet1988limit}.

The PW model, as a second-order system of hyperbolic conservation laws with a source term, can be numerically solved by Godunov methods \citep{jin2001PW}. With simulations it was shown that, in addition to modeling stable traffic like the LWR model, the PW model is capable of modeling the formation of vehicle clusters \citep{jin2003cluster}. However, the PW model has drawn some criticism since it allows wave solutions with speed higher than vehicle travel speed and may yield back-traveling (or negative-speed) results \citep{daganzo1995requiem}.

Another nonequilibrium traffic flow model is due to \citep{zhang1998theory, zhang1999analysis, zhang2000structural, zhang2001difference}. In this model, a modified momentum equation is included,
\bqn
v_t+vv_x+\frac{(\r \vs'(\r))^2}{\r} \r_x &=& \frac {V(\r)-v}{\tau}. \label {org:zhang2}
\eqn
This model bears shortcomings similar to the PW model, yet differs from the latter in that the sound speed $c=\r \vs'(\r)$ varies with respect to traffic density $\r$. Moreover, it is always stable and therefore acts like the LWR model \citep{li2001global}.

In order to avoid the negative-speed drawbacks of the PW model, \citet{aw2000arz} identified a number of principles and proposed a satisfactory momentum equation in the following form:
\bqn
v_t+(v-\r p'(\r))v_x&=&\frac{V(\r)-v}{\tau},\label{arz-ar}
\eqn
where $p(\r)$ is a pressure law and is increasing. Further, the Riemann problems were discussed for the model. However, how $p(\r)$ is related to driver-behavior was not specified in their study.

In the same spirit, \citet{zhang2002arz} derived a model similar to Aw and Rascle's model from a car-following model:
\bqn
v_t+(v+\r V'(\r))v_x&=&\frac{V(\r)-v}{\tau},\label{arz-z}
\eqn
which is in the framework of \refe{arz-ar} with $p(\r)=-V(\r)$. In this model, the definition of $p(\r)$ is derived from a car-following model. Both models, \refe{arz-ar} and \refe{arz-z}, no longer admit wave solutions faster than traffic and avoids back-traveling \citep{aw2000arz,zhang2002arz}. However, they are always stable \citep{li2001global} and, therefore, lose the PW model's ability to simulate unstable traffic and vehicle clusters.

Also in the framework of \refe{arz-ar}, another model was proposed in \citep{jiang2002JWZ} as
\bqn
v_t+(v-c_0)v_x&=&\frac{V(\r)-v}{\tau},\label{arz-pw}
\eqn
where the pressure function $p(\r)=c_0\ln \r$. Like the PW model, this model is unstable under some traffic conditions, but yields non-physical solutions when traffic is in unstable region \citep{jin2003arz}.

More complicated traffic systems, for example, with highway junctions \citep{liu1996michalopoulos,lee2000onramp}, have been studied with non-equilibrium models, primarily the PW model. Moreover, several simulation packages are also based on the PW model. In this dissertation, however, we focus on the LWR model, i.e., the first-order continuum model. 
%Wenlong Jin: 2003
\section{Fundamental diagrams}
In the kinematic wave models, fundamental diagrams capture constraints on a traffic system such as road characteristics, vehicle type, driver's behavior, weather conditions, and traffic rules. Therefore, the success of such models rely on the accuracy of the fundamental diagram.

One typical relationship that has been observed between flow rate and occupancy is shown in \reffig{occupancy_flow} from \citep{hall1986fd}. It's generally assumed that the equilibrium travel speed $V(\r)$ is decreasing with respect to traffic density; i.e., $V'(\rho)<0$, and the fundamental diagram $Q(\rho)\equiv \rho V(\rho)$ is concave; i.e., $Q''(\rho)<0$.

\bfg
\bc
\includegraphics[height=12cm, angle=-90]{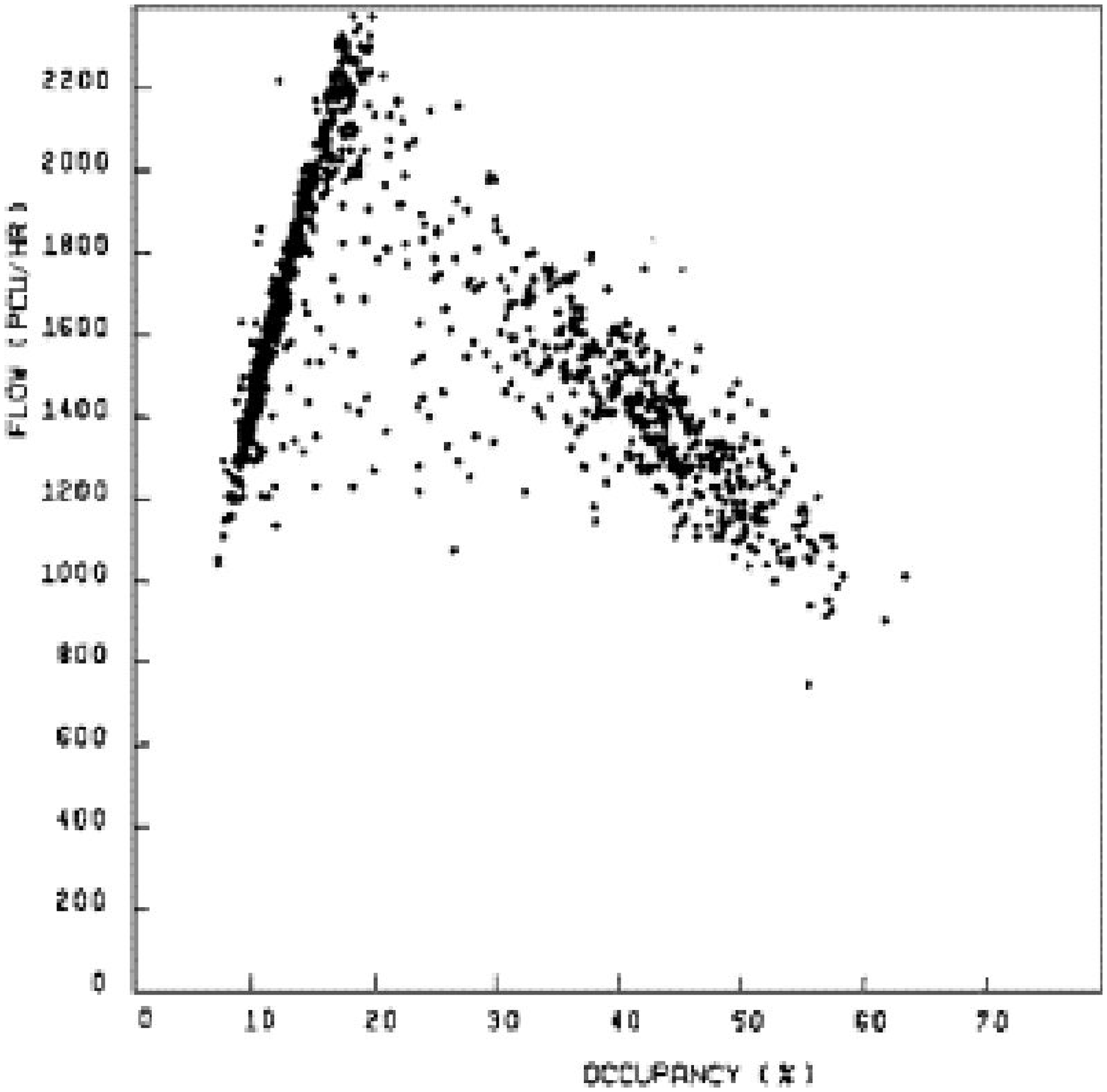}\caption{A typical occupancy-flow rate relationship \citep{hall1986fd}}\label{occupancy_flow}
\ec
\efg

Due to their importance, many fundamental diagrams have been proposed since the early era of traffic engineering practice. Some traditional speed-density relations are listed in \reft{table_fd}.

\btb
\begin {center}
\begin {tabular} {|c|c|}\hline
Functions&$V(\rho)$  \\\hline
\citet{greenshields1935capacity} & $v_f(1-\rho/\rho_j)$ \\\hline
\citet{greenberg1959analysis}& $v_0 \ln(\rho_j/\rho)$  \\\hline
\citet{newell1961nonlinear}& $v_f\left[1-\exp\left(\frac {|c_j|}{v_f}(1-\rho_j/\rho)\right) \right]$ \\\hline
\citet{underwood1961fd}& $v_f \exp (-\rho/\rho_j)$ \\\hline
\citet{drake1967fd}& $v_f \exp \left [-\frac 12 \left(\frac {\rho}{\rho_j}\right)^2 \right ]$\\\hline
\citet{pipes1967fd} & $v_f\left(1-\frac{\rho}{\rho_j} \right)^n,\: n>1$ \\\hline
\end {tabular}
\end {center}
\caption{Traditional speed-density relationship functions}\label{table_fd}
\etb

More recent examples are the following:  the triangular fundamental diagram \citep{newell1993sim},
\bqn
V(\r)&=&\cas{{ll}v_f &\m{ when under-critical},\\ v_f \frac{\r_c}{\r_j-\r_c} (\frac{\r_j}{\r}-1) & \m{ otherwise},
}
\eqn
a non-convex fundamental diagram \citep{kerner1994cluster, herrmann1998cluster},
\bqn
V(\r)&=&5.0461\left[\left(1+\exp\{[\rho/\r_j-0.25]/0.06\}\right)^{-1}-3.72\times 10^{-6}\right] l/\tau,
\eqn  and the exponential form \citep{delcastillo1995fd_theory,delcastillo1995fd_empirical},
\bqn
V(\r)&=&v_f\left\{ 1-\exp\left[1-\exp\left( \frac{|c_j|}{v_f}(\frac{\r_j}{\r}-1)\right)\right]\right\}.
\eqn

The aforementioned fundamental diagrams are all continuous functions. In addition, to capture ``two-capacity" phenomena \citep{banks1991twocapacity}, a discontinuous fundamental diagram in ``reverse-lambda" shape has been proposed \citep[e.g., ][]{koshi1983fd}. Finally, \citet{daganzo1997special} has proposed a novel, two-regime fundamental diagram for differentiated road and vehicle types.

In this dissertation, we focus on studying network inhomogeneities such as merges, diverges, and other junctions. Thus, we will consider continuous fundamental diagram, mostly of the triangular shape. 
%Wenlong Jin: 2003
\section{Motivation for the dissertation research}
While driving in ``stop and go" traffic, many people may have wondered how such congestion is formed, propagated, and diffused. As we know, the network structure is an important factor in determining the characteristics of traffic congestion. For example, congestion generally forms at merging junctions and propagates through diverges. Therefore, traffic conditions on a road network have to be considered as a whole.

The kinematic wave theory is a good candidate for studying vehicular traffic dynamics due to its theoretical rigor, numerical soundness, and computational efficiency. However, its potential has yet to be fully explored. Theoretically, there has been little progress in finding kinematic waves in network traffic systems since \citep{lighthill1955lwr,richards1956lwr}, and traffic congestion in road networks has not been well understood in the framework of kinematic waves. Numerically, little attention is paid to the stability and convergence of solution methods.  Computationally, many simulation models based on kinematic wave theory still keep track of individual vehicles and lose a part of the kinematic wave models' efficiency. Due to the aforementioned limitations, kinematic wave theories of network traffic flow have yet to be improved so that they can be applied in solving practical problems, e.g., in Intelligent Transportation Systems.

In the framework of kinematic wave theories, we study network traffic dynamics and provide a theoretically rigorous, numerically reliable, and computationally efficient simulation model for understanding and mitigating congestion in large-scale road networks. In order to fully explore the advantages of the kinematic wave theory, we consider the following four aspects of traffic modeling in this dissertation. First, we study traffic dynamics at crucial components of a road network, such as link bottlenecks, merges, and diverges, as well as different types of vehicles. In addition, we consider the dynamics of additive multi-commodity  traffic. These studies will offer a more rigorous understanding of traffic dynamics in terms of kinematic waves. Second, we discuss systematically with numerical methods for computing commodity-specified fluxes through junctions in the framework of supply-demand method \citep{daganzo1995ctm,lebacque1996godunov}. Third, we develop a macroscopic simulation platform of multi-commodity traffic on a road network. Finally, we apply this simulation model to study traffic on a road network and probe traffic phenomena related to network topology and route choice behavior.

The structure of this dissertation is as follows. In Chapter 2, we study the kinematic waves on inhomogeneous links. In Chapter 3, we discuss kinematic wave models of merging traffic, which are based on the supply-demand method. In Chapter 4, we propose a new kinematic wave theory for diverging traffic. In Chapter 5, we discuss theoretical and numerical solutions for mixed-type traffic. In Chapter 6, we develop a the multi-commodity kinematic wave simulation platform for network traffic flow. In Chapter 7, we apply this simulation model to study some traffic phenomena in road networks. Finally, we summarize our research results and present future research plans in Chapter 8.

\newpage
\pagestyle{myheadings}
\markright{  \rm \normalsize CHAPTER 2. \hspace{0.5cm}
 INHOMOGENEOUS LINK MODEL}
\chapter{Kinematic wave traffic flow model of inhomogeneous links}
%\thispagestyle{myheadings}
%Wenlong Jin: 2003
%Wenlong Jin: 2003
\section {Introduction}
The kinematic wave traffic flow model of LWR was introduced by
\citet{lighthill1955lwr} and \citet{richards1956lwr} for modeling dense traffic flow on crowded
roads, where the evolution of density
$\r(x,t)$ and flow-rate $q(x,t)$ over time is described by
equation,
\bqn
\r_t+q_x&=&0.\label{lwr_b}
\eqn
This equation follows conservation of traffic that vehicles are
neither
generated nor destroyed on a road section with no entries and exits.

The conservation equation alone is not sufficient to describe traffic
evolution, because it does not capture the unique character of
vehicular flow---drivers slow down when their front spacing is
reduced to affect safety.  The LWR model addresses this issue by
assuming a functional relationship between local flow-rate and
density,  i.e., $q=f(x,\r)$. This flow-density relation, also known as
the fundamental diagram of traffic flow, is often assumed to be
concave in $\rho$ and is a function of the local
characteristics of a road location, such as the
number of lanes, curvature, grades, and speed limit, as well as
vehicle and driver composition.  When a piece of
roadway is homogeneous; i.e., the aforementioned characteristics
of the road are uniform throughout the road section,
the fundamental diagram is invariant to location $x$ and the
LWR model becomes
\bqn
\r_t+f(\r)_x&=&0. \label{hom_LWR}
\eqn
In contrast, if a section of a roadway is inhomogeneous, the LWR
model is
\bqn
\r_t+f(x,\r)_x&=&0.\label{inh_b}
\eqn
Here we introduce a more explicit notation, an inhomogeneity
factor $a(x)$, into the flux function $f(x,\r)$ and obtain the
following equivalent  LWR model for an
inhomogeneous road
\bqn
\r_t+f(a,\r)_x&=&0. \label{inh_1st}
\eqn
This equation is particularly suited for our later analysis of the
 LWR model for inhomogeneous roads (We shall hereafter call
\refe{hom_LWR} the homogeneous LWR model and \refe{inh_1st} the
inhomogeneous LWR model).

Both the homogeneous and inhomogeneous LWR models have been
studied by researchers and applied by practitioners in the
transportation community. Note that the homogeneous
version  \refe{hom_LWR} is nothing more than a scalar conservation law.
Therefore, its wave solutions exist and are unique under the so-called
``Lax entropy condition" \citep{lax1972shock}. These solutions
are formed by basic solutions to the Riemann problem of
\refe{hom_LWR}, in which the initial conditions jump at a boundary
and are constant  both  upstream and downstream of the jump
spot.  Nevertheless, because analytical solutions are difficult
to   obtain  for \refe{hom_LWR} with arbitrary
initial/boundary  conditions, numerical solutions
have to be computed in most cases. The most often used approximation
of \refe{hom_LWR} is perhaps that of Godunov. In the Godunov method,
a roadway is partitioned into a number of cells; and the change
of the number of vehicles in each cell during a time
interval  is the net inflow of vehicles from its boundaries. The
rate of  traffic flowing through a boundary is obtained by solving
a  Riemann problem at this boundary. Besides the Godunov
method, there are other types of approximations of the homogeneous
LWR model, and some of them are shown to be variants of Godunov's
method \citep{lebacque1996godunov}.

In contrast to the well researched homogeneous LWR model,
the inhomogeneous model is less studied and less understood.
Of the few efforts to rigorously solve the inhomogeneous LWR model,
the works of  \citet{daganzo1995ctm} and \citet{lebacque1996godunov}
should be mentioned. In his cell transmission model, Daganzo
started with a discrete form of the conservation equation
and suggested that the flow through a boundary connecting two
cells of a homogeneous road is the minimum of the ``sending
flow" from the upstream cell and the ``receiving flow" of the
downstream cell. The ``sending flow" is equal to the upstream
flow-rate if the upstream traffic is undercritical (\uc) or the
capacity of the upstream section if the upstream traffic is
overcritical (\oc); on the other hand, the ``receiving flow" is equal to
the capacity of the downstream section if the downstream traffic is
\uc { or} the downstream flow-rate if the downstream traffic is \oc.
In the homogeneous case, the boundary flux computed from the
``sending flow" and the ``receiving flow" is the same as that
computed from solutions of the associated Riemann problem. Since the
definitions of ``sending flow" and  ``receiving flow" can be
extended to inhomogeneous sections, Daganzo's method can also be
 applied to  the inhomogeneous LWR model.
Different from Daganzo, Lebacque started his method with the
solution of  the ``generalized"  Riemann problem for \refe{inh_b}. In
this  work,  Lebacque came up with some  rules for solving the
``generalized"  Riemann problem. These rules play the same role as
entropy conditions.  Moreover, Lebacque found that the boundary flux
obtained from solving the Riemann problem is consistent with
that from Daganzo's method,  and he called Daganzo's ``sending flow"
 demand and ``receiving flow" supply.

The methods of Daganzo and Lebacque are streamlined versions of
Godunov's method for the inhomogeneous LWR model.
They hinge upon the definitions of the demand and supply functions, which
can be obtained unambiguously when $f(a,\r)$ is unimodal. When $f(a,\r)$
has multiple local maximum, or when the traffic flow model is of higher
order, it is yet to be determined if equivalent demand/supply
functions exist. Thus, these two methods may not
be applicable to solve the LWR model that has multiple critical
points on its fundamental diagram, nor higher-order models of traffic
flow, such as the Payne-Whitham \citep{payne1971PW,whitham1974PW} model and
a model by \citep{zhang1998theory,zhang1999analysis, zhang2000structural, zhang2001difference}.
Note that, however, these
higher-order models for homogeneous roads can still be
solved  with  Godunov's method \citep{zhang2001difference}.

In this chapter, we present a new
method for solving the Riemann problem for \refe{inh_1st}, which can
be extended to solve higher-order models.  By introducing an
additional  conservation law for $a(x)$, we consider
the inhomogeneous LWR model as a resonant nonlinear system and study
its properties (Section \ref{sec2}).  We also solve
the Riemann problem for \refe{inh_1st} and  show that the
boundary flux at the location of the inhomogeneity is consistent  with the one given by Lebacque and
Daganzo for the same initial condition (Section \ref{ri}). Finally, we
demonstrate our method through solving an initial value problem on a
ring road with a bottleneck, and draw some conclusions from our
analyses.

%Wenlong Jin: 2003
\section{Properties of the inhomogeneous  LWR model as a resonant
nonlinear system}\label{sec2}
Instead of directly study the inhomogeneous LWR
model described by \refe{inh_1st},  we augment \refe{inh_1st} into a
system of conservation laws through the introduction of an
additional conservation law  $a_t=0$ for {the} inhomogeneity
factor $a(x)$, which leads to
\bqn
U_t+F(U)_x&=&0,\label{system}
\eqn
where $U=(a,\r), F(U)=(0,f(a,\r)),x\in R,t\geq 0$. Without loss of
generality, we assume the inhomogeneity is the drop/increase of lanes
at a particular location, and write the fundamental diagram
as $f(a,\r)=\r\vs(\ra)$, where $v=\vs(\ra)$ is the
speed-density relation. The results obtained
hereafter apply to other types of inhomogeneities, such as changes
in grades.

The inhomogeneous LWR model \refe{system} can be linearized as
\bqn
U_t+\partial F(U) U_x&=&0, \label{linear}
\eqn
where the differential $\partial F(U)$ of the flux vector $F(U)$ is
\bqn
\partial F&=&\mat{{cc}0&0\\-\frac {\r^2}{a^2}
\vs'(\ra)&\vs(\ra)+\ra\vs'(\ra)}.
\eqn
The two eigenvalues of $\partial F$ are
\bqn
\l_0=0,\qquad \l_1=\vs(\ra)+\ra\vs'(\ra).
\eqn
The corresponding right eigenvectors are
\bqs
\vecr_0=\mat{{c}\vs(\ra)+\ra\vs'(\ra)\\(\ra)^2\vs'(\ra)},
\qquad \vecr_1=\mat{{c}0\\1},
\eqs
and the left eigenvector of $\partial f/\partial \r$ as $\vecl_1=1$.

System \refe{system} is a non-strictly hyperbolic system, since
it can happen that  $\l_1=\l_0$. We consider a traffic state
$U_{\ast}=(a_{\ast},\r_{\ast})$  in this system as  critical if
\bqn
\lambda_1(U_{\ast})=0\label{assum_1};
\eqn
i.e., at critical states, the two wave speeds are the same and system
\refe{system} is singular.
For a critical traffic state $U_{\ast}$ we also have
\bqn
\frac {\partial}{\partial \r} \l_1(U_{\ast})&=&f_{\r\r}<0
\label{assum_2}
\eqn
since flow-rate is concave in traffic density,
and
\bqn
\frac{\partial}{\partial a}
f(U_{\ast})&=&-(\ra)^2\vs'(\ra)\big|_{U_{\ast}}=\ra \vs(\ra)\big
|_{U_{\ast}}>0. \label{assum_3}
\eqn

A consequence of properties \refe{assum_2} and \refe{assum_3} is that
the linearized system \refe{linear} at $U_{\ast}$ has the following
normal form
\bqn
\mat{{c} \delta a\\\delta \r}_t+\mat{{cc}0&0\\1&0}\mat{{c} \delta
a\\\delta \r}_x&=&0. \label{nlinear}
\eqn
System \refe{nlinear} has the solution $\delta \r(x,t)=\delta
a'(x)t+c$, and the solution goes to infinity as $t$ goes to infinity.
Therefore \refe{nlinear} is a linear resonant system, and the
original inhomogeneous LWR model \refe{system} is a nonlinear
resonant system.

For \refe{system}, the smooth curve $\Gamma$ in $U$-space formed by
all critical states $U_{\ast}$ are  named a transition curve.
Therefore $\Gamma$ is defined as
\bqs
\Gamma&=&\left\{ U \big | \l_1(U)=0\right\}.
\eqs
Since $\l_1(U)=\vs(\ra)+\ra\vs'(\ra)$, we obtain
\bqn
\Gamma &=&\left\{(a,\r)\big|\ra=\alpha, \m{ where } \alpha \m{
uniquely solves }\vs(\alpha)+\alpha \vs'(\alpha)=0\right\};
\label{trancur}
\eqn
i.e., the transition curve for \refe{system} is a straight line
passing through the origin in $U$-space. In \refe{trancur}, $\alpha$
is unique since $f(a,\r)$ is concave in $\r$.

The entropy solutions to a nonlinear resonant system are different
from those to a {strictly} hyperbolic system of conservation laws.
\citet{isaacson1992resonance} proved that
solutions to the Riemann problem  for system \refe{system} exist and
are unique with the conditions \refe{assum_1}-\refe{assum_3}. \citet{lin1995resonant} presented solutions to a scalar
nonlinear resonant system, which is similar to our system
\refe{system} except that $f$ is convex in their study. In the next section we apply
those results to solve the Riemann problem for the inhomogeneous LWR
model.

%Wenlong Jin: 2003
\section{Solutions to the Riemann problem}\label{ri}
In this section we study the wave solutions to the Riemann problem
for \refe{system} with  the following jump initial conditions
\bqn
U(x,t=0)&=&\cas{{cc} U_L&\m{ if }x<0\\U_R &\m{ if }x>0},
\label{inh.ini}
\eqn
where the initial values of $U_L, U_R$ are constant. For
computational purposes, we are interested in  the average flux at the
boundary $x=0$ over a time interval $\dt$, which is denoted by
$f^{\ast}_0$ and defined as
\bqn
f^{\ast}_0&=&\frac 1 {\dt} \int_0^{\dt} f(U(x=0,t)) \m{dt}.
\eqn

The augmented inhomogeneous LWR model \refe{system} has two
families of basic  wave solutions associated with the two eigenvalues.
The solutions whose  wave speed is $\l_0$ are in the 0-family, and the
waves are called  0-waves. Similarly the solutions whose wave speed
is $\l_1$ are in  the 1-family, and the waves are called 1-waves.
{ The 0-wave is also called a}
standing wave  since its wave speed is always 0. The 1-wave
solutions are determined  by the solutions of the scalar conservation
law $\r_t+f(\bar  a,\r)_x=0$,  where $\bar a$ is a constant.
Corresponding to the two types of wave solutions, the integral curves
of the right  eigenvectors $\vecr_0$ and $\vecr_1$ in $U$-space are called 0- and 1-curves respectively. Hence the 0-curves are  given by $f(U)$=const, and the 1-curves are given by
$a=$const. {A 0-curve, a 1-curve}, and the transition curve $\Gamma$  passing through a critical state
$U_{\ast}$  are shown in
\reffig{F_waves}, where
$a$ is set as the vertical axis and $\r$ the horizontal
axis.

\bfg\bc\includegraphics
[height=12cm]{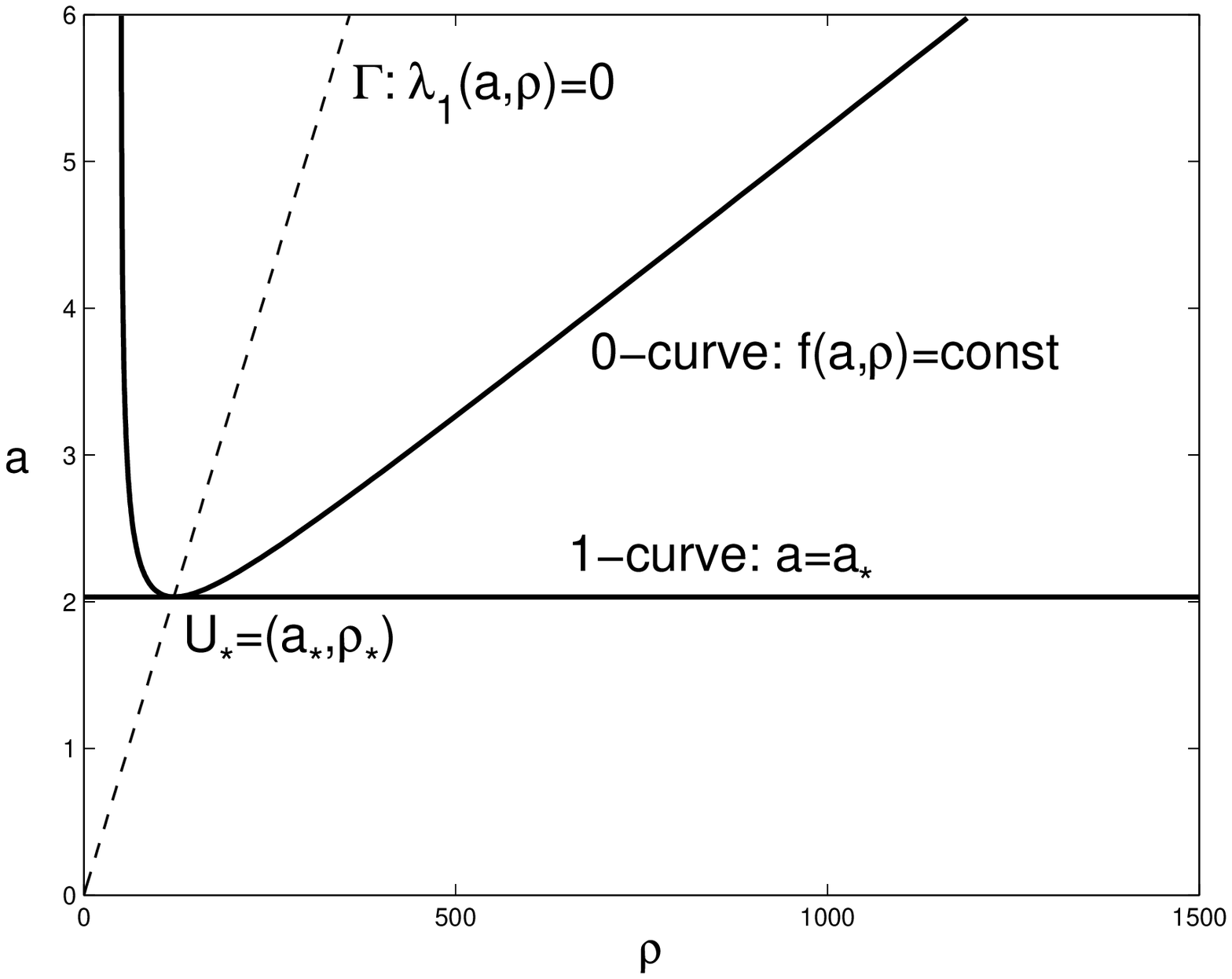}\ec\caption{Integral
curves}\label{F_waves}\efg

As shown in \reffig{F_waves}, the {0-curve} is convex, and the
{1-curve} is tangent to the {0-curve} at the critical state
$U_{\ast}$. The transition curve $\Gamma$ intersects the {0- and
1-curves}  at $U_{\ast}$, and there is only one
critical state on one {0-curve or 1-curve}. For any point $U$,
only one {0-curve and one 1-curve} pass it.
In \reffig{F_waves}, the states left to the transition curve are
undercritical ({\uc })  since $\r/a<\alpha$; and the states right to the
transition curve are  overcritical ({\oc })  since $\r/a>\alpha$.

 The wave solutions to the Riemann problem for \refe{system} are
combinations of basic 0-waves and 1-waves. Since \refe{system} is a
hyperbolic system of
conservation law,
its wave solutions must satisfy
Lax's entropy condition that the waves from left (upstream) to
right (downstream) should increase their wave speeds so that they
don't cross each other.
For \refe{system} as a resonant nonlinear system, an additional entropy
condition is introduced by Isaacson and Temple,
\bqn
\m{ The standing wave can NOT cross the transition curve } \Gamma.
\eqn
This entropy condition is equivalent to saying that, relative to the
apexes of the fundamental diagrams, traffic conditions upstream and
downstream of inhomogeneities
are on the same side. That is, they should be either both UC or both OC.

With the two entropy conditions, the solutions to the inhomogeneous
LWR model exist and are unique. The wave solutions for {\uc } left
state $U_L$ are shown in \reffig{F_riemann1}, and those for {\oc } left
state $U_L$ are shown in \reffig{F_riemann2}.

\bfg\bc\includegraphics [height=12cm]
{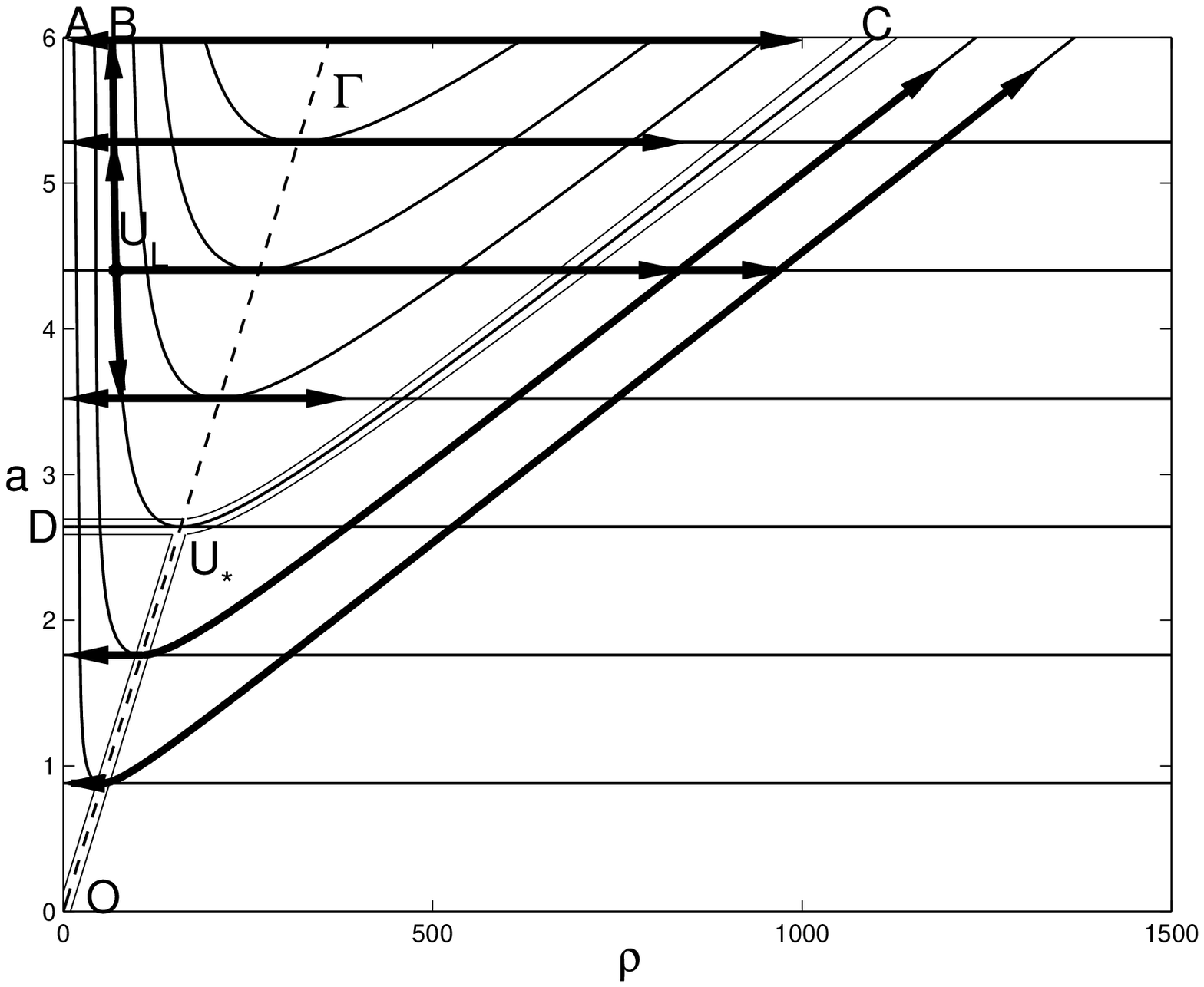}\ec\caption{The Riemann problem for $U_L$ left of
$\Gamma$}\label{F_riemann1}\efg
\bfg\bc\includegraphics [height=12cm]
{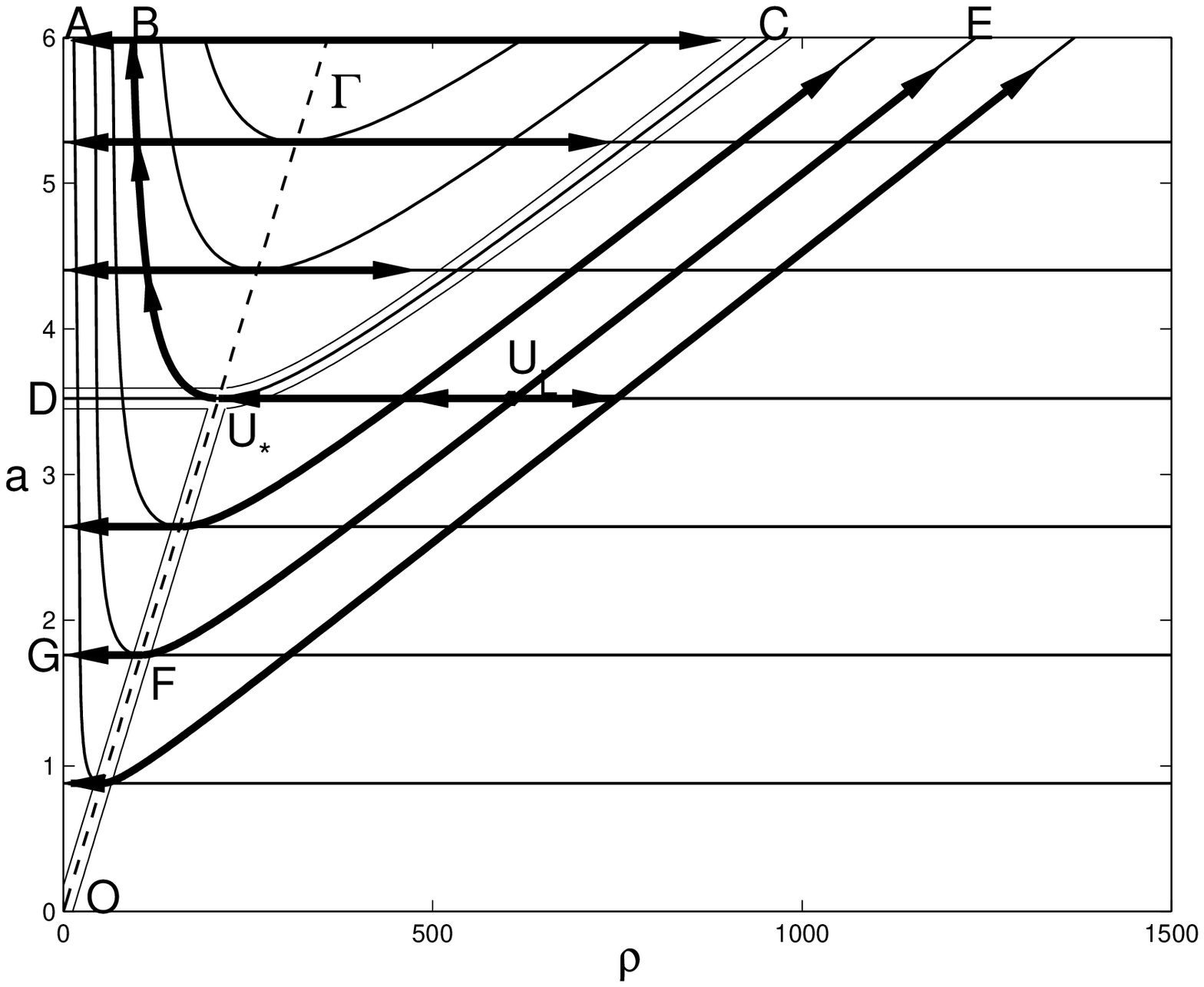}\ec\caption{The Riemann problem for $U_L$ right of
$\Gamma$}\label{F_riemann2}\efg

In the remaining part of this section, we discuss wave solutions to
the Riemann problem for \refe{system}, present the formula for the
boundary flux $f^{\ast}_0$ related to each type of {solution},
summarize our results and compare them with those found in
literature.

\subsection{Solutions of the boundary fluxes} \label{bf}
When $U_L=(a_L,\r_L)$ is {\uc }; i.e., $\r_L/a_L<\alpha$, where
$\alpha$ is defined in \refe{trancur}, we denote the special critical
point on standing wave passing $U_L$ as $U_{\ast}$. Thus, as shown in
\reffig{F_riemann1}, the $U$-space is partitioned into three regions by
$DU_{\ast}$, $OU_{\ast}$ and $U_{\ast}C$, where
$DU_{\ast}=\{(a,\r)|a=a_{\ast},\r<\r_{\ast}\}$,
$OU_{\ast}=\Gamma\cap\{0\leq \r\leq \r_{\ast}\}$ and
$U_{\ast}C=\{(a,\r)|f(a,\r)=f(U_L),\r>\r_{\ast}\}$.  Related to
different positions of the right state $U_R$ in the $U$-space, the
Riemann problem for \refe{system} with initial conditions
\refe{inh.ini} has the following four types of wave solutions. For each
type of solutions we  provide formula for
calculating the associated boundary flux $f^{\ast}_0$.

\bi
\item [Type 1] When $U_R$ is in region $ABU_LU_{\ast}DA$ shown in
\reffig{F_riemann1}; i.e.,
\bqn
f(U_R)<f(U_{\ast})=f(U_L),\quad \r_R/a_R<\alpha \m{ and } a_R\geq
a_{\ast},
\eqn
wave solutions to the Riemann problem are of type 1.  These solutions
consist of two basic waves with an intermediate state
$U_1=(a_R,\r_1|_{f(a_R,\r_1)=f(U_{\ast})=f(U_L)})$. Of these two
waves, the left one $(U_L, U_1)$ is a standing wave, and the right
 one $(U_1,U_R)$ is a rarefaction wave with characteristic velocity
$\l_1(a,\r)>0$.

\bfg\bc\includegraphics[height=12cm] {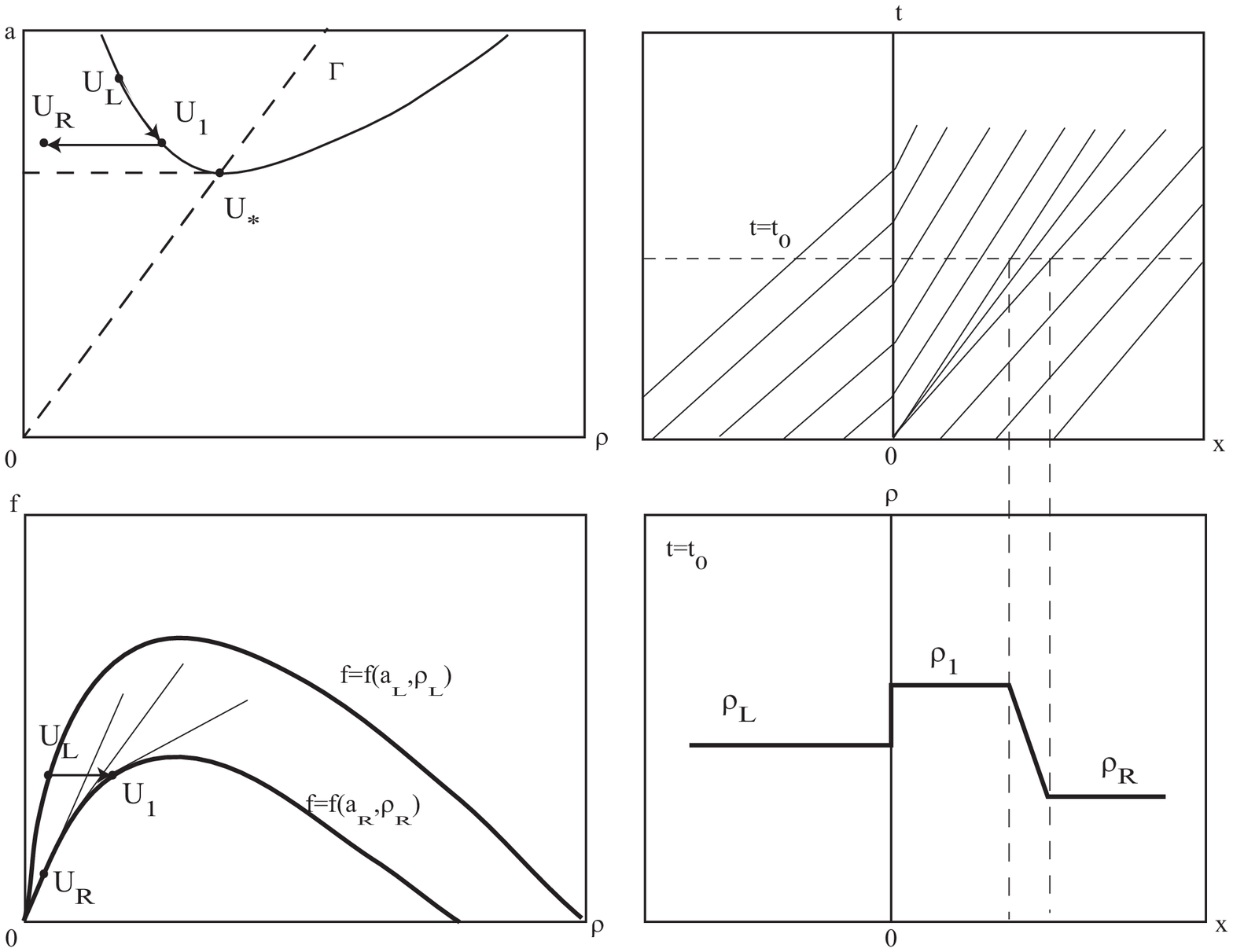}\ec\caption{An
example for wave solutions of type 1 for \refe{system} with initial
conditions \refe{inh.ini}}\label{F_case1}\efg

From \reffig{F_riemann1}, we can see that the Riemann problem may admit
this type of solutions when $a_L>a_R$ or $a_L\leq a_R$; i.e., when
the road merges or diverges at $x=0$. Here we present an example of
this type of solutions in \reffig{F_case1}, where the roadway merges at
$x=0$. In the case  when the roadway diverges at
$x=0$, we can find similar solutions.

From \reffig{F_case1},  we obtain the boundary flux
$f^{\ast}_0=f(U_L)=f(U_{\ast})$ for wave solutions of type 1.

\item [Type 2]  When $U_R$ is in region $BU_LU_{\ast}CB$ shown in
\reffig{F_riemann1}; i.e.,
\bqn
f(U_R)\geq f(U_{\ast})=f(U_L), \label{type2}
\eqn
wave solutions to the Riemann problem are of type 2.  These solutions
consist of two basic waves with an intermediate state
$U_1=(a_R,\r_1|_{f(a_R,\r_1)=f(U_{\ast})=f(U_L)})$. Of these two
waves, the left $(U_L, U_1)$ is a standing wave, and the right
$(U_1,U_R)$ is a shock wave with positive speed
$\sigma=\frac{f(U_R)-f(U_{\ast})}{\r_R-\r_1}>0$.

\bfg\bc\includegraphics[height=12cm] {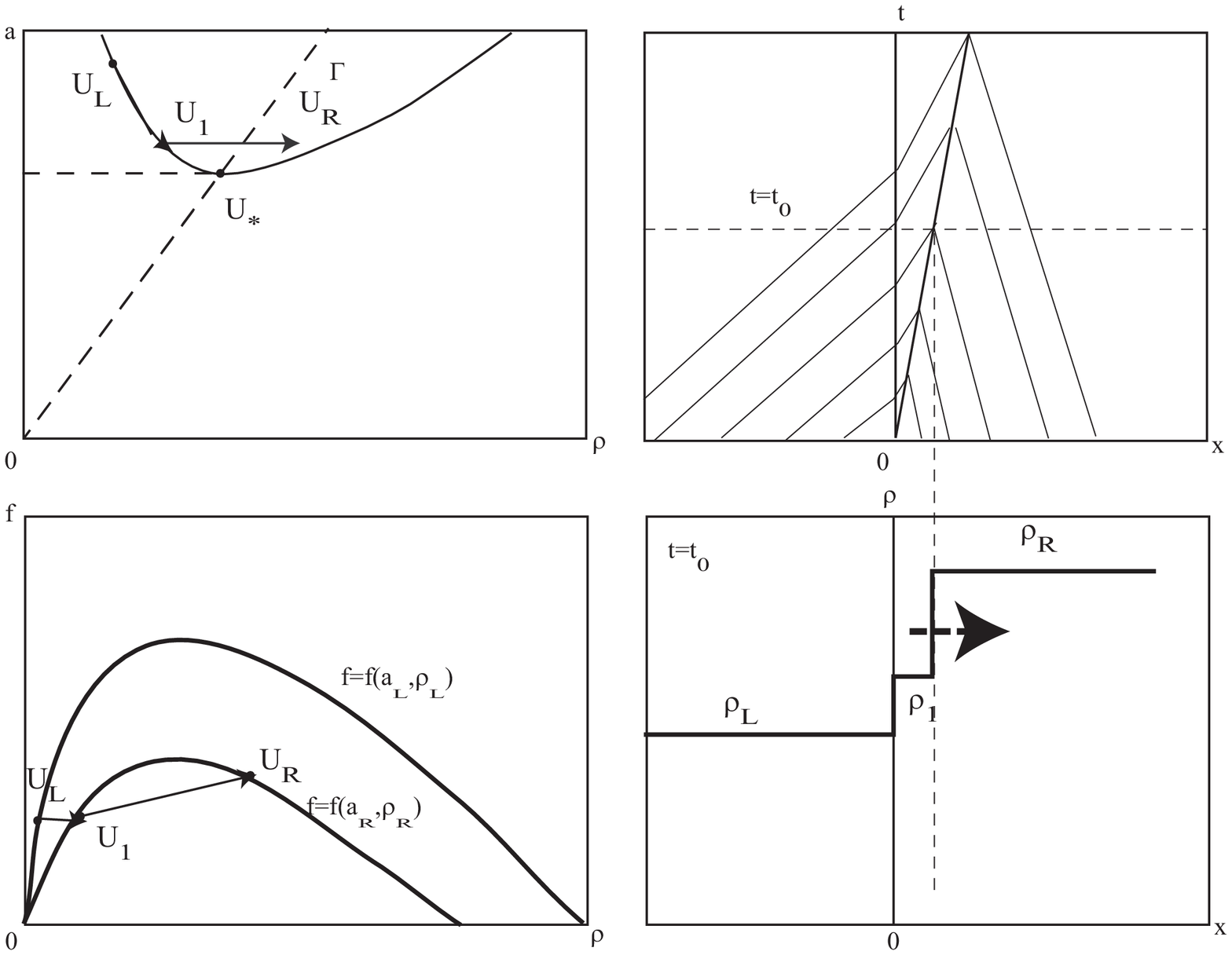}\ec\caption{An
example for wave solutions of type 2 for \refe{system} with initial
conditions \refe{inh.ini}}\label{F_case2}\efg

From \reffig{F_riemann1}, we can see that the Riemann problem  may
admit this type of solutions when the downstream traffic condition
$U_R$ is {\uc } or {\oc }, or the roadway merges or diverges at
$x=0$. Here we present an example of this type of solutions in
\reffig{F_case2}, where the downstream traffic condition is {\oc } and
the  roadway merges at $x=0$. Similar solutions can be found
for other situations that satisfy \refe{type2}.

From \reffig{F_case2},  we obtain the boundary flux
$f^{\ast}_0=f(U_L)=f(U_{\ast})$ for wave solutions of type 2. Here we
have the same formula as that for wave solutions of type 1.

\item [Type 3] When $U_R$ is in region $OU_{\ast}CO$ shown in
\reffig{F_riemann1}; i.e.,
\bqn
f(U_R)<f(U_{\ast})=f(U_L),\qquad \r_R/a_R \geq \alpha,
\eqn
wave solutions to the Riemann problem are of type 3. These solutions
consist of two basic waves with an intermediate state
$U_1=(a_L,\r_1|_{f(a_L,\r_1)=f(U_R)})$.  Of these two waves, the left
one $(U_L, U_1)$ is a shock wave with negative speed
$\sigma=\frac{f(U_1)-f(U_L)}{\r_1-\r_L}<0$, and the right one
$(U_1,U_R)$ is a standing wave.

\bfg\bc\includegraphics[height=12cm] {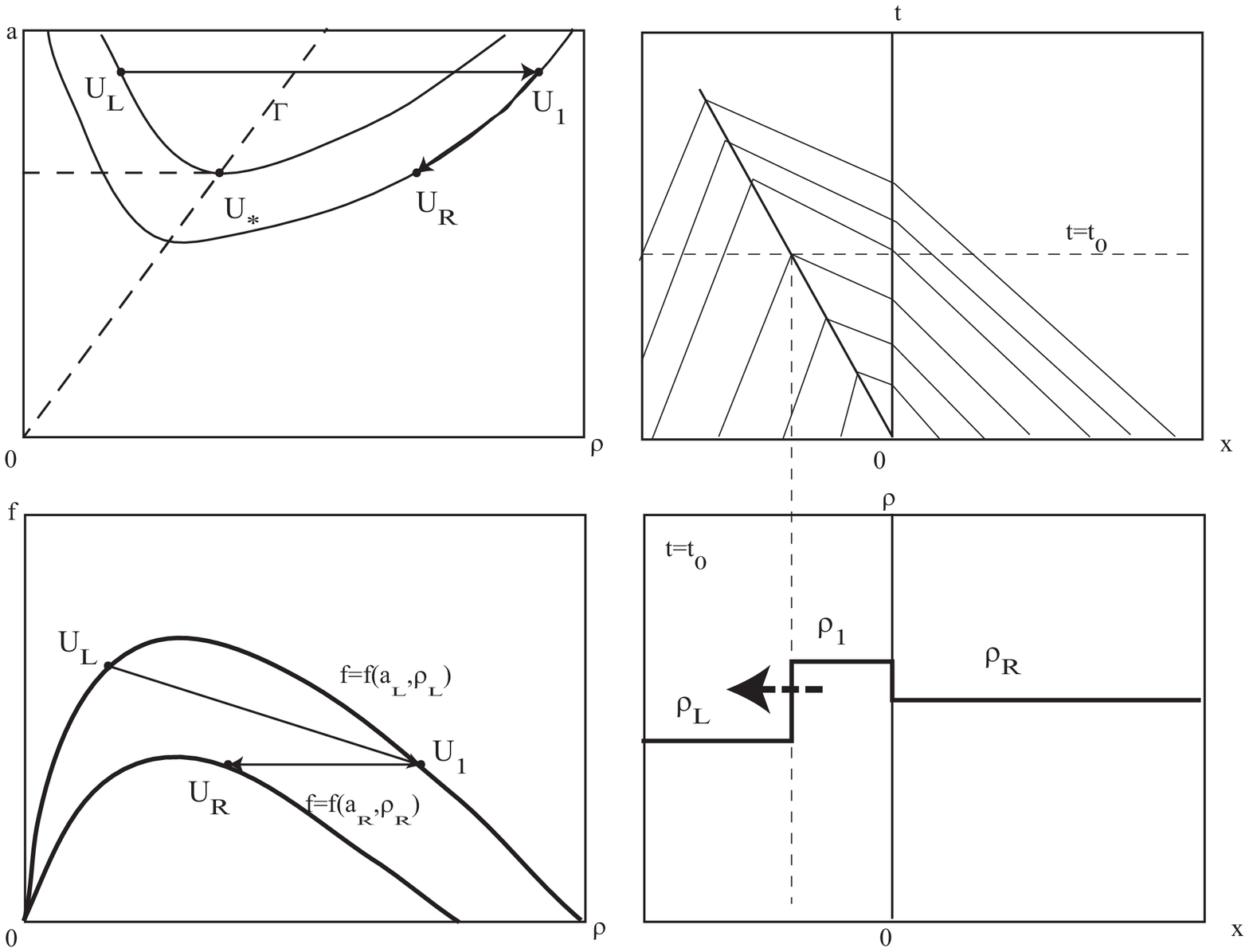}\ec\caption{An
example for wave solutions of type 3 for \refe{system} with initial
conditions \refe{inh.ini}}\label{F_case3}\efg

From \reffig{F_riemann1}, we can see that the Riemann problem may admit
this type of solutions when the roadway  merges or diverges at
$x=0$. Here we present an example of this type of solutions in
\reffig{F_case3}, where the roadway merges at $x=0$. In the case
when  the roadway diverges at $x=0$,  similar solutions can be found.

From \reffig{F_case3},  we obtain the boundary flux $f^{\ast}_0=f(U_R)$
for wave solutions of type 3.

\item [Type 4] When $U_R$ is in region $OU_{\ast}DO$ shown in
\reffig{F_riemann1}; i.e.,
\bqn
f(U_R)<f(U_{\ast})=f(U_L),\quad \r_R/a_R<\r_{\ast}/a_{\ast} \m{ and
}a_R<a_{\ast},
\eqn
wave solutions to the Riemann problem  are of type 4. These solutions
consist of three basic waves with two intermediate states:
$U_1=(a_L,\r_1|_{f(a_L,\r_1)=f(U_2)})$ and
$U_2=(a_R,\r_2|_{\r_2/a_R=\alpha})$. Of these three waves, the left
one $(U_L, U_1)$ is a shock wave with negative speed
$\sigma=\frac{f(U_1)-f(U_L)}{\r_1-\r_L}<0$, the middle one
$(U_1,U_2)$ is a standing wave with zero speed, and the right one
$(U_2,U_R)$ is a rarefaction wave with characteristic velocity
$\l_1(a,\r)>0$.

\bfg\bc\includegraphics[height=12cm] {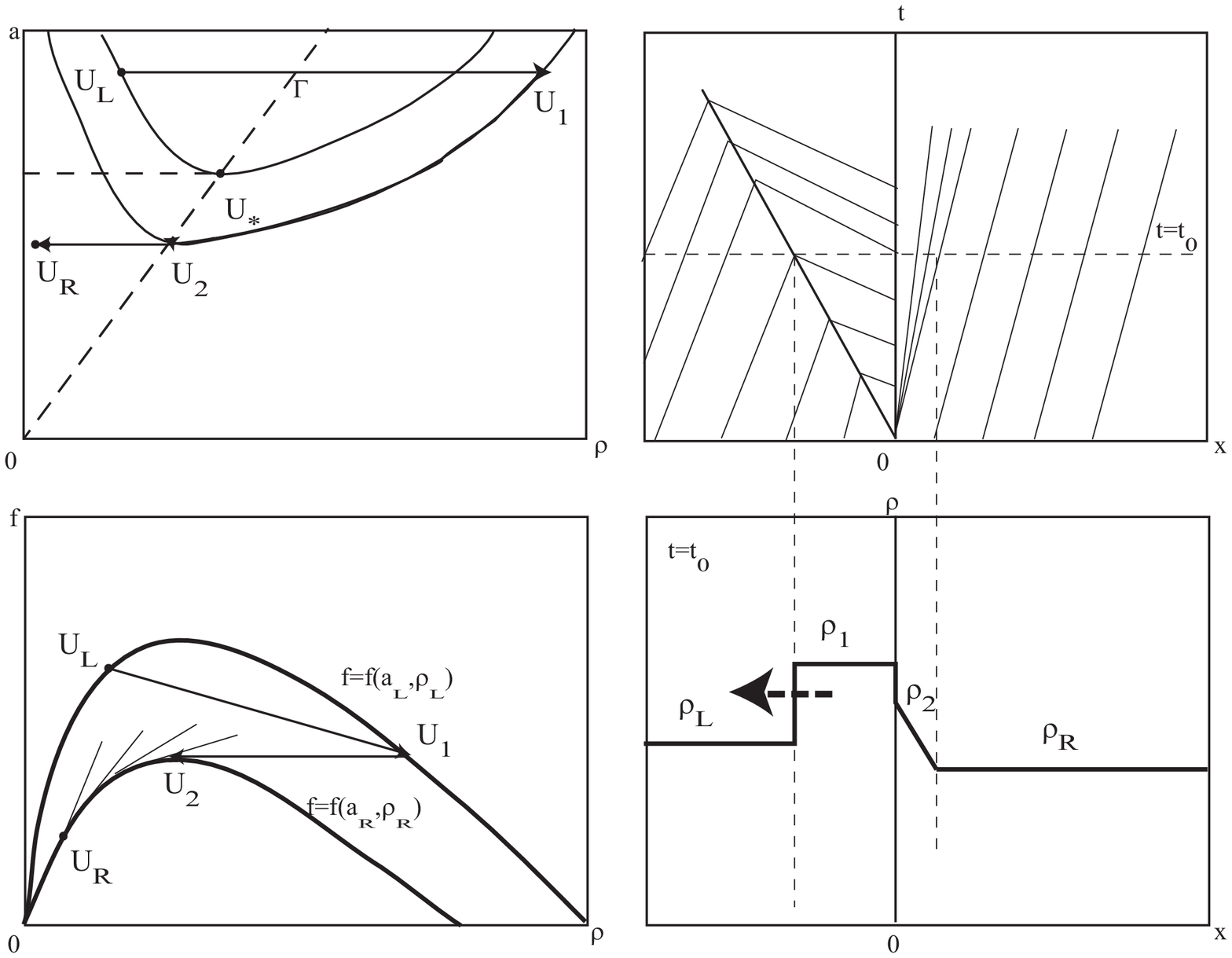}\ec\caption{An
example for wave solutions of type 4 for \refe{system} with initial
conditions \refe{inh.ini}}\label{F_case4}\efg

From \reffig{F_riemann1}, we can see that this type of solutions are
admitted only when the roadway merges at $x=0$. Here we present an
example of this type of  solutions in \reffig{F_case4}.

From \reffig{F_case4},  we obtain the boundary flux $f^{\ast}_0=f(U_2)$
for wave solutions of type 4.

\ei

When $U_L=(a_L,\r_L)$ is {\oc }; i.e., $\r_L/a_L>\alpha$, we denote
 the special critical
point on 1-wave curve passing $U_L$ as $U_{\ast}$; i.e.,
$U_{\ast}=(a_L,\r_{\ast}|_{\r_{\ast}/a_L=\alpha})$. Thus, as shown in
\reffig{F_riemann2}, the $U$-space is partitioned into three regions
by  three curves $DU_{\ast}=\{a=a_{\ast}=a_L,0\leq \r\leq
\r_{\ast}\}$, $OU_{\ast}=\{0\leq a\leq a_{\ast},\r=a\alpha\}$ and
$U_{\ast}C=\{a\geq a_{\ast},f(a,\r)=f(U_{\ast})$. Related to
different positions of the right state $U_R$ in the $U$-space, the
Riemann problem for \refe{system} with initial conditions
\refe{inh.ini} has the following six types of wave solutions. For each
type of solutions we  provide formula for
calculating the associated boundary flux $f^{\ast}_0$.

\bi
\item [Type 5] When $U_R$ resides in region $ABU_{\ast}DA$ shown in
\reffig{F_riemann2}; i.e.,
\bqn
f(U_R)<f(U_{\ast}), \quad \r_R/a_R<\alpha \m{ and } a_R \geq
a_{\ast}=a_L,
\eqn
wave solutions to the Riemann problem are of type 5. These solutions
consist of three basic waves with two intermediate states:
$U_1=U_{\ast}$ and $U_2=(a_R,\r_2|_{f(U_2)=f(U_{\ast})})$. Of these
three waves, the left one $(U_L,U_1)$ is a rarefaction wave with
negative characteristic wave velocity $\l_1(a,\r)$, the middle one
$(U_1,U_2)$ is a standing wave and the right one $(U_2,U_R)$ is a
rarefaction wave with positive characteristic velocity $\l_1(a,\r)$.

\bfg\bc\includegraphics[height=12cm] {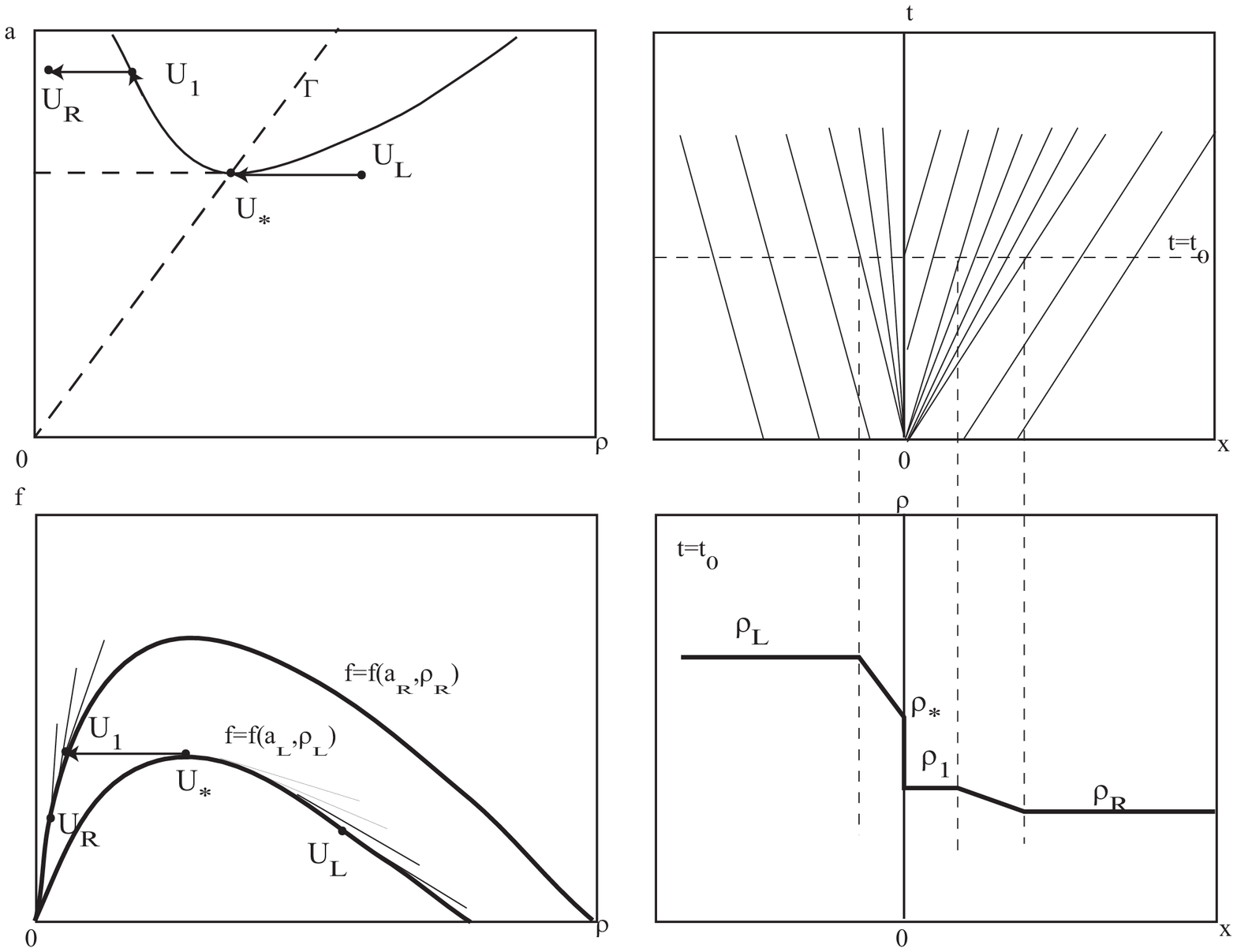}\ec\caption{An
example for wave solutions of type 5 for \refe{system} with initial
conditions \refe{inh.ini}}\label{F_case5}\efg

From \reffig{F_riemann2}, we can see that this type of solutions are
admitted only when the roadway diverges at $x=0$.
Here we present an example of this type of solutions in
\reffig{F_case5}.

From \reffig{F_case5}, we obtain the boundary flux $f^{\ast}_0=f(U_2)$
for wave solutions of type 5.

\item [Type 6] When $U_R$ resides in region $BU_{\ast}CB$  shown in
\reffig{F_riemann2}; i.e.,
\bqn
f(U_R)\geq f(U_{\ast}),
\eqn
solutions to the Riemann problem are of type 6. These solutions
consist of three basic waves with two intermediate states:
$U_1=U_{\ast}$ and $U_2=(a_R,\r_2|_{f(U_2)=f(U_{\ast})})$. Of these
three waves, the left one $(U_L,U_1)$ is a rarefaction wave with
negative characteristic velocity $\l_1(a,\r)$, the middle one
$(U_1,U_2)$ is a standing wave and the right one $(U_2,U_R)$ is a
shock wave with positive speed
$\sigma=\frac{f(U_R)-f(U_2)}{\r_R-\r_2}$.

\bfg\bc\includegraphics[height=12cm] {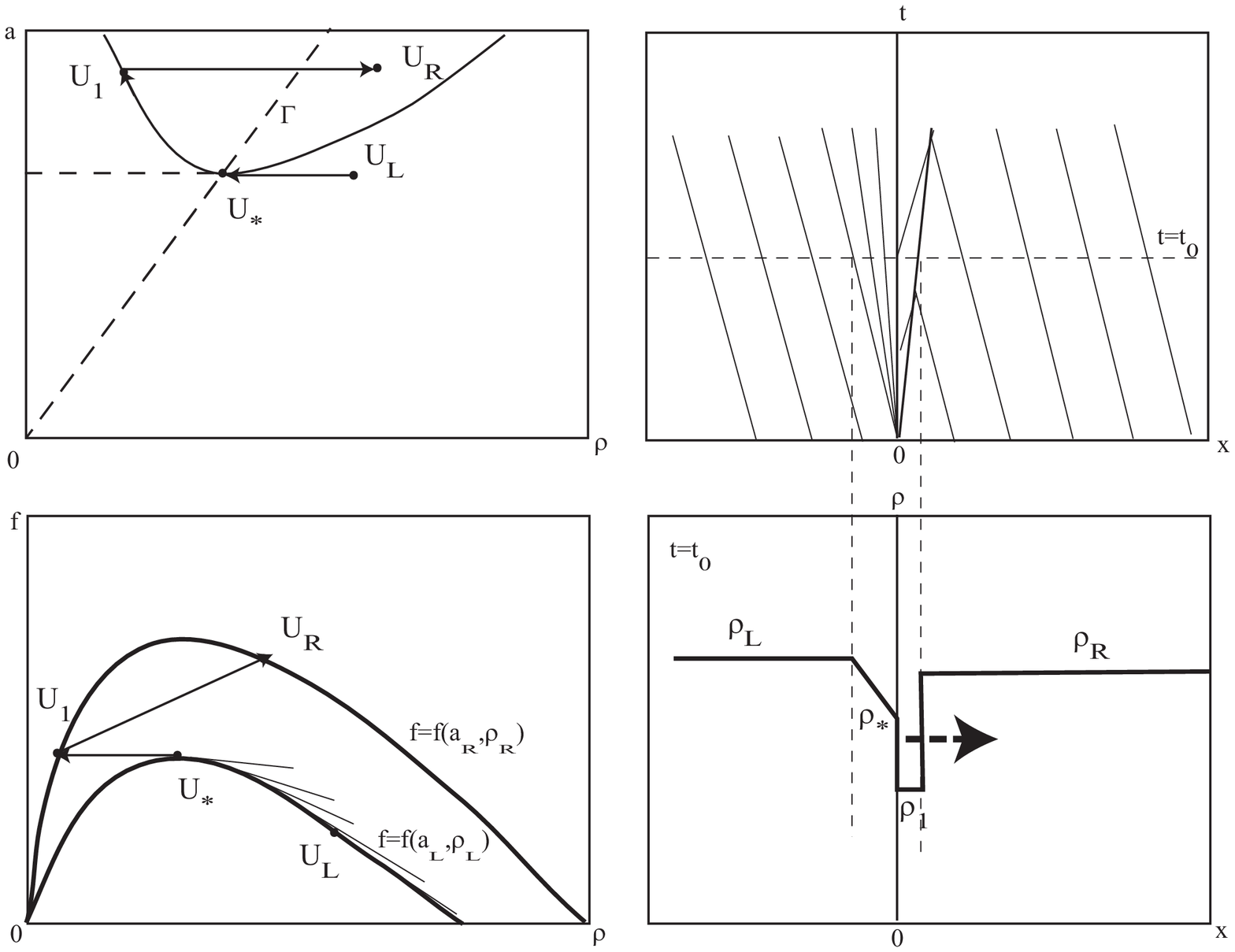}\ec\caption{An
example for wave solutions of type 6 for \refe{system} with initial
conditions \refe{inh.ini}}\label{F_case6}\efg

From \reffig{F_riemann2}, we can see that this type of solutions may be
admitted when the downstream traffic condition is {\uc } or {\oc };
However, they are admitted only when the roadway diverges at $x=0$.
Here we present an example of this type of solutions in
\reffig{F_case6}, where the downstream traffic condition is {\oc }. In
the case when the downstream traffic condition is {\uc }, we can find
similar solutions.

From \reffig{F_case6}, we obtain the boundary flux $f^{\ast}_0=f(U_2)$
for this type of wave solutions. Here we have the same formula as
that for wave solutions of type 5.

\item [Type 7] When $U_R$ resides in region $CU_{\ast}FU_LEC$ shown
in \reffig{F_riemann2}; i.e.,
\bqn
f(U_L)\leq f(U_R)<f(U_{\ast}) \m{ and }\r_R/a_R\geq\alpha,
\eqn
wave solutions to the Riemann problem are of type 7. These solutions
consist of two basic waves with an intermediate state
$U_1=(a_L,\r_1|_{f(U_1)=f(U_R)})$. Of these two waves, the left one
$(U_L,U_1)$ is a rarefaction with negative characteristic velocity
$\l_1(a,\r)$, and the right one $(U_1,U_R)$ is a standing wave.

\bfg\bc\includegraphics[height=12cm] {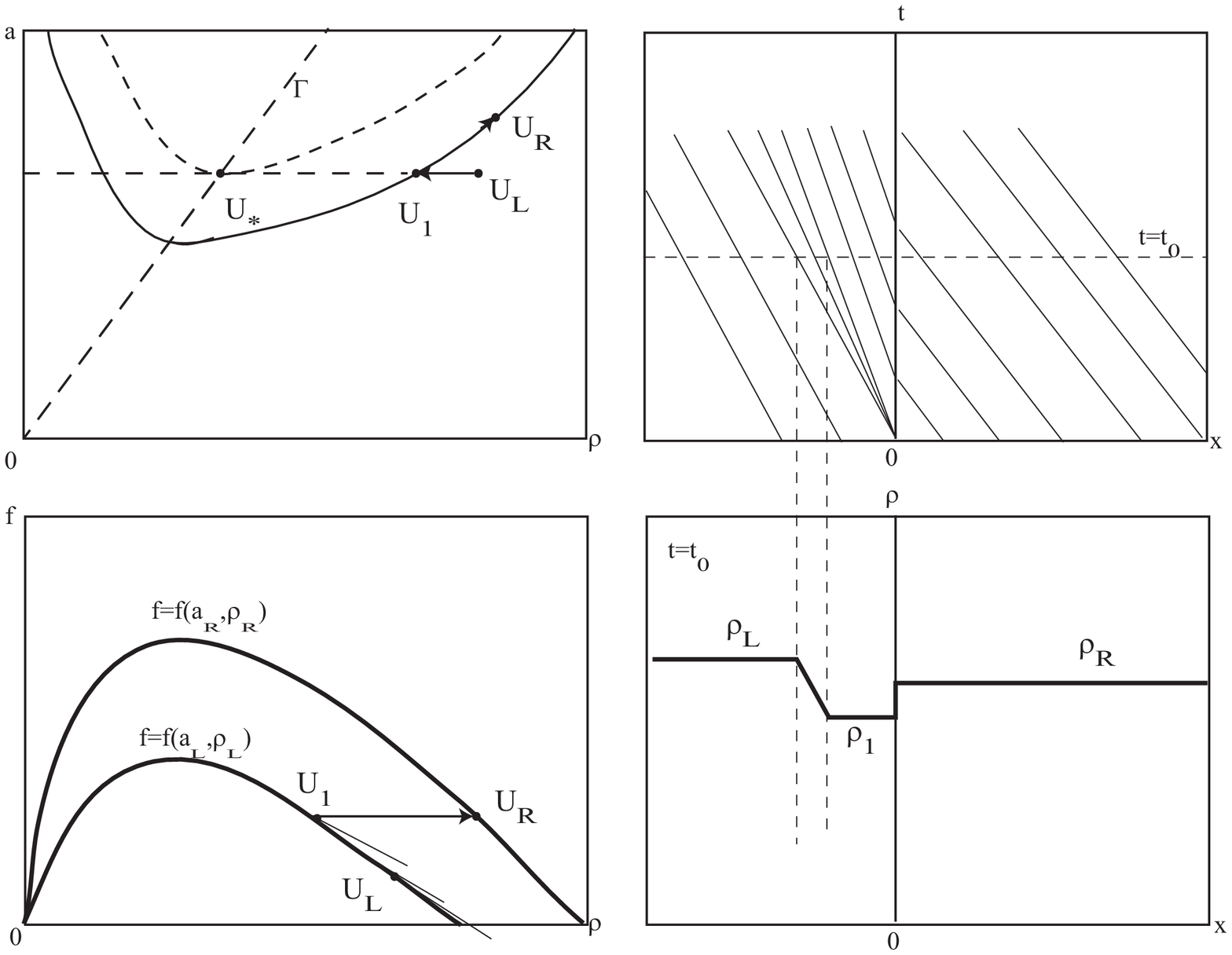}\ec\caption{An
example for wave solutions of type 7 for \refe{system} with initial
conditions \refe{inh.ini}}\label{F_case7}\efg

From \reffig{F_riemann2}, we can see that the Riemann problem may admit
this type of solutions when the roadway merges or diverges at $x=0$.
Here we present an example of this type of solutions in
\reffig{F_case7}, where the roadway diverges at $x=0$.
In the case when the roadway merges, we can find similar solutions.

From \reffig{F_case7}, we obtain the boundary flux $f^{\ast}_0=f(U_R)$
for wave solutions of type 7.

\item [Type 8] When $U_R$ locates in region $FU_LEOF$ shown in
\reffig{F_riemann2}; i.e.,
\bqn
f(U_R)<f(U_L)<f(U_{\ast}) \m{ and } \r_R/a_R \geq \alpha,
\eqn
wave solutions to the Riemann problem are of type 8. These solutions
consist of two basic waves with an intermediate state
$U_1=(a_L,\r_1|_{f(U_1)=f(U_R)})$. Of these two waves, the left one
$(U_L,U_1)$ is a shock with negative speed
$\sigma=\frac{f(U_L)-f(U_1)}{\r_L-\r_1}$, and the right one
$(U_1,U_R)$ is a standing wave.

\bfg\bc\includegraphics[height=12cm] {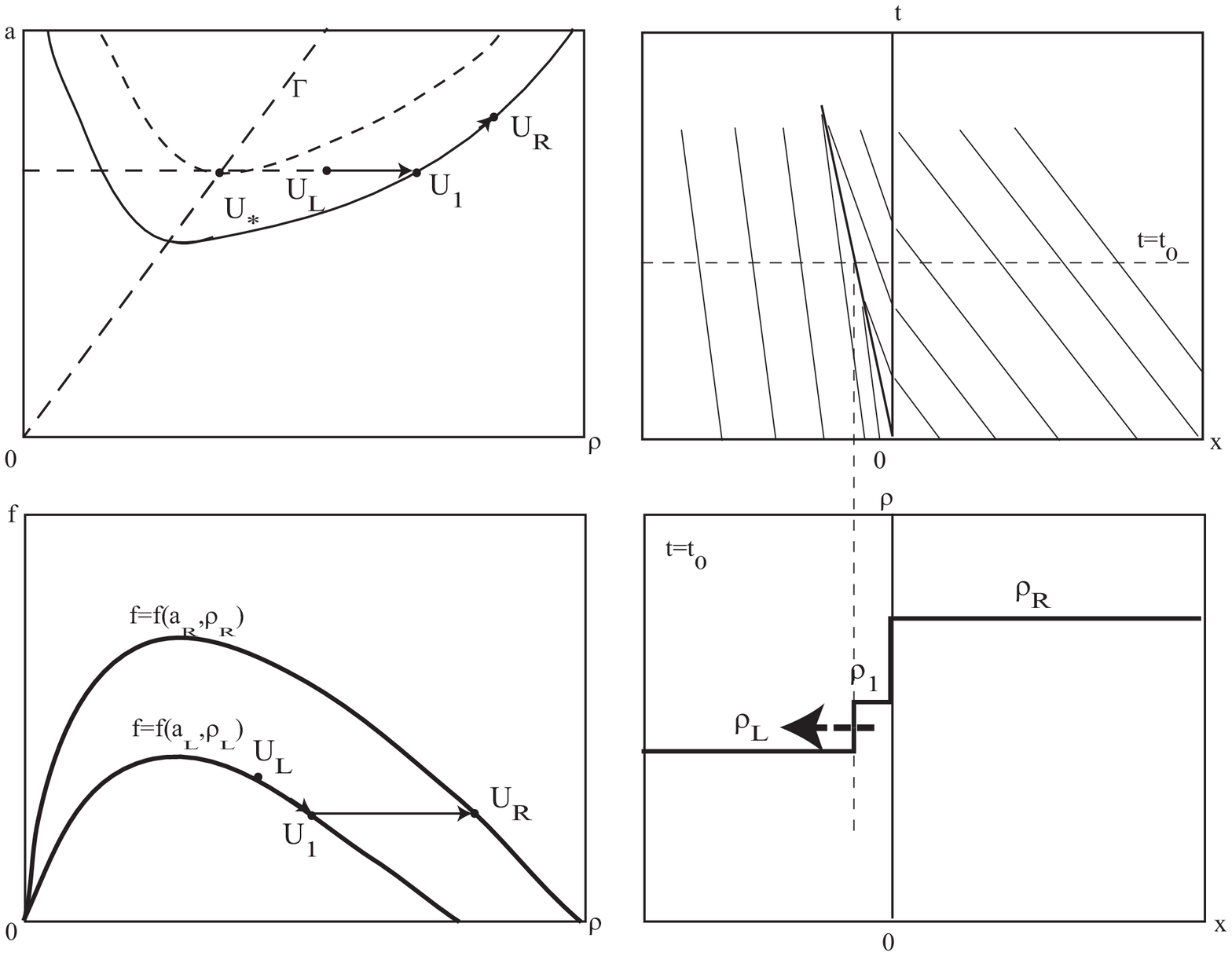}\ec\caption{An
example for wave solutions of type 8 for \refe{system} with initial
conditions \refe{inh.ini}}\label{F_case8}\efg

Like in the previous case, the Riemann problem may admit
this type of solutions when the roadway merges or diverges at $x=0$.
Here we present an example of this type of solutions in
\reffig{F_case8}, where the roadway diverges at $x=0$.
In the case when the roadway merges, we can find similar solutions.

From \reffig{F_case8}, we obtain the boundary flux $f^{\ast}_0=f(U_R)$
for wave solutions of type 8. The formula is the same as that
for wave solutions of type 7.

\item [Type 9] When $U_R$ resides in region $DU_{\ast}FGD$ shown in
\reffig{F_riemann2}; i.e.,
\bqn
f(U_L)\leq f(U_R)<f(U_{\ast}),\quad \r_R/a_R<\alpha \m{ and }
a_R<a_{\ast}=a_L,
\eqn
wave solutions to the Riemann problem are of type 9. These solutions
consist of three basic waves with two intermediate states:
$U_1=(a_L,\r_1|_{f(U_1)=f(U_2)})$ and
$U_2=(a_R,\r_2|_{\r_2/a_R=\alpha})$. Of these three waves, the left
one  $(U_L,U_1)$ is a rarefaction with negative characteristic velocity
$\l_1(a,\r)$, the middle one $(U_1,U_2)$ is a standing wave, and the
right one $(U_2,U_R)$ is a rarefaction with positive speed
$\l_1(a,\r)$.

\bfg\bc\includegraphics[height=12cm] {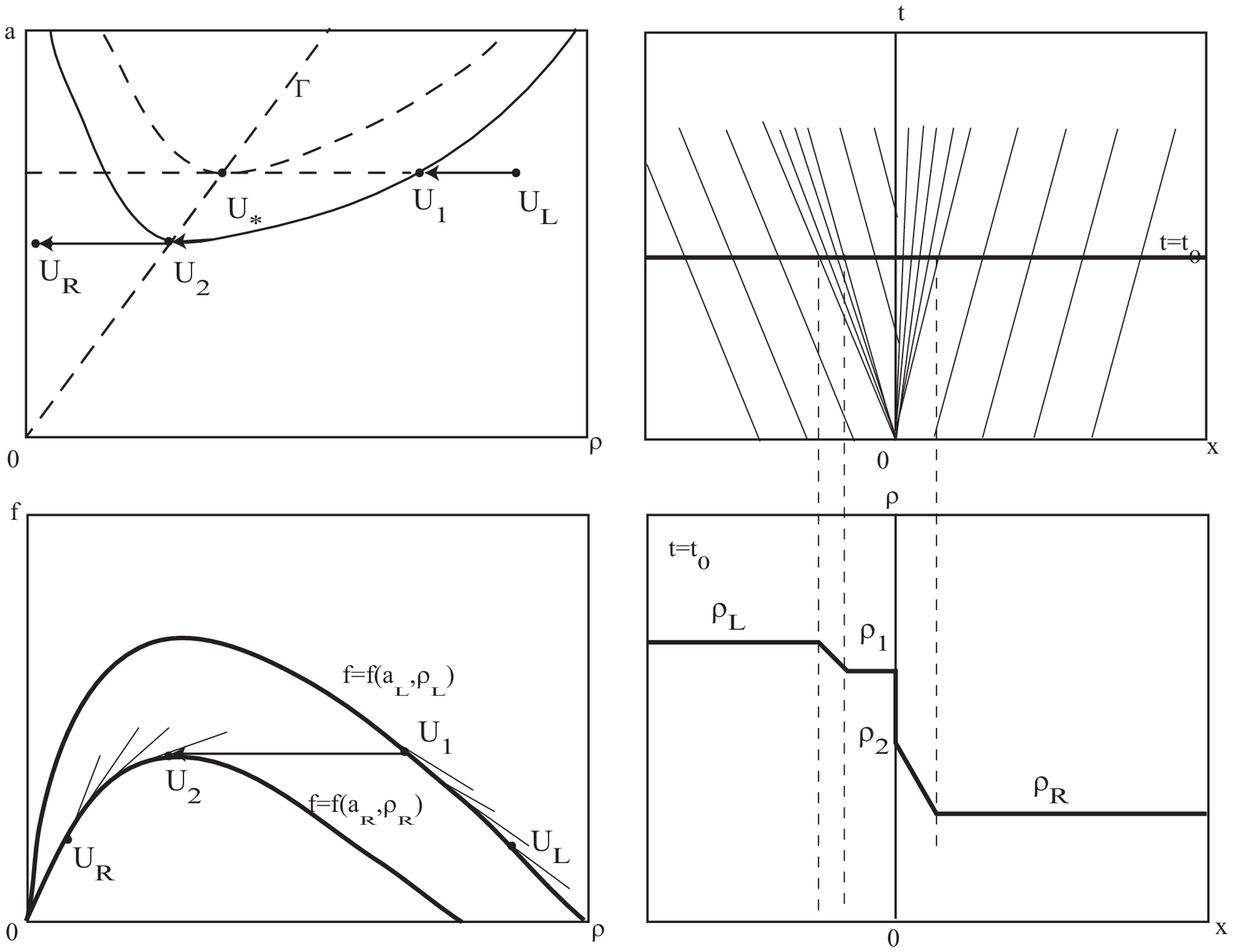}\ec\caption{An
example for wave solutions of type 9 for \refe{system} with initial
conditions \refe{inh.ini}}\label{F_case9}\efg

From \reffig{F_riemann2}, we can see that this type of solutions are
admitted only when the roadway merges at $x=0$. Here
we present an example of this type of solutions in \reffig{F_case9}.

From \reffig{F_case9}, we obtain the boundary flux $f^{\ast}_0=f(U_2)$
for wave solutions of type 9.

\item [Type 10] When $U_R$ resides in region $GFOG$ shown in
\reffig{F_riemann2}; i.e.,
\bqn
f(U_R)<f(U_L)<f(U_{\ast}),\quad \r_R/a_R<\alpha \m{ and
}a_R<a_{\ast}=a_L,
\eqn
wave solutions to the Riemann problem are of type 10. These solutions
consist of three basic waves with two intermediate states:
$U_1=(a_L,\r_1|_{f(U_1)=f(U_2)})$ and
$U_2=(a_R,\r_2|_{\r_2/a_R=\alpha})$. Of these three waves, the left
one $(U_L,U_1)$ is a shock with negative speed, the middle one
$(U_1,U_2)$ is a standing wave, and the right one $(U_2,U_R)$ is a
rarefaction wave with positive characteristic velocity $\l_1(a,\r)$.

\bfg\bc\includegraphics[height=12cm] {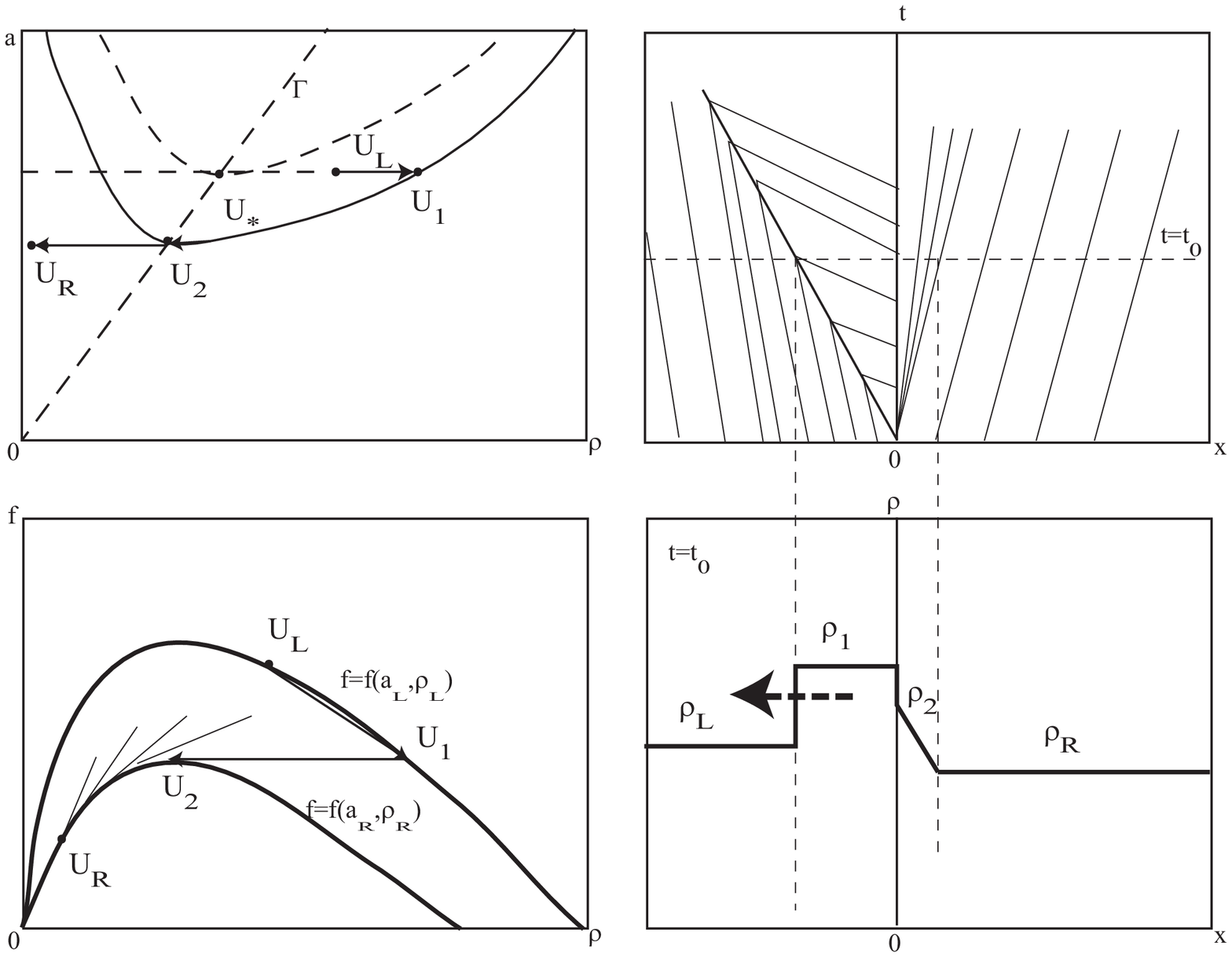}\ec\caption{An
example for wave solutions of type 10 for \refe{system} with initial
conditions \refe{inh.ini}}\label{F_case10}\efg

Like in the previous case, this type of solutions are
admitted only when the roadway merges at $x=0$. Here
we present an example of this type of solutions in \reffig{F_case10}.

From \reffig{F_case10}, we obtain the boundary flux $f^{\ast}_0=f(U_2)$
for wave solutions of type 10. Here we have the same formula as that
for wave solutions of type 9.
\ei

\subsection{Summary}
In each of the 10 cases discussed above, the boundary flux
$f^{\ast}_0$ is equal to one of the following four quantities: the
upstream flow-rate $f(U_L)$, the downstream flow-rate $f(U_R)$, the
capacity of the upstream roadway $f^{max}_L$ and the capacity of the
downstream roadway $f^{max}_R$.  For wave solutions of type 1 and 2,
the boundary flux is equal to the upstream traffic flow-rate; i.e.,
$f^{\ast}_0=f(U_L)$. For wave solutions of type 3, 7 and 8, the
boundary flux is equal to the downstream traffic flow-rate; i.e.,
$f^{\ast}_0=f(U_R)$. For wave solutions of type 4, 9 and 10, the
boundary flux is equal to the capacity of the downstream roadway;
i.e., $f^{\ast}_0=f^{max}_R$. For wave solutions of type 5 and 6, the
boundary flux is equal to the capacity of the upstream roadway; i.e.,
$f^{\ast}_0=f^{max}_L$. In Table \ref{Table0}, the boundary fluxes
are listed for the 10 types of wave solutions to the Riemann problem,
as well as the conditions when the Riemann problem admit those
solutions.

\btb
\begin{tabular}{|c||l|l||c|}\hline
No.&left state $U_L$ & right state $U_R$ &$f^{\ast}_0$
\\\hline
1&{\uc } &$f(U_R)<f(U_L)$, $a_R>a_{\ast}$, $\r_R/a_R<\alpha$ &
$f(U_L)$ \\\hline
2&{\uc } &$f(U_R)>f(U_L)$ & $f(U_L)$ \\\hline
3&{\uc } & $f(U_R)<f(U_L)$, $\r_R/a_R>\alpha$ &$f(U_R)$ \\\hline
4&{\uc } & $f(U_R)<f(U_L)$, $\r_R/a_R<\alpha$, $a_R<a_{\ast}$ &
$f^{max}_R$ \\\hline
5&{\oc } & $f(U_R)<f^{max}_L$, $a_R>a_L$, $\r_R/a_R<\alpha$ &
$f^{max}_L$ \\\hline
6&{\oc } & $f(U_R)>f^{max}_L$ & $f^{max}_L$ \\\hline
7&{\oc } & $f(U_L)<f(U_R)<f^{max}_L$, $\r_R/a_R>\alpha$ & $f(U_R)$
\\\hline
8&{\oc } & $f(U_R)<f(U_L)$, $\r_R/a_R>\alpha$ & $f(U_R)$ \\\hline
9&{\oc } & $f(U_L)<f(U_R)<f^{max}_L$, $\r_R/a_R<\alpha$, $a_R<a_L$ &
$f^{max}_R$ \\\hline
10&{\oc } & $f(U_R)<f(U_L)$, $\r_R/a_R<\alpha$, $a_R<a_L$ &
$f^{max}_R$ \\\hline
\end{tabular}
\caption {Solutions of the boundary fluxes $f^{\ast}_0$}\label{Table0}
\etb

Note that when $a_L=a_R$; i.e., when \refe{inh_1st} becomes a
homogeneous LWR model, wave solutions and the solutions of the
boundary fluxes provided here are the same as those for the
homogeneous LWR model.

\citet{lebacque1996godunov} studied the Riemann problem
of the inhomogeneous LWR for \refe{inh_b}.
He classified the problem according to two criteria. The first  criterion
is to compare
capacity of the upstream cell and that of the downstream cell. For the
roadway with variable  number of lanes, it is equivalently to
compare the number of lanes of the
upstream cell and that of the downstream cell.
The second criterion is to consider
whether the  upstream and downstream traffic conditions are {\uc } or
{\oc }. With  these criteria, he discussed 8 types of waves solutions to the
Riemann problem and obtained the formula for the boundary flux
related to each type of solutions. The conditions for those types of
wave solutions as well as the formulas related to those types of
solutions are listed in Table \ref{table1}. Under each of those
conditions, the Riemann problem may admit different types of
solutions discussed in Section \ref{bf}. The types of
solutions and our related formulas for the boundary flux are also
presented in Table \ref{table1}. From this table, we can see that our
results are consistent with those provided by Lebacque, although the
Riemann problem is solved through different approaches.

\btb
\begin{tabular}{|l|c||l|c|}\hline
Conditions& Solutions by Lebacque & Types &Our solutions\\\hline
$a_L\leq a_R$, $U_L$ \uc, $U_R$ \uc& $f(U_L)$& 1  & $f(U_L)$ \\\hline
$a_L\leq a_R$, $U_L$ \uc, $U_R$ \oc& $\min \{f(U_L),f(U_R)\}$ & 2 or
3 & $f(U_L)$ or $f(U_R)$\\\hline
$a_L\leq a_R$, $U_L$ \oc, $U_R$ \uc & $f^{max}_L$ & 5 or 6 &
$f^{max}_L$ \\\hline
$a_L\leq a_R$, $U_L$ \oc, $U_R$ \oc & $\min \{f^{max}_L, f(U_R)\}$
&6, 7 or 8&$f^{max}_L$ or $f(U_R)$\\\hline
$a_L\geq a_R$, $U_L$ \uc, $U_R$ \uc & $\min \{f^{max}_R, f(U_L)$ & 1
or 4 &$f(U_L)$, $f^{max}_R$ \\\hline
$a_L\geq a_R$, $U_L$ \uc, $U_R$ \oc & $\min \{f(U_L),f(U_R)\}$ & 2 or
3 & $f(U_L)$ or $f(U_R)$  \\\hline
$a_L\geq a_R$, $U_L$ \oc, $U_R$ \uc & $f^{max}_R$ & 9 or 10 &
$f^{max}_R$\\\hline
$a_L\geq a_R$, $U_L$ \oc, $U_R$ \oc & $f(U_R)$ & 7 or 8 & $f(U_R)$\\\hline
\end{tabular}
\caption{Comparison with Lebacque's results}\label{table1}
\etb

The consistency of our results with existing results can also be
shown by introducing a simple formula for the boundary flux.
If we define the upstream demand as
\bqn
f^{\ast}_L&=&\cas{{ll} f(U_L), & \r_L/a_L<\alpha \\f^{max}_L, &
\r_L/a_L\geq \alpha}
\eqn
and define the downstream supply as
\bqn
f^{\ast}_R&=&\cas{{ll} f^{max}_R, & \r_R/a_R<\alpha \\f(U_R),
&\r_R/a_R\geq \alpha}
\eqn
then the boundary flux can be simply computed as
\bqn
f^{\ast}_0&=&\min\{f^{\ast}_L, f^{\ast}_R\}. \label{ds}
\eqn
Note that $f^{\ast}_L=f(U_{\ast})$.
Formula \refe{ds} was also provided by \citet{daganzo1995ctm} and \citet{lebacque1996godunov}.

%Wenlong Jin: 2003
\section{Simulation of traffic flow on a ring road with a bottleneck}
\subsection{Solution method} \label{go}
The augmented inhomogeneous LWR model, expressed in
conservation  form \refe{system}, can be solved efficiently
with Godunov's method under general initial and boundary conditions.
In  Godunov's method, the roadway is partitioned
into $N$ cells and a duration of time is discretized into $M$ time
steps. In a cell $i$, we approximate the continuous equation
\refe{system} with a finite difference equation
\bqn
\frac {U_i^{m+1}-U_i^m} {\dt}+\frac
{F^{\ast}_{i-1/2}-F^{\ast}_{i+1/2}}{\dx}&=&0,
\eqn
whose component for $\r$ is
\bqn
\frac {\r_i^{m+1}-\r_i^m} {\dt}+\frac
{f^{\ast}_{i-1/2}-f^{\ast}_{i+1/2}}{\dx}&=&0, \label{dis.1}
\eqn
where $\r_i^m$ denotes the average of $\r$ in cell $i$ at time step
$m$, similarly $\r_i^{m+1}$ is the average at time step $m+1$;
$f^{\ast}_{i-1/2}$ denotes the flux through the upstream boundary of
cell $i$, and similarly $f^{\ast}_{i+1/2}$ denotes the downstream
boundary flux of cell $i$. In \refe{dis.1}, the boundary flux
$f^{\ast}_{i-1/2}$ is related to solutions to a Riemann problem for
\refe{system} with the following initial conditions:
\bqn
U(x=x_{i-1/2}, t=t_{m}) &=&\cas{{ll} U^m_{i-1} & x <x_{i-1/2}\\U^m_i
& x>x_{i-1/2}},
\eqn
which have been discussed in Section \ref{ri}.

%Wenlong Jin: 2003
\subsection {Numerical results}

We use the approximation developed earlier to
simulate traffic on a ring road. The  length of the ring road  is
$L=800l$ = 22.4 km. The simulation time is  $T=500 \tau$ = 2500 s =
41.7 min. We partition the road $[0, L]$ into
$N=100$ cells and the time interval $[0, T]$ into $K=500$ steps.
Hence, the length of each cell is $\dx=0.224$ km and the length of
each time step is $\dt=5$ s. Since $|\l_{\ast}|\leq v_f=5 l/\tau$, we
find the CFL \citep{courant1928CFL} condition number
\bqs
\max |\l_{\ast}| \frac {\dt}{\dx}&\leq& 0.625<1.
\eqs
Moreover,  we adopt in this simulation the fundamental diagram
used  in (Kerner and  Konh\"auser, 1994;
Herrmann and Kerner, 1998) with the following parameters: the
relaxation time
$\tau$ = 5 s;  the unit length $l$ = 0.028 km; the free flow speed
$v_f=5.0 l/\tau$ = 0.028 km/s = 100.8 km/h; the jam density of a
single lane $\r_j$ = 180 veh/km/lane.  The equilibrium speed-density relationship is
therefore
\bqs
v_{\ast}(\rho,a(x))=5.0461\left
[\left(1+\exp\{[\frac{\r}{a(x)\r_j}-0.25]/0.06\}\right)^{-1}-3.72\times
10^{-6}\right] l/\tau,
\eqs
where $a(x)$ is the number of lanes at location $x$.
The equilibrium functions $\vs(\r,a(x))$ and $f(\r,a(x))$ are
given in \reffig {fd_Kerner}.

\bfg
\bc\includegraphics[height=12cm] {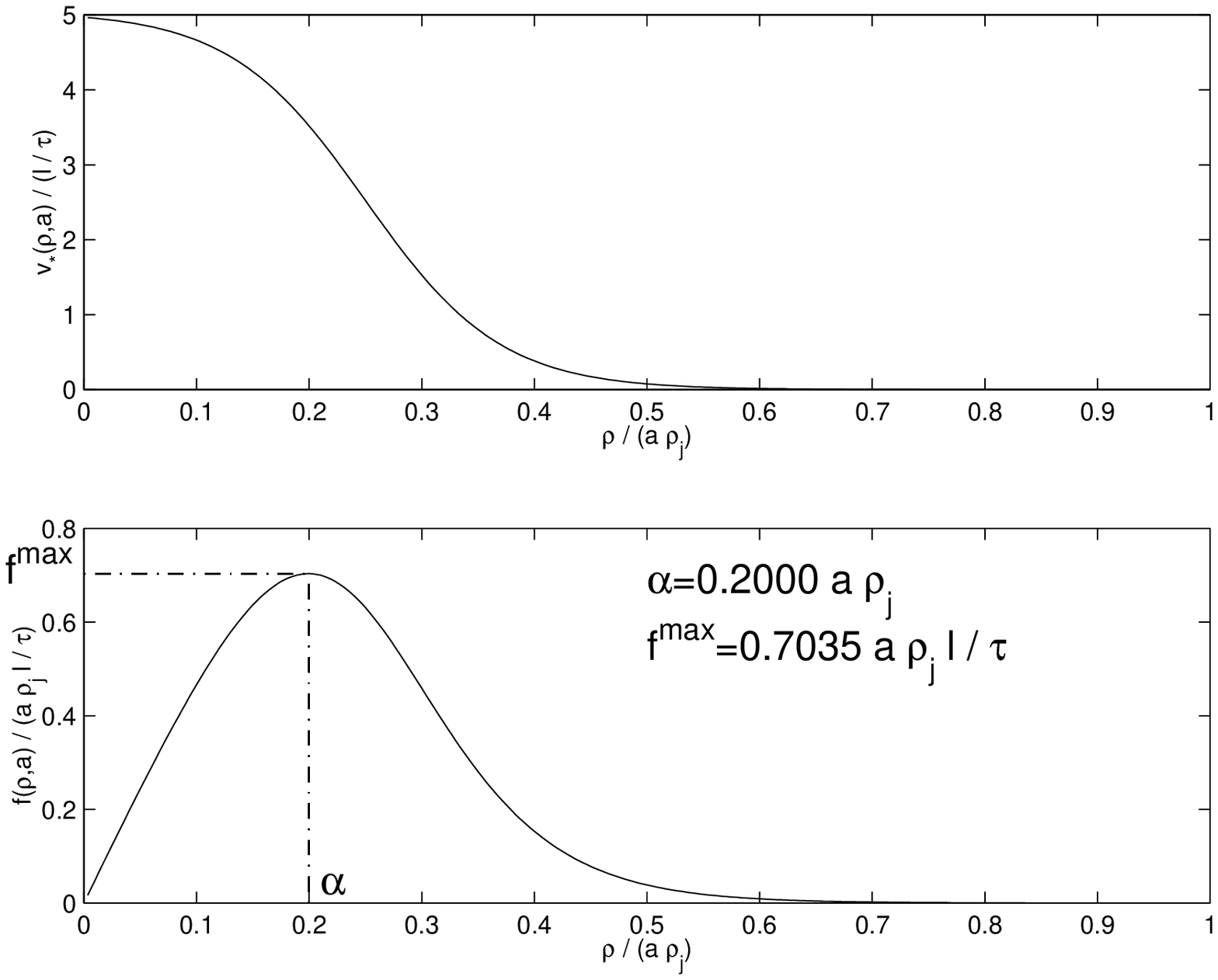}\ec
\caption {The Kerner-Konh\"auser model of speed-density and
flow-density relations} \label {fd_Kerner}
\efg

The first simulation is about the homogeneous LWR model. Here we
assume that the  ring road has single lane everywhere; i.e., $a(x)=1$
for any $x\in[0,L]$ , and  use a global perturbation as the initial condition
\bqn
\ba{ccll}
\rho(x,0)&=&\r_h+\Delta \r_0 \sin \frac {2\pi x}L, &x\in[0,L], \\
v(x,0)&=&v_{\ast}(\r(x,0),1),&x\in[0,L],
\ea\label{global}
\eqn
with  $\r_h$ = 28 veh/km and $\Delta
\r_0$ = 3 veh/km (the corresponding initial condition
\refe{global} is depicted in \reffig{F_ini1}).

\bfg
\bc\includegraphics[height=12cm] {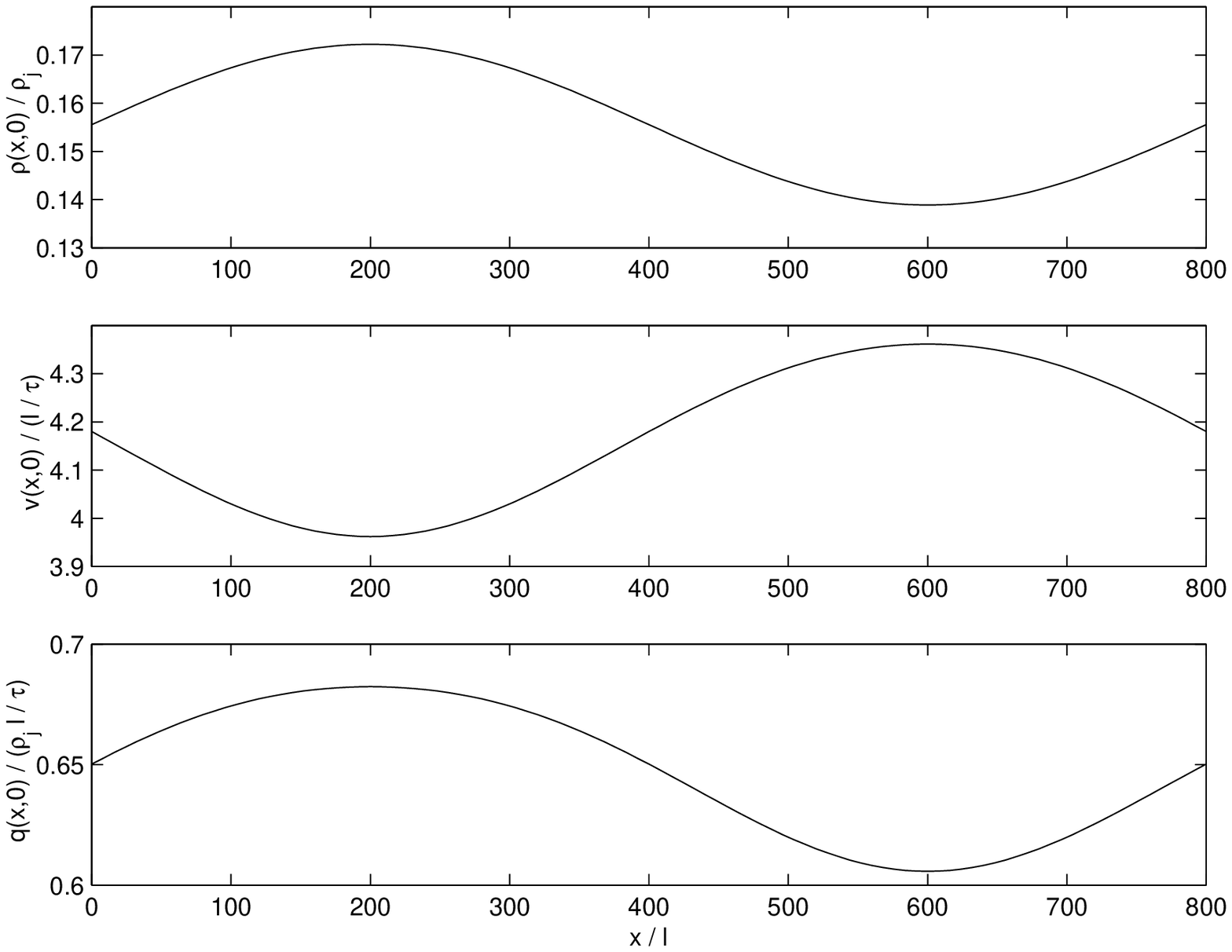}\ec
\caption {Initial condition \refe{global} with $\r_h$ = 28 veh/km and
$\Delta \r_0$ = 3 veh/km} \label {F_ini1}
\efg

The results are shown in \reffig{LWR_K_100}, from which we observe
that initially wave interactions are strong but gradually the bulge
sharpens from behind and expands from front to form a so-called
$N$-wave that travels around the ring with a nearly fixed profile.

\bfg
\bc\includegraphics[height=12cm] {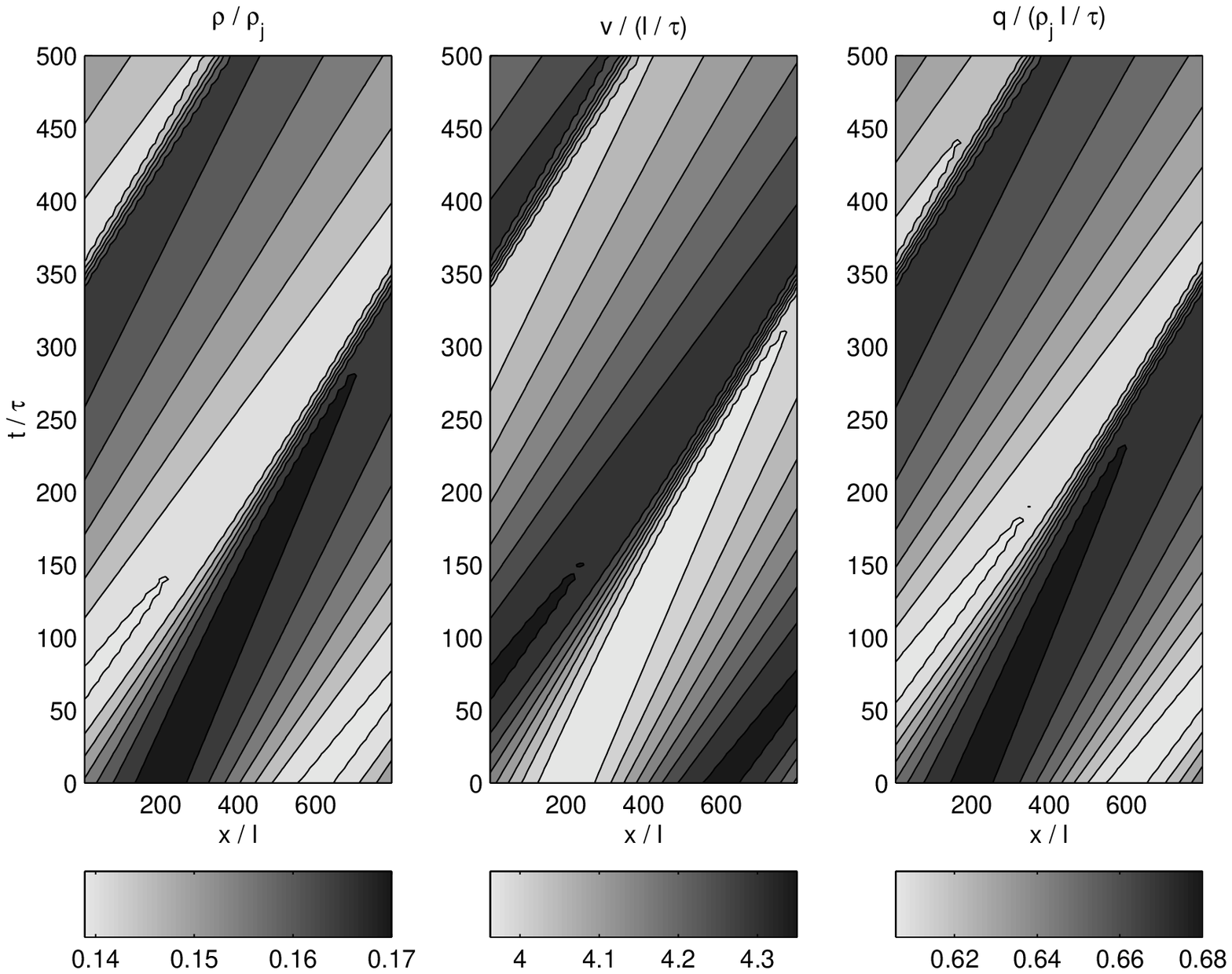}\ec
\caption{Solutions of the homogeneous LWR model with initial
condition in \reffig{F_ini1}}\label{LWR_K_100}
\efg

In the second simulation we created a bottleneck on the ring road
with the following lane configuration:
\bqn
a(x)&=&\cas{{ll} 1, & x\in[320 l, 400 l),\\2, & \m{ elsewhere }.}
\eqn
As before, we also use a global perturbation as the initial
condition
\bqn
\ba{ccll}
\rho(x,0)&=&a(x)(\r_h+\Delta \r_0 \sin \frac {2\pi x}L), &x\in[0,L],
\\
v(x,0)&=&v_{\ast}(\r(x,0),a(x)),&x\in[0,L],
\ea\label{global2}
\eqn
with $\r_h$ = 28 veh/km/lane and
$\Delta \r_0$ = 3 veh/km/lane (the corresponding initial
condition \refe{global2} is depicted in \reffig{F_ini2}).

\bfg
\bc\includegraphics[height=12cm] {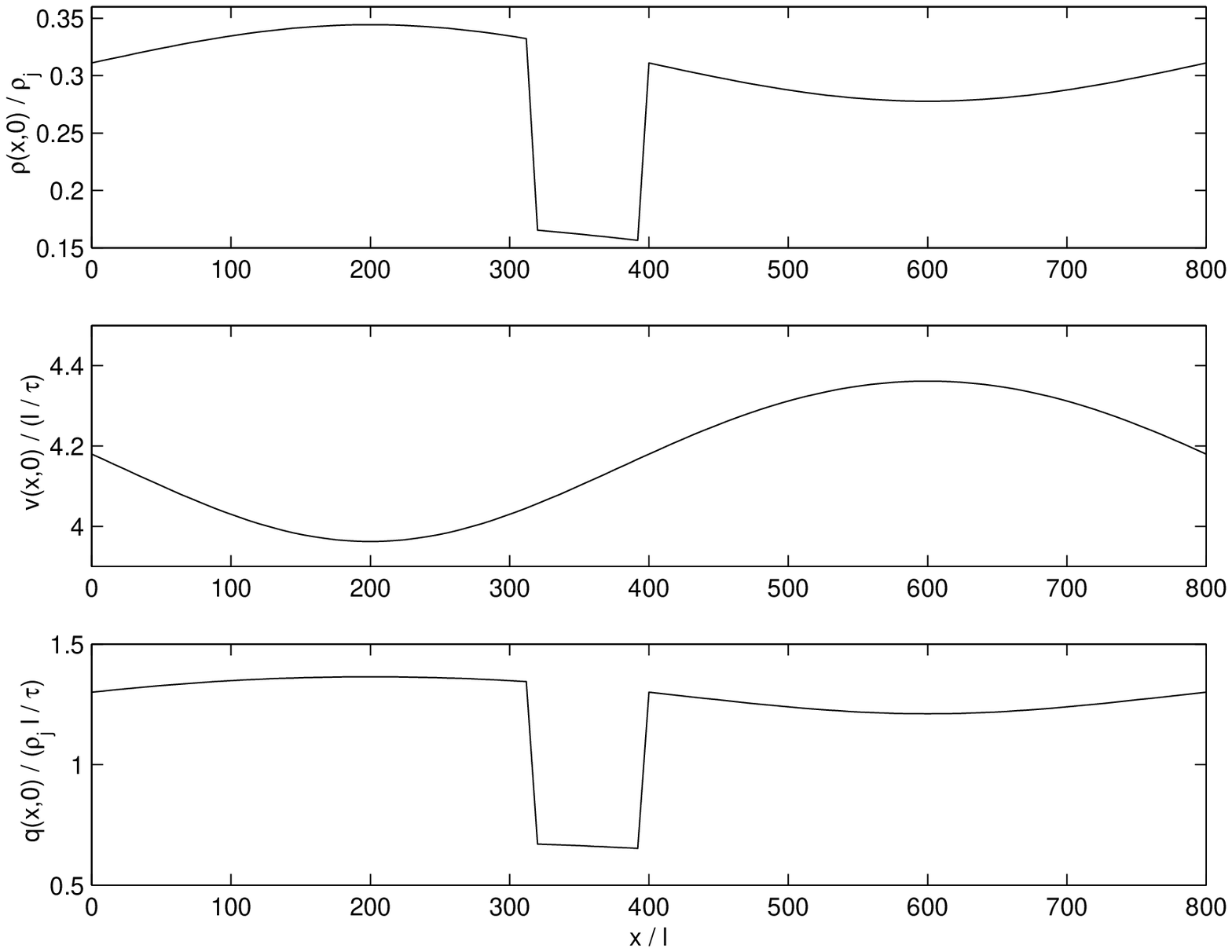}\ec
\caption {Initial condition \refe{global} with $\r_h$ = 28
veh/km/lane and $\Delta \r_0$ = 3 veh/km/lane} \label {F_ini2}
\efg

The results for this simulation are shown in \reffig{inhLWR_K_100}, and
are more interesting. We observe from this figure that at first flow
increases in the bottleneck to make the bottleneck saturated, then a
queue forms upstream of the bottleneck, whose tail propagates
upstream as a shock. At the same time, traffic emerging from the
bottleneck accelerates in an expansion wave. After a while, all the
commotion settles and an equilibrium state is reached, where a
stationary queue forms upstream of the bottleneck, whose in/out
flux equals the capacity of the bottleneck.
Similar situations can be observed in real world bottlenecks, although
queues formed at such bottlenecks rarely reach equilibrium because,
unlike in the ring road example, their traffic demands change over time.
Therefore, we observe queues forming, growing, and dissipating at
locations with
lane drops, upward slopes, or tight turns. Sometimes queues formed at a
bottleneck
can grow fairly long, to the extent that they entrap vehicles that do not
use the bottleneck. Under such situations, we can implement various
types of control strategies, such as ramp metering, to control the extent
of the bottleneck queues so that they do not block vehicles that wish to
exit upstream of the bottleneck. For this purpose the numerical method
presented here can be used to help model and design effective control.

\bfg
\bc\includegraphics[height=12cm] {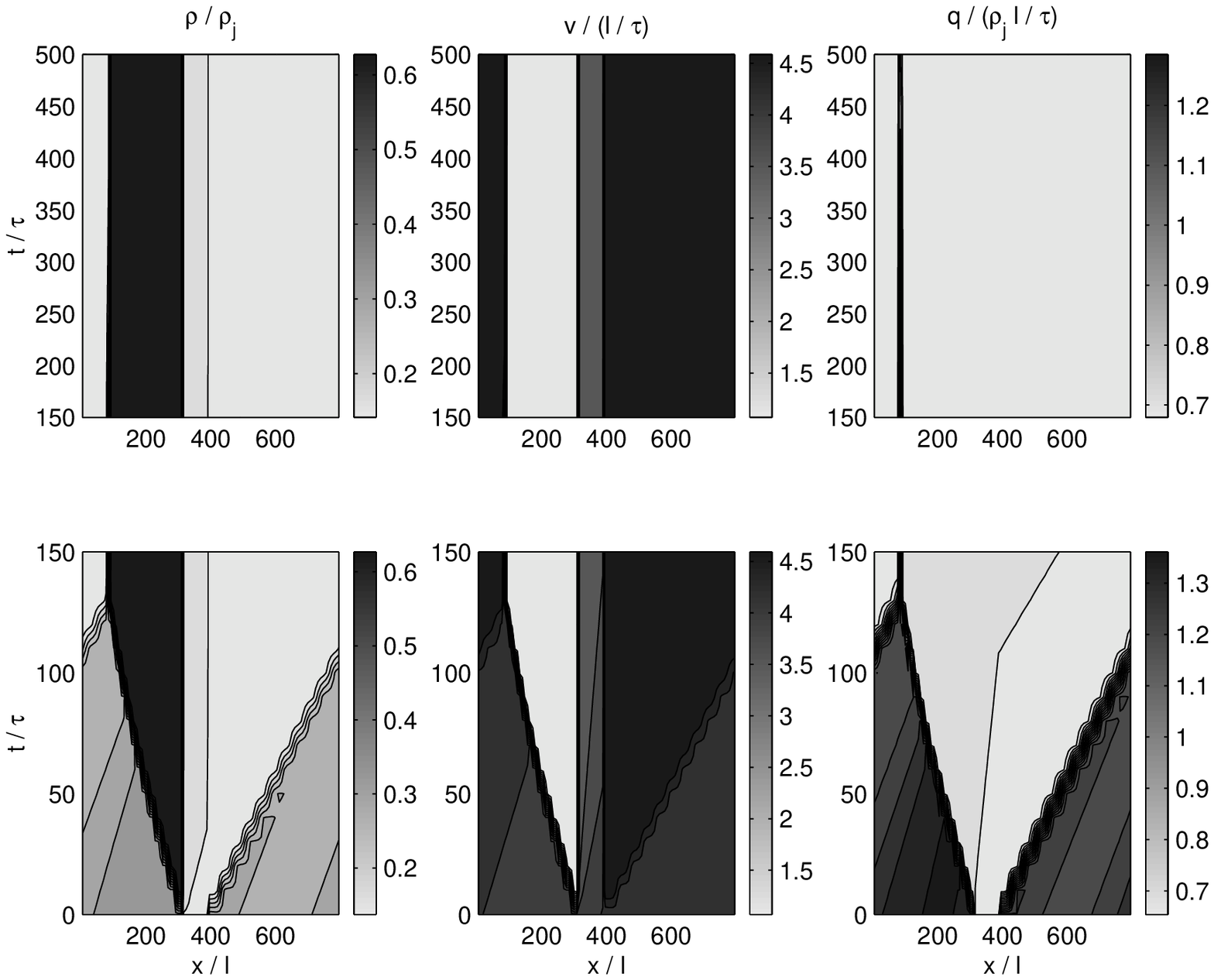}\ec
\caption{Solutions of the inhomogeneous LWR model with initial
condition in \reffig{F_ini2}}\label{inhLWR_K_100}
\efg

%Wenlong Jin: 2003
%last modified: 30Nov01
\section{Concluding remarks}
We studied the inhomogeneous LWR model as a nonlinear resonant
system.  The nonlinear resonance arises when the two
characteristics of the augmented LWR model coalesce.
Critical states and a transitional curve $\Gamma$  can be defined in
the $U$ space based on the behavior of these characteristics, which
are in turn used to solve the Riemann problem for the inhomogeneous
LWR model.  It is found that, under the entropy conditions of Lax and
of Isaacson and Temple, there exist ten types of wave solutions.
Formulas for computing the  boundary fluxes related to different types of
wave solutions were also obtained.  These formulas,  after translated
into the  supply/demand framework, are found to be consistent with
those found in  literature. For problems with general initial/boundary
conditions, the method of Godunov was applied to solve the
inhomogeneous model numerically.

The method presented here can be extended easily to model more
complicated situations, such as multiple inhomogeneities. Suppose at
location $x$, there are $i=1,\cdots,n$ types of inhomogeneities, such as
changes in number of lanes, grade, and curvature. We introduce an
inhomogeneity vector $\vec a(x)=(a_1(x),a_2(x),\cdots,a_n(x))^T$, and
express the flow-density function as $f(\vec a(x),\r)$. Then the
conservation law becomes
\bqs
\cas{{lcl}
\r_t+f(\vec a(x),\r)_x&=&0,\\
\vec a(x)_t&=&0,}
\eqs
and this higher-dimensional nonlinear resonant system can be solved in a
similar way.
 It is worth mentioning
that the augmentation approach taken in this chapter also applies to
higher-order traffic flow models for inhomogeneous roads.

\newpage
\pagestyle{myheadings}
\markright{  \rm \normalsize CHAPTER 3. \hspace{0.5cm}
 MERGING TRAFFIC MODEL}
\chapter{Kinematic wave traffic flow model of merging traffic}
%\thispagestyle{myheadings}
%Wenlong Jin: 2003

%Wenlong Jin: 2003
\section{Introduction}
For developing advanced traffic control strategies, dynamic traffic assignment (DTA)
algorithms, and other technologies in Advanced Traffic
Management Systems (ATMS) and Advanced Traveler Information Systems (ATIS), traffic
engineers need the assistance of network traffic flow models that can capture
system-wide features of traffic dynamics and are computationally efficient for a network
of realistic size. The kinematic wave model is a promising candidate for these tasks since it
provides a realistic description of dynamic traffic phenomena in the aggregate level in
terms of expansion and shock waves and as such is highly efficient for simulating traffic
dynamics in a large network.

In the seminal kinematic wave model by \citet{lighthill1955lwr} and \citet{richards1956lwr}, a.k.a. the LWR model, how a disturbance in traffic propagates through a link was thoroughly studied. To model traffic dynamics on a network with the kinematic wave
model, however, one needs to carefully study traffic dynamics at a merge, a diverge, or
other components of a network. The kinematic wave models of merging traffic have been
studied by \citet{daganzo1995ctm}, \citet{holden1995unidirection}, and \citet{lebacque1996godunov}. In the model by
Holden and Risebro, traffic flows through a merge are determined by an optimization
problem. However, the physical meaning and the objective function of the optimization
problem are not known or supported by observations. On the other hand, the models by
Daganzo and Lebacque are based on the definitions of the local traffic supply and demand
and can be considered as reasonable extensions of the kinematic wave model of link
traffic flow. In this chapter, we will examine the latter models so that they can be better
understood, more easily calibrated, and more efficiently applied in simulation.

As we know, the LWR model, in which the evolution of traffic density $\r(x,t)$, flow-rate
$q(x,t)$, and travel speed $v(x,t)$ is studied in space $x$ and time $t$, can be written as a
partial differential equation based on the fact of traffic conservation and the adoption of
a fundamental diagram. For the purpose of simulation, the LWR model is generally written
in a discrete form: a link and a duration of time are partitioned into a number of cells and
time steps respectively, and the increment of the number of vehicles in a cell at each time
step equals to the difference between the inflow into and outflow from that cell during the
time step. In the discrete LWR model, to solve the flow through a boundary (i.e., the
inflow into the downstream cell and the outflow from the upstream cell), two equivalent
approaches can be used: in the mathematical approach, one solves the Riemann problem at
that boundary \citep{lebacque1996godunov, jin2003inhLWR}; in the engineering
approach, a.k.a. the supply-demand method, the supply of the upstream cell and the
demand of the downstream cell are computed first and the boundary flow is taken as their
minimum. Here the concepts of supply and demand were first introduced by \citet{daganzo1995ctm}, but using the terms of ``sending flow" and ``receiving
flow" instead; the terms of ``supply" and ``demand" were first used by \citet{lebacque1996godunov}.
The definitions of demand and supply are as follows: the demand of a cell is equal to its
flow-rate when the traffic condition is under-critical (i.e., free flow) and its flow capacity
when overcritical (i.e., congested); the supply is equal to the flow capacity of the cell when
the traffic condition is under-critical and the flow-rate when overcritical.\footnote{From
the definitions of demand and supply, we can see that the flow through a boundary is
bounded by the capacity.}

 For computing flows through a merge, including the outflows from the
upstream cells and the inflow into the downstream cell, Daganzo (1995) extended the
supply-demand method as follows: the outflow from an upstream cell is smaller than or
equal to its demand, the inflow to the downstream cell is smaller than or equal to its supply,
and the inflow is equal to the sum of the outflows in order to preserve traffic conservation.
In this supply-demand method, the inflow is unique since it is equal to the minimum of the
supply and total demand. But the outflow from each upstream cell may not be unique. Thus
one has to find a way to distribute to each upstream cell a fraction of the total outflow,
which is equal to the inflow. Here we call such a way of determining the distribution
fractions \textit{the distribution scheme}.

Lebacque proposed another extension of the supply-demand method: the supply of the
downstream cell is first distributed as a virtual supply to each upstream cell, the outflow
from each upstream cell is the minimum of its demand and virtual supply, and the inflow
into the downstream cell is equal to the sum of the outflows. Thus the distribution scheme
in Lebacque's method is used to determine the fractions of virtual supplies, and is more
general since more feasible solutions of flows can be found in this method.

Both Daganzo and Lebacque provided general formulations of the kinematic wave model of merges. Here, we do not intend to extend these formulations. Rather, we are interested in the distribution schemes used in these models since a distribution scheme is the key to uniquely determine flows through a merge.
Since in possible applications of a merge model one wants to obtain unique flows under a
given situation, the distribution schemes are worth a thorough examination.

At a first glance, the determination of distribution fractions seems to be complicated since
they may be affected by travelers' merging behavior, the geometry of the studied merge,
traffic capacities, differences between the upstream cells, traffic conditions, and traffic
control. Considering part of these factors, both Daganzo (1995) and Lebacque
(1996) provided some suggestions on the distribution fractions: Lebacque suggested that
the distribution fraction of an upstream cell is proportional to its number of lanes; Daganzo
considered that upstream cells bear different priorities and hence introduced parameters
for priorities in his distribution fractions. Both suggestions have their limitations:
Lebacque's distribution scheme is very coarse and fails under certain situations, while
Daganzo's scheme becomes very complicated for a merge with three or more upstream
links. Moreover, priorities in Daganzo's distribution scheme vary with flow levels,
which seems to be counter intuitive. Therefore, we devote this study to the better
understanding of  various distribution schemes in a merge model, and  propose a
new distribution scheme which is well-defined, computationally efficient, and capable of
capturing the characteristic differences between different branches of the
merge.

In this chapter, we first review the discrete kinematic wave model of merges and discuss
different formulations of the supply-demand method for computing flows through a merge
(Section \ref{sec:merge2}). In Section \ref{sec:merge3}, after discussing existing distribution schemes, we propose a simple
distribution scheme, which incorporates the ``fairness" condition. In this scheme, the
distribution fractions are proportional to traffic demands of upstream cells. This scheme is
shown to work well in simulations due to its many merits: 1) it is capable of capturing the
characteristic differences between upstream cells (e.g. the speed difference between the
upstream freeway and on-ramps); 2) it is easy to calibrate because additional parameters
such as priorities do not need to be explicitly introduced; and 3) it is
computationally efficient. In Section \ref{sec:merge4}, we present an example of two
merging flows and demonstrate in numerical simulations that the discrete kinematic wave
merge model  incorporating the ``fairness" condition is well-defined and converges in first
order. In the conclusion part, we present the supply-demand method for computing flows
through a diverge and a general junction for single-commodity traffic flow, and discuss
related future research.
 
%Wenlong Jin: 2003
\section{The discrete kinematic wave model of merges with the supply-demand method}\label{sec:merge2}
In the kinematic wave traffic flow model of a road network with a
merge, the LWR model can be used to describe traffic dynamics of
each branch, for which flows through the merge can be considered as
boundary conditions. Thus, in this section, we first review the discrete
LWR model, the definitions of supply and demand, and the
supply-demand method for computing flows through link boundaries.
After reviewing the models of merges under the supply-demand
framework by Daganzo (1995) and Lebacque (1996), we then
demonstrate the importance of distribution schemes. At the end of
this section, we will discuss the properties of existing distribution
schemes.

\subsection{The discrete LWR model in  the supply-demand
framework}\label{godunovmethod} In the LWR model for each branch
of a merge, traffic dynamics are governed by a traffic conservation
equation,
\bqn
\r_t+q_x&=&0, \label{traf_cons}
\eqn
and an equilibrium relationship between $\r$ and $q$, also known as the fundamental diagram,
\bqn
q=Q(a,\r), \label{fund_diag}
\eqn
where $a(x)$ is an inhomogeneity factor, depending on road
characteristics, e.g., the number of lanes at $x$. Since $q=\r v$, we
also have a speed-density relation: $v=V(a,\r)\equiv Q(a,\r)/\r$. For
vehicular traffic, generally, $v$ is non-increasing and $q$ is concave in
$\r$. Examples of empirical models of speed- and flow-density
relations can be found in \citep{newell1993sim,kerner1994cluster}. Related to the
fundamental diagram, the following definitions are used in this chapter:
the maximum flow-rate at $x$ is called the traffic capacity, and the
corresponding density is called the critical density; traffic flow is
overcritical when its density is higher than the critical density, and
under-critical conversely.

From \refe{traf_cons} and \refe{fund_diag}, the LWR model can be
written as
\bqn
\r_t+Q(a,\r)_x&=&0, \label{LWR}
\eqn
where $0\leq\r\leq\r_j$ ($\r_j$ is the jam density).
When $a(x)$ is uniform with respect to location $x$, the LWR model is called homogeneous. Otherwise it is called inhomogeneous. Both the homogeneous and inhomogeneous models are hyperbolic systems of conservation laws. Actually the former, which is a strict hyperbolic conservation law \citep{lax1972shock}, is a special case of the latter, a non-strictly hyperbolic system of conservation laws and a resonant nonlinear system \citep{isaacson1992resonance}. Therefore, the following discussions for the inhomogeneous LWR model are valid for any kind of links.

With jump initial conditions, the LWR model \refe{LWR} is solved by shock waves,
expansion waves, and standing waves. These wave solutions are
unique under the so-called ``entropy" conditions. However, solutions
of the LWR model with general initial and boundary conditions can not
be expressed in analytical form, which calls for approximate solutions
with numerical methods. One efficient numerical method for solving
\refe{LWR} is due to \citet{godunov1959}. In the
Godunov method, the link is partitioned into $N$ cells, a duration of
time is discretized into $M$ time steps, and the discretization of
space and time satisfies the Courant-Friedrichs-Lewy \citep{courant1928CFL} (CFL) condition so that a vehicle is not
allowed to cross a cell during a time interval. Assuming that the
spacing $\dx$ and the time step $\dt$ are constant, $\r_i^m$ is the
average of $\r$ in the cell $i$ at time step $m$, $q^{m+1/2}_{i-1/2}$
and $q^{m+1/2}_{i+1/2}$ are the inflow into and the outflow from cell
$i$ from time step $m$ to $m+1$ respectively, the LWR model
\refe{LWR} for cell $i$ can be approximated with a finite difference
equation:
\bqn
\frac {\r_i^{m+1}-\r_i^m} {\dt}+\frac
{q^{m+1/2}_{i-1/2}-q^{m+1/2}_{i+1/2}}{\dx}&=&0. \label{dis.1m}
\eqn

In \refe{dis.1m}, the flow through the link boundary $x_{i-1/2}$, i.e., $q^{m+1/2}_{i-1/2}$, can be computed in two approaches.  One is from the wave solutions of the Riemann problem for \refe{LWR} with the following initial conditions \citep{jin2003inhLWR}:
\bqn
U(x=x_{i-1/2}, t=t_{m}) &=&\cas{{ll} U^m_{i-1} & x <x_{i-1/2}\\U^m_i
& x>x_{i-1/2}},
\eqn
where $U=(a,\r)$. Another is the supply-demand method \citep{daganzo1995ctm,lebacque1996godunov}, in which the flow through a link boundary is the minimum of the traffic demand of its upstream cell and the traffic supply of its downstream cell. The two approaches were shown to be equivalent \citep{jin2003inhLWR}. However, the method of solving the Riemann problem and the supply-demand method have different fates for studying traffic dynamics through a merge: there has been no formulation of the Riemann problem for merging dynamics in literature, but the supply-demand method has been extended and applied in the discrete kinematic wave models of merges.

In the following, we describe in detail the supply-demand method. Considering the link boundary at $x_{i-1/2}$, whose upstream and downstream cells are respectively denoted by cell $i-1$ and cell $i$, supposing that the traffic densities of the two cells are $\r_{i-1}^m$ and $\r_i^m$ at time step $m$, Daganzo (1995) and Lebacque (1996) suggested the following supply-demand method for computing $q^{m+1/2}_{i-1/2}$. First, traffic demand of the cell $i-1$ (called ``sending flow" by Daganzo), $D^{m+1/2}_{i-1}$ , and traffic supply of the cell $i$ (called ``receiving flow" by Daganzo), $S^{m+1/2}_i$, are defined by
\bqn
D^{m+1/2}_{i-1}&=&\cas{{ll} Q^m_{i-1}, & \m{cell } i-1 \m{ is under-critical}, \\Q^{max}_{i-1}, & \m{otherwise};}\label{def_demand}\\
S^{m+1/2}_i&=&\cas{{ll} Q^{max}_i, & \m{cell } i \m{ is under-critical},\\Q^m_i,
&\m{otherwise};}
\eqn
where $Q^{max}_i$ is the capacity of cell $i$, and $Q^m_i$ the flow-rate of cell $i$ at time step $m$. The demand can be considered as the maximum flow that can be discharged by the cell $i-1$ from time step $m$ to $m+1$; the supply $S^{m+1/2}_i$ is the maximum flow that can be received by the cell $i$. Thus, the boundary flow satisfies (all superscripts will be suppressed hereafter)
\bqn
\ba{ccc}
q_{i-1/2} &\leq &D_{i-1},\\
q_{i-1/2} &\leq &S_i.
\ea\label{constraint}
\eqn
Note that \refe{constraint} admits multiple solutions. To identify the
unique boundary flow, an additional ``optimality" condition, that the
actual boundary flow always reaches its maximum, is
assumed. Hence, the
boundary flow can be simply computed by
\bqn
q_{i-1/2}&=&\min\{D_{i-1}, S_i\}. \label{d-s}
\eqn
Here the ``optimality" condition can be considered as an entropy condition, which helps to choose a physical solution out of all feasible solutions.

\subsection{The kinematic wave model of merging traffic in the supply-demand framework}
In this subsection, we review the kinematic wave model of merging traffic in the supply-demand
framework. In this type of models, the supply-demand method is used
to compute flows through a merge.  Without loss of generality, we
consider a merge that connects two upstream cells to one down
stream cell. Furthermore, we assume that, at time step $m$, traffic
demands of the two upstream cells and the traffic supply of the
downstream cell are $D_1$, $D_2$, and $S_d$ respectively. We denote
the outflows from the upstream cells by $q_1$ and $q_2$ and the
inflow into the downstream cell by $q$ from time step $m$ to $m+1$.
Then, according to traffic conservation, we have $q=q_1+q_2$.

The basic assumption in the supply-demand method for computing the flows through a merge is that the flows, $q_1$, $q_2$, and $q$, are determined by traffic conditions $D_1$, $D_2$, $S_d$, and/or other characteristics of the merge.
Another assumption, as in the supply-demand method for computing
the flow through a link boundary, is the optimality condition. Two
types of optimality conditions have been proposed: one is that the
total flow $q$ reaches its maximum, and the other is that both $q_1$
and $q_2$ reaches their individual maximums. Following the first
assumption leads to Daganzo's merge model (1995), and following the
second leads to Lebacque's (1996).

In Daganzo's supply-demand method, we have the following optimization problem:
\bqn
\ba{lcl}
\max q&=&q_1+q_2\\
\m{s.t.}&&\\&&\ba{lcl}
q_1&\leq&D_1,\\
q_2&\leq&D_2,\\
q_1+q_2&\leq&S_d,\\
q_1,q_2&\geq&0,
\ea
\ea\label{optimality}
\eqn
from which we can find the total flow,
\bqs
q&=&\min\{D_1+D_2,S_d\}.
\eqs
However, $(q_1, q_2)$ may have multiple feasible solutions. This can be shown with \reffig{daganzo_form}: when $S_d\geq D_1+D_2$, the solution is unique and at point $Q$; i.e., $(q_1,q_2)=(D_1,D_2)$; but when $S_d<D_1+D_2$, the solution can be any point on the line segment $AB$. For the latter situation, Daganzo defined two (non-negative) distribution fractions $\alpha_1$ and $\alpha_2$, which satisfy $\alpha_1+\alpha_2=1$ and may be related to $D_1$, $D_2$, $S_d$, and other characteristics of the merge. Then, the total flow $q$ is distributed by $q_i=\alpha_i q$ ($i=1,2$). One example when $S_d<D_1+D_2$ is depicted in the figure, with given fractions $\alpha_1$ and $\alpha_2$. \reffig{daganzo_form} also shows that $\alpha_1$ or $\alpha_2$ are restricted by $D_1$, $D_2$, and $S_d$. For instance, for $S_d$ given in the figure, $\alpha_1$ can not be 1.

\bfg
\bc\includegraphics[height=12cm] {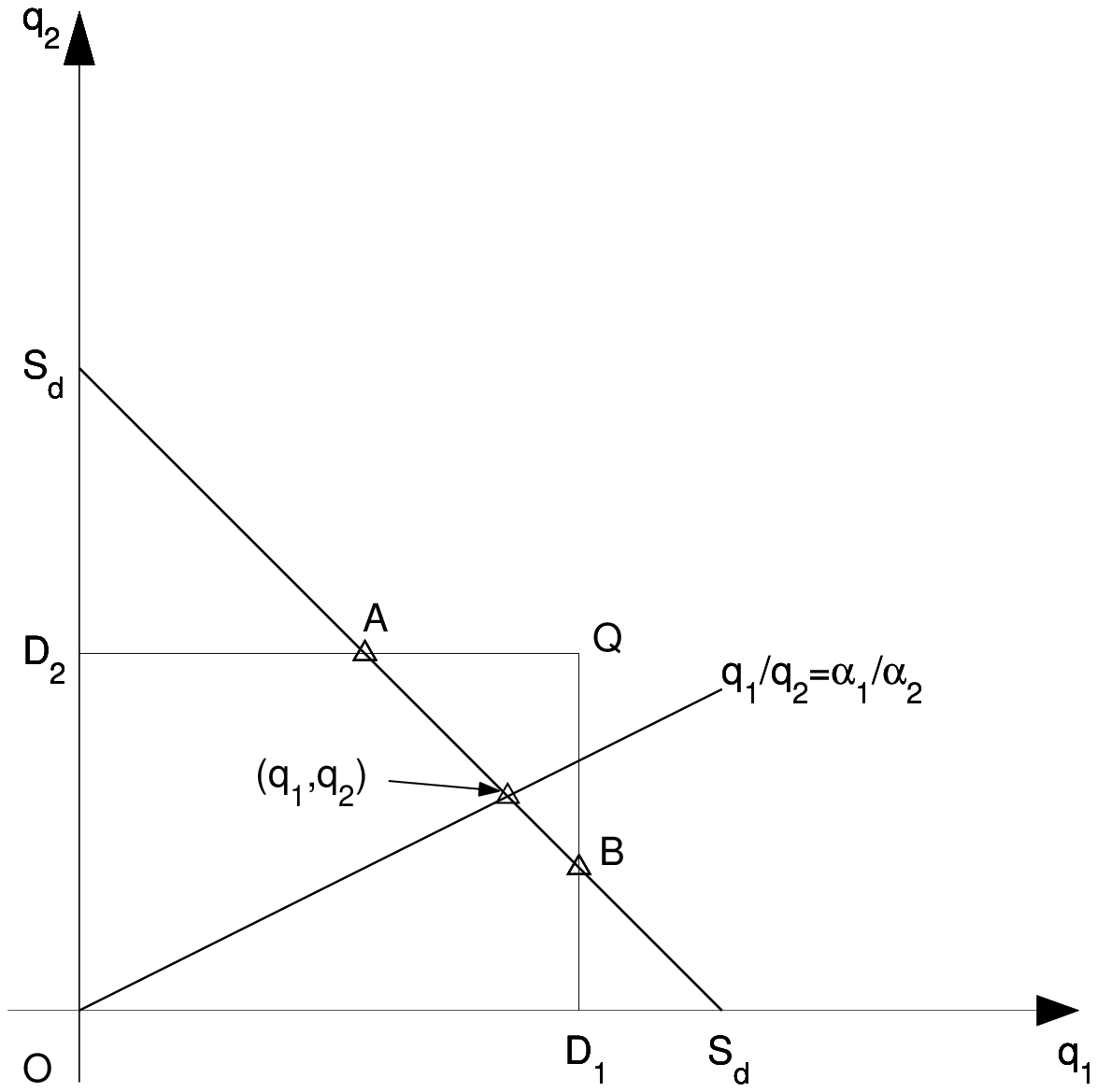}\ec
\caption{Feasible solutions in Daganzo's supply-demand method}\label{daganzo_form}
\efg

Lebacque suggested another supply-demand method: the supply of the downstream cell is first distributed to the two upstream cells with two fractions $\alpha_1$ and $\alpha_2$, and it is assumed that the flows $q_1$ and $q_2$ reaches their individual maximums.
i.e., we can compute the flow $q_i$ ($i=1,2$) as the following:
\bqn
\ba{lcl}
S_i&=&\alpha_i S_d,\\
q_i&=&\min\{D_i,S_i\}.
\ea\label{lebacque}
\eqn
The feasible solutions of Lebacque's method without fixed fractions are shown in \reffig{lebacque_form}. As shown, in this model, when $S_d\geq D_1+D_2$, $(q_1,q_2)$ can be any point on $D_1BQAD_2$; when $S_d<D_1+D_2$, $(q_1,q_2)$ can be any point on $D_1BAD_2$. In Lebacque's formulation, therefore, $\alpha_1$ and $\alpha_2$ are not restricted by $D_1$, $D_2$, or $S_d$, and the total flow $q$ may not reach its maximum $\min\{D_1+D_2,S_d\}$ in this method.

\bfg
\bc\includegraphics[height=12cm] {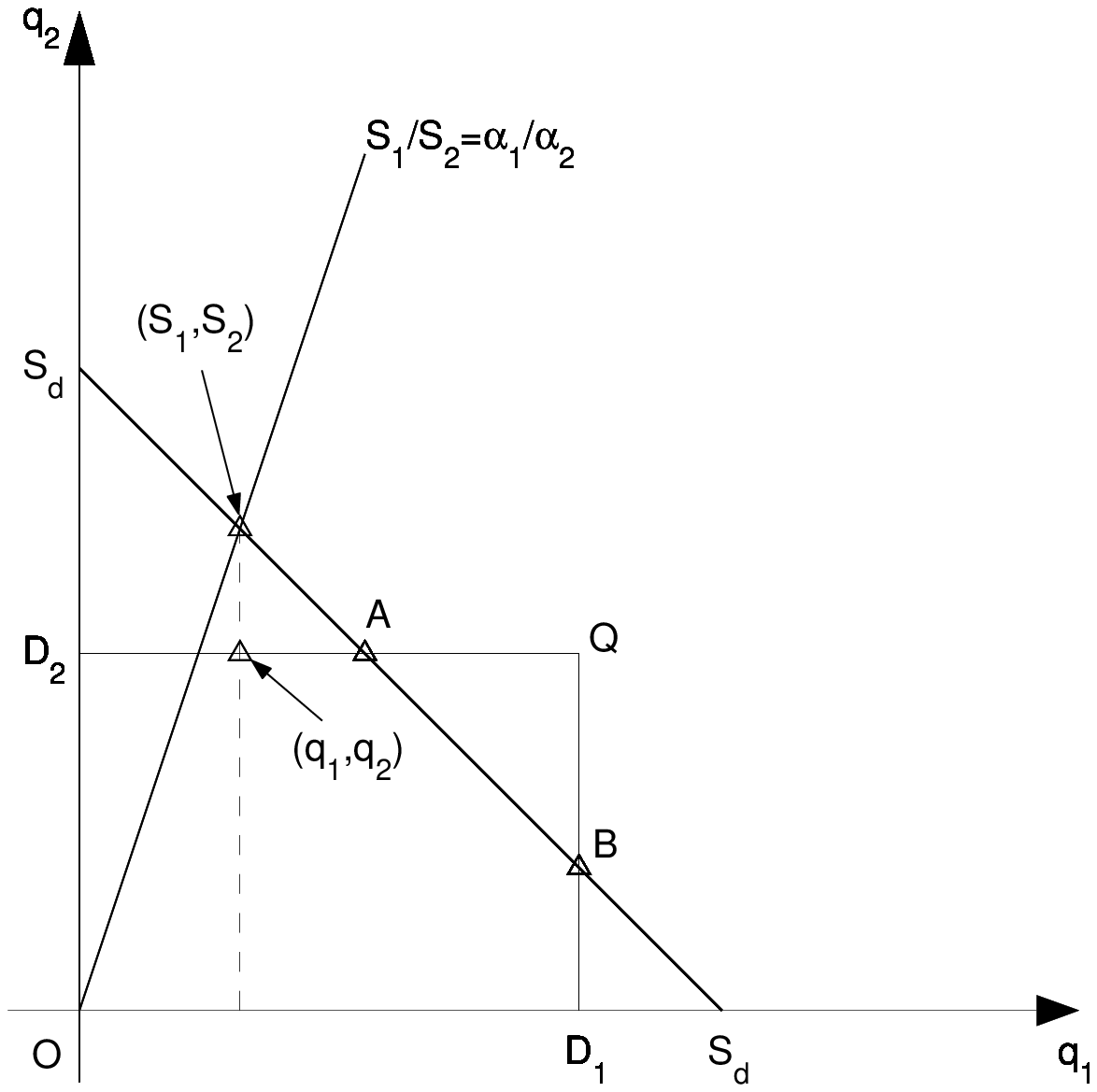}\ec
\caption{Feasible solutions in Lebacque's supply-demand method}\label{lebacque_form}
\efg

Comparing Daganzo's and Lebacque's methods, we can see that: 1)
when the fractions are the same, the two methods give the same
flows; 2) for given $D_1$, $D_2$, and $S_d$, the feasible solution
domain of Daganzo's method is contained by that of Lebacque's since
the distribution fractions in Daganzo's method (but not in Lebacque's)
are confined by the supplies and the demand.

From the above analysis, we can see that both Daganzo's and
Lebacque's models in the supply-demand framework are based on
reasonable assumptions, and Lebacque's method \refe{lebacque}
yields a larger set of feasible solutions than Daganzo's. In addition,
we think that both formulations are clear and general enough to
contain physical solutions. Thus, in this chapter, we do not intend to
investigate further the formulations. Instead, we are interested
in the distribution schemes used in these models.

The reasons why distribution schemes are worth further, deeper
discussions are as follows. First, we can see from
\reffig{daganzo_form} and \reffig{lebacque_form} that distribution
schemes play a key role in uniquely determining flows through a
merge. Thus whether solutions of flows are physical is highly
dependent on the distribution scheme used. Therefore, in order to
apply these models to simulate traffic dynamics at a merge, we need
a better understanding of their distribution schemes. Second, the
distribution fractions can be affected by travelers' merging
behaviors, the geometry of a merge, differences between the
upstream cells, traffic conditions, and possible control strategies
imposed on an on-ramp. On the surface, a distribution scheme that
models all these factors may be extremely complicated. A closer look
at various distribution schemes is needed to find a simple yet
physically meaningful one. Third, it is possible that many valid
distribution schemes are available. When this happens, a distribution
scheme that is easy to calibrate and computationally efficient is
always preferred.

\section{Investigation of various distribution schemes}
In this section, we take a closer look at various distribution schemes
and see how their distribution fractions are affected by traffic
conditions, i.e., $D_1$, $D_2$, and $S_d$, and other characteristics of a
merge. We start with a review and a discussion on the existing
distribution schemes of Daganzo and Lebacque, then propose a simple
distribution scheme and demonstrate that the supply-demand method
incorporating this scheme is capable of addressing all factors that we
concern about.

\subsection{Discussion of existing distribution schemes}
As we know, different types of links have different characteristics.  As
a result, even when an upstream highway and an on-ramp have the
same number of lanes and traffic density, the downstream link
usually receives more vehicles from the upstream highway than from
the on-ramp due to differences  in design speeds and geometry. For
example, when vehicles queues up on both a $L$-lane highway and
$1$-lane on-ramp that merge together, the ratio of flow from the
on-ramp to that from the highway is about $1/(2L-1)$
\citep{daganzo1996gridlock}. From these observations,
\citet{daganzo1995ctm} suggested that different upstream links bear
different priorities and proposed a distribution scheme including
parameters for priorities.

\reffig{daganzo_ds} shows how Daganzo's distribution scheme is
defined. In the figure, the priorities of the highway and the on-ramp
are denoted as $p_1$ and $p_2$ ($p_1+p_2=1$), respectively. Here the
upstream link $u_1$ is assumed to have higher priority than $u_2$;
i.e., $p_1/p_2>D_1/D_2$. Then the solution $(q_1,q_2)$ can be shown
to be one of three cases: i) when $S_d\leq D_1/p_1$; i.e., $S_d$ is
$x$- and $y$-intercept of line i, the solution $(q_1,q_2)=(p_1S_d,
p_2S_d)$ is at point $1$; ii) when $S_d\in(D_1/p_1, D_1+D_2)$; i.e.,
$S_d$ is the $x$- and $y$- intercept of line ii, the solution
$(q_1,q_2)=(D_1,S_d-D_1)$ is at point $2$; iii) when $S_d\geq
D_1+D_2$; i.e., $S_d$ is the $x$- and
$y$- intercept of line iii, the solution $(q_1,q_2)=(D_1,D_2)$ is at point
$3$.

\bfg
\bc\includegraphics[height=12cm] {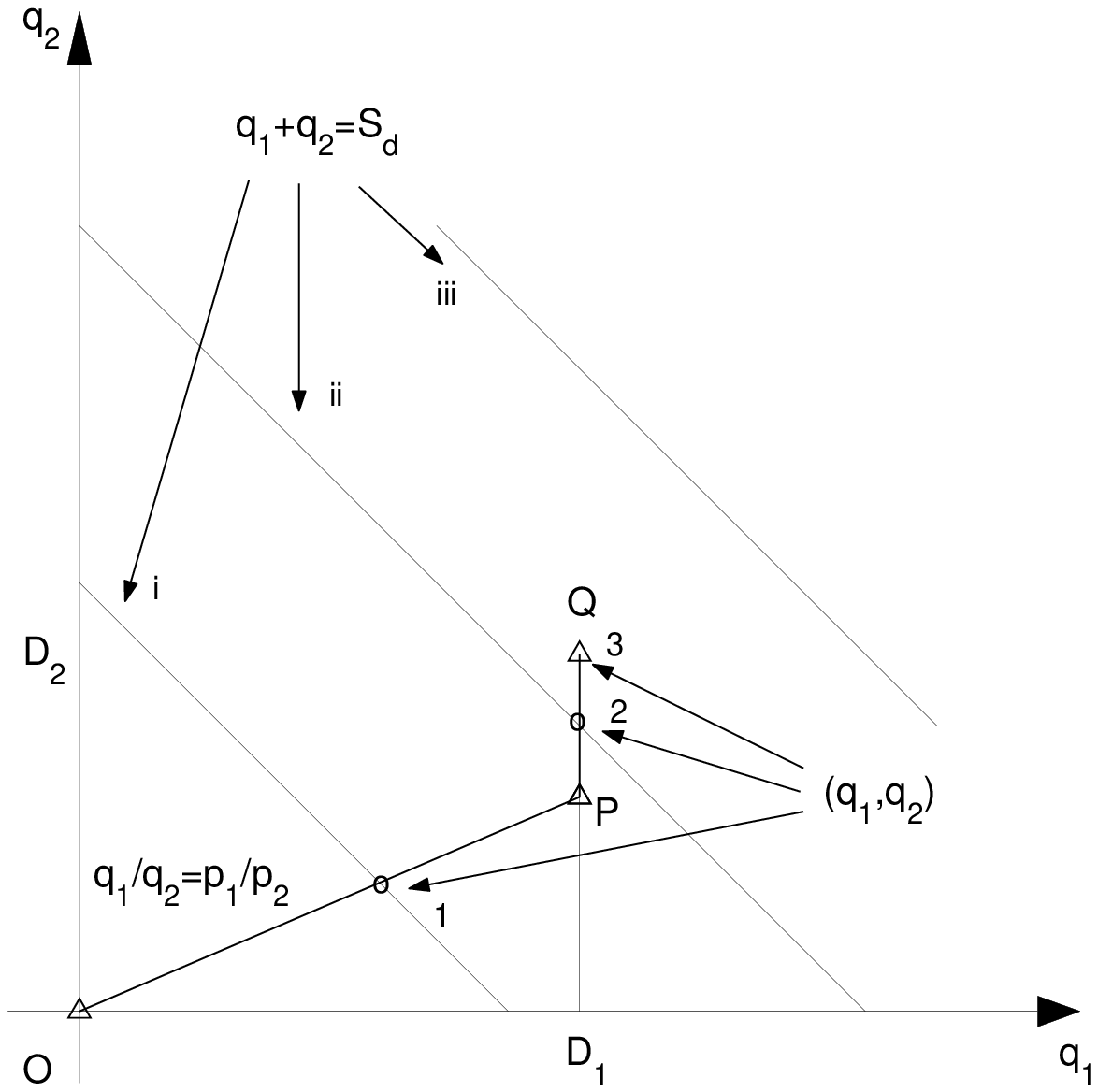}\ec
\caption{Solutions of flows in Daganzo's distribution scheme}\label{daganzo_ds}
\efg

Thus in Daganzo's distribution scheme, we can find that the  fraction
$\alpha_1$ is defined as
\bqn
\alpha_1&=&\cas{{ll} p_1,& S_d\leq D_1/p_1,\\\frac{D_1}{S_d},&
D_1/p_1<S_d\leq
D_1+D_2,\\\frac{D_1}{D_1+D_2},&S_d>D_1+D_2;}\label{d_frac}
\eqn
and $\alpha_2=1-\alpha_1$. From this definition and \reffig{daganzo_ds}, we can see that the priorities $p_1$ and $p_2$ have to satisfy $p_2/p_1 < D_2/D_1$.
i.e., they have to change with respect to traffic conditions $D_2$ and
$D_1$. Therefore, such ``priorities" are not uniquely determined by
road characteristics, as one would expect. This property makes this
distribution scheme less attractive. Moreover, even we allow
non-constant priorities, the distribution scheme with fractions
defined in \refe{d_frac} becomes quite complicated when considering
a merge with more than two upstream links.

Lebacque (1996) suggested another distribution scheme, in which
$\alpha_i$ equals to the ratio of the number of lanes of link $u_i$ to
that of link $d$. When all branches of a merge are highways with the
same characteristics, and the traffic conditions of upstream  links are
overcritical, the demand of each upstream link is equal to its
lane capacity  times the number of lanes. In this case, it is
reasonable that outflow from each upstream link is proportional to
the number of lanes; i.e., the distribution scheme by Lebacque works
well. However, when the upstream links are not similar, e.g., one is
highway and the other on-ramp, the fractions are obviously not
proportional to the number of lanes. Lebacque's distribution scheme
fails in another case when
$\alpha_1+\alpha_2>1$ and $S_i\leq D_i$, because it may yield invalid
solutions of $q=q_1+q_2>S_d$.

%Wenlong Jin: 2003
\subsection{A simple distribution scheme and its interpretation}\label{sec:merge3}
Our above analysis has revealed certain drawbacks of the existing distribution schemes, here we propose a simple distribution scheme, which, as we will see later, removes these drawbacks yet is capable of capturing characteristics of a merge.
In this distribution scheme, the distribution fractions are only related to the demands $D_1$ and $D_2$, as defined in \refe{fraction}.
\bqn
\ba{lcl}
\alpha_1&=&\frac{D_1}{D_1+D_2},\\
\alpha_2&=&\frac{D_2}{D_1+D_2}.
\ea\label{fraction}
\eqn
As shown in \reffig{simple_ds}, combining the distribution scheme \refe{fraction} with models \refe{optimality} or \refe{lebacque}, we are able to uniquely determine the flows: the solution $(q_1,q_2)$ with this distribution scheme is simply the intersection of $q_1+q_2=S_d$ and $q_1/q_2=D_1/D_2$ when $D_1+D_2\geq S_d$, and the point $Q$ otherwise.

\bfg
\bc\includegraphics[height=12cm] {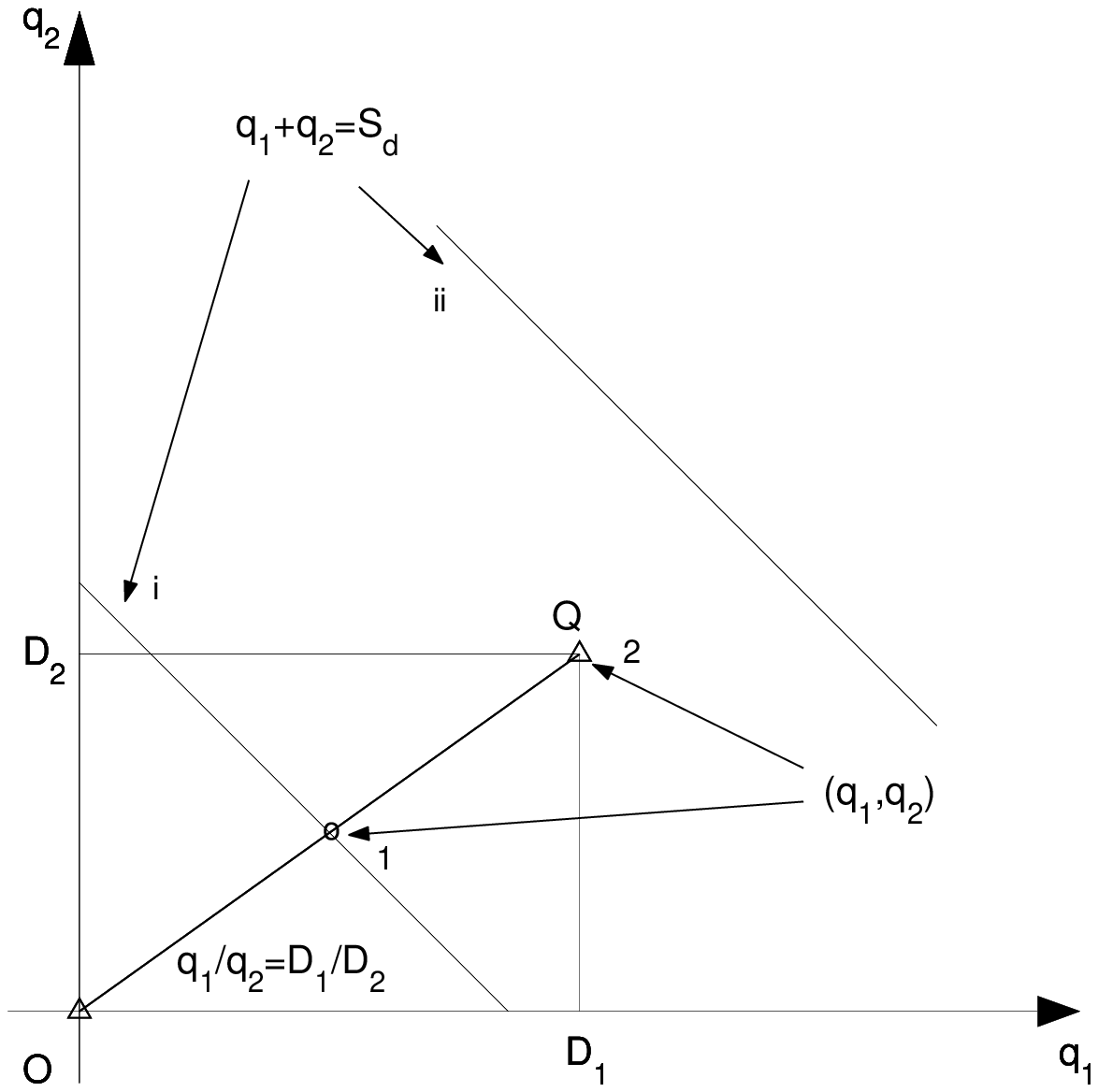}\ec
\caption{Solutions of flows in the simple distribution scheme}\label{simple_ds}
\efg

A distribution scheme is in fact equivalent to an additional entropy condition, which helps identify $q_1$ and $q_2$. Thus we also call the distribution scheme \refe{fraction}
 the ``fairness" condition, because the distribution fractions are
proportional to traffic demands of upstream links; i.e., the upstream
cell with more ``sending" flow is given more chances. This ``fairness"
condition, to some extent, is supported by observations at crowded
merges, e.g., vehicles from an on-ramp generally wait until there is a
big enough gap to merge when traffic is fluid. When many vehicles
from the on-ramp cannot merge and queue up, they may squeeze in
and force vehicles from the upstream mainline freeway to slow down
or switch lanes to give way to them \citep{kita1999merging}.
These observations show that vehicles from the upstream cells
compete ``fairly" with each other for admission into the downstream
cell.

From \refe{fraction}, we can see that the distribution scheme is not
\textit{directly} related to capacities, number of lanes, the difference
between upstream links, or control of on-ramps. Thus we can say that
this distribution scheme uses the fewest parameters. Therefore, it
will be easy to calibrate and efficient to compute. Indeed, it is the
simplest distribution scheme that we can have. Besides,  note that
the fractions  are independent of the downstream traffic supply
$S_d$, although the flows are related to $S_d$. As a mathematical
exercise, the following theorem shows the fractions in this scheme
are in fact the only fractions that are independent of $S_d$.
\begin{theorem}\label{unique}
Suppose that the fractions $\alpha_1$ and $\alpha_2$ are independent of the downstream supply $S_d$, then the fractions will be as in \refe{fraction}.
\end{theorem}
{\em Proof.} We have $q_1=\alpha_1 q$ and $q_2=\alpha_2 q$. From \refe{lebacque}, we then obtain
\bqn
\alpha_1 q&\leq&D_1,\label{e1}\\
\alpha_2 q&\leq&D_2.\label{e2}
\eqn
Since $\alpha_1$ and $\alpha_2$ are independent of $S_d$, set $S_d=D_1+D_2$, we obtain
$q=D_1+D_2$. Thus, both \refe{e1} and \refe{e2} have to take the ``=" sign, and we have \refe{fraction}. \eop

On the other hand, characteristics of a merge can be
\textit{indirectly} captured in the simplest distribution scheme
\refe{fraction} since capacities, number of lanes, design speeds, and
on-ramp control can be included in the definition of traffic
demands,
\refe{def_demand}.  For example, when upstream links have the same
per lane capacity and are congested,
\refe{fraction} will give fractions proportional to the number of lanes.
Thus this distribution scheme coincides with Lebacque's scheme based
on the number of lanes.  As to Daganzo's  ``priorities",   they are
embedded in the simple distribution scheme as follows:  when the
freeway and the on-ramp have the same number of lanes and density,
the freeway generally admits higher free flow speed, has higher
flow-rate and higher demand, and therefore has larger outflow that
reflects its higher priority. In addition, the resultant supply-demand
method can be applied to determine flows through a merge with a
controlled on-ramp
\citep{daganzo1995ctm}: when the metering rate of the on-ramp,
whose real traffic demand is $D_2$, is $r$, we can apply the controlled
traffic demand of the on-ramp $\min\{r,D_2\}$ in the supply-demand
method. Therefore, although the distribution scheme plays ``fairly",
the resultant supply-demand method of merges can address the
characteristics of a merge by incorporating them into the
computation of traffic demands and traffic supply.

From discussions above we can see that, in this simple distribution
scheme, characteristics of upstream links, the control of  on-ramps,
and other properties of a merge are captured in the definitions of
demand and supply. This is why distribution fractions
depends only on demands in this scheme. Although demand functions
of upstream cells are related to many factors and may not be easily
obtained, they have to be found in all supply-demand methods. In this
sense, the simple distribution scheme and, therefore, the
supply-demand method with this scheme, are the easiest to calibrate
and the most computationally efficient.

\subsection{Properties of the discrete kinematic wave model of merges with the simplest distribution scheme}
With the distribution scheme \refe{fraction}, the supply-demand method has the following further properties:
\begin{description}
\item [Equivalence of models by Daganzo and Lebacque:]
With the distribution scheme \refe{fraction}, as shown in
\reffig{simple_ds}, the solution $(q_1,q_2)$ will be on the line segment
$OQ$. We can see that Daganzo's model \refe{optimality} and
Lebacque's model\refe{lebacque} are equivalent with these
fractions. Both yield the following fluxes through the merge:
\bqn
\ba{lcl}
q&=&\min\{D_1+D_2,S_d\},\\
q_1&=&q\cdot\frac{D_1}{D_1+D_2},\\
q_2&=&q\cdot\frac{D_2}{D_1+D_2}.
\ea\label{newmodel}
\eqn

\item [Extensibility:]
The supply-demand method incorporating the simplest distribution scheme produces qualitatively similar results for merges with different number of upstream links.
When a merge has $U>2$ upstream links, the method \refe{newmodel} can be easily extended as
\bqn
\ba{lcl}
q&=&\min\{\sum_{i=1}^U D_i,S_d\},\\
q_i&=&q \frac{D_i}{\sum_{i=1}^U D_i}, \quad i=1,\cdots,U.
\ea\label{multi}
\eqn

\item [Convergence of the merge model:]
The discrete LWR model is considered as a good approximation to the continuous LWR model since it converges as $\dx\to0$ while $\dx/\dt$ is
constant.  Although analytical convergence analysis of the
merge model with the simplest distribution scheme has not yet been
performed, numerical results in section \ref{sec:merge2} do show that
it is convergent in the $L_1$ norm.

\item [Consistency of the merge model with the LWR model:]
Here we conceptually consider the consistency of the merge model
\refe{newmodel} with the LWR model for a link with multiple lanes. In
the LWR model for a multi-lane link, all lanes are assumed to be
identical; i.e., given the same initial and boundary conditions for each
lane, flows at the same location on each lane are identical and the
link's flow-rate or density at a location is simply the number of lanes
times the flow-rate or density at the location of each lane,
respectively.  i.e., the LWR model does not model lane-changing and
traffic in different lanes is treated as the same. Therefore a
multi-lane link can be considered as an artificial merge: for a
boundary inside the link, we separate its upstream part into two
links with identical flow characteristics, while keep the downstream
part intact. Next we check if traffic dynamics of this
artificial merge is indeed the same as those of the original link.
Assuming the two upstream links of the artificial merge have $N$ and
$M$ lanes, respectively, traffic demand of each lane is $D$, and traffic
supply of the downstream link is $S_d$. Since the lanes are identical in
the upstream links, traffic demands for the upstream  links are
$D_1=N D$,
$D_2=M D$; from \refe{newmodel}, we have
\bqs
q&=&\min\{ (N+M)D, S_d\},\\
q_1&=&q \frac{N}{N+M},\\
q_2&=&q\frac{M}{N+M}.
\eqs
Hence, as expected, flow from each upstream lane is  $\min\{ D,
S_d/(N+M)\}$, which is the same as the original boundary flow
computed from the LWR model.

\end{description}

The above analyses indicate that the  merge
model with the simplest distribution scheme is well-defined and
qualitatively sound, although the ultimate test of its validity rests on
empirical validation.

%Wenlong Jin: 2003
\section{Numerical simulations}\label{sec:merge4}
In this section, we present our numerical studies of the discrete kinematic wave model of merges using the simple distribution scheme. Here we apply Godunov's method discussed in Subsection \ref{godunovmethod} for each link, and the supply-demand method is used to find flows through link boundaries and merging boundaries. In particular, the simple distribution scheme is used for computing fluxes through the merge. The resulted numerical solution method is described as follows: in each cell, \refe{dis.1m} is used to update traffic density; we compute flows through link boundaries with \refe{d-s}; flows through a merge are computed from \refe{multi}.

In the numerical studies, we introduce a unit time $\tau$ = 5 s and  a
unit length $l$ = 0.028 km. Here we study a merge formed by a
two-lane mainline freeway and a one-lane on-ramp. The three
branches of the merge have the same length, $L=400l=11.2$ km.   The
upstream mainline freeway, the on-ramp, and the downstream
mainline freeway are labeled as links $u_1$, $u_2$, and
$d$, respectively. The simulation starts from $t=0$ and ends at
$t=500\tau=41.7$ min. In the following numerical simulations, we
partition each link into $N$ cells and the time interval into $K$ steps,
with $N/K=1/10$ always; e.g., if $N=50$  and $K=500$, the cell length
is $\dx=8 l$ and the length of each time step $\dt=1 \tau$.

For both the mainline freeway and the on-ramp, we use the triangular fundamental diagram; i.e., the flow-density relationships are
triangular. For the mainline freeway, the free flow speed is
$v_{f,m}=$65 mph=5.1877 $l/\tau$; the jam density is
$\r_{j,m}=2\r_j$=360 veh/km, where $\r_{j}$=180 veh/km/lane is the
jam density of a single lane; and the critical density
$\r_{c,m}=0.2\rho_{j,m}=0.4\r_j$=72 veh/km. Therefore, the speed-
and flow-density relationships can be written as follows:
\bqs
V_m(\r)&=&\cas{{ll} v_{f,m},& 0\leq\r\leq \r_{c,m};\\\frac{\r_{c,m}}{\r_{j,m}-\r_{c,m}}\frac{\r_{j,m}-\r}{\r} v_{f,m},& \r_{c,m}<\r\leq\r_{j,m}.}\\
Q_m(\r)&=&\r V_m(\r)=\cas{{ll} v_{f,m}\r,& 0\leq\r\leq \r_{c,m};\\\frac{\r_{c,m}}{\r_{j,m}-\r_{c,m}} v_{f,m} (\r_{j,m}-\r),& \r_{c,m}<\r\leq\r_{j,m}.}
\eqs
For the on-ramp, the free flow speed $v_{f,r}=$35 mph=2.7934 $l/\tau$; the jam density is $\r_j$; and the critical density $\r_{c,r}=0.2 \r_j$. Similarly, we can have the following speed- and flow-density relationships:
\bqs
V_r(\r)&=&\cas{{ll} v_{f,r},& 0\leq\r\leq \r_{c,r};\\\frac{\r_{c,r}}{\r_{j}-\r_{c,r}}\frac{\r_{j}-\r}{\r} v_{f,r},& \r_{c,r}<\r\leq\r_{j}.}\\
    Q_r(\r)&=&\r V_r(\r)=\cas{{ll} v_{f,r}\r,& 0\leq\r\leq \r_{c,r};\\\frac{\r_{c,r}}{\r_{j}-\r_{c,r}} v_{f,r} (\r_{j}-\r),& \r_{c,r}<\r\leq\r_{j}.}
\eqs
The above relationships are depicted in \reffig{fd_tri_merge}.

\bfg
\bc\includegraphics[height=12cm] {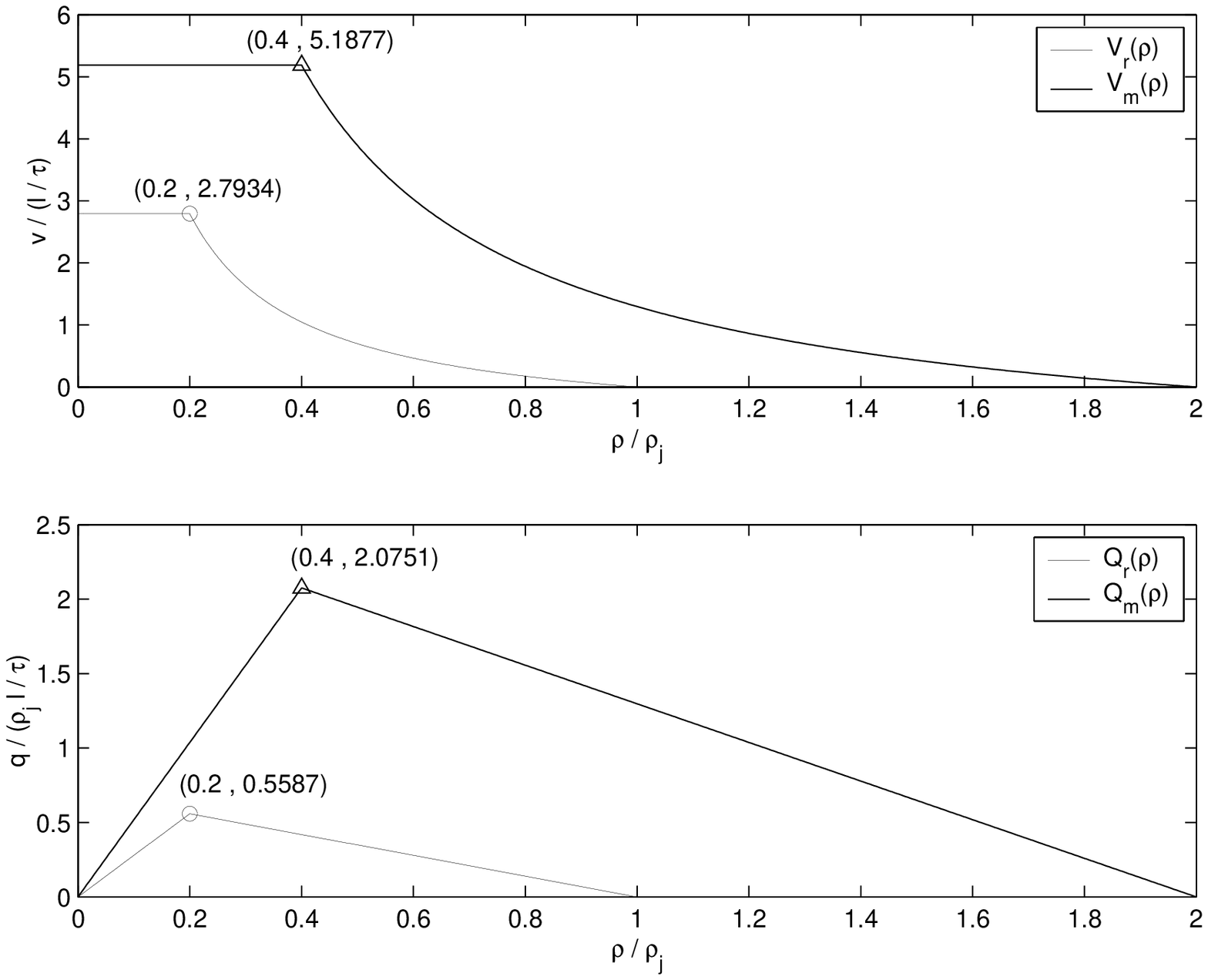}\ec
\caption {The triangular fundamental diagrams for the mainline freeway and the on-ramp} \label {fd_tri_merge}
\efg

Since $|\l_{\ast,m}|\leq v_{f,m}=5.1877 l/\tau$, where $\l_{\ast,m}(\r)=V_m(\r)+\r V_m'(\r)$ is the characteristic speed of \refe{LWR} on the mainline freeway, we find the CFL \citep{courant1928CFL} condition number
\bqs
|\l_{\ast,m}| \frac {\dt}{\dx}&\leq& 0.65<1.
\eqs
Since the characteristic speed on the on-ramp is smaller than that on  the
mainline freeway, the CFL condition is also satisfied for the on-ramp.

\subsection{Simulation of merging traffic without control}
In this subsection, we study the following merging traffic. Initially,  the
mainline freeway carries a constant flow: traffic
densities on the upstream and downstream freeway are the same,
$\r_{u_1}=\r_{u_2}=0.36 \r_j$, which is under-critical. After time $t=0$, a
constant flow with density $\r_{u_2}=0.175 \r_j$ arrives at the
on-ramp, and the on-ramp remains uncontrolled. In our simulation, the
Riemann boundary condition is imposed for the upstream boundaries of link
$u_1$ and $u_2$ and the downstream boundary of link $d$; i.e. the spatial
derivatives of traffic density at these boundaries are assumed to be zero.

After partitioning each of the three links into $N=500$ cells  and
discretizing the time duration into $K=5000$ steps, we obtain
simulation results as shown in \reffig{tri_cons_500}.
Figure 6(a) illustrates the evolution of traffic density on the freeway
upstream of the merge: at time $t=0\tau$, traffic density is uniformly at
$\r_A=0.36\r_j$; after the arrival of the on-ramp flow, freeway traffic
immediately upstream of the merging point\footnote{Traffic density at the merging
point is multiple-valued.} becomes congested and
reaches a new state
$\r_B=0.7394\r_j$; then a shock wave forms and travels upstream in a
constant speed $s_1\approx -0.61 l\tau=-7.6$ mph. Figure 6(b) shows the
evolution of traffic on the freeway downstream of the merge: initially,
traffic density is also uniformly at $\r_A$. After $t=0$, traffic immediately
downstream of the merging point reaches capacity flow at density
$\r_C=0.4\r_j$, and  a contact wave\footnote{A contact wave is a
flow/density discontinuity traveling at the same speed of traffic on both
sides of it.} appears since
$\r_C$ and $\r_A$ are both on the free-flow side of the triangular
fundamental diagram.  It travels downstream at the speed
$s_2=5.2 l/\tau=v_{f,m}$. Figure 6(c) shows that a backward
shock wave also forms on the on-ramp: traffic upstream of the shock
 has density $\r_D=0.175\r_j$ and downstream of the shock density
$\r_E=0.3697 \r_j$. This shock travels at
$s_3\approx-0.25 l/\tau=-3.1$ mph. The shock waves and the contact wave
observed on the three branches are shown in
Figure 6(d) on the $\r-q$ plane, in which line $AB$ represents
the shock wave on the freeway upstream of the merge and the slope of
$AB$ is its speed; line $CA$ represents the contact wave
on on the freeway downstream of the merge and its travel speed is the
slope of $CA$, which is also the free flow speed;
finally line $DE$ represents the shock wave on the on-ramp and the slope of
$DE$ is its travel speed.

\bfg
\bc\includegraphics[height=12cm] {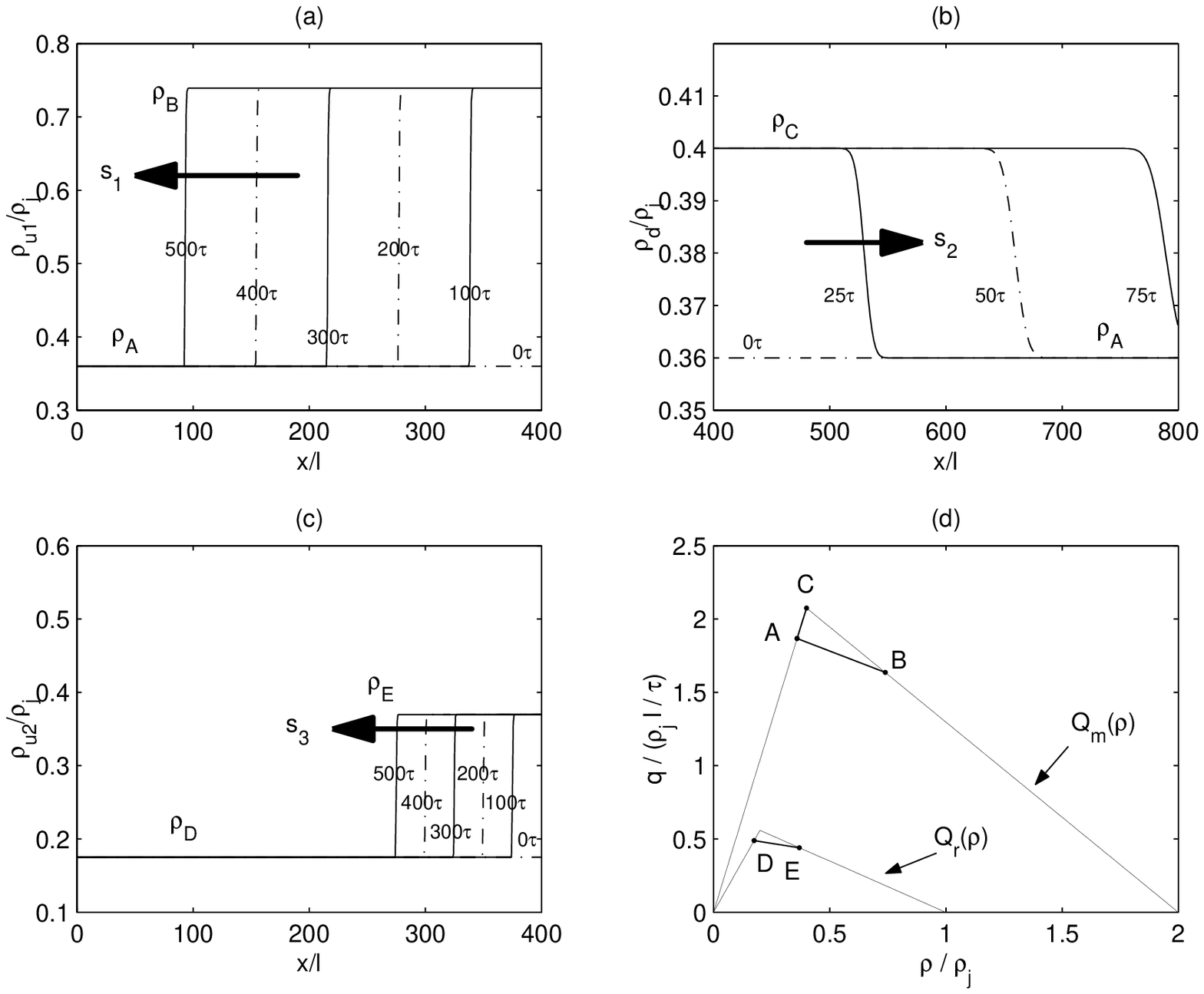}\ec
\caption {Simulation of merging traffic without control} \label {tri_cons_500}
\efg

Comparing the congested states on the upstream freeway (point $B$ in
\reffig{tri_cons_500}(d) ) and that on the on-ramp (point $E$ in
\reffig{tri_cons_500}(d)), we have the following observations: 1) the
ratio of the outflows from the upstream branches (freeway: $q_B=1.6349
\r_j l/\tau$=5933 veh/hr and on-ramp: $q_E=0.4402
\r_j l/\tau$=1597 veh/hr)  are not proportional to the lane ratio, owing to
the  different geometric and flow characteristics of the two
upstream branches, as reflected in their respective fundamental diagrams
and 2). as long as the freeway downstream of the merge is not
congested, and the total demand is greater than the total supply at the
merging point, traffic states surrounding the merge,
$\r_B$, $\r_E$, and $\r_C$,  are constant (i.e., stationary) states regardless
of the initial traffic conditions. This unique characteristics of the merge
model with our suggested distribution scheme offers a way to validate the
model and the fairness assumption.

\subsection{Simulation of merging traffic when the on-ramp is controlled}
In this subsection, we have the same initial/boundary conditions as in the
preceding subsection, but with the on-ramp controlled by a ramp meter.
For simplicity we take a constant metering rate $r=0.3445 \r_j
l/\tau$=1250 veh/hr. The simulation results are shown  in
\reffig{c_tri_cons_500}. From \reffig{c_tri_cons_500}(a), we can see
that a backward traveling shock wave forms on the freeway upstream  of
the merge, traveling at $s_1\approx -0.33
l/\tau=-4.1$ mph.  Traffic densities besides the shock are $\r_A=0.36 \r_j$
(upstream) and $\r_B=0.6278 \r_j$ (downstream), respectively.
\reffig{c_tri_cons_500}(b) shows the evolution of traffic on the
freeway downstream of the merge, which is identical to the case without
ramp control. From  \reffig{c_tri_cons_500}(c), we can see that a backward
traveling shock wave also forms on the on-ramp, traveling at
$s_3\approx-0.48 l/\tau= -6.0$ mph.  On the ramp, traffic densities besides
the shock are
$\r_D=0.175 \r_j$ (upstream) and $\r_E=0.577 \r_j$
(downstream), respectively. Again, \reffig{c_tri_cons_500}(d) shows the
the initial and congested states in the $\r-q$ plane.

\bfg
\bc\includegraphics[height=12cm] {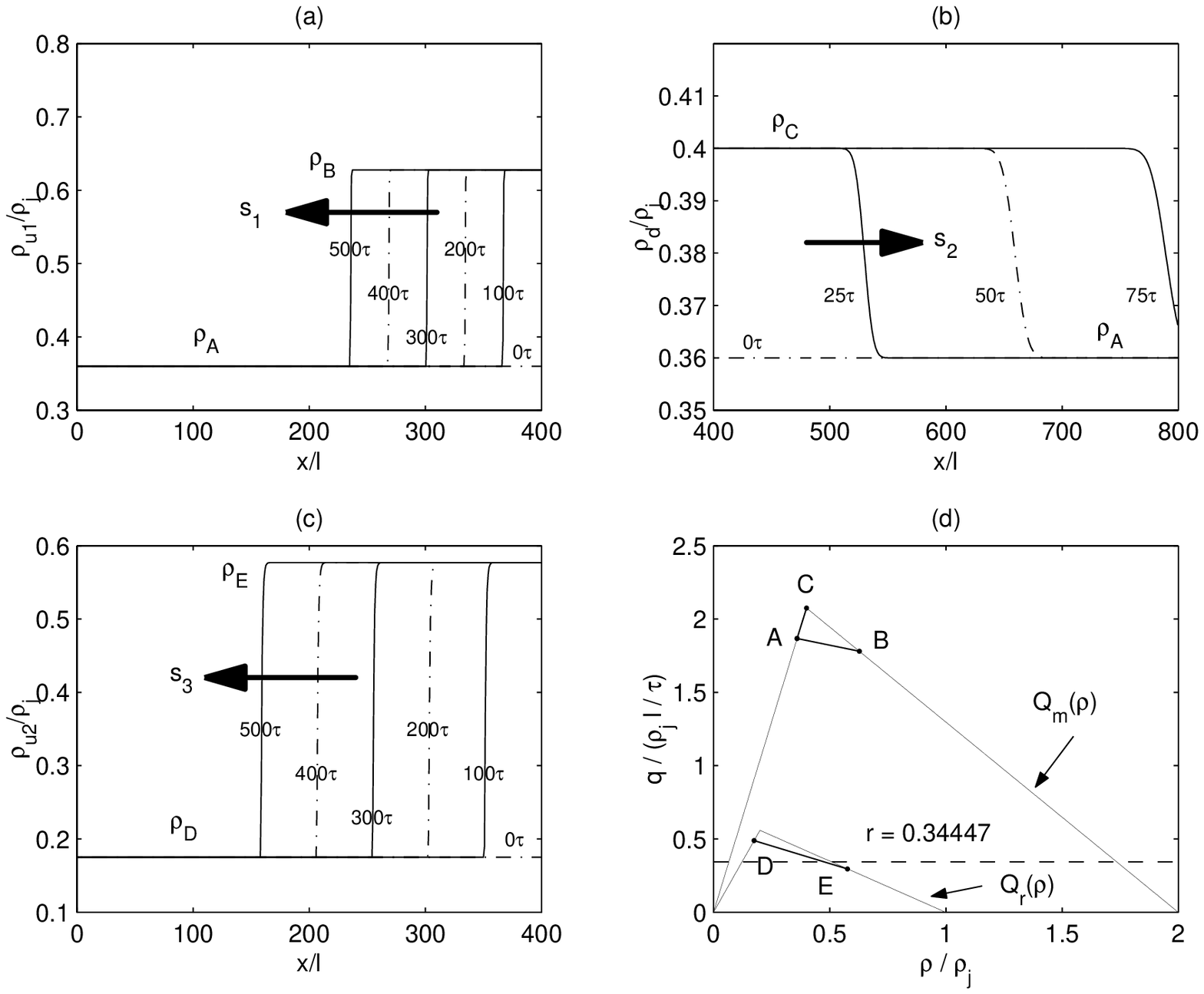}\ec
\caption {Simulation of merging traffic with on-ramp control} \label {c_tri_cons_500}
\efg

We note that with ramp control, the freeway upstream of the merge
becomes less congested ($\r_B$ is lower) while the on-ramp becomes more
congested ($\r_E$ is higher). Furthermore, the freeway queue grows
slower while the ramp queue grows faster with ramp control than without.
Finally,  the ramp control also affects the stationary states ($\r_B,
\r_C,\r_E$) and the distribution fractions, in favor of discharging higher
flow from the freeway.

\subsection{Computation of convergence rates}
In this subsection we will check the convergence of the merge model with
the distribution scheme \refe{fraction}, when the on-ramp is not
controlled. Here, this is done by computing convergence rates of traffic
density over the whole network.

First, we compare traffic density solutions at time $T_0=500\tau$ for two
different number of cells into which the network is partitioned and obtain
their difference. Denote solutions as $(U^{2N}_i)_{i=1}^{2N}$ for $2N$ cells
and $(U^N_i)_{i=1}^N$ for $N$ cells respectively, and define a difference
vector $(e^{2N-N})_{i=1}^N$ between these two solutions as
\bqn
\textbf e_i^{2N-N}&=&\frac12 (U^{2N}_{2i-1}+U^{2N}_{2i})-U^{N}_i, i=1,\cdots,N \label{def:error}.
\eqn
Then, the relative error with respect to $L^1$-, $L^2$- or $L^{\infty}$-norm can be computed as
\bqn
\epsilon^{2N-N}&=&\nm{{\textbf e}^{2N-N}}. \label{epsilon}
\eqn
Finally, a convergence rate is obtained when we compare two relative
errors:
\bqn
r&=&\log_2(\frac{\epsilon^{2N-N}}{\epsilon^{4N-2N}}).\label{rate}
\eqn
Here the vector of $U$ contains the densities of all three links which are
weighted by the number of lanes of each link.

We will compute convergence rate with the following conditions.
For link $u_1$, the number of lanes $a(u_1)=2$, we define its initial
condition as
\bqn
\ba{ccll}
\rho(x,0)&=&a(u_1)(0.18+0.05 \sin \frac {\pi x}L)\r_j, &x\in[0,L], \\
v(x,0)&=&V_m(\r(x,0)),&x\in[0,L].
\ea\label{ini_u1}
\eqn
For link $u_2$, the number of lanes $a(u_2)=1$, we define its initial
condition as
\bqn
\ba{ccll}
\rho(x,0)&=&a(u_2)(0.175+0.05\sin \frac {2\pi x}L)\r_j, &x\in[0,L], \\
v(x,0)&=&V_r(\r(x,0)),&x\in[0,L].
\ea\label{ini_u2}
\eqn
For link $d$, the number of lanes $a(d)=2$, we define its initial condition  as
\bqn
\ba{ccll}
\rho(x,0)&=&a(d)(0.18-0.05 \sin \frac {\pi x}L)\r_j, &x\in[L,2L], \\
v(x,0)&=&V_m(\r(x,0)),&x\in[L,2L].
\ea\label{ini_d}
\eqn
In addition, we impose the Riemann boundary condition for the upstream
boundaries of link $u_1$ and $u_2$ and the downstream boundary of link
$d$; i.e., $\r$ has zero derivative at these boundaries.

From \reft{conv_rate1}, convergence rates of the merge model with the
simple distribution scheme are of order one in $L^1$ norm.  The
convergence rate is of order one is  expected because the Godunov
method used here is a first-order method, in which traffic density on a link
is approximated by piece-wise constant functions. Unlike the Godunov
discretization of the LWR model, which we know is convergent to the LWR
model, we do not yet know what differential equations the discrete model
converges to. Nevertheless, the numerical convergence analysis gives us
comfort in applying the discrete merge model in the sense that we know
the properties of the computed solutions will not change over discretization
scales.

\btb
\bc
\begin{tabular}{lccccccc}\\\hline
$\rho/\rho_j$&128-64&Rate&256-128&Rate&512-256&Rate&1024-512\\\hline
$L^1$&3.31e-03&1.00&1.65e-03&1.00&8.27e-04&1.00&4.13e-04\\\hline
\end {tabular}
\ec
\caption{Convergence rates of the discrete merge model}\label{conv_rate1}
\etb 
%Wenlong Jin: 2003
\section{Discussions}\label{sec:merge5}
In this chapter, we studied the discrete kinematic wave model of merges
in the supply-demand framework, probed the supply-demand
merge models by Daganzo and Lebacque, and gained a better
understanding of various distribution schemes in merge models. More
importantly, we proposed the simplest distribution scheme, equivalent
to the so-called ``fairness" condition, in which the distribution
fractions are proportional to the upstream demands. We
demonstrated with both analytical discussions and numerical studies
that the discrete merge model with the simplest distribution scheme
is well-defined.

The simple distribution scheme can also be easily extended to
diverges in single-commodity traffic flow. In reality, a diverge is not
simply a reverse of a merge: in a merge, vehicles advance in the same
direction and once they travel into the downstream link, their origins
no longer affect traffic dynamics; while in a diverge, traffic dynamics
in the upstream link are affected by the combination of vehicles with
different destinations. However, for single-commodity traffic flows,
in which vehicles have no predefined routes and can choose any
downstream links when diverging, we can have similar
supply-demand method for computing fluxes through a diverge:
assuming that a diverge has $K$ downstream links, the upstream link
and the downstream links are properly discretized, the traffic demand
of the upstream link is $D_u$, and traffic supply of downstream link
$d_i$ ($i=1,\cdots,K$) is $S_i$, the supply-demand method with the
``fairness" condition, in which the influx into a downstream link is
proportional to its accommodation $S_i$, can be written as:
\bqn
\ba{lcl}
q&=&\min\{D_u,\sum_i^K S_i\},\\
q_i&=&p_i q,\\
p_i&=&\frac{S_i}{\sum_i^K S_i},\quad i=1,\cdots,K,
\ea\label{diverge}
\eqn
where $q_i$ is influx into the downstream link $d_i$.

Further, for a general junction with $J$ upstream links and $K$ downstream links, we can compute fluxes as
\bqn
\ba{lcl}
q&=&\min\{\sum_{j=1}^J D_j,\sum_{k=1}^K S_k\},\\
q^u_j&=&\alpha_j q,\\
q^d_k&=&\beta_k q, \\
\alpha_j&=&\frac{D_j}{\sum_{j=1}^J D_j},\quad j=1,\cdots,J\\
\beta_k&=&\frac{S_k}{\sum_{k=1}^K S_k},\quad k=1,\cdots,K,
\ea\label{junction}
\eqn
where $D_j$ is the traffic demand of the upstream link  $u_j$, $S_k$
is the traffic supply of the downstream link $d_k$, $q^u_j$ is the
out-flux from $u_j$, and $q^d_k$ the influx into $d_k$. Although this
model may not be proper for realistic traffic flow, it could be helpful
for the study of other unidirectional flows.

With a complete picture of the merge model, we can now simulate
how a traffic disturbance propagates through a merge. In particular,
since merges are well recognized as locations where congestion
often  initiates \citep{daganzo1999remarks}, this study will help
traffic engineers understand how congestion start and propagate at
a merge. Moreover, with the complete discrete kinematic wave merge
model, we are able to develop and evaluate local ramp metering
strategies, e.g., the ALINEA algorithm \citep{papageorgiou1997ALINEA}. In the future, we will be also
interested in the limit of the discrete kinematic wave model with the
simple distribution scheme. In particular, we will be interested in
kinematic wave solutions for a general Riemann problem of this
model. Finally, and more importantly, the supply-demand method and
the ``fairness" condition needs to be checked against observed
merging phenomena.

\newpage
\pagestyle{myheadings}
\markright{  \rm \normalsize CHAPTER 4. \hspace{0.5cm}
 DIVERGING TRAFFIC MODEL}
\chapter{Kinematic wave traffic flow model of diverging traffic}
%\thispagestyle{myheadings}
%Wenlong Jin: 2003
%Wenlong Jin: 2003
\section{Introduction}
Diverges have been well recognized as a major type of bottlenecks \citep{daganzo1999remarks,daganzo1999phase}, where congestion of one downstream branch can propagate to the upstream branch and further blocks vehicles traveling to other directions. Observations of diverging traffic \citep{munoz2002diverge, munoz2003structure} have shown that vehicles to different downstream branches tend to segregate themselves at a point, called the actual diverging point, as far as 2 km upstream of the physical diverging point. In addition, traffic  upstream of the actual diverging point is in the so-called ``1-pipe" regime in the sense that vehicles that appear at the same time and location have the same speed even they have different destinations; while between the actual and physical diverging points traffic is in a ``2-pipe" regime, since vehicles to different directions are separated and travel independently. In certain cases, however, the actual and physical diverging points may be the same \citep{newell1999exit}. Although it is an interesting problem to determine the location of an actual diverging point, here we only intend to study traffic dynamics at an actual diverging point or simply a highway diverge; i.e., we want to know, under certain traffic conditions, how much traffic will flow into each downstream branch.

Due to its simplicity in representation and efficiency in computation, the kinematic wave theory developed by \citep{lighthill1955lwr, richards1956lwr} has played an important role in understanding traffic dynamics on road networks. In the framework of the kinematic wave theory of diverging traffic, there have been two major lines of research in the framework of the kinematic wave theory. In the first line of research, the effect of diverging traffic to the mainline freeway is incorporated by considering an off-ramp as a sink \citep{kunhe1992continuum}. In these models, only part of the interaction between the diverging traffic and the main highway traffic, i.e., the impact of diverging traffic to the mainline traffic, is taken into account. In addition, the characteristics of sinks are not trivial to obtain. In contrast, in the second line of research, the complete dynamics at diverges is studied by treating each branch equally. When vehicles do not have predefined routes and can choose any downstream branches at a diverge, \cite{holden1995unidirection} determined the flow to each branch by solving an optimization problem, and \cite{jin2003merge} proposed a distribution scheme based on the definitions of traffic supply and demand. For more realistic situations when vehicles have predefined route choices, \cite{daganzo1995ctm} and \cite{lebacque1996godunov} proposed the so-called supply-demand method for computing traffic flows through a diverge. In this method, each branch is partitioned into a number of cells and a duration of time is discretized into time steps. Then a diverge connects an upstream cell and a number of downstream cells, and the flow to a downstream cell during a time period is assumed to be proportional to the fraction of vehicles in the upstream cell that wish to travel to this downstream cell. In addition, the flows to the downstream cells are assumed to reach their maxima subject to the constraint of traffic demand of the upstream cell and supplies of the downstream cells.

Although the supply-demand method is quite straightforward in computing flows through a diverge, it fails to provide a system-wide picture of how a small traffic perturbation on one branch will propagate through the diverge. That is, it gives no answer to what kind of kinematic waves will initiate in diverging traffic. In this chapter, we intend to study the kinematic waves (more specifically, the instantaneous kinematic waves) at a diverge. In our study, we still partition each branch into a number of cells and discretize a duration of time into some time steps. Then in a (short) time step\footnote{The length of this ``short" time step $\dt$, theoretically, can be infinitely small and, practically, can be obtained through a calibration exercise.}, we compute the in-flow to a downstream branch by holding traffic conditions of vehicles to the other destinations constant in the upstream cell and, after obtaining all flows, update traffic conditions in all the upstream and downstream cells. We can see that, in this theory, dynamics of traffic traveling to different destinations are first decoupled during a short time period and then combined together in an alternating manner. With this decoupling, corresponding to each downstream branch, we can obtain a hyperbolic system of conservation laws, which is independent of the systems corresponding to the other downstream branches. To the best of our knowledge, this decoupling approach in solving complicated dynamics of diverging traffic is the first of its kind in traffic flow.

Moreover, we will show that the decoupled systems are nonlinear resonant systems in the sense of \cite{isaacson1992resonance} and can be solved by combinations of shock, rarefaction, and standing waves. Here, these waves are considered instantaneous since they only exist in a short time interval, and waves corresponding to different downstream branches interact with each other at the end of each time interval. Thus, after a long time, the asymptotic kinematic waves arising from diverging traffic can be totally different from the instantaneous waves. In this chapter, we will focus on the instantaneous waves arising from each branch, and the interaction between these waves and the asymptotic kinematic waves will be further studied in a following chapter.

The rest of this chapter is organized as follows. In Section \ref{sec:diverge2}, we first formulate the kinematic wave theory of diverges in the framework of the LWR model, and introduce the instantaneous kinematic wave theory for dynamics of diverging traffic. In Section \ref{sec:diverge3}, we solve in detail for the instantaneous kinematic waves. In Section \ref{sec:diverge4}, we present a new definition of traffic demand for diverging vehicles and the supply-demand method based on this new definition. In Section \ref{sec:diverge5}, we carry out simulations of the new model for two cases. In Section \ref{sec:diverge6}, we make our conclusions and propose some follow-up research.

%Wenlong Jin: 2003
\section {A kinematic wave theory for diverges}\label{sec:diverge2}
As we know, in a traffic stream, traffic dynamics can be affected by characteristics such as vehicle types and destinations, and drivers' behavior. Here vehicles of the same characteristics are considered to belong to a commodity. In this study we are only interested in the impact of vehicles' destinations on traffic dynamics at a highway diverge. Therefore, for a diverge with $D$ downstream branches and one upstream branch, we can differentiate vehicles into $D$ commodities.

Here we consider an ideal diverge: on the upstream branch, traffic is considered in 1-pipe regime; once through the diverge, vehicles completely segregate.
Therefore, we have a $D$-commodity flow on the upstream link and a single-commodity flow on each downstream link, and traffic dynamics on a road network with a diverge consists of three parts: dynamics of the single-commodity flow on each downstream branch, dynamics of the $D$-commodity flow on the upstream branch, and dynamics at the diverge. In the following part, after reviewing the kinematic wave models of single-commodity and $D$-commodity traffic, we will focus on the study of traffic dynamics at the diverge.

\subsection{The kinematic wave theory of single-commodity traffic flow}
For dynamics of single-commodity traffic, the seminal LWR \citep{lighthill1955lwr,richards1956lwr} kinematic wave theory can be applied. In the LWR model for traffic flow on a long crowded road, two equations are used to describe traffic dynamics of aggregate quantities such as traffic density $\r$, flow-rate $q$ and travel speed $v$:
\bqn
\m{Traffic conservation:}& &\r_t+q_x=0,\label{trafficconservation}\\
\m{Fundamental diagram:}& &q=Q(a,\r).\label{fundamentaldiagram}
\eqn
In \refe{fundamentaldiagram}, $a(x)$ is an inhomogeneity factor depending on road characteristics such as the number of lanes at $x$. Since $q=\r v$, we also have a speed-density relation: $v=V(a,\r)=Q(a,\r)/\r$; for vehicular traffic, generally, $v$ is nonincreasing and the equilibrium flow-density relation is assumed to be concave. From \refe{trafficconservation} and \refe{fundamentaldiagram}, the LWR model can be written as
\bqn
\r_t+Q(a,\r)_x&=&0.\label{LWRd}
\eqn
When the road is homogeneous; i.e., $a(x)$ is constant, the LWR model is called homogeneous; otherwise, when $a(x)$ is dependent on location, \refe{LWRd} is called an inhomogeneous LWR model.

The homogeneous LWR model is nothing but a scalar conservation law and has been well studied both analytically and numerically. It is solved by combinations of shock and expansion (rarefaction) waves. In contrast, the inhomogeneous LWR model \refe{LWRd} can be considered as a nonlinear resonant system and has an additional type of kinematic waves - standing waves, which initiate and stay at the inhomogeneity \citep{jin2003inhLWR}.
Therefore, for a single-commodity flow on a road branch, traffic dynamics can be considered as the evolution of a combination of shock, rarefaction, and standing waves.

\subsection{The kinematic wave theory of multi-commodity traffic flow}
For traffic on the upstream link, it is assumed that vehicles with different destinations have no differences in dynamics; i.e., the speed-density relation is independent of commodities. Let $\r_i(x,t)$, $v_i(x,t)$, and $q_i(x,t)$ ($i=1,\cdots,D$) denote, at place $x$ and time $t$,  the density, speed, and flow-rate of commodity $i$, respectively. In contrast, $\r(x,t)$, $v(x,t)$, and $q(x,t)$ are the aggregate density, speed, and flow-rate. Then we have (with $(x,t)$ suppressed hereafter)
\bqn
\r&=&\sum_{i=1}^D \r_i,\\
v&=&v_i=V(a,\r),\qquad i=1,\cdots,D,\\
q&=&\sum_{i=1}^D q_i,
\eqn
where $q=\r v$, $q_i=\r_i v_i$ ($i=1,\cdots,D$), $V(a,\r)$ is the aggregate speed-density relation. We also call this kind of multi-commodity traffic as additive. In additive traffic, the flow of each commodity is called a partial flow and the flow containing all commodities the total flow.

Therefore, in additive traffic, dynamics of each partial flow can be written as
\bqn
(\r_i)_t+\left( \r_i V(a,\r)\right )_x&=&0, \qquad i=1,\cdots,D.
\eqn
Note that these $D$ equations are coupled through the aggregate speed-density relation.
By introducing the fraction of vehicles of commodity $i$, $\xi_i=\r_i/\r$, which satisfies $\sum_{i=1}^D \xi_i=1$, \cite{lebacque1996godunov} showed that traffic dynamics of additive traffic can also be written as
\bqn
\r_t+\left(\r V(a,\r)\right)_x&=&0,\\
(\xi_i)_t+V(a,\r)\: (\xi_i)_x&=&0, \qquad i=1, \cdots, D-1.
\eqn
We can see that the fractions travel in the same speed as vehicles. This shows that the First-In-First-Out principle is guaranteed in the kinematic wave theory of additive traffic \citep{lebacque1996godunov}.  In addition, we can see that, for multi-commodity traffic, the aggregate traffic dynamics can be described by the (inhomogeneous) LWR model and is still a combination of shock waves, rarefaction waves, and/or standing waves, and traffic dynamics of a commodity is of the similar pattern.

\subsection{A kinematic wave theory of diverging traffic}
To model traffic dynamics on a network with a diverge, the harder problem is to find the number of vehicles traveling from the upstream link to each downstream link during a time interval. This is an important problem to answer since flows from the upstream link to the downstream links are the necessary boundary conditions for the connected links. Here we propose a new approach to finding the approximate flows through a diverge. This new theory is different from the supply-demand method of \citep{daganzo1995ctm,lebacque1996godunov} in that this theory can be analytically solved by a combination of kinematic waves.

Like the supply-demand method of \citep{daganzo1995ctm,lebacque1996godunov}, the new kinematic wave theory is also in a discrete form: we partition each branch into a number of cells of length $\dx$, discretize the time duration into time steps of length $\dt$, and assume that the aggregate and partial traffic densities are constant in each cell at a time step. The discretization of space-time plane in this model is required to satisfy the CFL condition \citep{courant1928CFL}, which guarantees that the flows through the diverge are only dependent on traffic  densities in the adjacent cells and are independent of traffic conditions farther upstream or downstream.

In the discrete model, connected to the diverge with $D$ downstream branches are an upstream cell and $D$ downstream cells, denoted as downstream cell $i$ ($i=1,\cdots,D$). Assume that, at $t=0$, partial densities in the upstream cell are $\r_{i}^L$ and the aggregate density $\r^L=\sum_{i=1}^D \r_i^L$, and density in downstream cell $i$ is $\r_i^R$. With these initial traffic conditions, we will show how traffic dynamics will evolve from $t=0$ till $t=\dt$. In addition, we will be able to compute the in-flow to downstream cell $i$ ($i=1,\cdots,D$), $q_i$. Then from traffic conservation, the out-flow from the upstream cell is $\sum_{i=1}^D q_i$. With these flows, traffic conditions at the next time step $t=\dt$ can be obtained.

In this new model, to study the traffic dynamics of commodity $i$ through the diverge from $t=0$ to $t=\dt$, we hold partial densities of vehicles of the other commodities in the upstream cell constant during this time interval. That is, $\r_j^L$, where $j=1,\cdots,D$ and $j\neq i$, are all constant. Here we assume that traffic densities of commodities other than $i$ are constant during the short time interval. This assumption is equivalent to saying that the change of traffic conditions of commodity $i$ is, in a short time interval, independent of traffic dynamics in the other directions. This assumption appears to be a natural and reasonable way to decouple the interactions between traffic in different directions. This approach can be considered as an attempt to answer an important question: how do vehicles of different destinations segregate themselves at a diverge?

We set the diverging point at $x=0$ and vehicles travel in the positive direction of $x$-axis. Thus, $x<0$ for the upstream link and $x>0$ for the downstream link. Also, to simplify the notations, we rewrite $\r_i$ as $\r$ and the total traffic densities of other commodities as $k$. Therefore, traffic dynamics of vehicles of commodity $i$ in the upstream cell can be written as
$\r_t+\left (\r\: V(a(x),\r+k)\right)_x=0$ for $x<0$, and $\r_t+\left ( \r\: V(a(x),\r)\right)_x=0$ for $x>0$.
Since $k$ is held constant during the time interval $[0,\dt]$, and $k=0$ in downstream cell $i$, we have the following equation and conditions for traffic dynamics of commodity $i$:
\bqn
\r_t+\left ( \r\: V(a(x),\r+k)\right)_x&=&0, \label{model}
\eqn
with the following initial traffic conditions
\bqn
\ba{lcl}
\r&=&\cas{{ll} \r_L& x<0,\\\r_R&x>0;}\quad t=0,\\
k&=&\cas{{ll} k_L& x<0,\\0&x>0;}\quad t\in[0,\dt].
\ea\label{init}
\eqn
Here we extend our $x$-axis to both positive and negative infinity since, with the CFL condition, the traffic dynamics only depend on the conditions above.

Thus, \refe{model} with \refe{init} is a description of traffic dynamics of commodity $i$ during a short time interval. This new model will be shown to be solved by instantaneous kinematic waves in the following section.

%Wenlong Jin: 2003
\section{The instantaneous kinematic waves}\label{sec:diverge3}
In this section, we discuss the instantaneous kinematic waves of \refe{model} with initial conditions \refe{init}. For the purpose of exposition, we here assume that both the upstream and downstream branches share the same speed-density relation; i.e., we exclude the road characteristics $a(x)$. We also assume that $q$ is concave and $v$ non-increasing in $\r\in[0,\r_j]$, where $\r_j$ is the jam density. With modifications, the results here are applicable to the case when the upstream and downstream branches have different characteristics.

Since we assume that $k$ is time independent, we have another conservation equation: $k_t=0$. Therefore \refe{model} with \refe{init} can be rewritten as a $2\times2$ system:
\bqn
\ba{lcl}
\r_t+\left(\r V(\r+k)\right)_x&=&0,\\
k_t&=&0,
\ea\label{CL}
\eqn
together with initial conditions
\bqn
\ba{lcl}
\r(x,t=0)&=&\cas{{ll}
\r_{L} & x<0,\\
\r_{R} & x>0,}\\
k(x,t=0)&=&\cas{{ll}
k_L &x<0,\\
0&x>0.}
\ea\label{initCL}
\eqn

Next, we will show that \refe{CL} is indeed a nonlinear resonant system studied by \citet{isaacson1992resonance}. Consequently, we can use their technique to solve this system, which describes traffic dynamics of one commodity during a short time period. Note that here we do not restrict our time domain to be $[0,t_1]$. This is because \refe{CL} with \refe{initCL} has self-similar solutions, which means that the flow $f(x=0,t)$ is constant for any time $t>0$.

\subsection{The properties of \refe{CL} as a nonlinear resonant system}
Denote the state $U=(k,\r)$, the flux vector $F(U)=\left(0,\r V(\r+k)\right)$, \refe{CL} can be written as a hyperbolic system of conservation laws:
\bqn
U_t+F(U)_x&=&0,\label{sys_cons}
\eqn
where $x\in (-\infty, \infty)$, $t\geq 0$. This system can be written in the following form of a quasi-linear system
\bqn
U_t+\partial F(U)U_x&=&0,
\eqn
where the Jacobian matrix $\partial F(U)$ is
\bqn
\partial F&=&\mat{{cc}0&0\\\r V'(\r+k)&V(\r+k)+\r V'(\r+k)}.
\eqn
Thus the two eigenvalues of $\partial F(U)$ are
\bqs
\lambda_0=0,\quad \lambda_1=V(\r+k)+\r V'(\r+k),
\eqs
and their corresponding right eigenvectors are
\bqs
\vec R_0=\mat{{c}V(\r+k)+\r V'(\r+k)\\-\r V'(\r+k)},\quad\vec R_0=\mat{{c}0\\1}.
\eqs

System \refe{sys_cons} is a non-strictly hyperbolic system, since it may happen that
\bqn
\l_1(U_\ast)=\l_0(U_\ast)=0. \label{cond1}
\eqn
When \refe{cond1} is satisfied, we say that traffic state $U_\ast$ is partially critical. If denoting the partial fundamental diagram as $Q(\r;k)=\r V(\r+k)$, we can define other partial quantities of the commodity that we are interested in, partial capacity $Q^{\max}(k)$ and partial critical density $\gamma(k)$, as follows:
\bqn
Q(\r;k)&\leq&Q^{\max}(k), \quad \forall \r\in[0,\r_j-k],\\
Q^{\max}(k)&=&Q(\gamma(k);k).
\eqn
Since $V'\leq0$ and $d^2 (\r V(\r))/ d \r^2<0$, we have
\bqs
\frac{\partial^2 Q(\r;k)}{\partial \r^2} &<&0,\\
\frac{\partial Q(\r;k)}{\partial k} &\leq&0;
\eqs
therefore, $Q^{\max}(k)$ and $\gamma(k)$ are unique for given $k\in[0,\r_j]$. Thus the partially critical state $U_\ast=(k_\ast,\r_\ast)=(k_\ast,\gamma(k_\ast))$. We call the collection of the partially critical states, $\Gamma=\left\{ U_\ast | k\in[0,r_j]\right\}$, the transition curve.

Moreover, at a partially critical state $U_\ast$, \refe{CL} is genuinely nonlinear; i.e.,
\bqn
\frac{\partial}{\partial \r} \l_1(U_\ast)&=&\frac{\partial^2}{\partial \r^2} (\r V(\r+k))<0,\label{cond2}
\eqn
since $q=\r V(\r+k)$ is concave in $\r$, and $\partial F(U)$ is nondegenerate; i.e.,
\bqn
\frac{\partial}{\partial k} (\r_\ast V(\r_\ast+k_\ast))&=&\r_\ast V'(\r_\ast+k_\ast)=-V(\r_\ast+k_\ast)<0. \label{cond3}
\eqn

With conditions \refe{cond1}-\refe{cond3} satisfied, \refe{CL} is a nonlinear resonant system in the sense of \citep{isaacson1992resonance}, and it is guaranteed that, in a neighborhood of the state $U_\ast$, the Riemann problem for \refe{CL} with \refe{initCL} has a unique solution with a canonical structure.

\subsection{The instantaneous kinematic waves of \refe{CL} with \refe{initCL}}
The Riemann problem for \refe{CL} with initial conditions \refe{initCL} is solved by two families of basic waves, associated with the two eigenvalues. The 0-waves with wave speed $\l_0=0$, also called standing waves, are the integral curves of $\vec R_0$ in $U$-space, and therefore are given by $Q(\r;k)$=const. Similarly, the 1-waves with wave speed $\l_1(U)$ are integral curves of $\vec R_1$ in $U$-space, and are given by $k$=const. As shown in \reffig{basic_waves}, the 0-curve is concave and tangent to the 1-curve at the critical state $U_\ast$; states left to the transition wave $\Gamma(k)$ are undercritical (UC), and overcritical (OC) right to the transition wave.

\bfg
\bc\includegraphics[height=12cm] {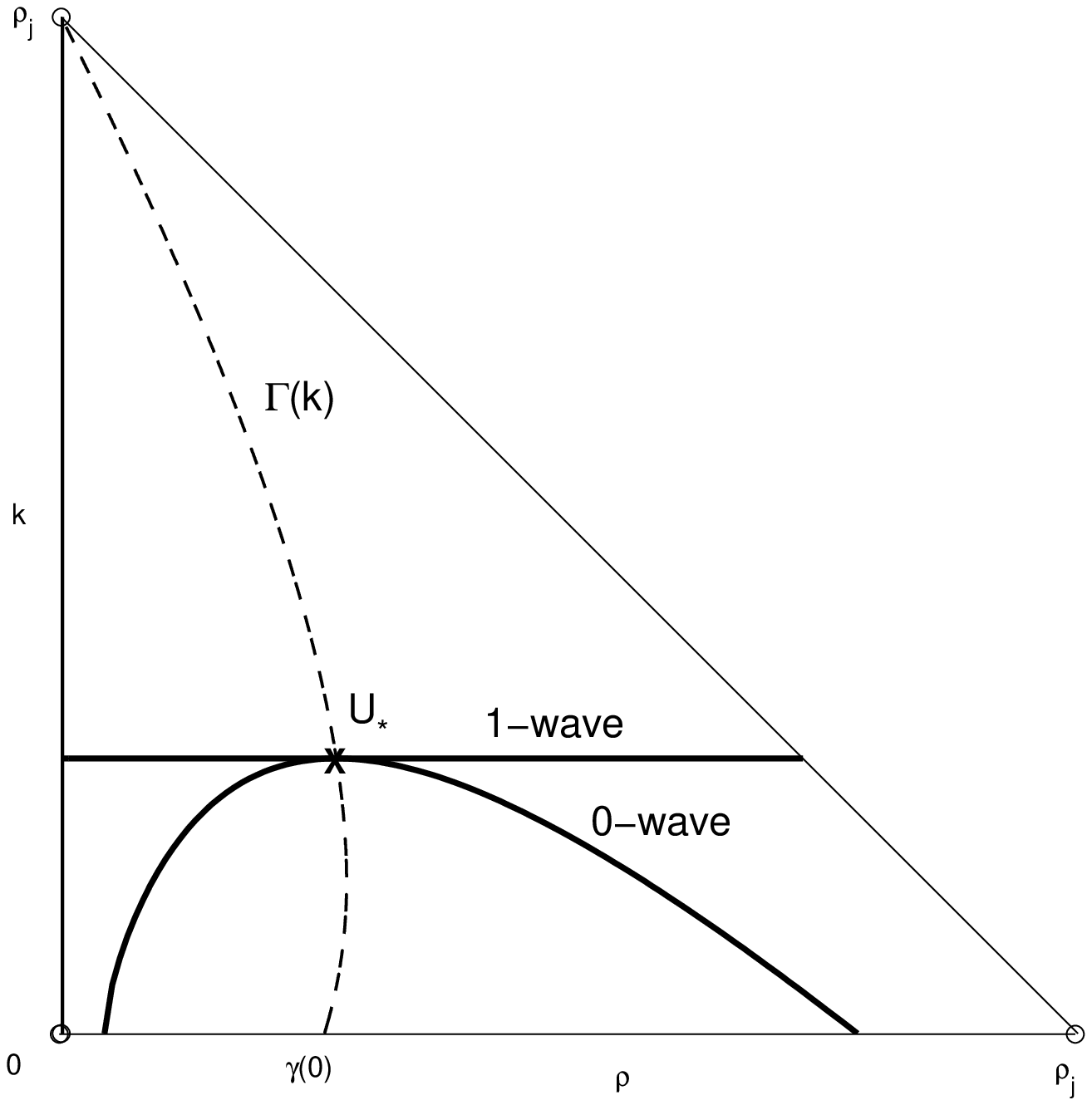}\ec
\caption{Integral curves in $(\r,k)$-space}\label{basic_waves}
\efg

\citet{isaacson1992resonance} showed that wave solutions to the Riemann problem are formed by no more than three basic waves, including standing waves, shock waves, and rarefaction waves, and these waves satisfy two entropy conditions: Lax's entropy condition (1972), i.e., that these waves increase their wave speeds from left (upstream) to right (downstream), and an additional condition that an UC state and an OC state can not be connected by a standing wave. With these conditions, wave solutions exist and are unique: when $U_L$ is UC, wave solutions in $U$-space are shown in \reffig{left_waves}, and in \reffig{right_waves} when $U_L$ is OC. Note that the solutions can also be categorized according to whether $\r_R$ is OC or UC. Also recall that $U_R$ always lies on $k=0$.

\bfg
\bc\includegraphics[height=12cm] {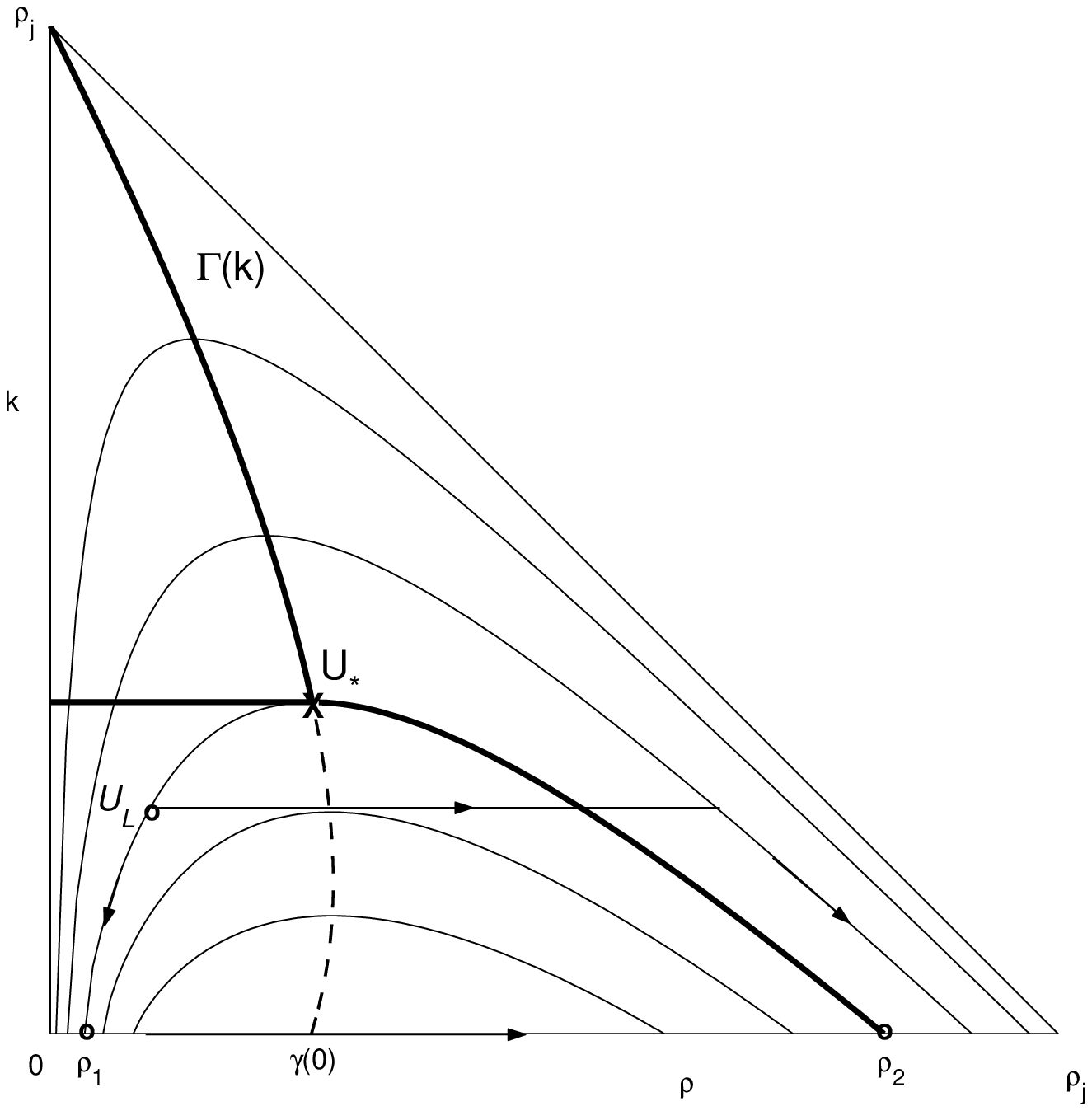}\ec
\caption{The Riemann solutions when $U_L$ is UC}\label{left_waves}
\efg
\bfg
\bc\includegraphics[height=12cm] {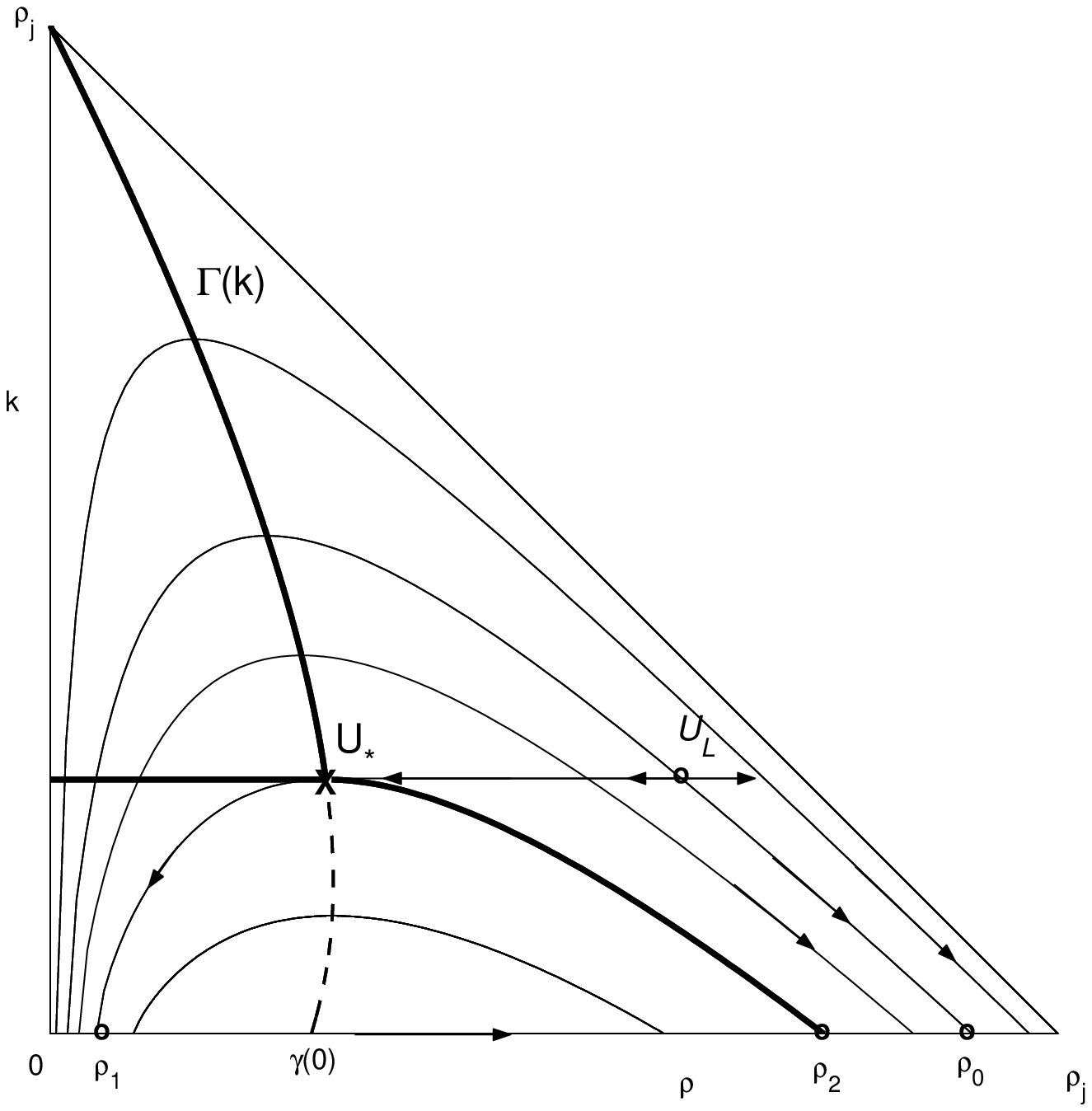}\ec
\caption{The Riemann solutions when $U_L$ is OC}\label{right_waves}
\efg

In the remaining part of this subsection, we will discuss the Riemann solutions of \refe{CL} with \refe{initCL}, and compute the boundary flux $q(x=0,t)$.

\ben
\item When $U_L$ is UC; i.e., $\r_L<\gamma(k_L)$, we denote the intersection between $Q(U_L)$ and $k=0$ by $\r_1$ and $\r_2$, where $\r_1\leq \gamma(0)\leq \r_2$. Hence, we obtain three types of solutions when $\r_R\in [0,\r_1]$, $(\r_1,\r_2)$, or $[\r_2,\r_j]$ as shown in \reffig{left_waves}.
\bi
\item [Type 1] When $\r_R\in[0,\r_1]$; i.e., $Q(\r_R;0)\leq Q(U_L)$ and $U_R$ is UC, wave solutions to the Riemann problem consist of two basic waves with the intermediate state, $U_1=(\r_1,0)$: the left wave $(U_L,U_1)$ is a standing wave and the right wave $(U_1,U_R)$ is a rarefaction wave with characteristic velocity $\l_1(U)>0$. As shown in \reffig{type_1}, the boundary flux $q(x=0,t>0)=Q(U_L)$.

\bfg
\bc\includegraphics[height=12cm] {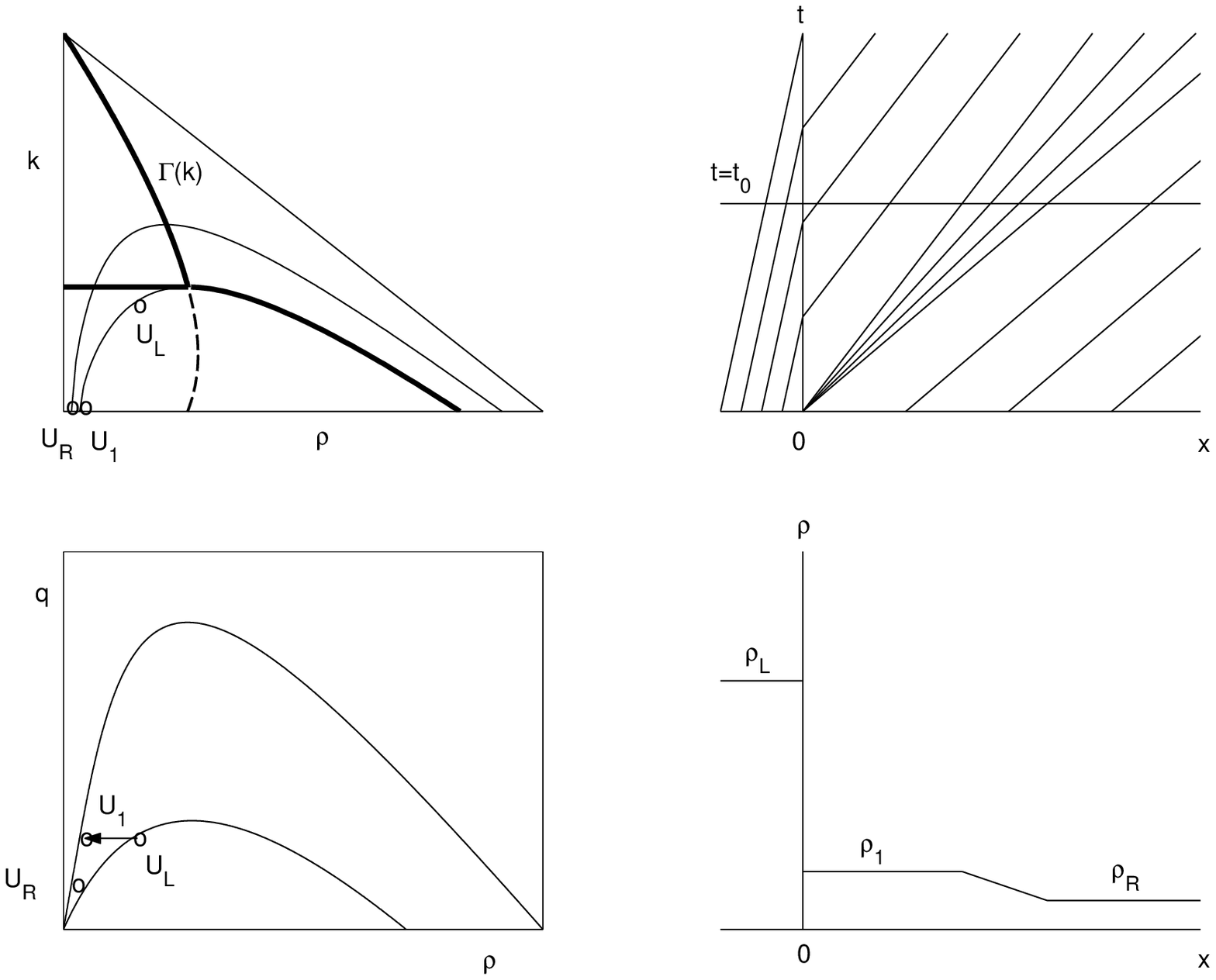}\ec
\caption{An example for wave solutions of type 1 for \refe{CL} with initial conditions \refe{initCL}}\label{type_1}
\efg

\item[Type 2] When $\r_R\in(\r_1,\r_2)$; i.e., $Q(\r_R;0)> Q(U_L)$, wave solutions to the Riemann problem consist of two basic waves with the intermediate state, $U_1=(\r_1,0)$: the left wave $(U_L,U_1)$ is a standing wave and the right wave $(U_1,U_R)$ is a shock wave with positive wave speed $s(U_1,U_R)=(Q(\r_R;0)-Q(U_L))/(\r_R-\r_1)$. As shown in \reffig{type_2}, the boundary flux $q(x=0,t>0)=Q(U_L)$.

\bfg
\bc\includegraphics[height=12cm] {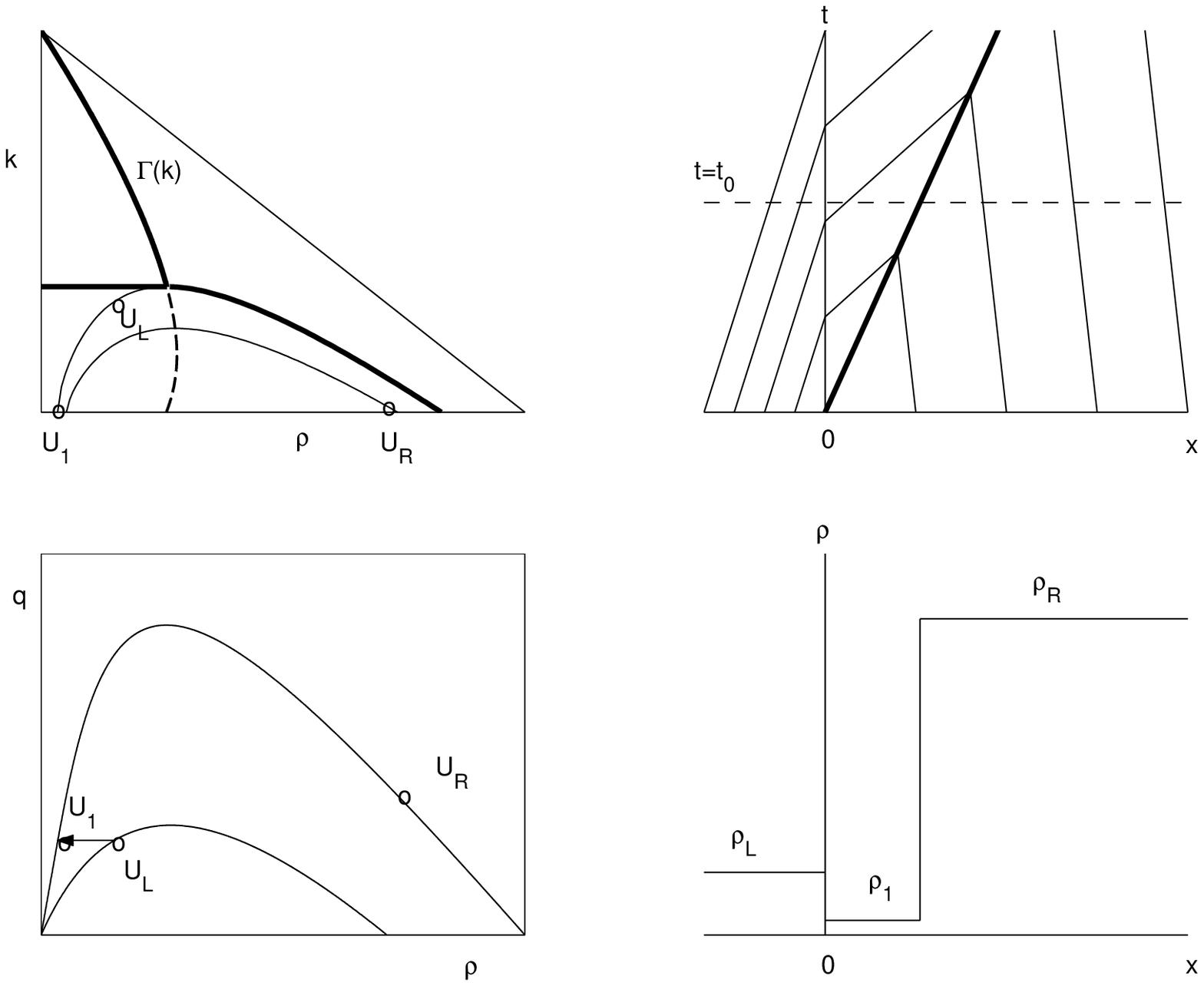}\ec
\caption{An example for wave solutions of type 2 for \refe{CL} with initial conditions \refe{initCL}}\label{type_2}
\efg

\item[Type 3] When $\r_R\in[\r_2,\r_j]$; i.e., $Q(\r_R;0)\leq Q(U_L)$ and $U_R$ is OC, wave solutions to the Riemann problem consist of two basic waves with the intermediate state, $U_1=(\r_1,\r_L)$ and $Q(U_1)=Q(\r_R;0)$: the left wave $(U_L,U_1)$ is a shock wave with negative wave speed $s(U_L,U_1)=(Q(U_L)-Q(U_1))/(\r_L-\r_1)$ and the right wave $(U_1,U_R)$ is a standing wave. As shown in \reffig{type_3}, the boundary flux $q(x=0,t>0)=Q(U_L)$.
\ei

\bfg
\bc\includegraphics[height=12cm] {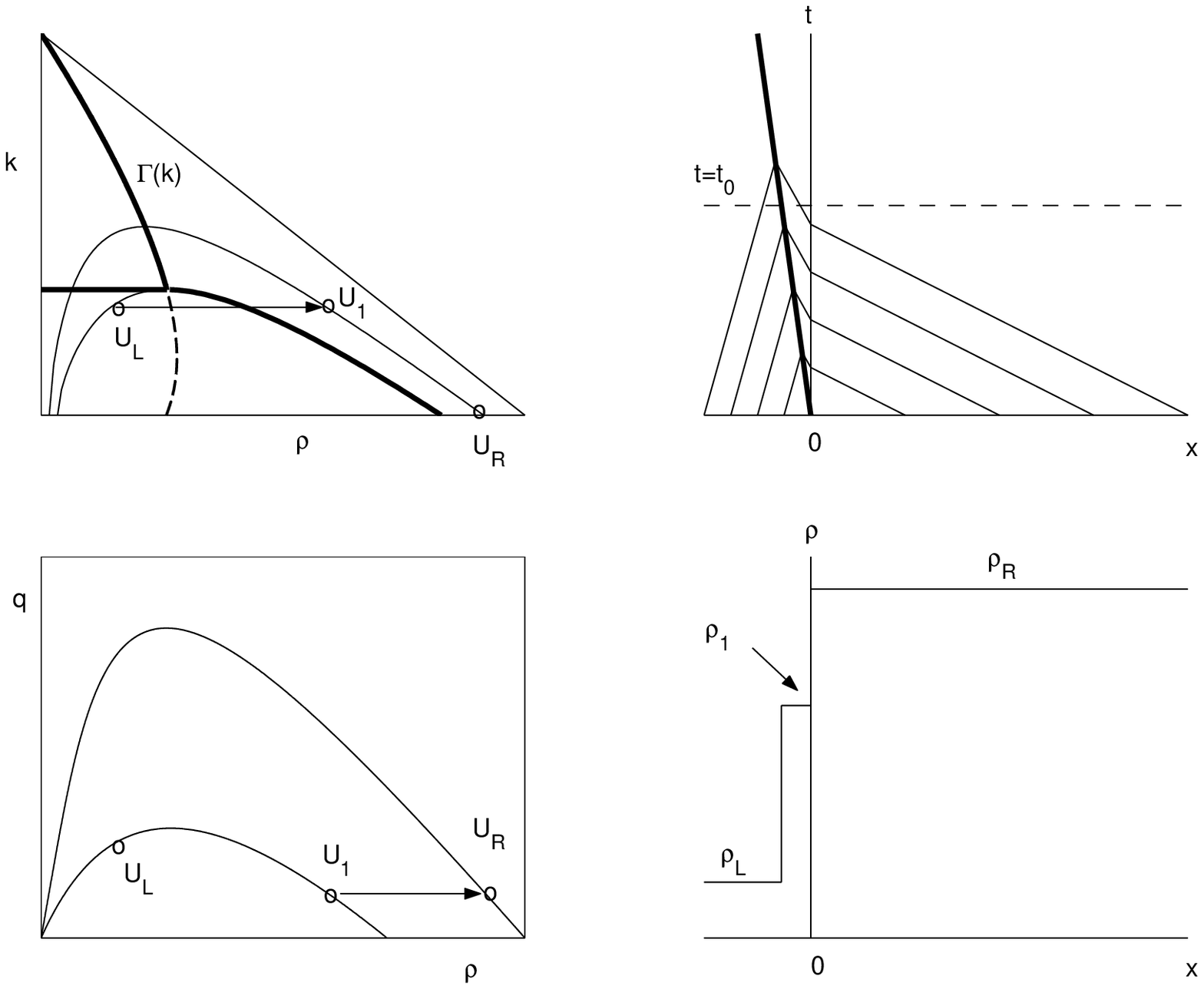}\ec
\caption{An example for wave solutions of type 3 for \refe{CL} with initial conditions \refe{initCL}}\label{type_3}
\efg

\item When $U_L$ is OC; i.e., $\r_L>\gamma(k_L)$, we define $U_\ast$ by $U_\ast\in\Gamma(k)$ and $k_\ast=k_L$. Therefore, $Q(U_\ast)$ is the partial capacity $Q^{\max}(k_L)$. Denoting the intersection between $Q(\r;k)=Q(U_L)$ and $(\r\geq \gamma(0), k=0)$ by $\r_0$ and the intersections between $Q(\r;k)=Q(U_\ast)$ and $k=0$ by $\r_1$ and $\r_2$, where $\r_1\leq \gamma(0)\leq \r_2$, we can obtain four types of solutions when $\r_R\in [0, \r_1]$, $(\r_1,\r_2)$, $[\r_2,\r_0]$, or $(\r_0,\r_j]$, as shown in \reffig{right_waves}.
\bi
\item [Type 4] When $\r_R\in [0, \r_1]$; i.e., $Q(\r_R;0)\leq Q(U_\ast)$ and $U_R$ is UC, wave solutions to the Riemann problem consist of three basic waves with two intermediate states, $U_\ast$ and $U_1$, where $Q(U_1)=Q(U_\ast)$ and $k_1=0$: the left wave $(U_L,U_\ast)$ is a rarefaction wave with non-positive characteristic velocity, the center wave $(U_\ast,U_1)$ is a standing wave, and the right wave $(U_1,U_R)$ is a rarefaction wave with positive characteristic velocity. As shown in \reffig{type_4}, the boundary flux $q(x=0,t>0)=Q(U_\ast)=Q^{\max}(k_L)$.

\bfg
\bc\includegraphics[height=12cm] {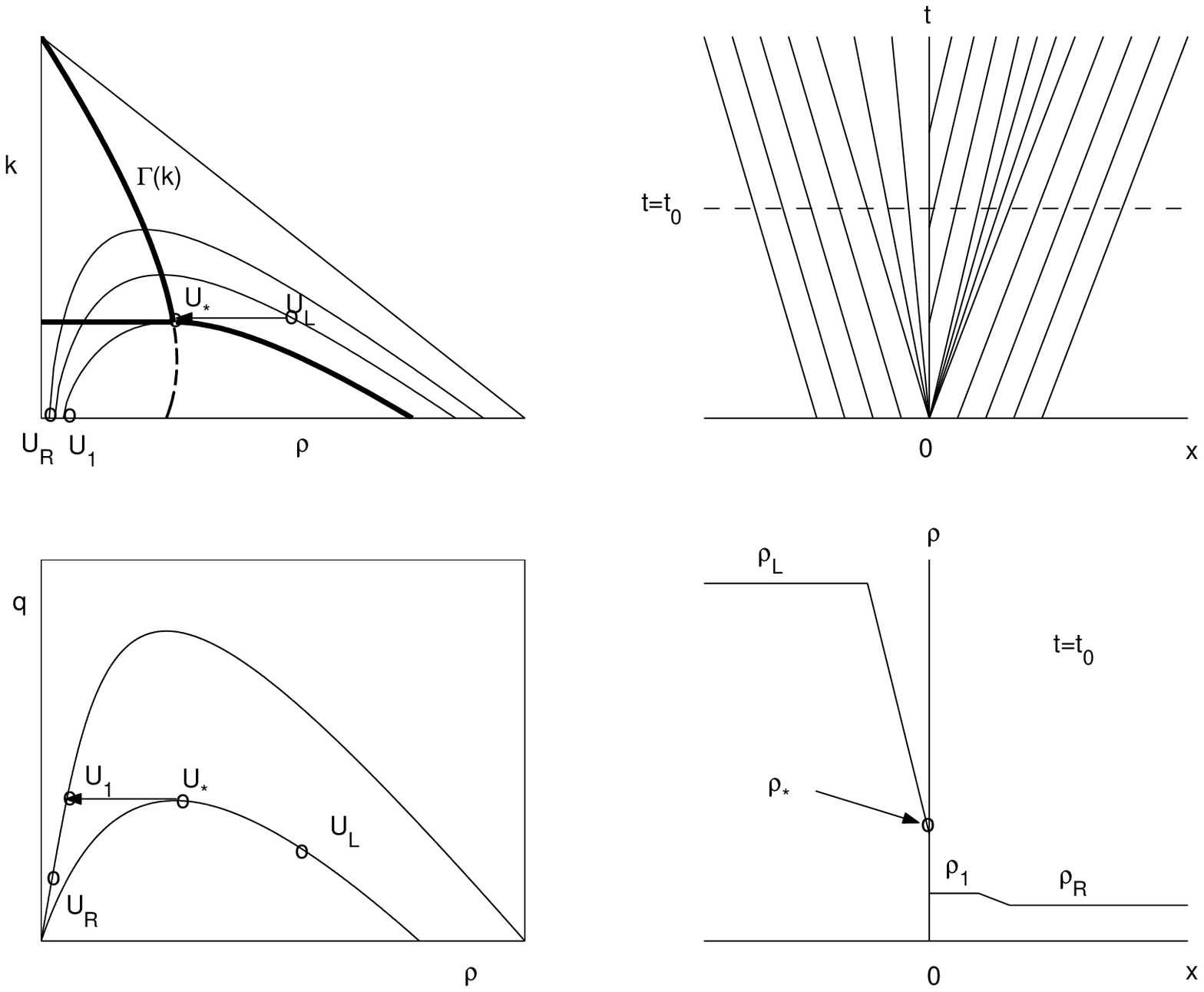}\ec
\caption{An example for wave solutions of type 4 for \refe{CL} with initial conditions \refe{initCL}}\label{type_4}
\efg

\item [Type 5] When $\r_R\in (\r_1, \r_2)$; i.e., $Q(\r_R;0)> Q(U_\ast)$, wave solutions to the Riemann problem consist of three basic waves with two intermediate states, $U_\ast$ and $U_1$, where $Q(U_1)=Q(U_\ast)$ and $k_1=0$: the left wave $(U_L,U_\ast)$ is a rarefaction wave with non-positive characteristic velocity, the center wave $(U_\ast,U_1)$ is a standing wave, and the right wave $(U_1,U_R)$ is a shock wave with positive speed. As shown in \reffig{type_5}, the boundary flux $q(x=0,t>0)=Q(U_\ast)=Q^{\max}(k_L)$.

\bfg
\bc\includegraphics[height=12cm] {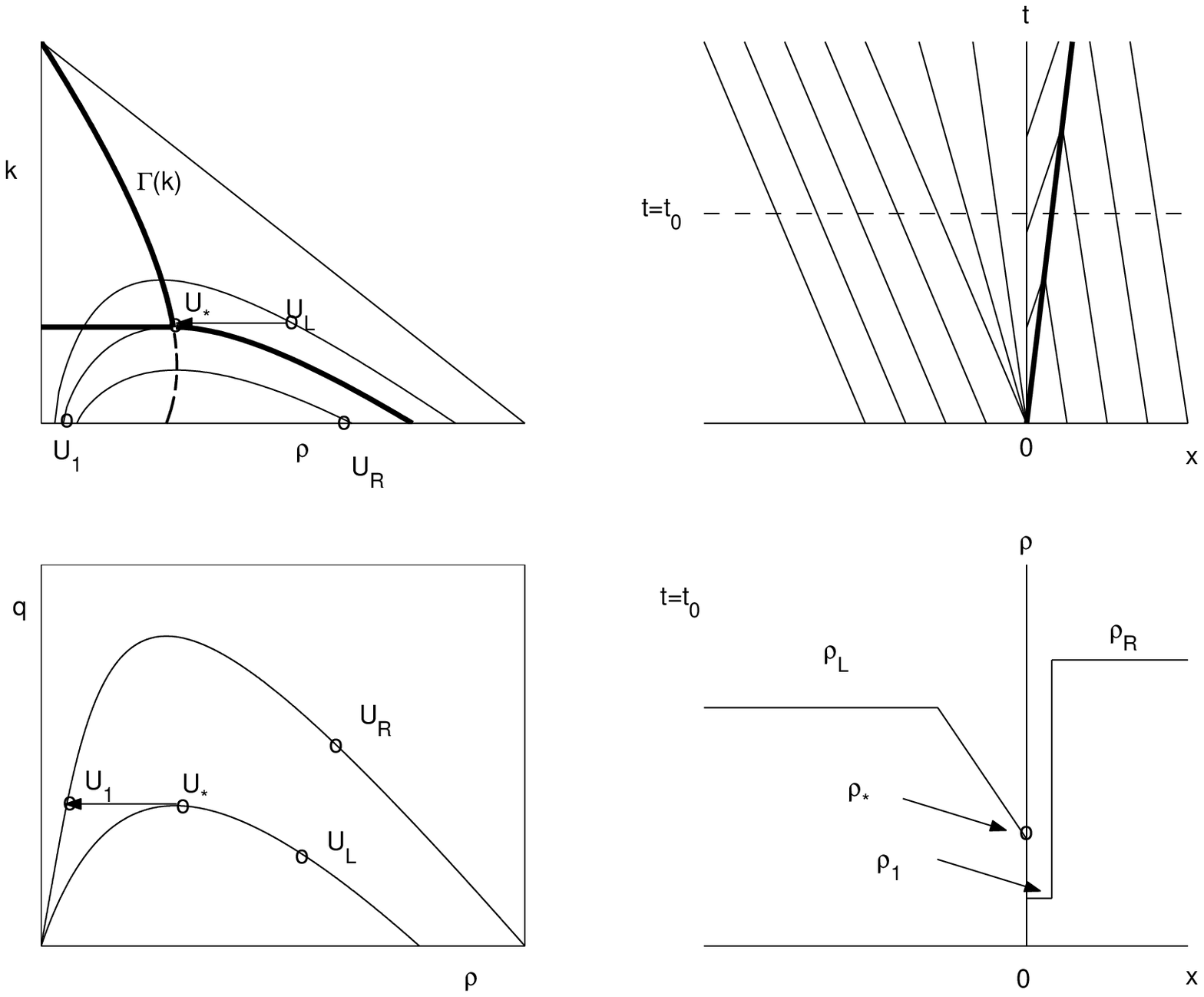}\ec
\caption{An example for wave solutions of type 5 for \refe{CL} with initial conditions \refe{initCL}}\label{type_5}
\efg

\item [Type 6] When $\r_R\in [\r_2, \r_0]$; i.e., $Q(U_L)\leq Q(\r_R;0)\leq Q(U_\ast)$ and $U_R$ is OC, wave solutions to the Riemann problem consist of two basic waves with the intermediate state $U_1$, where $Q(U_1)=Q(U_R)$ and $k_1=k_L$: the left wave $(U_L,U_1)$ is a rarefaction wave with negative characteristic velocity and the right wave $(U_1,U_R)$ is a standing wave. As shown in \reffig{type_6}, the boundary flux $q(x=0,t>0)=Q(U_R)$.

\bfg
\bc\includegraphics[height=12cm] {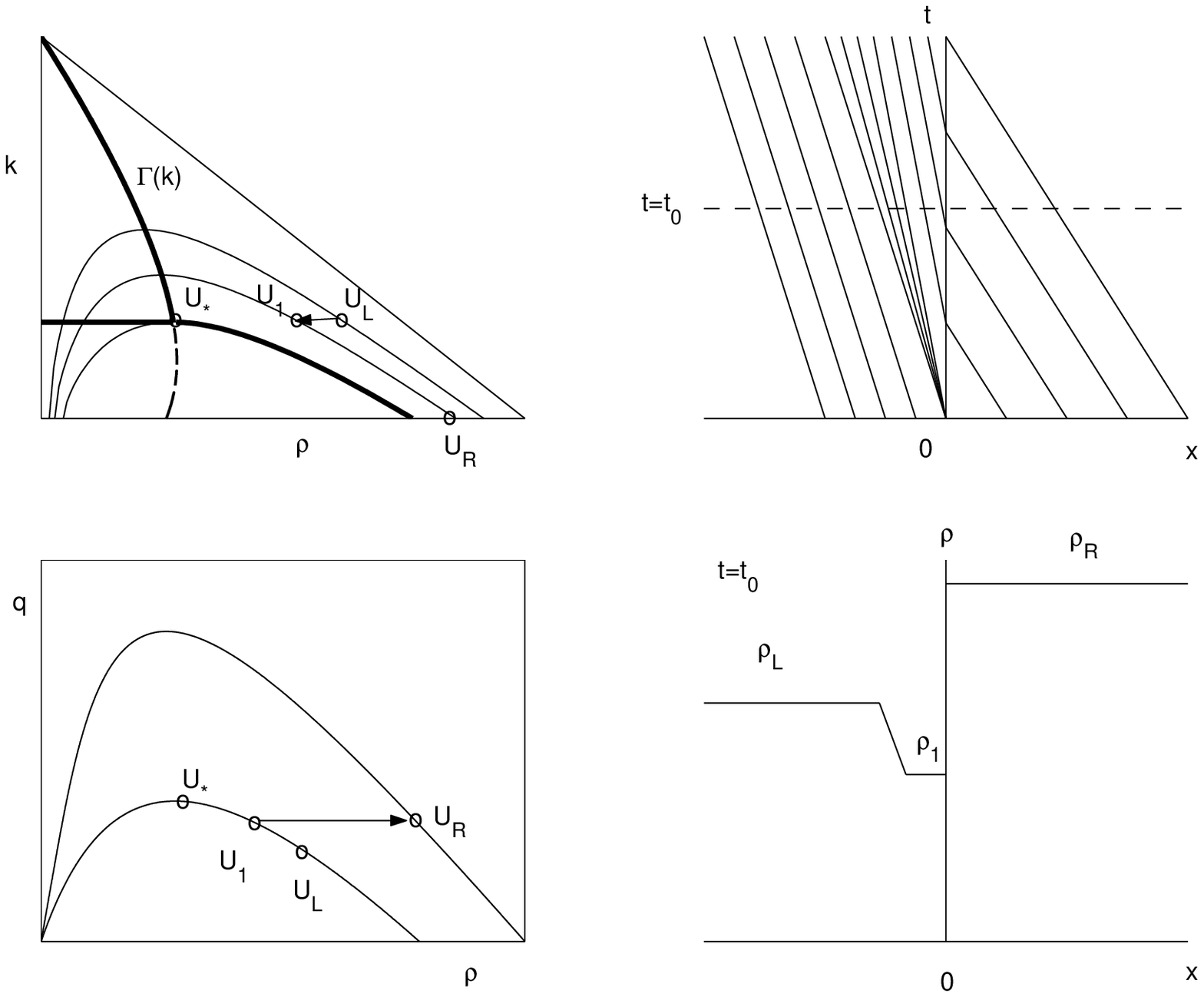}\ec
\caption{An example for wave solutions of type 6 for \refe{CL} with initial conditions \refe{initCL}}\label{type_6}
\efg

\item [Type 7] When $\r_R\in (\r_0, \r_j]$; i.e., $Q(\r_R;0)<Q(U_L)$ and $U_R$ is OC, wave solutions to the Riemann problem consist of two basic waves with the intermediate state $U_1$, where $Q(U_1)=Q(U_R)$ and $k_1=k_L$: the left wave $(U_L,U_1)$ is a shock wave with negative speed and the right wave $(U_1,U_R)$ is a standing wave. As shown in \reffig{type_7}, the boundary flux $q(x=0,t>0)=Q(U_R)$.
\ei

\bfg
\bc\includegraphics[height=12cm] {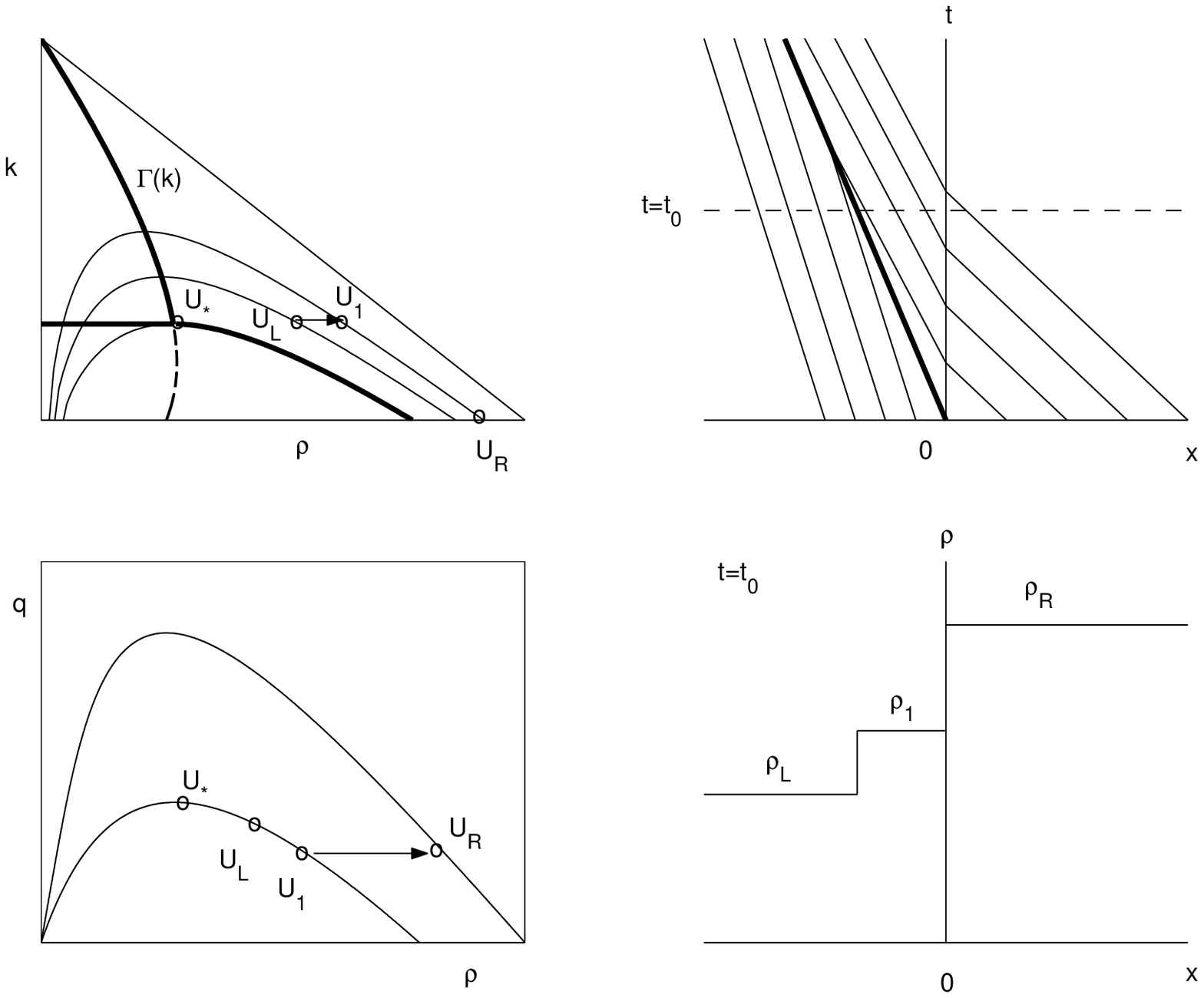}\ec
\caption{An example for wave solutions of type 7 for \refe{CL} with initial conditions \refe{initCL}}\label{type_7}
\efg

\een

%Wenlong Jin: 2003
\section{The supply-demand method with a new definition of traffic demand}\label{sec:diverge4}
Based on the discussions in the previous section, we summarize solutions of boundary flux $q(x=0,t>0)$ in \reft{table}.

\btb
\begin{tabular}{|c||l|l||c|}\hline
Solution type&left state $U_L$ & right state $\r_R$ &$q(x=0,t>0)$
\\\hline
1&{UC } &$Q(\r_R;0)\leq Q(U_L)$, $\r_R<\gamma(0)$ & $Q(U_L)$ \\\hline
2&{UC } &$Q(\r_R;0)> Q(U_L)$ & $Q(U_L)$ \\\hline
3&{UC } & $Q(\r_R;0)\leq Q(U_L)$, $\r_R>\gamma(0)$ & $Q(U_R)$ \\\hline
4&{OC } & $Q(\r_R;0)\leq Q(U_L)$, $\r_R>\gamma(0)$& $Q^{\max}(k_L)$ \\\hline
5&{OC } & $f(U_R)<f^{max}_L$, $a_R>a_L$, $\r_R/a_R<\alpha$ &$Q^{\max}(k_L)$ \\\hline
6&{OC } & $f(U_R)>f^{max}_L$ & $Q(U_R)$ \\\hline
7&{OC } & $f(U_R)<f(U_L)$, $\r_R/a_R<\alpha$, $a_R<a_L$ &
$Q(U_R)$ \\\hline
\end{tabular}
\caption {Solutions of the boundary flux $q(x=0,t>0)$}\label{table}
\etb

Further, if we introduce a new definition of partial traffic demand for commodity $1$ in the upstream link as
\bqn
D(\r_L;k_L)&=&\cas{{ll} Q(U_L)&U_L \m{ is UC} \\Q^{\max}(k_L) &  \m{otherwise}},\label{def_pd}
\eqn
the boundary flux $q(x=0,t>0)$ can then be computed by
\bqn
q(x=0,t>0)&=&\min\{S_1,D(\r_L;k_L)\},\label{demand-supply}
\eqn
where the supply of the downstream link, $S_1$, is the same as in \citep{daganzo1995ctm,lebacque1996godunov}; i.e.,
\bqn
S_1&=&\cas{{ll}Q(U_R)& U_R\m{ is OC}\\Q^{\max}(0) &\m{otherwise}}.
\eqn

The supply-demand method \refe{demand-supply} can also be used to calculate the boundary fluxes of other commodities. It yields the same solutions to the Riemann problem of \refe{CL} with \refe{initCL} as the analytical solution method. Moreover, it is much simpler in computation and can be easily extended to the more complicated cases when the upstream and downstream branches have different road characteristics.

%Wenlong Jin: 2003
\section{Numerical simulations}\label{sec:diverge5}
In this section, we carry out numerical simulations of the instantaneous kinematic wave model presented in this chapter. We will study a small diverging network consisting of two downstream links and one upstream link: the lengths of three links are the same, $L=400l=11.2$ km, with unit length $l=0.028$ km; one downstream link, labeled as link $d_1$, has $a(d_1)=2$ lanes, another downstream link, link $d_2$, has $a(d_2)=1$ lane, and the upstream link, link $u$, has $a(u)=2$ lanes. The simulation starts from $t=0$ and ends at $t=500\tau=41.7$ min, with unit time $\tau=5$ s. In the following simulations, we partition each link into $N$ cells and the time interval into $K$ steps, with $N/K=1/10$ always; e.g., if $N=50$  and $K=500$, the cell length  is $\dx=8 l$ km and the length of each time step $\dt=1 \tau$ s.

We will use the exponential fundamental diagram \citep{delcastillo1995fd_empirical}, and the parameters are given as follows: the free flow speed $v_f=5.0 l/\tau$ = 0.028 km/s = 100.8 km/h; the jam density of a single lane $\r_j$ = 180 veh/km/lane; the wave velocity for jam traffic $c_j=-1.0 l/\tau$ = -0.0056 km/s = -20.16 km/h;  and the equilibrium speed-density relationship
\bqs
V(a,\rho)=5\left [1- \exp\left\{\frac{1}{5} (1-\frac{a\rho_j}{\r})\right \}\right] l/\tau,
\eqs
where $a(x)$ is the number of lanes at location $x$.
The equilibrium functions $V(a,\r)$ and $Q(a,\r)$ are shown in \reffig {fd_Newell}, in which the critical traffic density is $\alpha=0.259 a \r_j$.

\bfg
\bc\includegraphics[height=12cm] {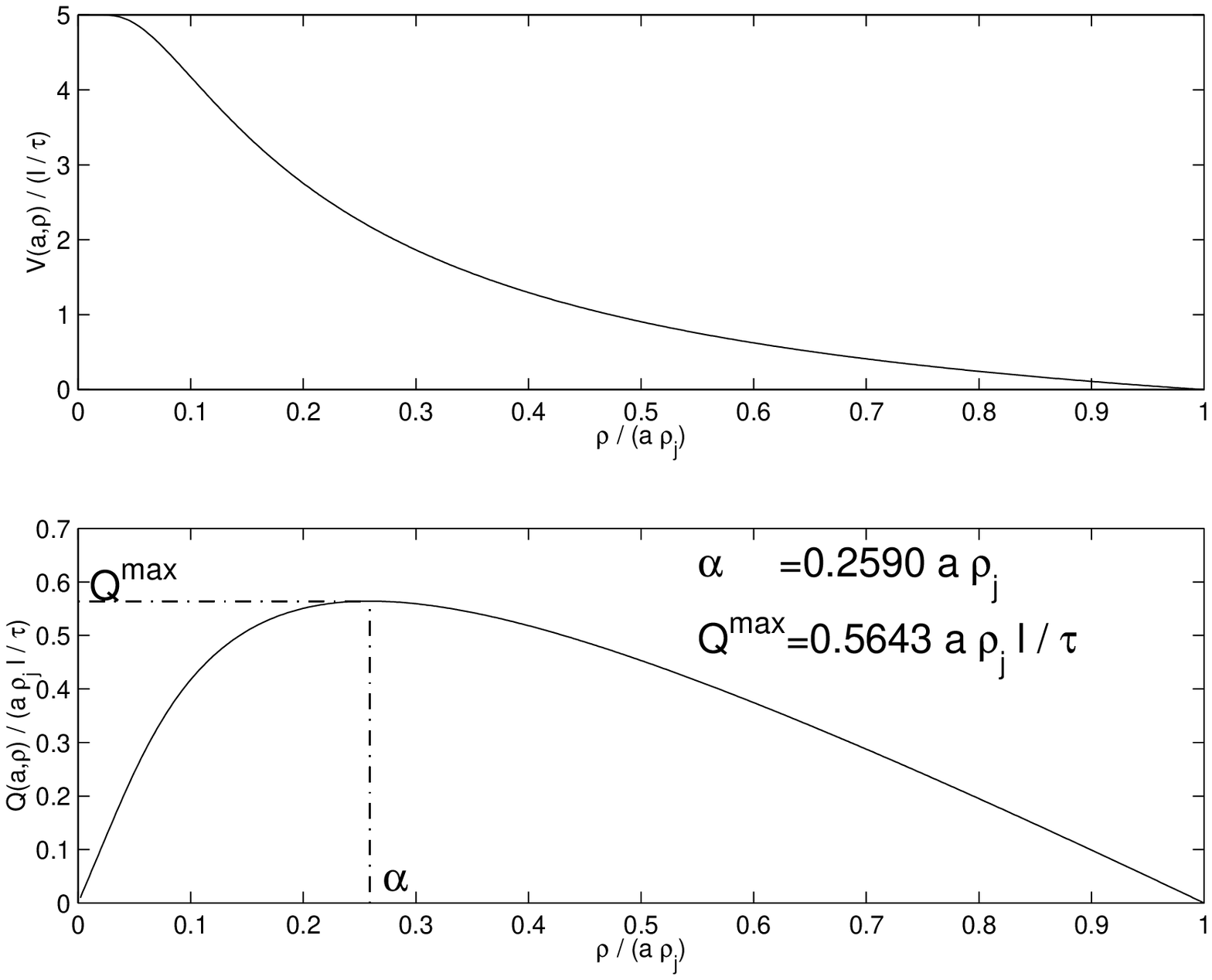}\ec
\caption {The Newell model of speed-density and flow-density relations} \label {fd_Newell}
\efg

Here we apply the first-order Godunov method \citep{godunov1959}, in which the traffic conditions are updated for each cell based on traffic conservation, and boundary fluxes are computed with the supply-demand method described in Section 4. Since $\l_{\ast}\leq v_f=5 l/\tau$, we find the CFL condition number
\bqs
\l_{\ast} \frac {\dt}{\dx}&\leq& 0.625<1.
\eqs
Therefore, the first-order Godunov method can solve the instantaneous kinematic wave model efficiently.

\subsection{Simulation I: A general case}
In this simulation, we study a general case of diverging traffic. Initially, the upstream branch carries a constant flow with traffic density $\r_u=1.1111 \r_j=200$ veh/km and the proportion of vehicles traveling to downstream branch $d_1$ is $80\%$; downstream link $d_1$ is empty; traffic density on downstream link $d_2$ is $\r_{d_2}=0.5556\r_j=100$ veh/km. In addition, we impose the Neumann boundary conditions on the boundary of this  diverging  network; i.e., the spatial derivatives of traffic densities at the boundary are set to be zero.

With $N=500$ and $K=5000$, we obtain simulation results as shown in \reffig{sim_1}. Fig. \ref{sim_1}(a) illustrates the evolution of traffic on the upstream link: at time $t=0\tau$, traffic density is uniformly at $\r_A=1.11\r_j$; after the beginning of diverging process, traffic immediately upstream of the diverging point reaches a new state $\r_B=0.69\r_j$, which keeps propagating on the upstream link; as a result, an expansion wave forms and travels upstream. Fig.\ref{sim_1}(b) shows the evolution of traffic on downstream link $d_1$: initially, this link is empty; after $t=0$, traffic immediately downstream of the diverging point reaches state $\r_C=0.22\r_j$; along with the propagation of $\r_C$, on this link, another expansion wave forms and travels downstream; after around $200\tau$, traffic density on this link is uniformly $\r_C$. Fig. \ref{sim_1}(c) presents the evolution of traffic on downstream link $d_2$: initially, traffic density is $\r_D=0.56\r_j$; after vehicles start to diverge at $t=0$, traffic immediately downstream of the diverging point reaches $\r_E=0.04\r_j$; then a shock forms and travels downstream in a constant speed $s=0.38 l/\tau$. The expansion waves and the shock wave observed on the three branches can be shown on the $\r-q$ plane as in Fig. \ref{sim_1}(d), in which dashed line $AB$ represents the expansion wave on the upstream link $u$, dashed line $OC$ represents the expansion wave on downstream link $d_1$, and solid line $DE$ represents the shock wave on downstream link $d_2$, whose slope gives the shock speed.

\bfg
\bc\includegraphics[height=12cm]{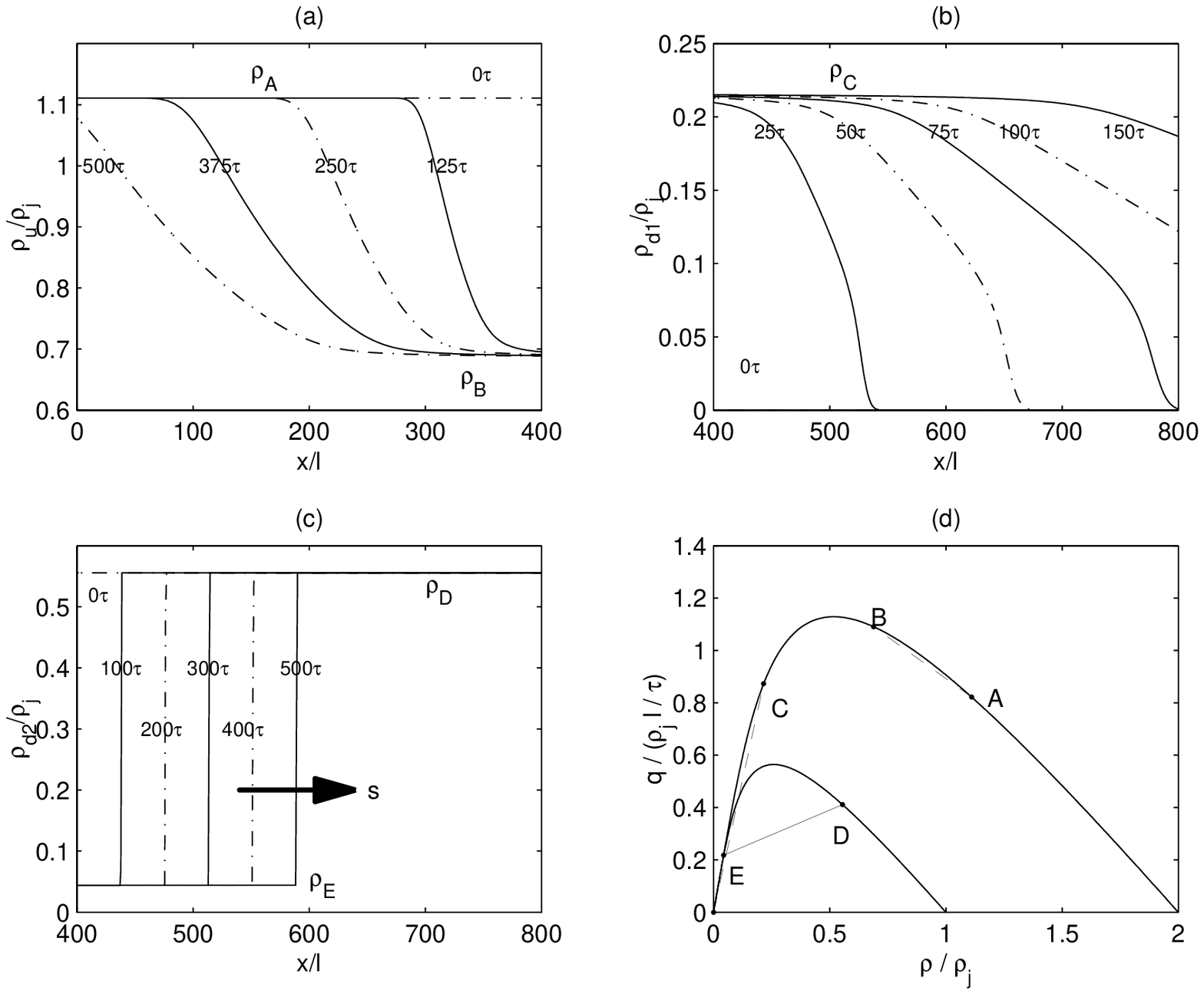}\ec
\caption{Simulation I: A general case}\label{sim_1}
\efg

As shown in \reffig{sim_1}, at the diverging point, there are three traffic states, represented by $B$, $C$, and $E$ in Fig. \ref{sim_1}(d). Flow rates at these states are $q_B=1.09 \r_jl/\tau$, $q_C=0.87\r_jl/\tau$, and $q_E=0.22 \r_jl/\tau$ respectively. We can see that $q_B=q_C+q_E$. Note that $q_B$ is the outflow of link $u$, $q_C$ inflow to link $d_1$, and $q_E$ inflow to link $d_2$. Thus traffic is conserved at the diverge. Further, $q_C/q_B=80\%$, which is the proportion of vehicles on the upstream link traveling to $d_1$. Therefore, this is consistent with a general observation that diverging flows are proportionally determined by the composition of traffic on the upstream link \citep{papageorgiou1990assignment}. This property of diverging flows guarantees that the composition of vehicles on the upstream link never changes as observed in our simulation.

\subsection{Simulation II: An extreme case}
In this subsection, an extreme case is studied. Initially, like in the previous simulation, the upstream branch carries a constant flow with traffic density $\r_u=1.1111 \r_j=200$ veh/km, the proportion of vehicles on the upstream link traveling to downstream branch $d_1$ is $80\%$, and downstream link $d_1$ is empty. However, traffic on downstream link $d_2$ is jammed; i.e., $\r_{d_2}=\r_j=180$ veh/km. Still, we impose the Neumann boundary conditions on the boundary of the diverging  network and have the same discretization to the three links and the time duration: $N=500$ and $K=5000$.

Simulation results are shown in \reffig{sim_2}. In this simulation, traffic density on link $d_2$ is uniformly $\r_j$ as expected and not shown. Fig. \ref{sim_2}(a) demonstrates traffic evolution on the upstream link $u$: initially traffic density is $\r_A=1.11 \r_j$; after the beginning of the diverging process, traffic density immediately upstream of the diverging point reaches the jammed density $\r_B=2\r_j$; then, jammed traffic propagates upstream as a back-traveling shock at a speed $s=-0.92 l/\tau$. Fig. \ref{sim_2}(b) shows traffic dynamics on link $d_1$. In this figure, we can observe small spikes of density, which travel downstream and shrink along time. After around $t=200\tau$, link $d_1$ is almost empty as initially. Fig. \ref{sim_2}(d) presents traffic states on the diverging network on the $\r-q$ plane. Here line $AB$ represents the shock wave on the upstream link, and its slope gives the shock wave speed.

\bfg
\bc\includegraphics[height=12cm]{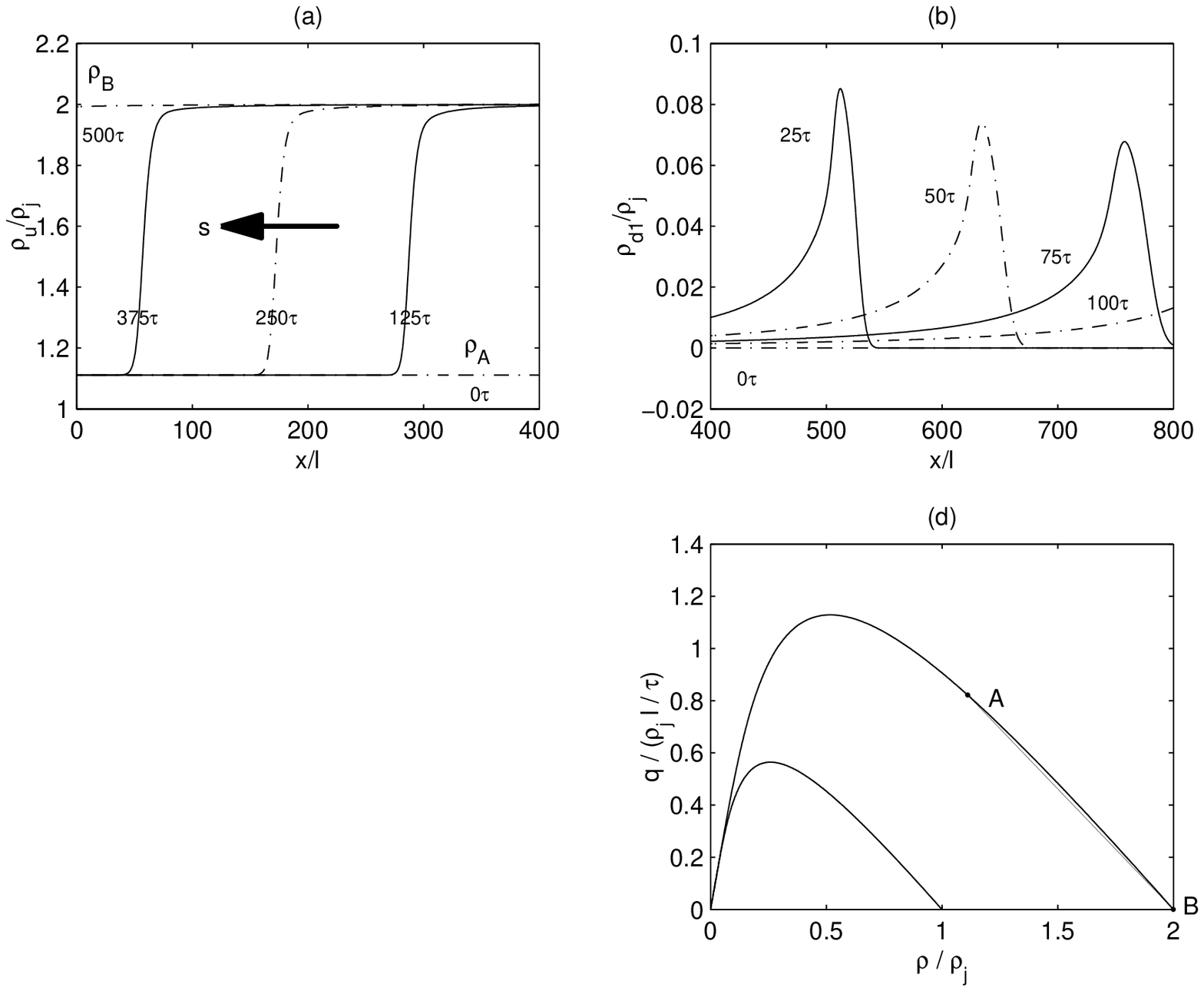}\ec
\caption{Simulation II: An extreme case}\label{sim_2}
\efg

From these figures, we can see that, if a downstream link is blocked, the upstream link and the other downstream link(s) will be blocked, as expected. However, the small spikes in Fig. \ref{sim_2}(b) suggests there is still a small number of vehicles get out of the upstream link at the very beginning. Although this phenomenon seems interesting, the appearing of small spikes in our model is due to numerical error of finite difference, caused by the finite number of cells, and are expected to disappear if we partition the links into fine enough cells.

%Wenlong Jin: 2003
\section{Discussions}\label{sec:diverge6}
In this chapter, we have introduced a new approach to modeling diverging traffic dynamics, which is represented by  instantaneous kinematic waves. We presented analytical solutions to these waves and found these solutions are consistent with the supply-demand method if the traffic demand is modified as in Section \ref{sec:diverge4}. With numerical simulations, we assert that this model satisfies two important properties: (1) diverging flow to a downstream link is proportional to the fraction of vehicles traveling to this link, and (2) the upstream link is blocked after a downstream link is blocked.

The instantaneous kinematic wave theory of diverging traffic is different from existing models since it provides more details on diverging traffic dynamics and sheds more light at the aggregate level how vehicles segregate themselves at a diverge. In the follow-up studies, we would like to discuss the convergent kinematic wave solutions when all branches are partitioned into infinite number of cells. We will also be interested in finding the location of actual diverging point and how it is related to traveler's behavior, the geometry of a diverge, and traffic conditions. Finally we will test this model with field data and discuss possible applications.

\newpage
\pagestyle{myheadings}
\markright{  \rm \normalsize CHAPTER 5. \hspace{0.5cm}
 MIXED TRAFFIC MODEL}
\chapter{Kinematic wave traffic flow model for mixed traffic}
%\thispagestyle{myheadings}
%Wenlong Jin: 2003
%Wenlong Jin: 2003
\section{Background}
\label{intro}

Vehicular traffic on highways often comprises different types of
vehicles with varying driving performances. This heterogeneity
affects traffic flow characteristics in significant ways, a fact that
has long been recognized by the transportation engineering
profession. For example, in the computation of flow capacity on a
highway or at a signalized intersection, the Highway Capacity Manual
recommends a series of adjustments to take account of the capacity
reduction caused by heavy vehicles (i.e., trucks/buses/recreational
vehicles). If one is interested in the effects of heavy vehicles on
traffic flow over space and time, however, the Highway Capacity
Manual procedures are not adequate. For this one needs a
dynamic model for mixed traffic.

Mixed traffic can be modeled at three different levels---microscopic,
mesoscopic and macroscopic. It is perhaps most straightforward to
model mixed traffic on a microscopic level---one simply endow, at one
extreme, each individual vehicle with different performance and
behavior characteristics. Many commercially available simulation
packages, such as
CORSIM, PARAMICS, and VISSIM allow the specification of
multiple vehicle classes. Major challenges arise when one models mixed
traffic on a mesoscopic level, mainly due to the correlation between
various probability distributions of vehicular speeds. Nevertheless,
a number of mesoscopic models of mixed traffic have been developed in
recent years \citep{helbing1997,hoogendoorn2000multiclass}. Aggregation of mesoscopic
models of mixed traffic through expectation operations lead to
multi-class traffic flow models in the macroscopic level. There is,
however, another approach to develop macroscopic mixed traffic
flow  models. This is the approach of continuum modeling. It is this
approach that we
shall adopt in developing a traffic flow model for mixed traffic.

In the continuum description of traffic flow, vehicular traffic is
described as a special kind of fluid that are characterized by its
concentration (density, $\rho$),  mean velocity ($v$), and
vehicle flux (flow rate $q$), all are functions of space ($x$) and
time
($t$).
The starting point of any  continuum model of traffic flow is the
conservation
of vehicles
    \begin{equation}
            \frac{\partial}{\partial
t}\int_{x_{1}}^{x_{2}} \rho(x,t)dx =q
            (x_{1},t)-q (x_{2},t),
            \label{conserv1}
    \end{equation}
and the relation between flow, density, and mean velocity $ q=kv$.

\refe{conserv1} is an integral form
of traffic conservation. When the road segment $[x_{1}, x_{2}]$
shrinks to a point in space, one obtains the familiar
differential  form of traffic conservation:
\begin{equation}
\rho_{t}+q_{x}=0, \; \hbox{or} \; \rho_{t}+(\rho v)_{x}=0.
\label{differential}
\end{equation}
If one introduces a relation between vehicle concentration and traffic
speed $v=V(\rho)$, one obtains the classic kinematic wave model
developed by Lighthill, Whitham \citep{lighthill1955lwr} and Richards
\citep{richards1956lwr}:
\beq
\rho_{t}+ (\rho V(\rho))_{x}=0, \; \rho V(\rho) \equiv Q(\rho).
\label{kw}
\eeq

The classic kinematic wave model of Lighthill and Whitham  was
formulated
for homogeneous flows on a long crowded road.  It does not consider
the
effects of performance differences among different types of vehicles.
Recently,
Daganzo extended the theory to treat a freeway system with two
types of lanes, special lanes and regular lanes, and two types of
vehicles, priority vehicles and regular
vehicles~\citep{daganzo1997special,daganzo1997it,daganzo2002behavior}. Priority
vehicles are
allowed to travel on either regular lanes or special lanes, whereas
regular
vehicles can only travel on the regular lane. The two types of
vehicles in Daganzo's special lane model have different vehicle
performances in free-flow traffic, but are indistinguishable in heavy
traffic, where both types of vehicles travel at the same speed.

In this chapter, we extend the kinematic wave model to vehicular
traffic with a mixture of vehicle types. In the mixed flow each
vehicle type is conserved and travels at the group velocity, but the
differences among vehicle types are accounted for in determining the
states of the collective flow.  This model can be used to study
traffic
evolution on long crowded highways where low performance vehicles
entrap high performance ones. It can also give a more accurate
description of the I-pipe state in Daganzo's special lane model.

The remaining parts of the chapter are organized as follows. In Section
\ref{sec:mixed2} we give the extended KW model and its basic properties. In Section
\ref{sec:mixed3} we analyze the Riemann problem for this model. In Section \ref{sec:mixed4} we
propose a fundamental diagram of mixed traffic and discuss its
properties. In Section \ref{sec:mixed5} we provide numerical examples and in Section
\ref{sec:mixed6} we conclude the chapter.

%Wenlong Jin: 2003
\section{The extended KW model for mixed traffic}\label{sec:mixed2}

Let us assume that there are $i=1, \cdots, n$ types of
vehicles in the traffic stream ($n \ge 2$), each type has
concentration $\rho_{i}(x,t)$ and velocity $v_{i}(x,t)$. By
conservation
of each vehicle type, we have
\beq
(\rho_{i})_{t}+(\rho_{i}v_{i})_{x}=0, \; i=1, \cdots, n;
\label{gmkw}
\eeq

As in the development of the classic KW model, we postulate that
equilibrium relations exist between vehicular speeds and traffic densities:
\beq
     v_{i}=V_{i}(\rho_{1}, \cdots, \rho_{n}),
\label{kv}
\eeq
with $v_{i}(0) =v_{fi}$, the free-flow speed of each vehicle type and
$\partial_{\r_{j}}V_{i} <0, \; i=1, \cdots, n; \; j=1, \cdots, n$.

\refet{gmkw}{kv} are the general governing
equations of mixed traffic flow without special lanes. Note that if
one adopts $v_{i}=V(\sum_{i=1}^{n}\rho_{i})$, one recovers the I-pipe
state in Daganzo's special lane model. In this chapter we study a
special case of the general equations for mixed traffic in which we
consider two types of vehicles---one represents passenger cars
($\rho_{1}$) and the other represents heavy vehicles such as trucks
($\rho_{2}$), and two traffic flow regimes---free-flow and  congested traffic.

When traffic is light and there are adequate opportunities for
passing, different classes of vehicles would travel at their own
free-flow speeds $v_{fi}$. The traffic flow in this case can
be described by
\beq
(\rho_{i})_{t}+ v_{fi} (\rho_{i})_{x}=0, \; i=1, 2;
\label{gmkw-2}
\eeq
By defining $\rho=\sum_{i=1}^{n} \rho_{i}$, and
$v_{f}=\frac{\sum_{i=1}^{n} \rho_{i}v_{fi}}{\sum_{i=1}^{n}
\rho_{i}}$, we can use
\beq
\rho_{t} +  (\rho v_{f})_{x}=0,
\label{average}
\eeq
to approximately model the average behavior of light traffic.

When traffic concentration reaches a critical value $\rho_{c}$,
passing opportunities diminish and vehicles of lower performance (e.g.
trucks) start to entrap vehicles of higher performance (e.g.,
passenger cars).  Under such conditions it is assumed that
the various classes of vehicles are completely mixed and move
at the group velocity $V$. That is, mixed traffic flow  in this
regime is described by
\beq
(\rho_{i})_{t}+  (\rho_{i} V)_{x}=0, \; i=1, 2.
\label{gmkw-3}
\eeq

Through the definition of a proper average free-flow speed and the
selection of a suitable critical density, we  combine \refe{gmkw-2}
(for free-flow traffic) and
\refe{gmkw-3} (for congested traffic) into one modeling equation:
\beq
\left(
\ba{l}
       \rho_{1} \\ \rho_{2}
\ea
\right)_{t} +
\left(
\ba{l}
       \rho_{1} V(\rho_{1},\rho_{2}) \\  \rho_{2} V(\rho_{1},\rho_{2})
\ea
\right)_{x}  =0,
\label{mkw}
\eeq
where
\[
V(\rho_{1},\rho_{2}) =\left\{
\ba{ll}
\frac{\sum_{i=1}^{2} \rho_{i}v_{fi}}{\sum_{i=1}^{2}
\rho_{i}}, &   \gamma_{1} \rho_{1} +  \gamma_{2}\rho_{2}  < 1\\
V_{*} (\rho_{1},\rho_{2}), &   \gamma_{1} \rho_{1} +  \gamma_{2}\rho_{2}\ge 1
\ea
\right. .
\]
Here $\gamma_{1}$ and $\gamma_{2}$ are parameters that determine the
critical density in $(\rho_{1},\rho_{2})$ coordinates. $V_{*}$ is a
two-dimensional speed-density relation for congested traffic. It is
understood that $\partial_{\r_{i}}V_{*} <0, \; i=1, 2$.

\refe{mkw} is a system of conservation laws with characteristic velocities:
\[
\lambda_{1}=V + \rho_{1} V_{\dot{1}} + \rho_{2}V_{\dot{2}}, \;
\lambda_{2} =V(\rho_{1}, \rho_{2}).
\]
Here we used a special notation for partial derivatives of $V$ with
respect to $\rho_{1}$ and $\rho_{2}$ : $\partial_{\rho_{1}} V \equiv
V_{\dot{1}}$ and $\partial_{\rho_{2}} V \equiv V_{\dot{2}}$. Owing to
the nature of $V(\rho_{1},\rho_{2})$, we have $\lambda_{1} \le
\lambda_{2} = V$, that
is, both characteristics travel no faster than \textit{average
traffic}. In fact, the second
characteristic travels at precisely the speed of traffic.
When the free-flow speeds of both types of vehicles are identical,
the extended KW model preserves the anisotropic property of the KW
model. Otherwise, the extended model is not anisotropic in light
traffic (the nature of this violation of anisotropy is explained in
detail in \citep{zhang2000anisotropic}). Moreover, it can be shown that
\[
\left(\frac{\rho_2}{\rho_1}\right)_t + V
\left(\frac{\rho_2}{\rho_1}\right)_x=0,
\]
from which one obtains
$\frac{d}{dt}\left(\frac{\rho_2}{\rho_1}\right)=0$, that is, the
level curves of $\left(\frac{\rho_2}{\rho_1}\right)$ in the $t-x$
plane coincide with
vehicle trajectories.  The separation of
$\left(\frac{\rho_2}{\rho_1}\right)$ level
curves therefore implies first-in-first-out traffic flow behavior
between vehicle
classes.

Furthermore, the corresponding eigenvectors of the flow Jacobian
matrix are
\[
r_{1} = \left( \ba{c} \frac{\rho_{1}}{\rho_{2}} \\ 1 \ea \right), \;
r_{2} = \left( \ba{c} - \frac{V_{\dot{2}}}{V_{\dot{1}}} \\ 1 \ea
\right)
\]
and the Riemann invariants $(w,z)$, defined as $\bigtriangledown w
\bullet r_{1}
=0$, $\bigtriangledown z \bullet r_{2} =0$, are
\[
w=\frac{\rho_{2}}{\rho_{1}}, \; z=V.
\]
They are used here to obtain the expansion wave solutions of a
Riemann  problem (see next section. For more details on Riemann
problems and Riemann invariants, refer to \citep{whitham1974PW}).

It can be shown that the first characteristic field is nonlinear and
the second characteristic field is linearly degenerate. We therefore
have both shock and smooth expansion waves in the first field and
contact waves (or slips) in the second field. We shall derive the
expressions for these waves related to Riemann data in the next
section.

%Wenlong Jin: 2003
\section{The Riemann problem and basic wave solutions}\label{sec:mixed3}
\label{Riemann}
In this section we discuss the solutions of the extended KW model,
\refe{mkw}, with the following so-called Riemann data:
\beq
\rho(x,0)=\left\{
    \ba{ll}
     \rho^{l}, &   x<0 \\
     \rho^{r}, &   x>0
    \ea
\right.
\rho=\left(
   \ba{l}
       \rho_{1} \\ \rho_{2}
   \ea
   \right)
\eeq
To solve the above Riemann problem, we first study the right
(downstream) states that can be connected to the left (upstream)
states by an elementary wave, i.e., a smooth expansion (rarefaction)
wave, a contact, or a shock (readers are referred to
\citep{leveque2002fvm} and \cite{aw2000arz,zhang2000structural,zhang2002arz,zhang2001perspective} for a
more detailed discussion of Riemann problems
related to systems of conservation laws in general and traffic flow
in particular). Throughout the remaining sections, we assume
that $v_{f1}=v_{f2}=v_f$. This assumption ensures that our proposed
$V(\rho_{1},\rho_{2})$ function is continuous over the entire feasible
$(\rho_{1},\rho_{2})$ region. The Riemann problem of \refe{mkw} with
discontinuous $V(\rho_{1},\rho_{2})$ is more involved and will be
discussed elsewhere. Nevertheless, the
analysis of the model, \refe{mkw}, with $v_{f1}=v_{f2}=v_f$  still reveals many key
features of mixed traffic flow.

\vspace{2mm}
\noindent {\bf The 1-expansion waves}:
An upstream state $\rho^{l}$ can be connected to a downstream state
$\rho^{r}$ by a 1-expansion wave if and only if the downstream state
satisfies
\[
w(\rho^{l})=w(\rho^{r}), \; \rho^{l} > \rho^{r},
\]
i.e.,
\beq
\frac{\rho_{2}^{l}}{\rho_{1}^{l}} = \frac{\rho_{2}^{r}}{\rho_{1}^{r}}.
\label{1-r}
\eeq
This means that in the $\rho-$plane the two states are
on a ray from the origin. Clearly across an expansion wave traffic
composition does not change, that is, vehicles observe the
first-in-first-out rule.

\vspace{2mm}
\noindent {\bf The contact waves}:

A contact wave is a slip that separates two traffic regions of
different traffic densities and vehicle compositions but the same travel speed.
That is,
\beq
V(\rho^{l})=V(\rho^{r}).
\label{contact}
\eeq
In the $\rho-$plane, all states on a level curve of
$V(\rho)$ are connected by a contact wave.

\vspace{2mm}
\noindent {\bf The shock waves}:

The shock waves in the extended KW model are given by the jump
condition:
\bqn
    s (\rho_{1}^{l} -\rho_{1}^{r}) &=& \rho_{1}^{l}
    V(\rho_{1}^{l},\rho_{2}^{l}) - \rho_{1}^{r}
    V(\rho_{1}^{r},\rho_{2}^{r })  \\
      s (\rho_{2}^{l} -\rho_{2}^{r}) &=& \rho_{2}^{l}
    V(\rho_{1}^{l},\rho_{2}^{l}) - \rho_{2}^{r}
    V(\rho_{1}^{r},\rho_{2}^{r })
\eqn
After elimination of $s$ from the equations and some algebraic
manipulations one obtains
\[
(\rho_{1}^{l}\rho_{2}^{r} - \rho_{1}^{r} \rho_{2}^{l})(
V(\rho_{1}^{l},\rho_{2}^{l}) - V(\rho_{1}^{r},\rho_{2}^{r}))=0.
\]
Two possibilities exist: $ V(\rho_{1}^{l},\rho_{2}^{l}) -
V(\rho_{1}^{r},\rho_{2}^{r}) =0 $ which gives the contact waves that
we
have already discussed, or
\beq
\rho_{1}^{l}\rho_{2}^{r} - \rho_{1}^{r} \rho_{2}^{l}=0,
\label{1-s}
\eeq
this gives the downstream states $\rho^{r}$ that can be connected to
the
upstream state $\rho^{l}$ by a shock. Note that all these states also
fall on a ray originating from the origin of the $\rho-$plane. This
implies that across a shock vehicle composition also does not change,
that
is, vehicles observe first-in-first-out rule.  Moreover, we have the
following entropy conditions
\[
\rho^{l} < \rho^{r}.
\]
to ensure the stability of the shock.

Now we can state the procedure to solve a Riemann problem for the
extended KW model. Note that for any state $\rho^{l}$, the two
curves/lines given by \refet{1-r}{contact} divide the
feasible $\rho-$plane into four regions (\reffig{phase}). If the downstream
state $\rho^{r}$ falls on any of these two curves/lines, it can be
connected to the upstream state by an elementary wave. If it falls on
any of the four regions, however, an intermediate state $\rho^{m}$ is
generated on the line given by \refe{1-r}, which is connected with the
upstream state by a 1-wave (i.e., an expansion
or shock wave) and with the
downstream state by a contact (\reffig{phase}. \reffig{wave_sol} shows a few examples
of Riemann
solutions.

\bfg
\begin{center}
\includegraphics[height=12cm]{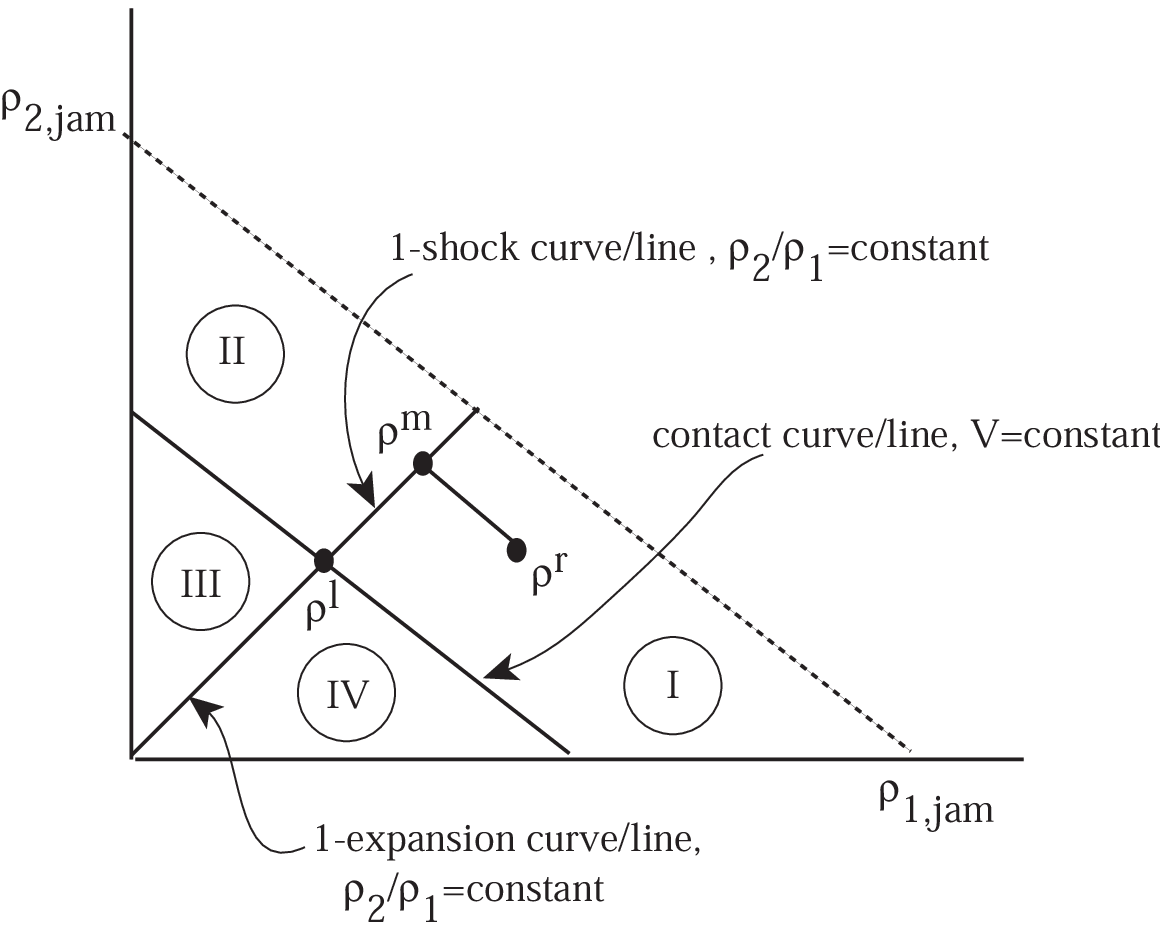}
\caption{Phase diagram for determining elementary and simple
waves}
\label{phase}
\end{center}
\end{figure}

\bfg
\begin{center}
\includegraphics[height=12cm]{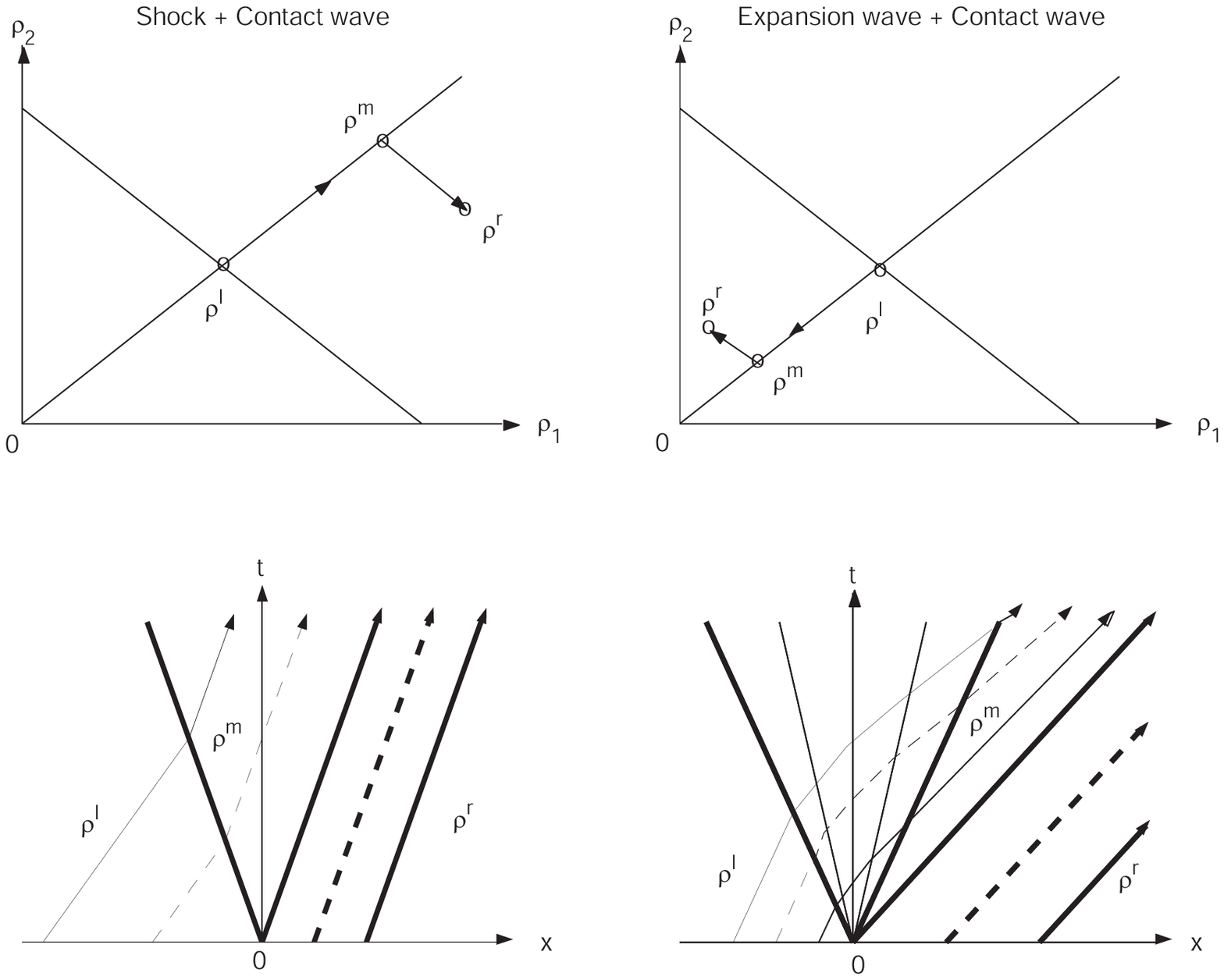}
\caption{Wave solutions to the Riemann problem: Shock + Contact wave
(left)
and Expansion wave + Contact wave (right) (In the bottom figures,
thick
(dashed) lines are characteristics, and lines with arrows are
vehicles'
trajectories.)} \label{wave_sol}
\end{center}
\end{figure}

%Wenlong Jin: 2003
\section{Fundamental diagrams for mixed traffic}\label{sec:mixed4}
We propose the following $\rho-V$ relation, which can be derived from
a car-following model under steady-state conditions
\citep{zhang2000mixed}, to be used in the mixed traffic flow model.
\[
V=\left\{
\ba{ll}
\frac{ \rho_{1}v_{f1}+\rho_{2}v_{f2}}{\rho_{1}+\rho_{2}}, &   (l_{1}
+ \tau_{1}v_{f1}) \rho_{1} +   (l_{2} +
\tau_{2}v_{f2}) \rho_{2}  < 1\\
\frac{1-\rho_{1} l_{1} - \rho_{2}l_{2}}{\rho_{1} \tau_{1} +
           \rho_{2}\tau_{2}},  & (l_{1} + \tau_{1}v_{f1})
           \rho_{1} +  (l_{2} + \tau_{2}v_{f2}) \rho_{2}  \ge 1
\ea
\right.
\]
This fundamental diagram is shown in \reffig{fd_Zhang}.

\bfg
\begin{center}
\includegraphics[height=12cm]{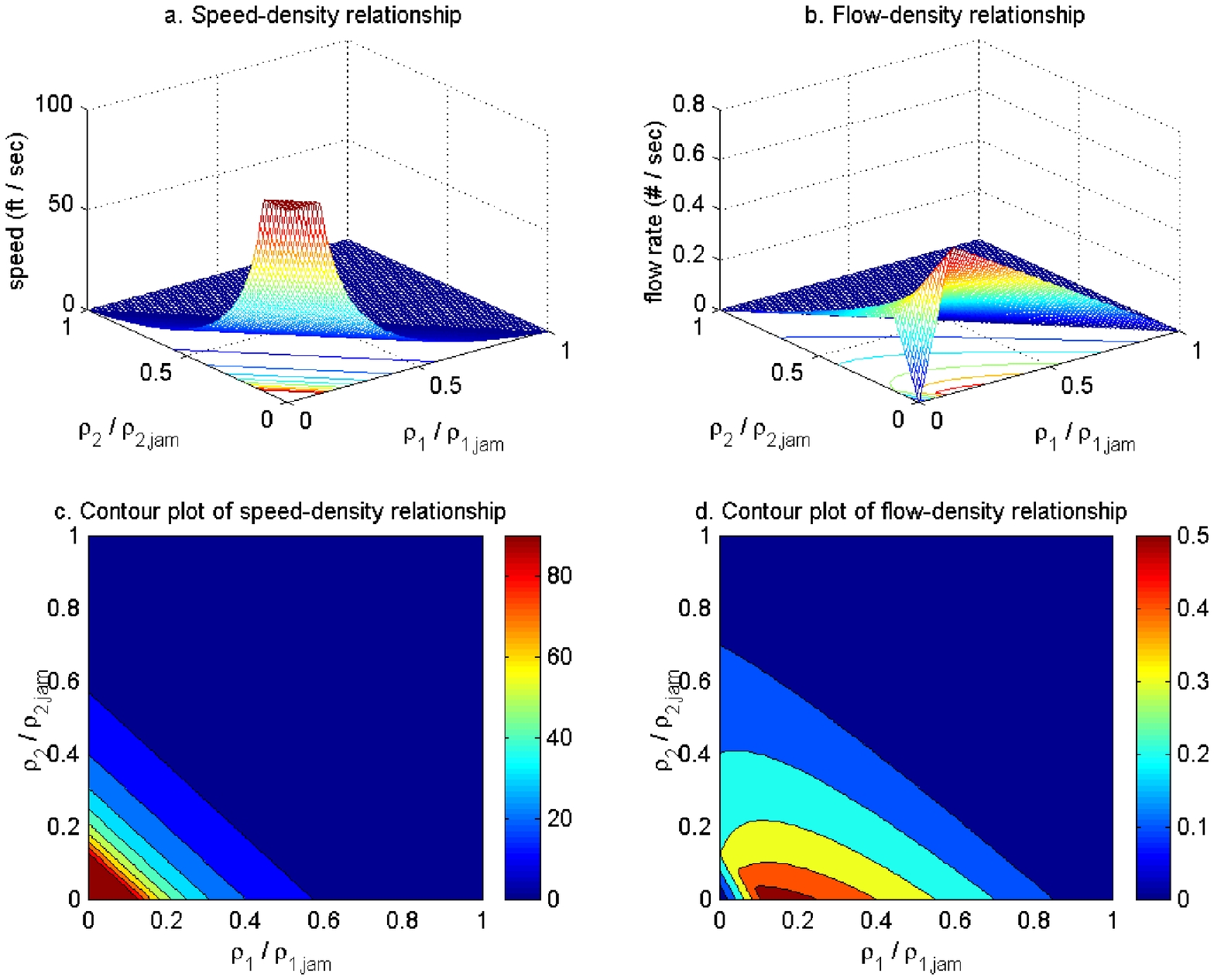}
\caption{The extended triangular fundamental diagram}
\label{fd_Zhang}
\end{center}
\efg

We call the above relation the extended speed-density relation for the
triangular fundamental diagram. The  parameters are: free flow speeds  for
both types of vehicles
$v_{f1}$ and $v_{f2}$,  effective vehicle lengths for  type-1 and
type-2 vehicles
$l_{1}$ and
$l_{2}$,    and response times of type-1 and type-2 vehicles
$\tau_{1}$ and $\tau_{2}$. The last two parameters capture, to a
certain degree, the  acceleration/deceleration differences between the two
classes of vehicles.

The capacity of mixed flow depends on vehicle composition. For
example, in the case of $v_{f1}=v_{f2}=v_{f}$, let
$\frac{\r_{2}}{\r_{1}}=p < \infty$, then the critical densities of the
proposed speed-density relation are
\[
\r_{1c}=\frac{1}{(l_{1}+pl_{2})+v_{f}(\tau_{1}+p \tau_{2})}, \;\;
\r_{2c}=p\r_{1c}
\]
Note that when $p=0$, i.e., there are no type-2 vehicles in the
traffic stream, we recover the critical density for type-1 vehicles
\[
\r_{1c}=  \frac{1}{l_{1} + v_{f} \tau_{1}}
\]
and when  $p=\infty$, i.e., no type-1 vehicles present in the traffic
stream, we can switch the positions of $p$ in relation to $l$'s and
$\tau$'s in the above formulas and replace it with $\frac{1}{p}=0$.
Again we recover the critical density
for type-2 vehicles
\[
\r_{2c}= \frac{1}{l_{2} + v_{f} \tau_{2}}
\]
For any vehicle  mixture (i.e. $0<p<\infty$), the capacity flow is
$\left( \frac{v_{f}}{l_{2} +
v_{f}\tau_{2}} ,  \frac{v_{f}}{l_{1} + v_{f}\tau_{1}} \right)$.

%Wenlong Jin: 2003
\section{Numerical solution method and simulations}\label{sec:mixed5}
\subsection{The Godunov method}
We approximate the mixed traffic KW model with the Godunov
method \citep{godunov1959}:
\bqs
\frac{\r_{1,i}^{j+1}-\r_{1,i}^j}{\dt}+\frac{\r_{1,
i+1/2}^{*j}V(\r_{1, i+1/2}^{*j},\r_{2,
i+1/2}^{*j})-\r_{1,i-1/2}^{*j}V(\r_{1,i-1/2}^{*j},\r_{2,i-1/2}^{*j})}{
\dx}&=&0,\\
\frac{\r_{2,i}^{j+1}-\r_{2,i}^j}{\dt}+\frac{\r_{2,
i+1/2}^{*j}V(\r_{1, i+1/2}^{*j},\r_{2,
i+1/2}^{*j})-\r_{2,i-1/2}^{*j}V(\r_{1,i-1/2}^{*j},\r_{2,i-1/2}^{*j})}{\dx}&=&0,
\eqs
in which $\r_{k,i}^{j}$ is the average of $\r_k$ over cell $i$ at
time $t_j$; i.e., $\r_{k,i}^j=\int_{x=x_{i-1/2}}^{x_{i+1/2}}
\r_k(x,t_j)\m{ dx}/\dx$, and $\r_{k,i-1/2}^{*j}$ is the average over
time interval $[t_j, t_{j+1}]$ at the boundary $x_{i-1/2}$ between
cells $i$ and $i-1$; i.e., $\r_{k,i-1/2}^{*j}=\int_{t=t_j}^{t_{j+1}}
\r_{k}(x_{i-1/2},t) \m{ dt}/\dt$. Given $(\r_1,\r_2)$ at $t_j$, we
hence can compute traffic states at the following time $t_{j+1}$.

In the above equations, the boundary average $\r_{k,i-1/2}^{*j}$ can
be found by solving the Riemann problem for the extended KW model,
\refe{mkw}, with the following initial condition
($\r^l=(\r_{1,i-1}^j,\r_{2,i-1}^j)$ and
$\r^r=(\r_{1,i}^j,\r_{2,i}^j)$)
\bqs
\r&=&\cas{{ll}\r^l,& \m{ if } x-x_{i-1/2}<0, \\ \r^r, &\m{ if }
x-x_{i-1/2}\geq 0.}
\eqs
As shown in section 2, the Riemann solutions consist of a shock or
rarefaction wave with an intermediate state $\r^m$ and a contact
wave. Since all the waves are self-similar, we have
$\r^{*j}_{i-1/2}$=$\r(x_{i-1/2},t)$= const for all $t>0$, which is
determined by the shock or rarefaction wave connecting $\r^l$ and
$\r^m$ since the contact wave has non-negative velocity and is not
involved.

 From \refet{1-r}{1-s}, we have
\bqn
\frac{\r_2^l}{\r_1^l}&=&\frac{\r_2^m}{\r_1^m},\label{rm1}
\eqn
and from \refe{contact}
\bqn
V(\r^m)=V(\r^r).\label{rm2}
\eqn
Combining \refet{rm1}{rm2}, we can find the intermediate
state $\r^m$, from which we can compute $\r^{*j}_{i-1/2}$ as
described in the following cases.

\bi
\item [Case 1] When $\r^l<\r^m$, they are connected by a shock, and
the shock speed is
\bqn
s&=&\frac{\r_1^lV(\r_1^l,\r_2^l)-\r_1^mV(\r_1^m,\r_2^m)}{\r_1^l-\r_1^m}.
\label{defs}
\eqn
In this case, solutions of $\r^{*j}_{i-1/2}$ are summarized in the
\reft{mixed_shock}.

\btb
\begin {center}
\begin {tabular} {||c||c|c|c||}\hline
&$s$ given in \refe{defs} &$\r_{1,i-1/2}^{*j}$&$\r_{2,i-1/2}^{*j}$
\\\cline{2-4}
Shock&$s>0$ & $\r_1^l$ &$\r_2^l$ \\\cline{2-4}
&$s\leq0$ & $\r_1^m$ & $\r_2^m$ \\\hline
\end {tabular}
\end {center}
\caption{Shock wave solutions in mixed traffic}\label{mixed_shock}
\etb

\item [Case 2] When $\r^l>\r^m$, they are connected by a rarefaction
wave, in which the characteristic velocity is $\la_1(\r_1,\r_2)$,
and $\lambda_1(\r^l)<\lambda_1(\r^m)$.
If $\lambda_1(\r^l)\geq 0$, $\r^{*j}_{i-1/2}$ is the same as the left
state $\r^l$; if $\lambda_1(\r^m)\leq 0$, it is the same as the
intermediate state $\r^m$. Otherwise, $\r^{*j}_{i-1/2}$ satisfies
\bqn
\ba{lcl}
\lambda_1 (\r^{*j}_{1,i-1/2},\r^{*j}_{2,i-1/2})&=&0,\\
{\r^{*j}_{2,i-1/2}}/{\r^{*j}_{1,i-1/2}}&=&{\r_2^l}/{\r_1^l},
\ea\label{maximum}
\eqn
which implies that $\r^{*j}_{i-1/2}$ maximizes the total flow
$(\r_1+\r_2)V(\r_1,\r_2)$ along the line $\r_2/\r_1=\r_2^l/\r_1^l$.

In this case, therefore, solutions of $\r^{*j}_{i-1/2}$ are
summarized in \reft{mixed_rarefaction}.

\ei

\btb
\begin {center}
\begin {tabular} {||c||c|c|c||}\hline
&$\lambda_1$ &$\r_{1,i-1/2}^{*j}$&$\r_{2,i-1/2}^{*j}$ \\\cline{2-4}
&$\lambda_1(\r^l)\geq0$ & $\r_1^l$ &$\r_2^l$ \\\cline{2-4}
Rarefaction&$\lambda_1(\r^m)\leq0$ & $\r_1^m$ & $\r_2^m$
\\\cline{2-4}
&o.w. & \multicolumn {2} {c||} {given in \refe{maximum} } \\\hline
\end {tabular}
\end {center}
\caption{Rarefaction wave solutions in mixed traffic}\label{mixed_rarefaction}
\etb

\subsection{Numerical simulations}

In our simulations, we will use the
extended triangular fundamental diagram (\reffig{fd_Zhang}), in
which
the parameter values are: free flow speed
for both types of vehicles $v_{f1}=v_{f2}=v_{f}=65$ mph = 95.3333 ft/sec,
effective vehicle lengths for
type-1 and type-2 vehicles: $l_{1}=20$ ft, $l_{2}=40$ ft,
and response times of type-1 and type-2 vehicles:  $\tau_{1} =1.5$ s,
$\tau_{2}=3$ s. Therefore, we have $\r_{1,jam}=1/l_1=0.05$ veh/ft and
$\r_{2,jam}=1/l_2=0.025$ veh/ft. Moreover, since $l_1/l_2=\tau_1/\tau_2$, we have that in the extended triangular fundamental diagram $V(\r_1,\r_2)$ is a function of $\r_1/\r_{1,jam}+\r_2/\r_{2,jam}$. Thus, as we will see later, the evolution pattern of travel speed is the same as that of $\r_1/\r_{1,jam}+\r_2/\r_{2,jam}$.

We will conduct numerical simulations for a ring road, whose length
$L=2000l_1$=40,000 ft, during a time interval from $t=0$ to $T=100
\tau_1$=150 s. In order to apply the Godunov method, we partition
the ring road into $N=1000$ cells with length $\dx=L/N$=40 ft, and
discretize the time interval to $M=N/2$ steps with the duration of
each time step $\dt=T/M$=0.3 s. Since the CFL condition number
\citep{courant1928CFL}
\bdm
\max \{|\la_1|,|\la_2|\}\frac\dt\dx\leq v_f \frac\dt\dx=0.715 <1,
\edm
the Godunov method yields convergent solutions.

For the numerical simulations, we use the following global
perturbation as initial traffic conditions,
\bqn
\ba{lcl}
\r_1(x,0)&=&(0.2+0.16\sin (2\pi x/L))\r_{1,jam},\\
\r_2(x,0)&=&(0.15+0.1\sin (2\pi x/L))\r_{2,jam},
\ea
\eqn
in which the density of the 1-type vehicles, i.e., has higher average
but smaller oscillation.

With these conditions, the solutions of the mixed traffic flow
model are depicted as contour plots and shown in \reffig{Z_g_contour}. The horizontal axis in each of
the sub figures represents space and the vertical axis time. These
figures depict the traffic conditions (speed and density).
 As shown in the contour plots of $v$
and $\r_1/\r_{1,jam}+\r_2/\r_{2,jam}$, there are roughly two traffic regions along the ring road
initially. In one region  $(l_1+\tau_1v_f)\r_1+(l_2+\tau_2v_f)\r_2\geq 1$, waves initiated from this
region travel backward in the  same speed, which is
$\la_1(\r_1,\r_2)=-(\r_1l_1+\r_2l_2)/(\r_1\tau_1+\r_2\tau_2)=-l_1/\tau_1=-40/3$
ft/sec. In another region, $(l_1+\tau_1v_f)\r_1+(l_2+\tau_2v_f)\r_2 < 1$, waves initiated from this region
travel forward at free-flow speed. Two waves separate the two regions: an expansion wave on the left and
a shock wave on the right. The shock wave travels forward initially but
eventually travel at $- l_1/\tau_1$
as traffic density increases to critical density in the second (free-flow) region. The patterns of the
solutions of $\r_1$ and $\r_2$, however, are not the same as that of $v$ because the change of the
overall traffic conditions affect each vehicle class differently.

\bfg
\begin{center}
\includegraphics[height=12cm]{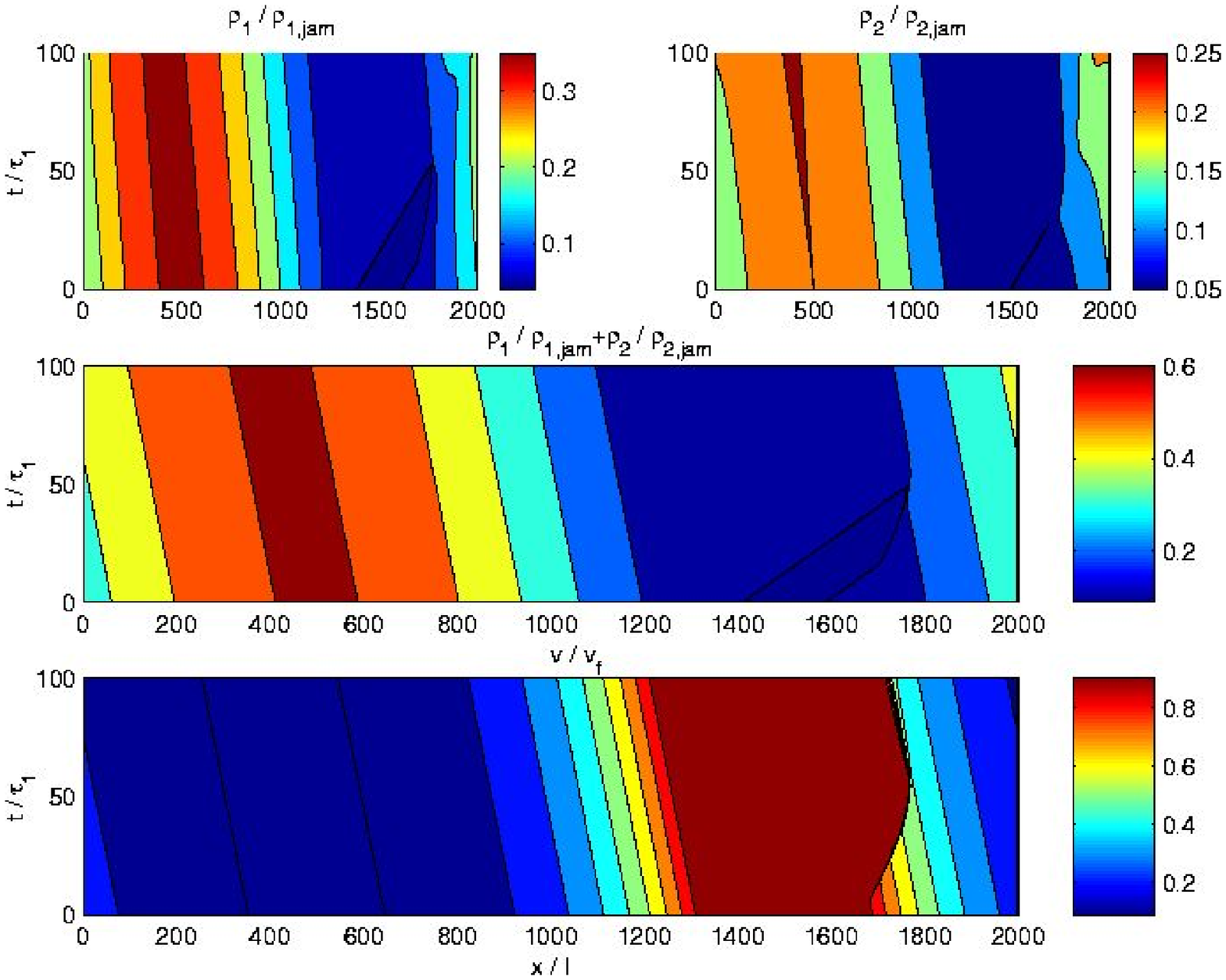}
\caption{Contour plots of solutions on the $x-t$ space with the
extended triangular fundamental diagram}
\label{Z_g_contour}
\end{center}
\end{figure}

The contour plot of $\r_2/\r_1$ is given in \reffig{Z_fifo}, from which we can see that the level curves
do not intersect. Remembering that the level curves of $\r_2/\r_1$ coincide with vehicle trajectories,
the solutions shown here indicate that under the given conditions
first-in-first-out rule is respected by the mixed flow model. From this figure we can also see the
expansion and shock waves as they move through traffic, which are marked by changes in the slopes of the
level curves.

\bfg
\begin{center}
\includegraphics[height=12cm]{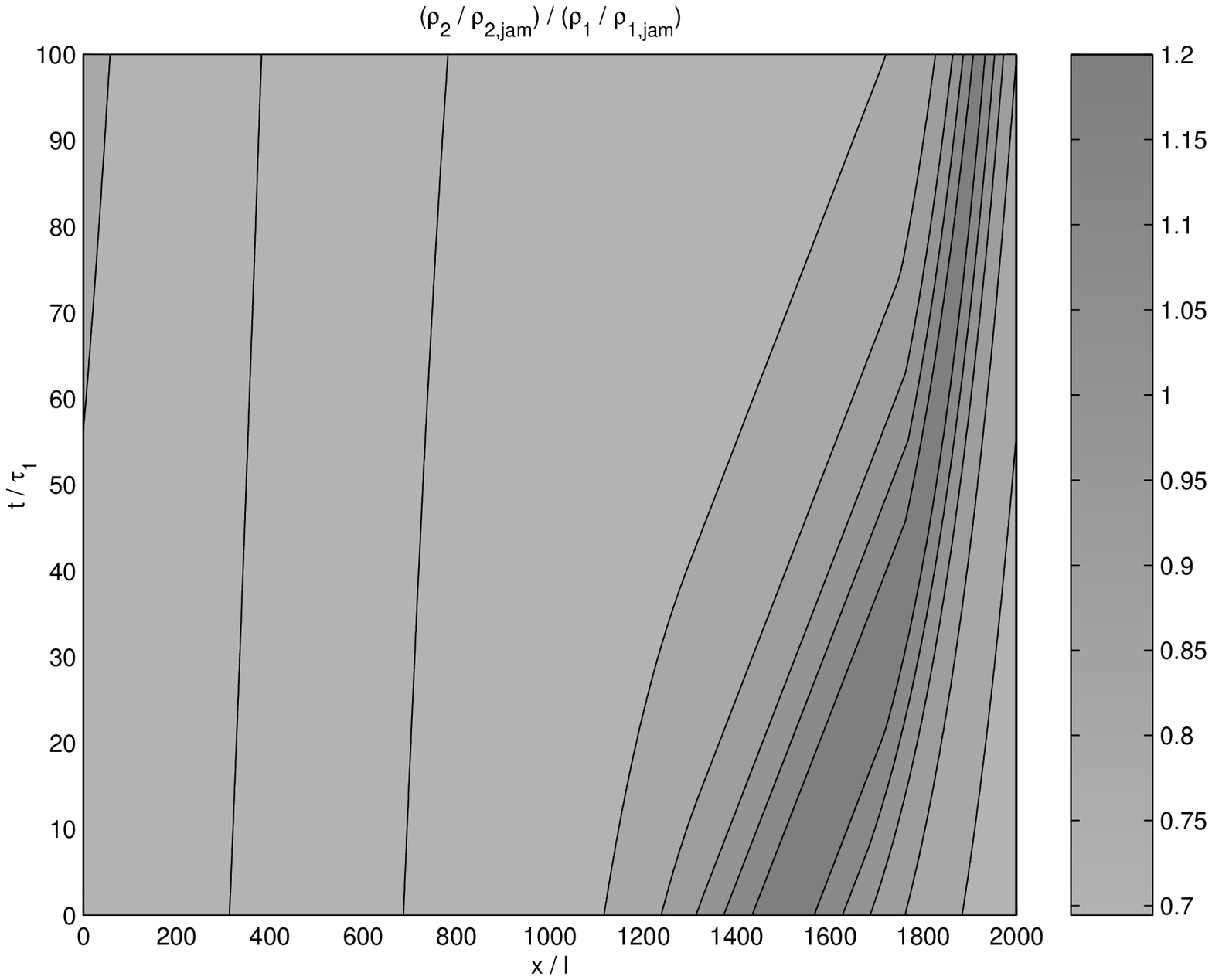}
\caption{Contour plot of $\r_2/\r_1$ on the $x-t$ space with the
extended triangular fundamental diagram}
\label{Z_fifo}
\end{center}
\end{figure}

%Wenlong Jin: 2003
\section{Concluding remarks}\label{sec:mixed6}
In this chapter we extend the Lighthill-Whitham-Richards kinematic
wave traffic flow model to describe traffic with different types of
vehicles, where all types of vehicles are completely mixed and travel
at the same group velocity. A study of such a model with two vehicle
classes (e.g., passenger cars and trucks) shows that,when
both classes of traffic have identical free-flow speeds, the model 1)
satisfies first-in-first-out rule, 2) is anisotropic, and 3) has the
usual shock and expansion waves, and a family of contact waves.
Different compositions of vehicle classes in this model propagate
along contact waves. Such models can be used to study traffic
evolution on long crowded highways where low performance vehicles
entrap high performance ones.

\newpage
\pagestyle{myheadings}
\markright{  \rm \normalsize CHAPTER 6. \hspace{0.5cm}
 MULTI-COMMODITY NETWORK TRAFFIC MODEL}
\chapter{Kinematic wave simulation model for multi-commodity network traffic flow}
%\thispagestyle{myheadings}
%Wenlong Jin: 2003
%Wenlong Jin: 2003
\section{Introduction}
Recurrent or non-recurrent traffic congestion in many major metropolitan areas have seriously deteriorated the mobility, convenience, and productivity of highway transportation. To tackle the congestion problem, traffic engineers and scientists are facing a number of challenges, including: 1) the evaluation of the performance of a road network (e.g. total travel time or the level of service (LOS)), 2) the prediction of occurrence of incidents, 3) the development of traffic control schemes (e.g. ramp metering) and management strategies (e.g. traffic guidance) in Advanced Transportation Management and Information Systems (ATMIS) or other Intelligent Transportation Systems (ITS), and 4) the estimation of traffic demand associated with a given origin/destination (i.e. O/D information). As we all know, a fundamental and essential issue in addressing all these challenges is the understanding of traffic dynamics on a road network, i.e., the evolution of traffic on a road network under initial and boundary conditions, for which traffic flow models of road networks play an important role.

Among many traffic flow models, the kinematic wave models are advantageous in studying traffic dynamics in a large-scale road network, due to their inherit compliance with ITS applications, theoretical rigor, and computational efficiency. Practically, people are more interested in aggregate-level traffic conditions such as average travel speeds, densities, flow-rates, and travel times, which are directly concerned or can be easily derived in the kinematic wave models. Theoretically, the evolution of traffic conditions on a link can be studied as hyperbolic conservation laws, given the fundamental diagram, which defines a functional relationship between density and flow-rate or travel speed, and traffic conservation, which means that the change of the number of vehicles in a section in a time interval equals to the number of vehicles entering the section minus the number of vehicles leaving it. Computationally, the kinematic wave models can be solved with the Godunov method, in which in-flux and out-flux of a cell can be computed through kinematic wave solutions or the supply-demand method \citep{daganzo1995ctm,lebacque1996godunov}.

In literature, the following approaches can be considered as typical when modeling network traffic flow in the framework of kinematic wave theory. First, \citet{vaughan1984od} studied network traffic flow with two continuous equations: a ``local equation", which ensures traffic conservation and is consistent with the traditional LWR model on a link at the aggregate level, and a ``history equation", which computes the trajectory of each vehicle at the disaggregate level. Since the trajectories of all vehicles are not always required for many applications, this model tends to waste significant amount of computational resources. Second, \citet{jayakrishnan1991dissertation} introduced a discrete network flow model, in which each link is partitioned into a number of cells, vehicles adjacent to each other and with the same O/D and path are called a ``macroparticle", and the position of a macroparticle at each time step is determined by its travel speed and the cell's length. However, with the given mechanism, this model may not be consistent with the LWR model since traffic conservation can be violated. Third, \citet{daganzo1995ctm} introduced a network flow model based on his Cell Transmission Model \citep{daganzo1994ctm}, which is a numerical approximation of the LWR model with a special fundamental diagram --- a triangle. In this discrete model, macroparticles in the sense of \citet{jayakrishnan1991dissertation} in a cell are ordered according to their waiting times and moved to the downstream cell when their waiting times are greater than a threshold minimum waiting time, which is computed from traffic flow at the aggregate level. However, the determination of the threshold minimum waiting time is quite tedious. Fourth, in the KWave98 simulation package \citep{leonard1998kwave}, which is based on the simplified kinematic wave theory by \citet{newell1993sim}, vehicles on a link are considered to be ordered as a queue. Then, with the in-queue and out-queue of the link determined, one can easily update the link queue. However, in-queues and out-queues at typical highway junctions such as merges and diverges create difficulties for this model.

In many applications, such as dynamic traffic assignment (DTA), First-In-First-Out (FIFO) principle is a key concern. Consequently, all the aforementioned models order vehicles but in different fashions. Vehicles are ordered according to their trajectories in \citep{vaughan1984od}, locations in \citep{jayakrishnan1991dissertation}, waiting times in \citep{daganzo1995ctm}, and positions in a link queue in \citep{leonard1998kwave}. That is, these simulation models track either vehicles' trajectories, or positions, or waiting times, or queue orders. \footnote{In \citep{jayakrishnan1991dissertation}, although the introduction of macroparticles save a certain amount of memory in computation, this saving is diminished when vehicles of different classes are evenly mixed.} In these models, therefore, traffic conditions are also considered at the vehicle level. As a result, the computational efficiency of the kinematic wave models is not fully utilized.

In order to fully explore the computational efficiency inherited in  the kinematic wave theory, we propose a new network traffic flow simulation model, in which traffic dynamics are studied at the aggregate level with commodity specified densities and flow-rates. Traffic of a commodity has the same characteristics, which can be vehicle type, destination, path, or any classification criterion. We hereafter refer to this model as  multi-commodity kinematic wave (MCKW) simulation model. In this dissertation research, we start with a simple traffic system, where vehicles are categorized into multiple commodities according to their paths and no differentiations are made in vehicle types, driver classes, or lane types, such as high-occupancy-vehicles (HOV's) or HOV lanes. In the MCKW simulation platform, besides link characteristics such as free flow speed, capacity, and number of lanes, traffic dynamics are highly related to geometrical characteristics of a road network including link inhomogeneities, merges, diverges, and other junctions. Although it does not give a complete, detailed picture of traffic dynamics, this model is still of importance for many applications, in which the difference in vehicles, drivers, or lanes are negligible not interested or can be averaged out without major loss of accuracy.

In the MCKW simulation platform, the Godunov-type approximation is applied, and fluxes through boundaries inside a link and junctions such as merges and diverges are computed systematically in the framework of the supply-demand method. In addition, since vehicles are categorized by vehicles' origins, destinations, and paths, traffic dynamics are also differentiated for different commodities in the MCKW simulation platform. We will show that FIFO is observed in the commodity specific kinematic wave theory.

In the rest of this chapter, we will discuss theories underlying the MCKW model in Section \ref{sec:sim2}, and the data, network, and program structures in the MCKW model in Section \ref{sec:sim3}. In Section \ref{sec:sim4}, we will discuss the process of output from the MCKW model to obtain interested data of a road network, such as travel times. In Section \ref{sec:sim5}, we carry out some numerical simulations of a simple road network. In Section \ref{sec:sim6}, we will draw some conclusions and provide further discussions about the calibration of the MCKW model.

%Wenlong Jin: 2003
\section{Underlying theories of the MCKW simulation model}\label{sec:sim2}
In the multi-commodity kinematic wave (MCKW) simulation of network traffic flow, dynamics of total traffic, i.e., the evolution of traffic conditions of all commodities, are studied at the aggregate level and governed by the kinematic wave theories. These theories have been studied in the previous chapters for fundamental network components such as inhomogeneous links, merges, and diverges. They provide the underlying algorithms and building blocks for the MCKW simulation platform. At the disaggregate level, traffic of each commodity in the form of proportions is studied with its proportion, and First-In-First-Out property on a link will be discussed.

\subsection{Kinematic wave theories at the aggregate level}
In the MCKW model, we use the discrete form of the kinematic wave theories, which can be obtained through the first-order\footnote{A second-order method was discussed in \citep{daganzo1999lagged}.} Godunov method \citep{godunov1959} of the continuous versions. In the discrete form, each link is partitioned into $N$ cells, of equal length or not, and the time interval is discretized into $K$ time steps. Then, we obtain the Godunov-type finite difference equation for total flow in cell $i$ from time step $j$ to time step $j+1$ as follows:
\bqn
\frac{\r_i^{j+1}-\r_i^j}{\dt}+\frac{f^{j\ast}_{i-1/2}-f^{j\ast}_{i+1/2}}{\dx}&=&0, \label{fi}
\eqn
where $\dx$ is the length of cell $i$, $\dt$ is the time from time step $j$ to time step $j+1$, and the choice of $\frac {\dt}{\dx}$ is governed by the CFL \citep{courant1928CFL} condition. In \refe{fi}, $\r_i^j$ is the average of traffic density $\r$ in cell $i$ at time step $j$, similarly $\r_i^{j+1}$ is the average of $\r$ at time step $j+1$; $f^{j\ast}_{i-1/2}$ is the flux through the upstream boundary of cell $i$ from time step $j$ to time step $j+1$, and similarly $f^{j\ast}_{i+1/2}$ is the downstream boundary flux. Given traffic conditions at time step $j$, we can  calculate the traffic density in cell $i$ at time step $j+1$ as
\bqn
\r_i^{j+1}&=&\r_i^j+\dxt(f^{j\ast}_{i-1/2}-f^{j\ast}_{i+1/2}). \label{fi2}
\eqn

Defining ${\bf N}_i^j=\r_i^j \dx$ as the number of vehicles in cell $i$ at time step $j$, ${\bf N}_i^{j+1}=\r_i^{j+1}\dx$ as the number of vehicles at time step $j+1$, $F_{i-1/2}^j=f(\r^{j\ast}_{i-1/2})\dt$ as the number of vehicles flowing into cell $i$ from time step $j$ to $j+1$, and $F_{i+1/2}^j$ as the number of vehicles flowing out of cell $i$,  \refe{fi2} can be written as:
\bqn
{\bf N}_i^{j+1}&=&{\bf N}_i^j+F_{i-1/2}^j-F_{i+1/2}^j, \label{diffe}
\eqn
which is in the form of traffic conservation.

Given the initial and boundary conditions, we will use the supply-demand method \citep{daganzo1995ctm,lebacque1996godunov} for computing fluxes through cell boundaries: $F_{i-1/2}^j$ or $f^{j\ast}_{i-1/2}$. In a general road network, there are the following types of boundaries: boundaries inside a link, merges, diverges, and more complicated intersections.
\ben
\item When the boundary at $x_{i-1/2}$ is a boundary inside a link, whose upstream cell is denoted as $u$ and downstream cell $d$, we follow the supply-demand method discussed in \citep{daganzo1995ctm,lebacque1996godunov,jin2003inhLWR}. I.e., if we define the upstream demand as
\bqn
D_u&=&\cas{{ll} Q(U_u), & \m{when } U_u \m{ is under-critical} \\Q^{max}_u, & \m{when } U_u \m{ is over-critical}} \label{demands}
\eqn
and define the downstream supply as
\bqn
S_d&=&\cas{{ll} Q^{max}_d, & \m{when } U_d \m{ is under-critical} \\Q(U_d), & \m{when } U_d \m{ is over-critical}} \label{supplys}
\eqn
then the boundary flux can be simply computed as
\bqn
f^{j\ast}_{i-1/2}&=&\min\{D_u, S_d\}, \label{sd_boundary}
\eqn
where $U_d$ and $U_u$ are traffic conditions including density $\r$ and road inhomogeneity $a$ at $j$th time step, of the downstream and upstream cells, respectively. As discussed in \citep{jin2003inhLWR}, this method is consistent with analytical solutions of the Riemann problem for inhomogeneous roadway.

\item When $x_{i-1/2}$ is a merging junction with $P$ upstream merging cells, which are denoted as
$u_p$ ($p=1,\cdots,P$), and a downstream cell $d$. The demand of upstream cell $u_p$, $D_p$, is defined in \refe{demands}, and the supply of the downstream cell, $S_d$, is defined in \refe{supplys}. Then, we apply the simplest distribution scheme \citep{jin2003merge} and compute the boundary fluxes as
\bqn
\ba{lcl}
f^{j\ast}_{i-1/2,d}&=&\min\{\sum_{p=1}^P D_p,S_d\},\\
f^{j\ast}_{i-1/2,u_p}&=&q \frac{D_p}{\sum_{p=1}^P D_p}, \quad p=1,\cdots,P,
\ea\label{sd_merge}
\eqn
where $f^{j\ast}_{i-1/2,d}$ is the in-flow of the downstream cell, and $f^{j\ast}_{i-1/2,u_p}$ is the out-flow of upstream cell $u_p$.

In addition, if an upstream cell, e.g. $u_p$, is signalized and denote $r$ as the proportion of green light in a cycle (i.e. green ratio), then we can apply the controlled traffic demand of $u_p$, $\min\{r Q^{max}_p,D_p\}$, in the supply-demand
method above \citep{daganzo1995ctm, jin2003merge}. Note that $r$ can be a continuous function when considering the average effect or a piece-wise constant function when the simulation interval $\dt$ is smaller than a signal cycle.

\item When $x_{i-1/2}$ is a diverging junction with $P$ downstream cells, which are denoted as $d_p$ ($p=1,\cdots,P$), and an upstream cell $u$. In the model proposed in \citep{jin2002diverge}, we introduced a new definition of partial traffic demand of vehicles traveling to $D_p$ in cell $u$ as follows,
\bqn
D_p&=&\cas{{ll} Q(\r_p;\hat\r_p)&\r_p \m{ is UC} \\Q^{\max}(\hat \r_p) &  \m{otherwise}}, \label{diverge2.1}
\eqn
where $\hat \r_p$ is equal to the density of vehicles not traveling to $d_p$, and at critical density $Q(\r;\hat\r_p)$ reaches its maximum. The traffic supply for $d_p$, $S_p$, is defined by \refe{supplys}. Then, the boundary flux to $d_p$, $f^{j\ast}_{i-1/2,d_p}$, can be computed by
\bqn
f^{j\ast}_{i-1/2,d_p}&=&\min\{S_p,D_p\},
\eqn
and the out-flow of $u$, $f^{j\ast}_{i-1/2,u}$, is the sum of these flows,
\bqn
f^{j\ast}_{i-1/2,u}&=&\sum_{p=1}^P  f^{j\ast}_{i-1/2,d_p}. \label{diverge2.3}
\eqn

Another model of traffic diverging  to $D$ downstream links we will implement in the MCKW simulation was proposed in \citep{daganzo1995ctm,lebacque1996godunov}:
\bqn
\ba{lcl}
f^{j\ast}_{i-1/2,u}&=&\min_{d=1}^D\{D_u,S_d / \xi_d\},\\
f^{j\ast}_{i-1/2,d}&=&\xi_d f^{j\ast}_{i-1/2,u}, \qquad d=1,\cdots ,D,
\ea\label{diverge1}
\eqn
where $\xi_d$ is the proportion of commodity $d$ in total traffic, and here $D_u$ is the demand of the upstream cell as defined in \refe{demands}.

When vehicles have no predefined route choice and can choose every downstream link at a diverge, we use the model proposed in \citep{jin2003merge}:
\bqn
\ba{lcl}
f^{j\ast}_{i-1/2,u}&=&\min\{D_u,\sum_{d=1}^D S_d\},\\
f^{j\ast}_{i-1/2,d}&=&\frac{S_d}{\sum_{d=1}^D S_d}f^{j\ast}_{i-1/2,u}, \qquad d=1,\cdots, D.
\ea\label{diverge3}
\eqn

\item For intersections with two or more upstream and downstream links, we can combine the merge and diverge models together. Note that only the computation of demands and supplies may change, and the supply-demand method is still the same.

For example, when we combine the supply-demand methods in \refe{sd_merge} and \refe{diverge1} for an intersection with $U$ upstream branches and $D$ downstream branches, we can compute fluxes by
\bqn
\ba{lcl}
f^{j\ast}_{i-1/2}&=&\min_{d=1}^D\{\sum_{u=1}^U D_u,S_d \big / \left(\frac{\sum_{u=1}^U D_u \xi_{u,d}}{\sum_{u=1}^U D_u} \right) \},\\
f^{j\ast}_{i-1/2,d}&=&\frac{\sum_{u=1}^U D_u \xi_{u,d}}{\sum_{u=1}^U D_u} f^{j\ast}_{i-1/2}, \qquad d=1,\cdots ,D,\\
f^{j\ast}_{i-1/2,u}&=&\frac{D_u}{\sum_{u=1}^U D_u} f^{j\ast}_{i-1/2}, \qquad u=1,\cdots ,U,
\ea\label{sd_general}
\eqn
where $\xi_{u,d}$ is the proportion of traffic heading downstream link $d$ in upstream link $u$, $f^{j\ast}_{i-1/2}$ is the total flux through the boundary, $f^{j\ast}_{i-1/2,d}$ flux heading downstream link $d$, and $f^{j\ast}_{i-1/2,u}$ flux from upstream link $u$. In this model, the intersection is considered as a combination of a merge with $U$ upstream branches and a diverge with $D$ downstream links in the fashion of \citep{daganzo1995ctm}. Note that the merge model, \refe{sd_merge}, and the diverge model, \refe{diverge1}, are specific cases of \refe{sd_general}.
\een

\subsection{Commodity-based kinematic wave theories}
In the MCKW simulation platform, commodities are differentiated by origin/destination or path.
We assume that a road network has $P'$ origin/destination (OD) pairs and $P$ paths ($P\geq P'$). When vehicles have predefined paths, we then have a $P$-commodity traffic flow on the road network and label vehicles taking $p$th path as $p$-commodity. When vehicles of an O/D have no predefined paths, we have $P'$-commodity traffic flow.

In the kinematic wave theories of multi-commodity traffic, we denote total traffic density, travel speed, and flow-rate respectively by $\r$, $v$, and $q$, which are all functions of location $x$ and time $t$. In contrast, these quantities for $p$-commodity vehicles are $\r_p$, $v_p$, and $q_p$ respectively.  The fundamental diagram of total traffic  defines a functional relationship between density and travel speed or flow-rate: $q=Q(a,\r)$ and $v=V(a,\r)\equiv Q(a,\r)/\r$, where $a(x)$ stands for road inhomogeneities at location $x$ such as changes in the number of lanes, curvature, and free flow speed. Further, we assume traffic on all links is additive in the following sense \citep{jin2002diverge}:
\bqn
\r&=&\sum_{p=1}^P \r_p,\\
v&=&v_p=V(a,\r),\qquad p=1,\cdots,P,\\
q&=&\sum_{p=1}^P q_p,
\eqn

The kinematic wave theory of additive multi-commodity traffic on a link can be described by the following theory \citep{jin2003inhLWR},
\bqn
\ba{lcl}
\r_t+Q(a,\r)_x&=&0,\\
(\r_p)_t+(\r_p V(a,\r))_x&=&0, \qquad p=1,\cdots, P.
\ea\label{lwrs}
\eqn
If denoting the local proportion of $p$-commodity ($p=1,\cdots,P$) by $\xi_p=\r_p/\r$, we then have the following advection equations \citep{lebacque1996godunov}
\bqn
(\xi_p)_t+V(a,\r) (\xi_p)_x&=&0, \qquad p=1,\cdots, P. \label{FIFO}
\eqn
From \refe{FIFO}, we can see that proportions of all commodities travel forward in a link along with vehicles in traffic flow, as the change of $\xi_p$ in material space, $(\xi_p)_t+V(a,\r) (\xi_p)_x$, equals to zero. This is also true for all kinds of junctions, in particular diverges, in their supply-demand models in the preceding subsection \footnote{That traffic is anisotropic is believed to regulate this property.} Therefore, \refe{FIFO} also means that the profile of proportions coincides with vehicles' trajectories on a link. That is, if two or more commodities initially completely are divided by an interface,
this interface will move forward along with vehicles on both sides of, and these commodities will never mix. Since each single vehicle can considered as a commodity, all vehicles' trajectories keep disjoint in the commodity-based kinematic wave models. Therefore, FIFO principle is respected in this continuous model.

In the previous subsection, we studied the discrete kinematic wave theory for total traffic. Here, we will present the discrete kinematic wave theory for each commodity. Given traffic conditions of $p$-commodity at time step $j$, i.e., $\r_{p,i}^j$ in all cells, we can  calculate the traffic density of $p$-commodity in cell $i$ at time step $j+1$ as
\bqn
\r_{p,i}^{j+1}&=&\r_{p,i}^j+\dtx (f^{j\ast}_{p,i-1/2}-f^{j\ast}_{p,i+1/2}), \label{ceCommodity}
\eqn
where $f^{j\ast}_{p,i-1/2}$ is the in-flux of $p$-commodity through the upstream boundary of cell $i$ during time steps $j$ and $j+1$, and $f^{j\ast}_{p,i+1/2}$ out-flux. Furthermore, since the profile of the proportion of a commodity always travels forward at traffic speed, the proportion of a commodity in out-flux of cell $i$ compared to all commodities is equal to the proportion of the commodity in the  cell. I.e. \citep{lebacque1996godunov},
\bqn
f^{j\ast}_{p,i+1/2}\cdot \r_{i}^j = f^{j\ast}_{i+1/2} \cdot \r_{p,i}^j, \qquad p=1,\cdots, P. \label{FIFO_discrete}
\eqn
This is true for cells right upstream of merging junctions \citep{jin2003merge} and  diverging junctions \citep{papageorgiou1990assignment,daganzo1995ctm,lebacque1996godunov,jin2001diverge2}.

During time steps $j$ and $j+1$ at a boundary $x_{i+1/2}$, which has $U$ upstream cells  and $D$ downstream cells, if we know the out-flux from upstream cell $u$ ($u=1,\cdots,U$), $f^{j\ast}_{p,u,i+1/2}$ ($p=1,\cdots,P$), we can obtain the in-flux of downstream cell $d$ ($d=1,\cdots,D$), $f^{j\ast}_{p,d,i+1/2}$ ($p=1,\cdots,D$), from traffic conservation in $p$-commodity:
\bqn
\sum_{u=1}^U f^{j\ast}_{p,u,i+1/2} = \sum_{d=1}^D f^{j\ast}_{p,d,i+1/2}.\label{tc_commodity}
\eqn
However, when $p$-commodity vehicles can take more than one downstream cells, we have
\bqn
\sum_{p=1}^P f^{j\ast}_{p,d,i+1/2}=f^{j\ast}_{d,i+1/2}.
\eqn

Note that, in \refe{lwrs}, the kinematic wave solutions are determined by those of total traffic, which are obtained by the first-order convergent Godunov method. Also from \refe{FIFO}, we can see that \refe{FIFO_discrete} yields an up-wind method for $\xi_p$ in \refe{ceCommodity}. Therefore, the discrete model for the commodity-based kinematic wave model, \refe{lwrs}, converges in first order to continuous version, whose solutions observe FIFO principle. That is, in numerical solutions, error in travel time of any vehicle is in the order of $\dt$. That is, in the MCKW simulation, FIFO is accurate to the order of $\dt$ and $\dx$. Therefore, when we decrease $\dt$, this approximation becomes more accurate.

%Wenlong Jin: 2003
\section{Network structure, data structure, and program flow-charts in the MCKW simulation platform}\label{sec:sim3}
In this section, we will discuss the programming details of the MCKW simulation of an illustrative road network.

\subsection{Network structure}
For the purposes of exposition, a simple traffic network, shown in \reffig{F_network}, is considered, and these discussions can be extended to more general road networks.

\bfg
\bc\includegraphics[height=6cm]{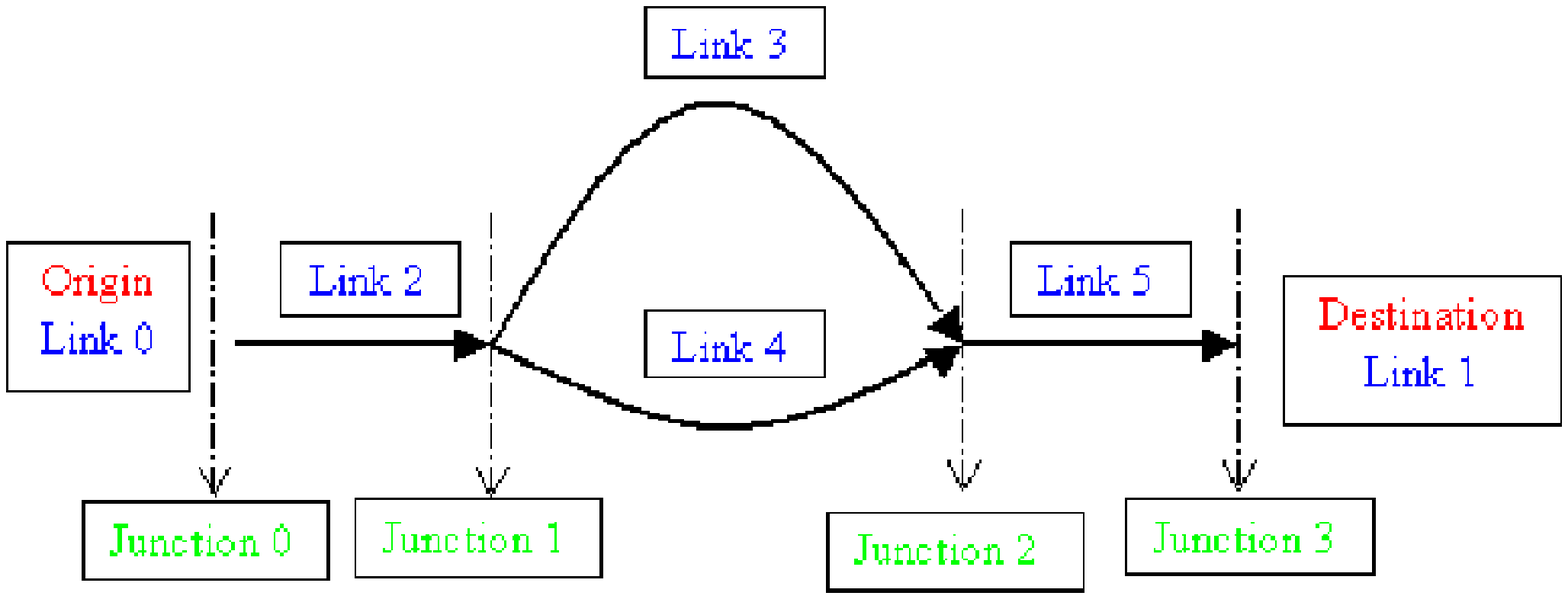}\ec
\caption{A demonstration road network}\label{F_network}
\efg

The road network in \reffig{F_network}, where the arrows show traffic direction, consists of one origin/destination (O/D) pair and four links. In this network, there are two paths. Furthermore, we assume vehicles have predefined paths. \footnote{When vehicles have no predefined paths, these discussions are also applicable.} Thus, traffic flow on this road network consists of two commodities.

In the MCKW simulation, origins and destinations have the same data structure as regular links and are treated as links. For a road network with \texttt{num\_origin} origins, \texttt{num\_od} origins and destinations, and \texttt{num\_link} links, all links are numbered: origins from \texttt{0} \footnote{The numbering starts from 0 rather than 1 according to C/C++ conventions.} to \texttt{num\_origin}-1, the number of origins; destinations from \texttt{num\_origin} to \texttt{num\_od}-1; and regular links from \texttt{num\_od} to \texttt{num\_links}-1. In the sample network, links are numbered as shown in \reffig{F_network}.

With these numbers, the paths or commodities are denoted as follows: commodity 0 takes links 0, 2, 3, 5, and 1, and commodity 1 takes links 0, 2, 4, 5, and 1. This is equivalent to saying that links 0, 2, 5, and 1 carry 0-commodity and 1-commodity flows, link 3 carries only 0-commodity flow, and link 4 carries only 1-commodity flow. Further, the network structure and traffic flow direction is represented by the upstream and downstream links of each link. For example, the upstream links of link 5 and the downstream links of link 2 are links 3 and 4, respectively.

In the MCKW simulation, each link is partitioned into a number of cells.  \footnote{Origins and destinations have only one cell.} Since fluxes through cell boundaries are important in computation in the kinematic wave theories, cell boundaries are also included in the structure of a link. \footnote{Note that there are one more cell boundaries than cells.} In the MCKW simulation platform, cells and boundaries are ordered according to traffic direction: adjacent cells and boundaries are either upstream or downstream to a cell.

In the MCKW simulation, therefore, network can be constructed if we know the structures of all links, commodities on a link, and the upstream and downstream links of all links. That is, junctions are not used to store network structure although they are also numbered in \reffig{F_network}. The representation of network structure in the MCKW simulation largely simplifies the data structure, in which only links are used.

\subsection{Data structure}
In the MCKW simulation platform, the structure and characteristics of a road network and traffic conditions are all represented by links as well as cells and boundaries inside a link. Therefore, the major data structure is \texttt{linkType}, through which network and traffic conditions are dealt with, together with sub-structures for cells, \texttt{cellType}, and boundaries, \texttt{boundaryType}. The data structures are shown in \reffig{F_data} and explained in detail as below.

\begin{figure}[p]
\bc\includegraphics[height=15cm]{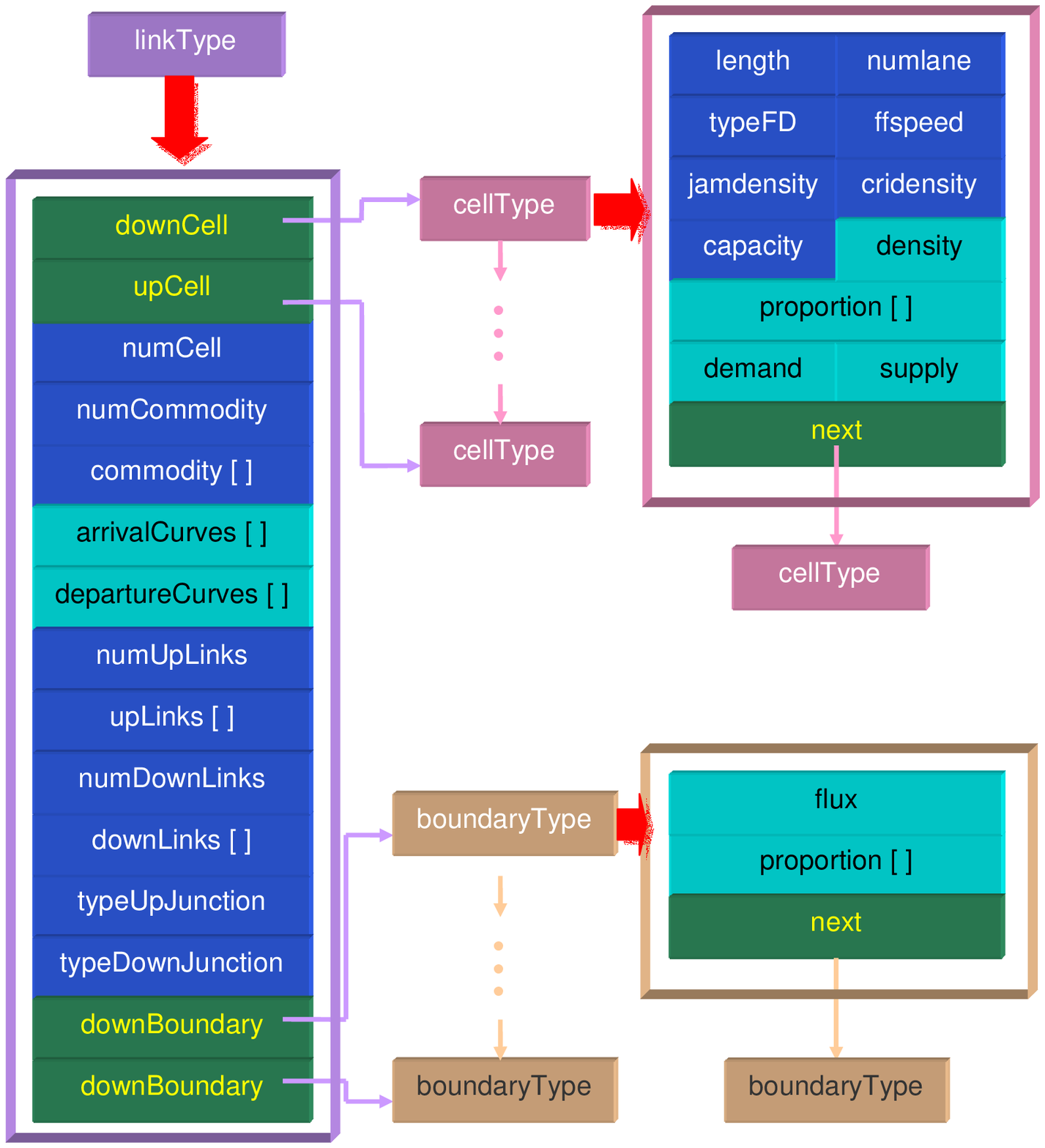}\ec
\caption{Data structure in the MCKW mckw platform}\label{F_data}
\efg

In \reffig{F_data}, the data structure for a link, \texttt{linkType}, is shown in the left box. Its fields are explained as follows:
\bds
\item[\texttt{downCell}] is an address pointing to the furthest downstream cell of a link.
\item[\texttt{upCell}] is an address pointing to the furthest upstream cell of a link.
\item[\texttt{numCell}] is the number of cells in a link.
\item[\texttt{numCommodity}] is the number of commodities traveling on the link.
\item[\texttt{commodity [ ]}] is an array of commodities, whose length is \texttt{numCommodity}. Commodities are ordered increasingly for all links.
\item[\texttt{arrivalCurves [ ]}] is an array of accumulative flow entering the link. Its length is \texttt{numCommodity+1}. The first \texttt{numCommodity} entries store the cumulative flows of corresponding commodities, and the last entry stores total cumulative flow.
\item[\texttt{departureCurves [ ]}] is the same as \texttt{arrivalCurves [ ]} except that we consider exiting flows.
\item[\texttt{numUpLinks}] is the number of upstream links adjacent to a link.
\item[\texttt{upLinks [ ]}] stores all the adjacent upstream links. Its length is \texttt{numUpLinks}.
\item[\texttt{numDownLinks}] is the number of downstream links adjacent to a link.
\item[\texttt{downLinks [ ]}] stores all the adjacent downstream links. Its length is \texttt{numDownLinks}.
\item[\texttt{typeUpJunction}] denotes the type of the upstream junction incident to a link. Here, type $0$ stands for a linear junction connecting one upstream link and one downstream link, type $1$ for a merging junction, type $2$ for a diverging junction as described by \refe{diverge1}, type $3$ for a diverging junction by \refe{diverge2.1}-\refe{diverge2.3}, and type $4$ for a diverging junction by \refe{diverge3}. For different types of junctions, traffic flow models are different, as shown in Section \ref{sec:sim2}.
\item[\texttt{typeDownJunction}] denotes the type of the downstream junction incident to a link. The definition of junction types are the same as in \texttt{typeUpJunction}.
\item[\texttt{downBoundary}] is an address pointing to the furthest downstream boundary of a link.
\item[\texttt{upBoundary}] is an address pointing to the furthest upstream boundary of a link.
\eds
In these cells, four fields in dark green are pointers with no physical meaning, two fields in cyan are time-dependent quantities, and the rest in blue represent quantities that determine network structure and are time invariant.

As exhibited by the top right box of \reffig{F_data}, \texttt{cellType}, has the following fields, which characterize a cell.
\bds
\item[\texttt{length}] is the length of a cell.
\item[\texttt{numlane}] is the number of lanes of a cell.
\item[\texttt{typeFD}] denotes the type of fundamental diagrams in a cell. Type $0$ stands for the triangular fundamental diagram \citep{newell1993sim}. For other types of fundamental diagrams, refer to \citep{delcastillo1995fd_empirical,kerner1994cluster}. With the number of lanes considered, we can have the fundamental diagram for the cell.
\item[\texttt{ffspeed}] is the cell free flow speed.
\item[\texttt{jamdensity}] is the cell jam density of each lane.
\item[\texttt{cridensity}] is the cell critical density of each lane, at which traffic flow reaches capacity.
\item[\texttt{capacity}] is the cell maximum flow-rate of a lane.
\item[\texttt{density}] is the cell total density.
\item[\texttt{proportion [ ]}] is an array of proportions of commodities in a cell. Its length is \texttt{numCommodity} of the corresponding link.
\item[\texttt{demand}] is the cell traffic demand during a time interval, as defined in \refe{demands}.
\item[\texttt{supply}] is the cell traffic supply, as defined in \refe{supplys}.
\item[\texttt{next}] is the address pointing to the upstream cell.
\eds
As in \texttt{linkType}, the dark green field is a pointer, seven blue fields are for time-independent quantities, which denote major characteristics of a cell, and four cyan fields are for time-dependent quantities. Note that the number of lanes, free flow speed, and, therefore, the critical density and capacity may change when accidents occur. Besides, if there are signals on the boundaries of a cell, the demand and supply may be restricted \citep{daganzo1995ctm,jin2003merge}.

The data structure for cell boundaries, \texttt{boundaryType}, is illustrated by the bottom right box in \reffig{F_data} and has the following fields:
\bds
\item[\texttt{flux}] is the flux through the boundary during a time interval.
\item[\texttt{proportion [ ]}] stores the proportions of all commodities in the total flux. Its length is \texttt{numCommodity} of the corresponding link.
\item[\texttt{next}] is the address pointing to the upstream boundary.
\eds
As we can see, both \texttt{flux} and \texttt{proportion [ ]} are time dependent.

In the MCKW simulation, cells and boundaries are ordered in the direction opposite to traffic's. However, it is also straightforward to order them in the traffic direction. Besides, ordering will not affect computation efficiency significantly.

\subsection{Program flow-chart}
The program flow chart in the MCKW simulation platform is shown in \reffig{F_program}. The modules in the program are explained below in the same order as they appear in the chart.

\begin{figure}[p]
\bc\includegraphics[height=15cm] {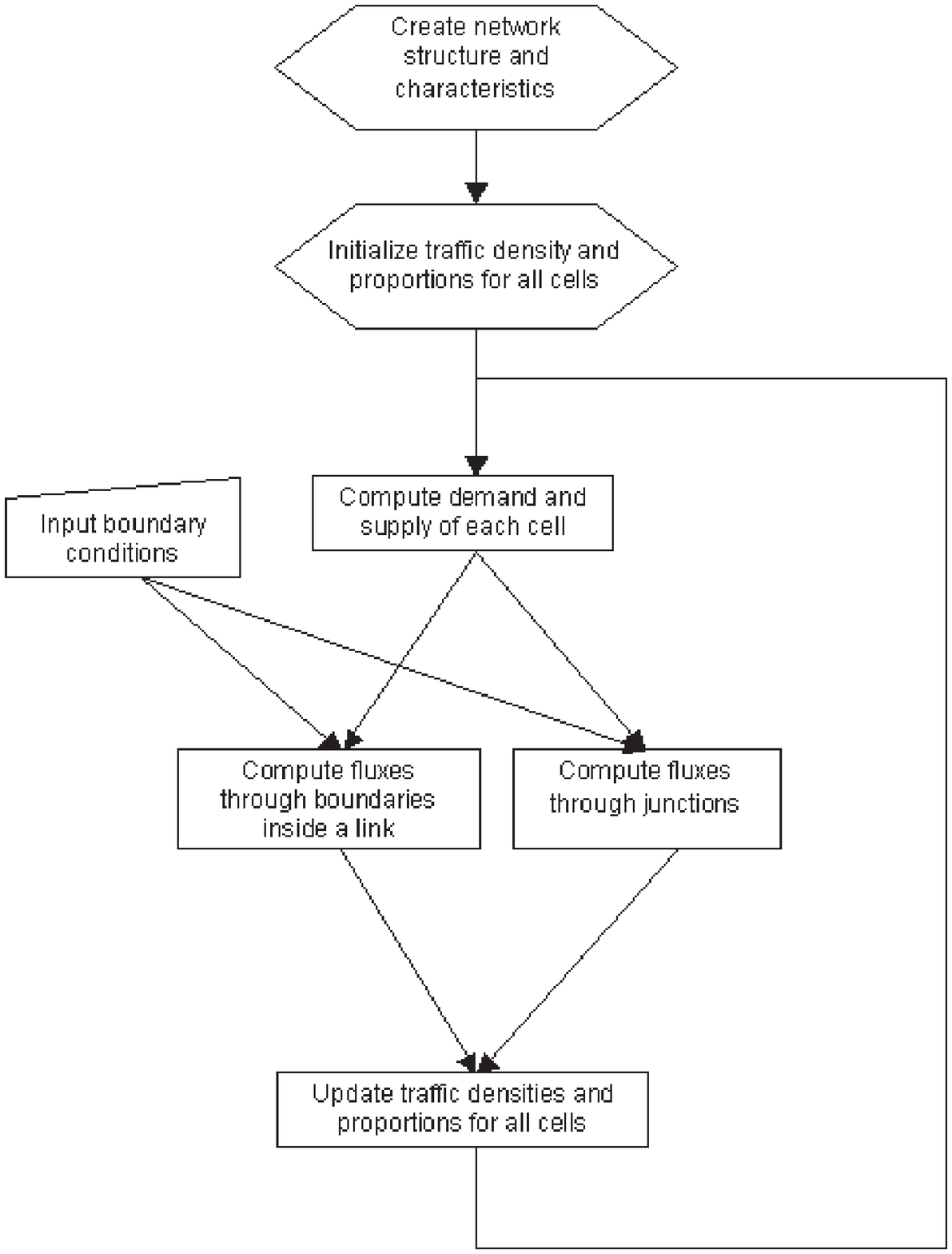}\ec\caption{The program flow chart in the MCKW simulation}\label{F_program}\efg

\ben
\item {\bf Create network.} Network structure and characteristics are created. That is, we provide values for the blue fields of each link and cell in \reffig{F_data}. We also assign all locations of the pointers.
\item {\bf Initialize traffic.} Traffic density and proportions of all commodities are initialized for each cell. One typical initialization is to set network empty; i.e., traffic density of each cell is zero.
\item {\bf Compute supply/demand.} Given traffic density of a cell, we are able to compute traffic demand and supply according to \refe{demands} and \refe{supplys}.
\item {\bf Resolve boundary conditions.} Several types of boundary conditions can be used. The first and most important type of boundary conditions is conditions for origins and destinations. In the MCKW simulation, we use traffic demand, specified for all commodities, for origins and traffic supply for destinations. This is different from the boundary conditions used in the previous chapters, where the Dirichlet, Neumann, and periodic boundary conditions are generally used. However, from those boundary conditions, we can easily compute the supply and demand according to \refe{demands} and \refe{supplys}. Moreover, these boundary conditions are different from the O/D flow matrix in dynamic loading studies in the sense that we cannot determine in-flow, which is also affected by current traffic conditions on a link incident to an origin. Therefore, with the same pattern in origin demand, we may have totally different O/D flows. This observation suggests a criteria for evaluating the level of service of a road network: the amount of time for loading a number of vehicles. In addition, the effect of assignment algorithms can be studied through the proportions of origin demands. Second, signals and accidents are considered as boundary conditions in the MCKW simulation: signals acting at cell boundaries put a constraint on supplies of the downstream cells and demands of the upstream cells \citep{daganzo1995ctm,jin2003merge}, and accidents will change the number of lanes and free flow speed in a cell. From the discussions above, we can see that the influence of incidents and accidents can be studied through imposing different boundary conditions.
\item {\bf Compute link flows.} From \refe{sd_boundary}, we are able to compute fluxes through boundaries inside a link. From the FIFO principle of traffic flow in the kinematic wave theories, the proportion of a commodity in fluxes is equal to that in density in the upstream cells.
\item {\bf Compute junction flows.} From \refe{sd_merge}, \refe{diverge2.1}-\refe{diverge2.3}, \refe{diverge1}, and \refe{diverge3}, we can compute fluxes through different types of junctions. The proportions of different commodities can be obtained from the FIFO property and traffic conservation of each commodity, \refe{tc_commodity}. Since links share junctions, we only need compute junction fluxes for a set of links. For example, for the road network in \reffig{F_network}, we only compute junction fluxes at the upstream and downstream junctions of links $2$ and $5$.
\item {\bf Update traffic conditions.} Traffic densities of all commodities in a cell can be updated by \refe{ceCommodity}. However, one has to be careful when computing the proportions since total density, as a divider, may be very small. Although the proportions may not be accurate for very small densities, it rarely matters.
\een

%Wenlong Jin: 2003
\section{Cumulative flow, travel time, and other properties of a road network}\label{sec:sim4}
In the MCKW simulation, we keep track of the change of traffic densities of all cells and fluxes through all boundaries. Besides, these quantities are specified for commodities. In this section, we will discuss how to obtain other traffic information from these quantities.

\subsection{Cumulative flow and vehicle identity}
Cumulative flow at a boundary $x_{i-1/2}$ from time $t_0$ to $t$, $N(x_{i-1/2};[t_0,t])$, is the total number of vehicles passing the spot during the time interval. If the flux is $f^{\ast}(x_{i-1/2},s)$ at time $s$, then we have
\bqn
N(x_{i-1/2};[t_0,t])&=&\int_{s=t_0}^t \: f^{\ast} (x_{i-1/2},s) \: \m{ds}.
\eqn
Correspondingly, the discrete cumulative flow, $N(x_{i-1/2};[J_0,J])$, which is from time steps $J_0$ to $J$, is defined as
\bqn
N(x_{i-1/2}; [J_0, J]) &=& \sum_{j=J_0} ^ {J-1}  \: f^{j\ast}_{i-1/2} \: \dt,
\eqn
where $f^{j\ast}_{i-1/2}$ is the flux at $x_{i-1/2}$ during time steps $j$ and $j+1$.

A curve of cumulative flow versus time is also known as a Newell-curve or simply N-curve \citep{daganzo1994ctm}, since  \citet{newell1993sim} developed a simplified version of the LWR kinematic wave theory based on this concept.

From the definition of cumulative flow, we can see that an N-curve is non-decreasing in time. Further, it is increasing when passing flow is not zero.

Although densities and fluxes are quantities at the aggregate level, the MCKW model is capable of tracking traffic information at the disaggregate level. This can be done also with cumulative flows: a vehicle passing a cell boundary at a time step can be labeled by the corresponding cumulative flow. If all cumulative flows are synchronized; for example, when the initial traffic in a road network is empty, then the same cumulative flow of a commodity refer to the same vehicle. This fact is due to the FIFO property in all commodities. \footnote{When type $4$ diverge appears, this has to be checked.}

Therefore, in the MCKW simulation, with curves of cumulative flows as a bridge between the aggregate and disaggregate quantities, we are able to keep track of vehicle trajectories, accurate to the order of $\dx$ and $\dt$, from cumulative flows at all cell boundaries. Further, with finer partition of each link, we can obtain more detailed information at the disaggregate level.

\subsection{Travel time}
For a vehicle, which can be identified by its commodity cumulative flow number under FIFO, its travel time across a link or from the origin to the destination can be inferred from N-curves. For example, when we know its arrival and departure times to a link from the corresponding N-curves, we can easily compute its travel time across the link.

This can be demonstrated in \reffig{F_traveltime}. In this figure, the left curve is the N-curve at location $x_1$, and the right curve at $x_2$. These two curves are synchronized in the sense that the vehicles between $x_1$ and $x_2$ at $t=0$ are not counted in $N(x_2; [0, t])$. Therefore, from FIFO, we can see that the $N_0$ vehicle on the left N-curve is the same as the $N_0$ vehicle on the right N-curve. Then, from the curve, we know that the times of the $N_0$ passing $x_1$ and $x_2$ are $t_1$ and $t_2$ respectively. Thus, its travel time from $x_1$ to $x_2$ is $t_2-t_1$.

\bfg\bc\includegraphics[height=12cm] {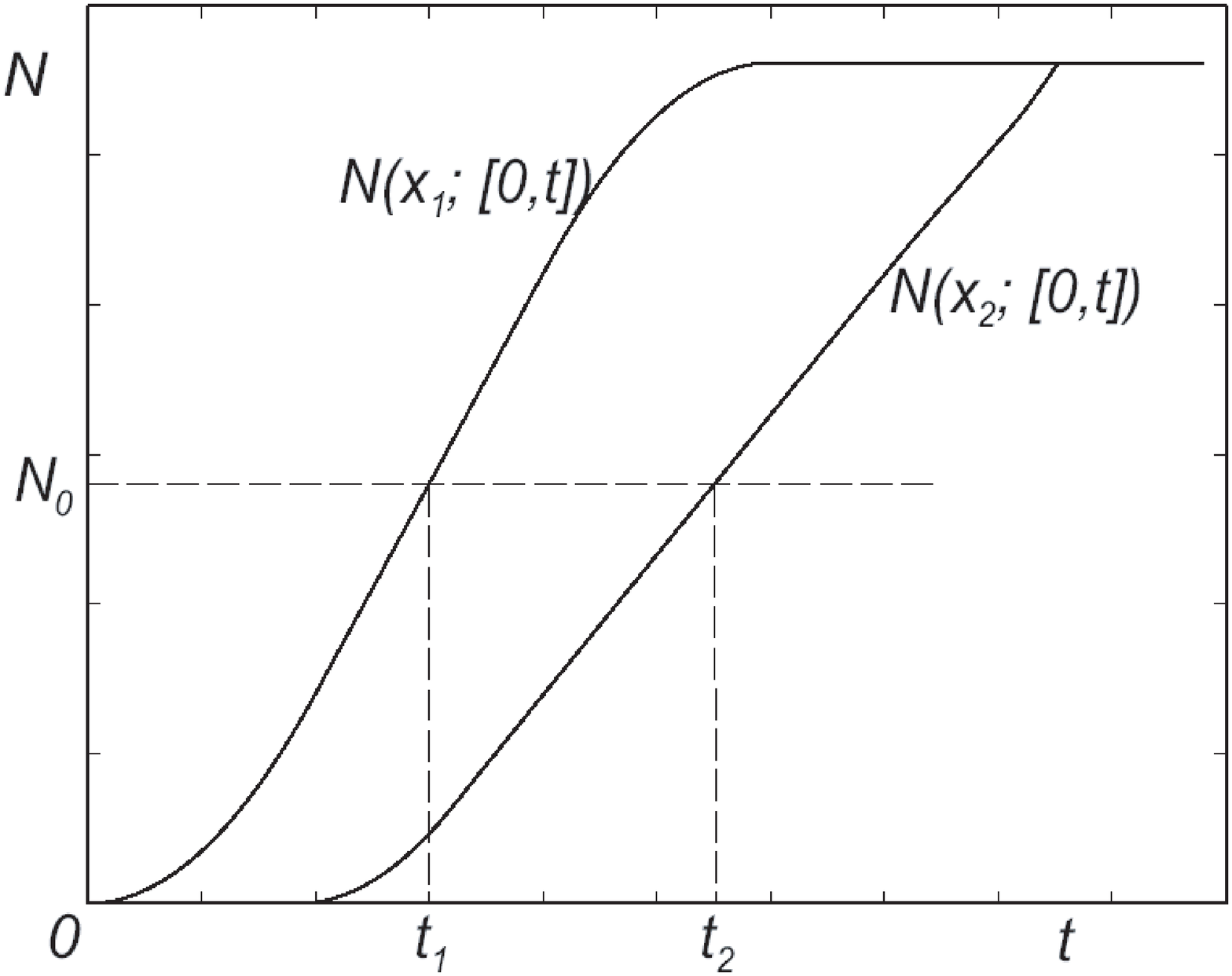}\ec\caption{Cumulative flows and travel time}\label{F_traveltime}\efg

In \reffig{F_traveltime}, the left N-curve reaches a maximum at some time and stop increasing after that. This means that no flow passes $x_1$ after that time. The right N-curve has the same pattern. In such cases, one has to be cautious when computing travel time for the last vehicle, identified by the maximum cumulative flow, which corresponds to  multiple values in time. Rigorously, therefore, the time for a vehicle $N_0$ passing a location $x$, where the N-curve is $N(x; [t_0, t])$, can be defined by
\bqn
T(N_0; x)&=&\min_s \{s \big | \m{when } N(x; [t_0,s])=N_0 \}. \label{def_passingtime}
\eqn
Further, the travel time for the $N_0$ vehicle from $x_1$ to $x_2$ is
\bqn
T(N_0; [x_1, x_2]) &=& T(N_0; x_2)-T(N_0; x_1). \label{ind_tt}
\eqn

With the definition of passing time in \refe{def_passingtime}, at $x$, the vehicle identity $N_0$ has a one-to-one relationship with its passing time $T(N_0; x)$. Therefore, the passing time can be considered as another identity of a vehicle. For a vehicle $N_0$, if we know its passing time at any location in a road network, we then obtain its trajectory.

From the travel times of individual vehicles, we are able to compute the total travel time between two locations, in particular between an O/D pair, as follows:
\bqn
T([N_1,N_2]; [x_1,x_2]) &=& \sum_{M=N_1}^{N_2} \: T(M; [x_1, x_2]), \label{total_tt}
\eqn
where $N_1$ is the first vehicle and $N_2$ the last. We can see that, in \reffig{F_traveltime}, the total travel time is equal to the area between the two N-curves. Then the average travel time for each vehicle will be
\bqn
\bar T([N_1,N_2]; [x_1,x_2]) &\equiv& \frac{T([N_1,N_2]; [x_1,x_2])}{N_2-N_1}=\frac {\sum_{M=N_1}^{N_2} \: T(M; [x_1, x_2])}{N_2-N_1}.\label{ave_tt}
\eqn

Moreover, for a road network, we can integrate travel times for all O/D pairs and, therefore, obtain the total travel time and the average travel time for the whole road network. These quantities are important indicators of the performance of a road network. Besides, we consider the loading time for an amount of flow to be released from an origin as another performance indicator.

Hence, the MCKW simulation platform can be applied to evaluate traffic management and control strategies, such as route assignment and ramp metering algorithms. These applications will be discussed in the next chapter. 
%Wenlong Jin: 2003
\section{Numerical simulations}\label{sec:sim5}
In this section, we investigate the properties of the MCKW simulation model through numerical simulations. We will show the evolution of traffic and examine the convergence of solutions. In this section, the diverge connecting links 2, 3, and 4 is modeled by \refe{diverge1}.

\subsection{Simulation set-up}
For these simulations, the network has the structure as shown in \reffig{F_network}. In this network, links 2, 3, and 5 have the same length, 20 miles, and the length of link 4 is 40 miles (not drawn to proportion); link 2 has three lanes, and the other links has two lanes; all links have the same triangular fundamental diagram \citep{newell1993sim}:
\bqn
Q(a,\r)&=& \cas{{ll} v_{f}\r,& 0\leq\r\leq a \r_{c};\\\frac{\r_{c}}{\r_{j}-\r_{c}} v_{f} (a\r_{j}-\r),& a\r_{c}<\r\leq a\r_{j};}
\eqn
where $\r$ is the total density of all lanes, $a$ the number of lanes, the jam density $\r_j$=180 vpmpl, the critical density $\r_c$=36 vpmpl, the free flow speed $v_f$=65 mph, the capacity of each lane $q_c=\r_c v_f$=2340 vphpl, and the corresponding shock wave speed of jam traffic is $c_j=-\r_c/(\r_j-\r_c) v_f\approx-17$ mph.

Initially, the road network is empty. Boundary conditions are defined as follows. Traffic supply at the destination is always $2q_c$. At origin, traffic demand at the origin is $3q_c$ during $[0, 6.0]$ and zero after that, and the proportion of commodity 0, which takes link 3 instead of 4, is always $\xi=70\%$.

Links 2, 3, and 5 are partitioned into $N$ cells each, and link 4 into $2N$ cells, with each cell of the same length, $\dx=20/N$ miles. The total simulation time of 8.4 hours is divided into $K$ time steps, with the length of a time step $\dt=8.4/K$. In our simulations, we set $N/K=1/30$. Thus the CFL \citep{courant1928CFL} number is no larger than $v_f \dt/dx=0.91$, which is valid for Godunov method \citep{godunov1959}.

\subsection{Traffic patterns on the road network}\label{C6S5_pattern}
We let $N=400$ and $K=12000$. Hence $\dx=0.05$ miles=80 meters, and $\dt=0.0007$ hours=2.52 seconds. Here the sizes of road cells and time steps are relatively small, in order for us to obtain results closer to those of the kinematic wave theories with the Godunov method.

The contour plots of the solutions are shown in \reffig{contour_C6S5}. From these figures, we can divide the evolution of traffic dynamics on the road network into three stages.

\bfg\bc\includegraphics[height=12cm] {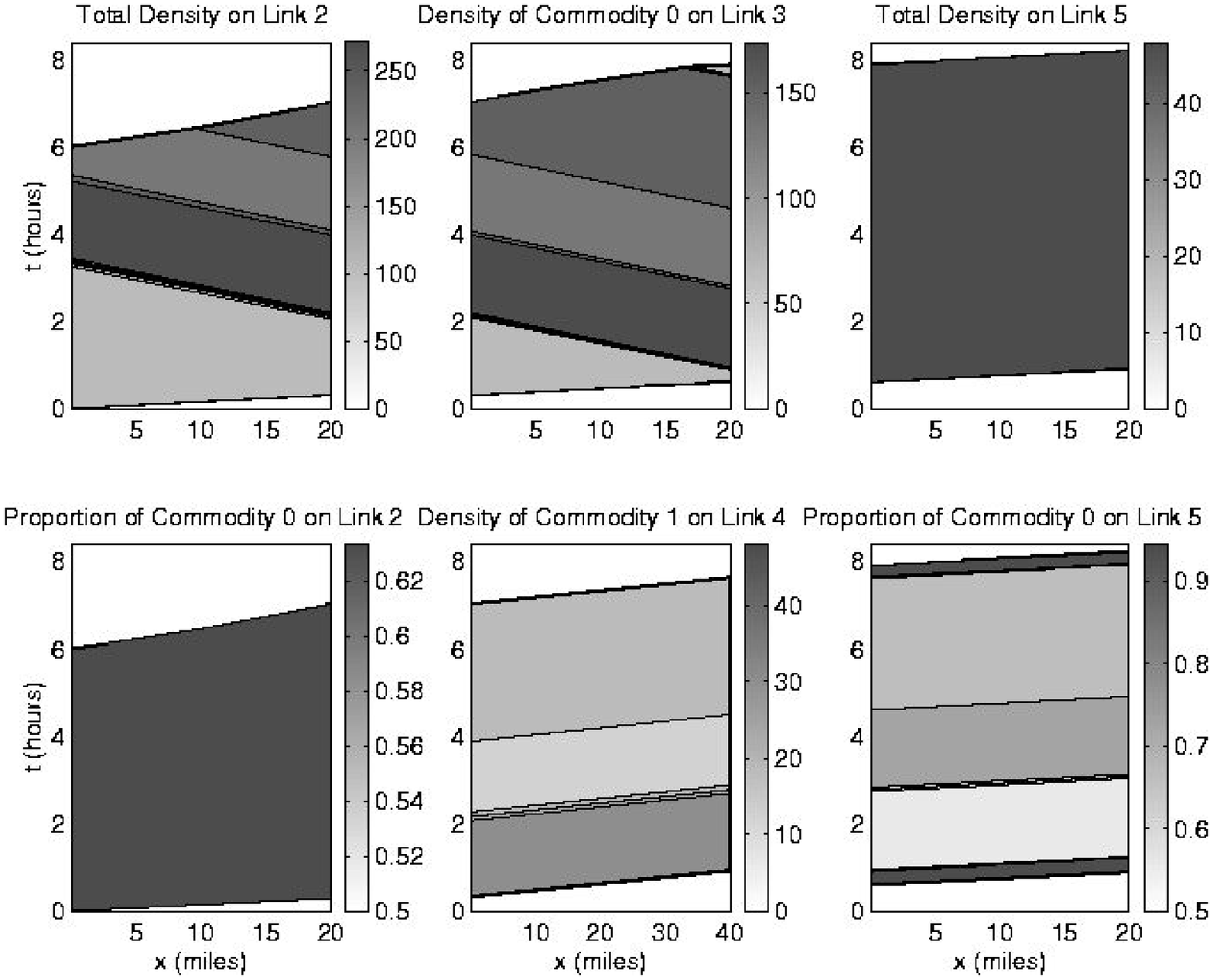}\ec\caption{Contour plots of network traffic flow}\label{contour_C6S5}\efg

In the first stage starting from 0, vehicles embark link 2 with the free flow speed, prevail the link in its critical density $3\r_c$, and arrive junction 1 at $t_1=20/65$ hr. At the diverge, junction 1, fluxes are computed from \refe{diverge1}: out-flux of link 2 is
\bqs
f_{2,out}&=&\min\{3q_c,\frac{2q_c}{0.7},\frac{2q_c}{0.3}\}=\frac{20}7 q_c,
\eqs
which is slightly smaller than the in-flux of link 2,
in-flux of link 3 is $f_3=0.7 f_2=2 q_c$, which is its capacity, and in-flux of link 4 is $f_4=0.3 f_2=\frac 67 q_c$, which is less than half of its capacity. After $t_1$, two streams of free flow form on links 3 and 4, and a backward travelling shock wave forms on Link 2, and the shock wave speed is
\bqs
v_j&=&-\frac 14 v_f.
\eqs
At $t_2=t_1+20/v_f$, the first vehicle on link 3 reaches junction 2, which is a merge. At $t_2$, the first vehicle on link 4 is half way back since the length of link 4 is double of link 3's. From the merge traffic flow model, \refe{sd_merge}, we have the in-flux of link 5 as
\bqs
f_{5,in}&=&\min\{2q_c, 2q_c\}=2q_c,
\eqs
which is also the out-flux of link 3. After $t_2$, the proportion of commodity 0 on link 5 is 1, as we can see on the bottom right figure.

The second stage starts at $t_3=t_1+40/v_f=60/65$ hr, when the first vehicle on link 4 reaches junction 2. After that, the in-flux of link 5 is still $2q_c$, but the proportion of commodity 0 reduces to 0.5714 since commodity 1 also contributes; on link 3, a new state forms at $\r=195.4290$ vpm, which is over-critical, and a shock wave travels upstream at the speed of $|c_j|\approx 17$ mph; on link 4, $\r=30.8571$ vpm, which is under-critical. At $t_4=t_3+20/|c_j|=140/65$ hr, the back-traveling shock on link 3 hits junction 1, and the traffic supply on link 3 is reduced. Therefore, the out-flux of link 2 is further reduced, and link 2 becomes more congested, as shown in the top left picture. This also reduces traffic flow on link 4, and the reduced flow reaches junction 2 at $t_5=t_4+40/v_f=180/65$ hr. After $t_4$, link 3 becomes less congested, and a rarefaction wave travels backward on it at $|c_j|\approx17$ mph; the proportion of commodity on link 5 gets higher. From the bottom middle figure, we can see that at $t_5=t_4+20/|v_j|$, traffic density on link 4 swings back a little due to the back traveling rarefaction on link 3. This shift is transported to junction 2 at $t_6=t_5+40/v_f$ and oscillates traffic density on link 3 and the proportion of commodity 0 on link 5.

The third stage starts at $t_7=6$, when traffic demand from origin subsides to zero. After that, a shock forms on link 2 and travels forward, and propagates to link 4 and link 3. On link 4 the shock travels at $v_f$, and on link 3 it travels slower. This is why the proportion of commodity 0 on link 5 becomes 1 before it is emptied.

This simulation indicates that oscillation of traffic conditions can be caused by network merges and diverges even the initial and boundary conditions are very nice. The traffic flow pattern on this road network suggests that, if we keep the same demand from the origin, an equilibrium state will be reached after some time. This equilibrium state will be further investigated in the following chapter.

The traffic patterns in \reffig{contour_C6S5} can be partially observed from the top left figure in \reffig{ncurves_C6S5}, where the thicker four curves give cumulative flows for commodity 0, the thinner for commodity 1, and solid, dashed, dotted, and dash-dot curves are for cumulative flows at junction 0, 1, 2, and 3, respectively. In the bottom left figure, the solid lines are cumulative flows of commodity 0 at origin/destination, i.e., junction 0 and 3, the dashed curve is the number of commodity 0 vehicles in the network at a time, the dashed line shows the average number, the dotted curve is the travel time of a commodity 0 vehicle identified by its cumulative flow, and the dotted line is the average travel time. The bottom right figure has the same curves and lines as the bottom left figure, except that they are for commodity 1. Here travel time of individual vehicle is computed by \refe{ind_tt}, total travel time by \refe{total_tt}, and average travel time by \refe{ave_tt}. Although the formulas were developed for link travel times, it is valid for a network as long as vehicles of each commodity observe FIFO principle. Total (TTT) and average (ATT) travel times are listed in \reft{t_tt} (unit=hours).

\btb
\bc\begin{tabular}{|l||c|c|c||}\hline
& Total Number of Vehicles & TTT& ATT\\\hline
Commodity 0& 23,8585&4.7291$\times 10^4$& 1.9822\\\hline
Commodity 1& 10,2251&1.7372$\times 10^4$& 1.6989\\\hline
\end{tabular}\ec
\caption{Total travel time (TTT) and average travel time (ATT) for two commodities}\label{t_tt}
\etb

\bfg\bc\includegraphics[height=12cm] {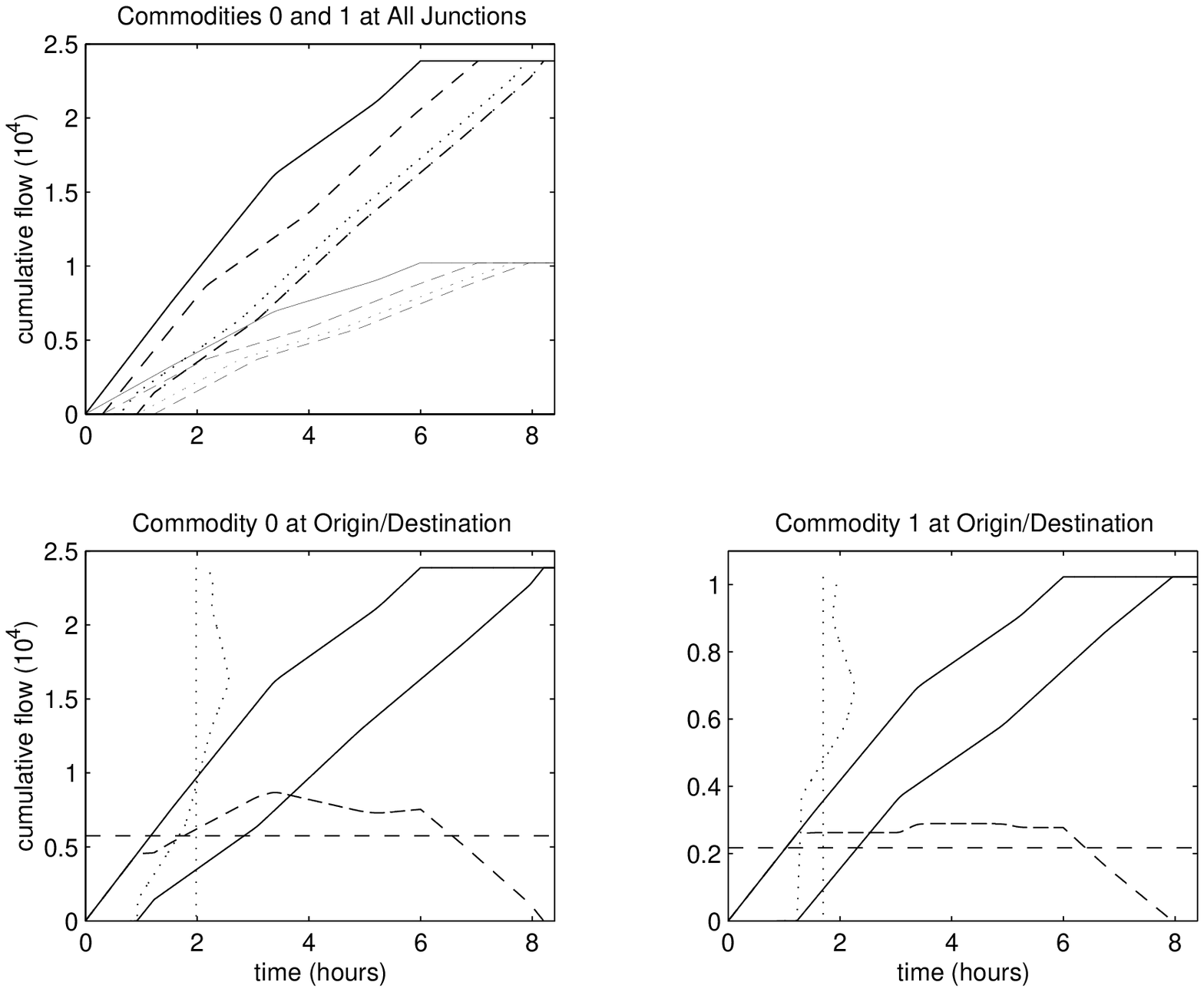}\ec\caption{N-curves and travel times of each commodity in the road network}\label{ncurves_C6S5}\efg

Thus, vehicles that take link 4 on average use shorter time, which is still longer than the free flow travel time, $80/65$ hr. Obviously, if all vehicles at origin decide to take link 3, the travel time, $60/65$ will be the shortest possible travel time between the origin and destination. Therefore, this assignment fraction, $70\%$, is not an optimum one. More detailed analysis of the influence of the assignment fraction on travel times will be engaged in the following chapter.

\subsection{Convergence of the MCKW simulation model}\label{convergence_MCKW}
In this subsection, we study the convergence of the MCKW simulation platform with increasing number of cells. Here we use the same road network, initial and boundary conditions as in the preceding subsection.

As we show in the previous chapters, the discrete forms of kinematic wave theories converge to their continuous counterparts as we partition each link into finer cells \citep{jin2003merge}. In those studies, we generally check convergence in terms of traffic densities. Here we intend  to show convergence in average travel times of both commodities.

Denoting the average travel time of a commodity, $T$, as a function of the number of cells, $N$; i.e., $T=T(N)$, we can define the relative error, from $N$ to $2N$, by
\bqn
\epsilon^{2N-N}&=&|T(2N)-T(N)|.
\eqn
Then a convergence rate is computed by
\bqn
r&=&\log_2(\frac{\epsilon^{2N-N}}{\epsilon^{4N-2N}}).
\eqn
The convergence rates of the average travel times are given in \reft{convergence}.

\btb
\bc\begin{tabular}{|l||c|c|c|c|c||}\hline
Commodity 0& N=200 & N=400& N=800&N=1600&N=3200\\\hline
ATT& 1.98189893& 1.98215215&1.98227240&1.98234941&1.98239377\\\hline
Error $[10^{-3}]$& &0.2532  & 0.1202&0.07700&   0.0444\\\hline
Rate&&&{\bf 1.074}  & {\bf 0.6430}&{\bf 0.7958}\\\hline\hline
Commodity 1& N=200 & N=400& N=800&N=1600&3200\\\hline
ATT& 1.69922958& 1.69892887&1.69877593&1.69871236&1.69868722\\\hline
Error $[10^{-3}]$& &-0.3007&  -0.1529&  -0.0636&  -0.0251\\\hline
Rate&&&{\bf 0.9755}  & {\bf 1.2664}&   {\bf 1.3384}\\\hline
\end{tabular}\ec
\caption{Convergence rates for the MCKW simulation platform}\label{convergence}
\etb

From the table, we can see that average travel times are also convergent in first order. Note that this convergence is different from the aforementioned traffic conditions converging to certain equilibrium state. Moreover, we can see that the results with $N=200$ is already accurate enough in this case. Since the computation time of the MCKW simulation platform is quadrupled when $N$ is doubled, in later simulations, we use $\dx=0.1$ mile=160 meters and $\dt=0.0014$ hours= 5.04 seconds with the same simulation period. 
%Wenlong Jin: 2003
\section{Discussions}\label{sec:sim6}
In this chapter, we proposed the Multi-Commodity Kinematic Wave (MCKW) simulation model. In this simulation model, we integrated the kinematic wave theories studied in the previous three chapters to form the foundation of the algorithms, carefully discussed commodity-based kinematic wave theories, and presented the data structure and program structure for implementation. We further demonstrated how to obtain cumulative flows and travel times from outputs of the MCKW simulation model. Simulations show that numerical results converge to FIFO solutions although the FIFO condition is not strictly enforced in the discrete form of commodity-based kinematic wave theories.

Different from many existing simulation packages, where traffic is tracked down to vehicle level, the MCKW simulation concerns traffic conditions down to commodity level. The simulation model is designed for handling very large road networks and can be applied in studies of Intelligent Transportation Systems, such as dynamic traffic assignment, dynamic O/D estimation, and so on. However, as pointed out earlier, the effects of ``departure from FIFO" should be carefully considered in these applications.

In the future, the MCKW simulation model can be enhanced in three aspects. Theoretically, vehicle types and special lanes can be incorporated \citep{daganzo2002behavior}, and nonequilibrium continuum models \citep{jin2003arz} may also be integrated. Numerically, parallel algorithms can be applied to improve computational speed since traffic conditions on different links can be updated simultaneously, and consumption of computer memory will be checked. Finally, for different applications, we also plan to design different input/out interfaces. For example, the network structure can be imported from GIS (Geographic Information System) data, and boundary conditions and output can be manipulated for different applications.  

\newpage
\pagestyle{myheadings}
\markright{  \rm \normalsize CHAPTER 7. \hspace{0.5cm}
 STUDIES OF NETWORK VEHICULAR TRAFFIC}
\chapter{Studies of network vehicular traffic with kinematic wave simulations}
%\thispagestyle{myheadings}
%Wenlong Jin: 2003
%Wenlong Jin: 2003
\section{Introduction}
As discussed in the preceding chapter, the Multi-Commodity Kinematic Wave (MCKW) simulation model is based on the kinematic wave theories of traffic dynamics at various road network components, including link bottlenecks, merges, and diverges. After vehicles are differentiated into commodities by their paths or origin/destination (O/D) pairs, we can then keep track of the evolution of traffic densities of all commodities on a road network over time. Further, information of individual vehicles can be obtained through cumulative flows. Therefore, compared to some other simulation models based on the kinematic wave theory (e.g., \cite{daganzo1995ctm, vaughan1984od}), the MCKW simulation  is a pure macroscopic simulation program. As such, the MCKW simulation model is expected to be less costly in computation but adequate for many applications in Advanced Transportation Information Systems (ATIS) and Advance Transportation Management Systems (ATMS).

In this chapter, we will explore traffic dynamics in a road network with the MCKW simulation model and its implications in applications. First in Section \ref{sec:app2}, we are interested in equilibrium states of a road network under a given O/D demand pattern and the influence of different assignment strategies. Then in Section \ref{sec:app3}, we study the formation mechanism of periodic oscillations in a network and their properties.

These simulations are studied on a small, demonstration network, where traffic dynamics at diverges obey \refe{diverge1}. The major purpose of these studies is to show the capabilities of the MCKW simulation model. Moreover, these studies can be considered as initial steps to the understanding of more complicated traffic phenomena in a large-scale road network and  practical applications.

%Wenlong Jin: 2003
\section{Equilibrium states of a road network and preliminary examination of traffic assignment}\label{sec:app2}
In Subsection \ref{C6S5_pattern}, we show that traffic flow on a road network can reach a certain pattern after sufficient amount of time. We call these time-invariant traffic patterns as equilibrium states.

For the equilibrium states, we are interested in the following questions: 1) how are the equilibrium states related to the boundary conditions? 2) what is the effect of changing proportions in origin demand? 3) do these equilibrium states have anything to do with the so-called User Equilibrium States?

\subsection{The simulated network}
In this section, we study a demonstration network shown in \reffig{F_assignmentNetwork}. In this network, links 2, 3, and 5 have the same length, 20 miles, and the length of link 4 is 40 miles; link 2 has three lanes, and the other links has two lanes; all links have the same fundamental diagram:
\bqs
Q(a,\r)&=& \cas{{ll} v_{f}\r,& 0\leq\r\leq a \r_{c};\\\frac{\r_{c}}{\r_{j}-\r_{c}} v_{f} (a\r_{j}-\r),& a\r_{c}<\r\leq a\r_{j};}
\eqs
where $\r$ is the total density of all lanes, $a$ the number of lanes, the jam density $\r_j$=180 vpmpl, the critical density $\r_c$=36 vpmpl, the free flow speed $v_f$=65 mph, and the capacity of each lane $q_c=\r_c v_f$=2340 vphpl.

\bfg\bc\includegraphics[height=3cm]{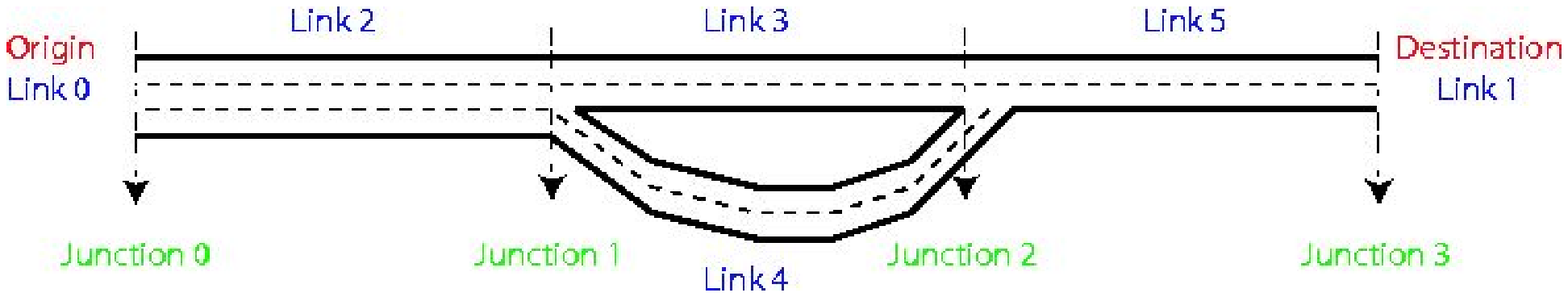}\ec\caption{Network for studying equilibrium state and traffic assignment}\label{F_assignmentNetwork}\efg

Initially, the road network is empty. Traffic supply at the destination is always the capacity of two lanes. For studying equilibrium states, we assume traffic demand at the origin is always the capacity of link 2, $3q_c$.

Links 2, 3, and 5 are partitioned into 200 cells each, and link 4 into 400 cells. Therefore, all cells have the same length, 0.1 miles. The total simulation time of 8.4 hours is divided into $6000$ time steps, with the length of a time step $\dt=0.0014$ hours=5.04 seconds. Thus the CFL \citep{courant1928CFL} number is no bigger than $v_f \dt/dx=0.91$, which is valid for Godunov method \citep{godunov1959}. As we show in subsection \ref{convergence_MCKW}, the simulation results are numerically convergent.

In the following simulations, we will study the equilibrium states and the performance of the road network for a proportion of commodity 0, $\xi$.

\subsection{Equilibrium states}\label{ana_es}
When $\xi=0.5$, and commodity 1 vehicles have the same proportion, $1-\xi=0.5$, solutions of traffic dynamics on the road network are shown in \reffig{network55}. The traffic pattern evolves as follows. From $t_0=0$ to $t_1=20/65$ hr \footnote{Hereafter, the unit of time is always hour except if otherwise mentioned.} , the first vehicle traverses link 2 and reaches junction 1, where fluxes can be computed from \refe{diverge1} as
\bqs
f_{2,out}&=&\min\{3q_c, 2q_c/0.5, 2q_c/0.5\}=3q_c,\\
f_{3,in}&=&1.5q_c,\\
f_{4,in}&=&1.5q_c,
\eqs
where subscript {\em in} denotes in-flux, and {\em out} out-flux. After $t_1$, links 2, 3, and 4 all carry free flow. At $t_2=t_1+20/65$, traffic on link 3 arrives at the merge, junction 2, where from \refe{sd_merge} fluxes are the following:
\bqs
f_{5,in}&=&\min\{2q_c, 2q_c\}=2q_c,\\
f_{3,out}&=&2q_c,\\
f_{4,out}&=&0.
\eqs
At $t_3=t_2+20/65$, vehicles on link 4 also reaches junction 2, and fluxes become
\bqs
f_{5,in}&=&\min\{2q_c+2q_c, 2q_c\}=2q_c,\\
f_{3,out}&=&q_c,\\
f_{4,out}&=&q_c.
\eqs
Thus, after $t_3$, shock waves form on both links 3 and 4 and travel upstream in the speed of $-v_f/4$. At $t_4=t_3+80/65$, the shock wave on link 3 reaches junction 1, and fluxes become
\bqs
f_{2,out}&=&\min\{3q_c, q_c/0.5, 2q_c/0.5\}=2q_c,\\
f_{3,in}&=&q_c,\\
f_{4,in}&=&q_c.
\eqs
After this, link 3 equilibrates at $q_3=q_c$ and $\r_3=1.2\r_j$, a shock travels backward on link 2 at $v_f/4$, and a shock travels forward at the free flow speed. At $t_5=t_4+20/65$, the two shock waves on link 4 meet at the middle of the link and form a stable zero-speed shock wave, where the upstream half is at $0.2\r_j$ and the downstream half at $1.2\r_j$. After $t_6=t_4+80/65$, traffic on link 2 is uniformly at $q_2=2q_c$ and $\r_2=1.4 \r_j$.

\bfg\bc\includegraphics[height=12cm]{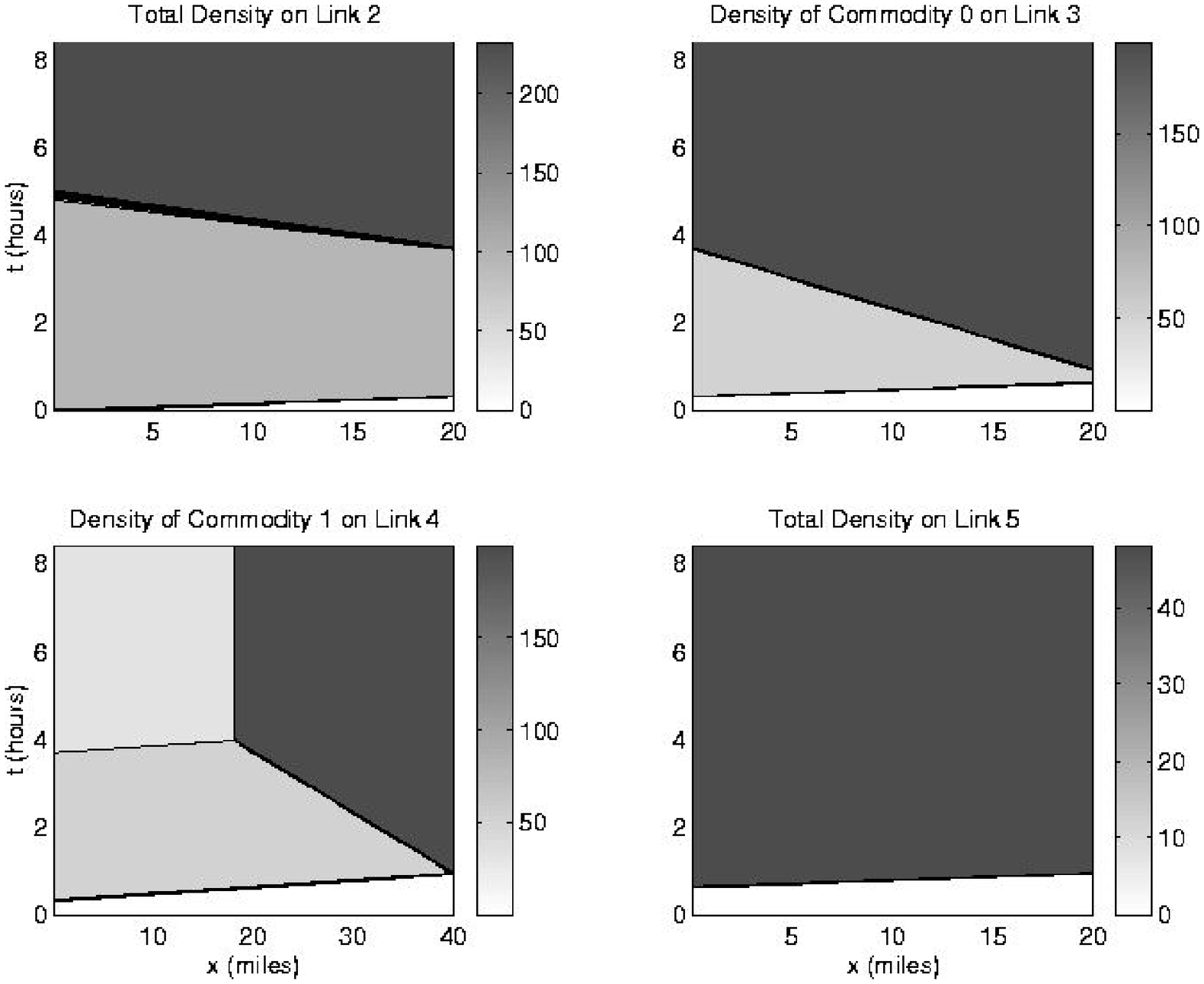}\ec\caption{Solutions when $\xi=0.5$}\label{network55}\efg

Therefore, the whole road network reaches an equilibrium state after $t_6$. In the equilibrium state, links 2 and 5 both carry flow-rate $2q_c$, which is the capacity of the road network, links 3 and 4 carry flow rate $q_c$ corresponding to the proportion of each commodity. Besides, links 2 and 3 are congested, and link 4 has a zero-speed shock, which, however, is unstable in the sense that a small oscillation in the upstream or downstream flow will make the zero-shock disappear. The instability of this equilibrium state on the road network can also been seen later when we consider $\xi$ away from $0.5$.

When $\xi=0.6$, the contour plots of traffic densities on four links are given in \reffig{network64}. The traffic evolution pattern is similar to that in subsection \ref{C6S5_pattern}. Here we directly go into the discussion of the equilibrium state, in which, from \reffig{network64}, we can see that all links carry uniform flows. The equilibrium flow-rate is determined by the network bottleneck, link 5, and is $2q_c$. From observations for $\xi=0.5$, the equilibrium density on link 5 is always $0.4\r_j$, and traffic density on link 2 is $1.4\r_j$, which is brought by a back-traveling shock wave. Then, determined by traffic proportions, the flow-rate on links 3 and 4 are $1.2q_c$ and $0.8q_c$, respectively. The remaining task is to determine the traffic densities on links 3 and 4.

\bfg\bc\includegraphics[height=12cm]{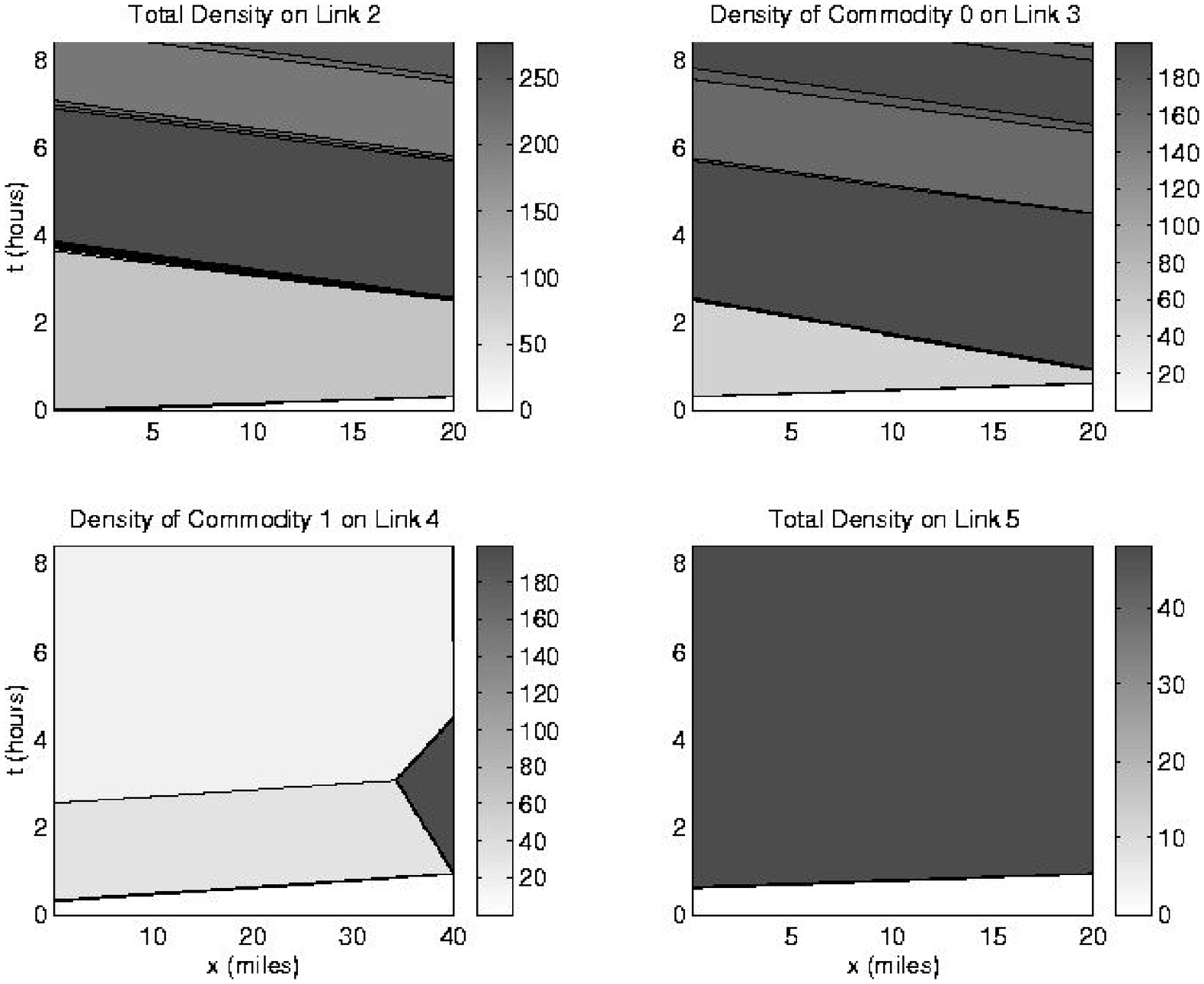}\ec\caption{Solutions when $\xi=0.6$}\label{network64}\efg

We denote traffic demand and supply on these links as $D_3$, $S_3$, $D_4$, and $S_4$. Then $(D_3, S_3)=(1.2 q_c, 2q_c)$ when link 3 is under-critical, and $(D_3, S_3)=(2 q_c, 1.2q_c)$ otherwise. Similarly, $(D_4,S_4)=(0.8 q_c, 2q_c)$ when link 4 is under-critical, and $(D_4,S_4)=(2 q_c, 0.8 q_c)$ otherwise. Since at junction 1
\bqs
f_{2,out}&=&\min\{3q_c, S_3/0.6, S_4/0.4\}=2q_c,
\eqs
links 3 and 4 cannot be under-critical at the same time. Moreover, since at junction 2
\bqn
\ba{lcl}
f_{3,out}&=&\frac{D_3}{D_3+D_4} 2 q_c=1.2 q_c,\\
f_{4,out}&=&\frac{D_4}{D_3+D_4} 2 q_c=0.8 q_c,
\ea\label{network64eq}
\eqn
links 4 cannot be over-critical at the same time. Therefore, the only possible case is that link 3 is over-critical and link 4 under-critical. Thus, the equilibrium densities on links 3 and 4 are $\r_3=1.04\r_j$ and $\r_4=0.32\r_j$. The solutions in \reffig{network64} support these solutions.

However, $D_3=2 q_c$ and $D_4=0.8 q_c$ do not satisfy \refe{network64eq} exactly. In fact, there is an intermediate state at the downstream boundary of link 4 \citep{jin2002continuousmerge}. At the intermediate state, also under-critical, the  demand is $\bar D_4=\frac 43 q_c$, which satisfies \refe{network64eq}, and $\bar \r_4=\frac {0.8} 3 \r_j$. This intermediate state, theoretically, only exists at a point, but is stable since its in-flux and out-flux is equal. This kind of intermediate state unlikely exists on a single link \citep{jin2003inhLWR}.

When $\xi=0.4$, the solutions are shown in \reffig{network46}. Although the evolution process is different from that for $\xi=0.6$, the equilibrium state is exactly the same as before if links 3 and 4 are switched.

\bfg\bc\includegraphics[height=12cm]{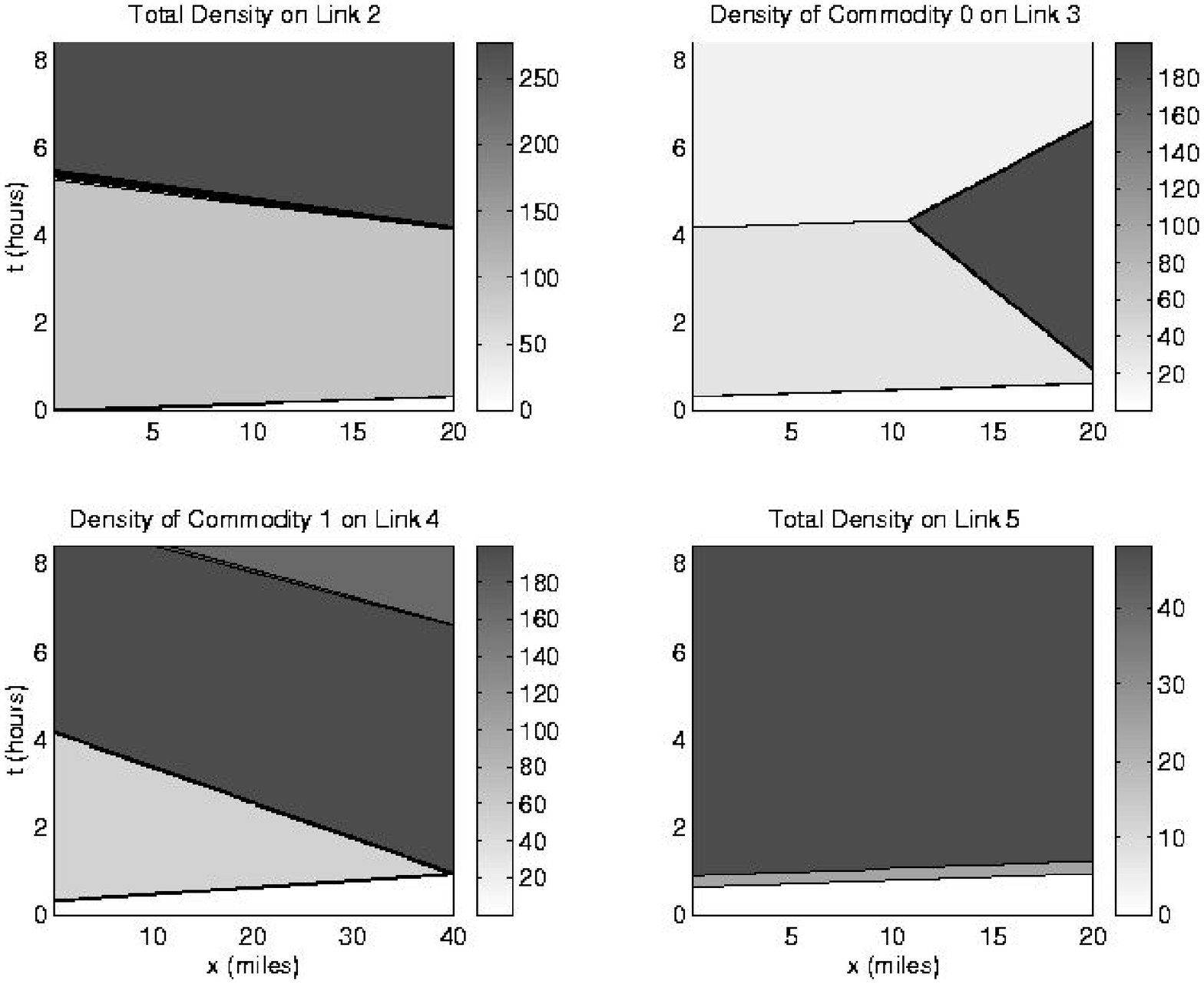}\ec\caption{Solutions when $\xi=0.4$}\label{network46}\efg

For a general proportion $\xi\in[0,1]$, the equilibrium states are listed in \reft{tab_equsta}. In this table, we omit the tenuous (unstable) intermediate states, on link 3 when $0<\xi<0.5$ and on link 4 when $0.5<\xi<1$.

\btb
\bc\begin{tabular}{|c||c|c|c|}\hline
$(\r/\r_j, q/q_c)$ & $0\leq \xi<0.5$ &$\xi=0.5$& $0.5<\xi\leq 1$\\\hline
Link 2& $(1.4, 2)$ &$(1.4, 2)$&$(1.4, 2)$\\\hline
Link 3& $(0.4\xi, 2\xi )$& $(1.2, 1)$  & $(2-1.6\xi, 2\xi)$\\\hline
Link 4& $(0.4+1.6\xi, 2(1-\xi))$ &$(0.2, 1)$ and $(1.2, 1)$ & $(0.4(1-\xi), 2(1-\xi))$\\\hline
Link 5&$(0.4, 2)$&$(0.4, 2)$&$(0.4, 2)$\\\hline
\end{tabular}
\ec\caption{Equilibrium density and flow-rate v.s. $\xi$}\label{tab_equsta}
\etb

\subsection{Travel times at equilibrium states}
For equilibrium states of different proportion $\xi$, \reft{tab_estt} shows travel speed and travel time on each link as well as average travel times of commodity 0 ($ATT_0$), commodity 1 ($ATT_1$), and all commodities ($ATT=\xi ATT_0+(1-\xi) ATT_1$). In the table, travel times during the establishment of equilibrium states are not considered. The relationship between these travel times and proportion $\xi$ is illustrated in \reffig{F_estt}.

\btb
\bc\begin{tabular}{|c||c|c|c|}\hline
$(v/v_f, T/hours)$ & $0\leq \xi<0.5$ &$\xi=0.5$& $0.5<\xi\leq 1$\\\hline
Link 2& $(\frac 27, 1.0769)$ &$(\frac 27, 1.0769)$&$(\frac 27, 1.0769)$\\\hline
Link 3& $(1, 0.3077 )$& $(\frac 16, 1.8462)$  & $(\frac{\xi}{5-4\xi}, \frac{20-16\xi}{13\xi})$\\\hline
Link 4& $(\frac{1-\xi}{1+4\xi}, \frac{8(1+4\xi)}{13(1-\xi)})$ &(1 or $\frac 16$, 2.1538) & $(1, 0.6154)$\\\hline
Link 5&$(1, 0.3077)$&$(1, 0.3077)$&$(1, 0.3077)$\\\hline\hline
$ATT_0$& 1.6923&3.2308 & $\frac{20-16\xi}{13\xi}$+1.3846\\\hline
$ATT_1$&$\frac{8(1+4\xi)}{13(1-\xi)}$+1.3846&3.5384&2\\\hline
$ATT$ &$2+2.7692\xi$&3.3846&$3.5385-1.8462\xi$\\\hline
\end{tabular}
\ec\caption{Equilibrium speed and travel times v.s. $\xi$}\label{tab_estt}
\etb

\bfg\bc\includegraphics[height=12cm]{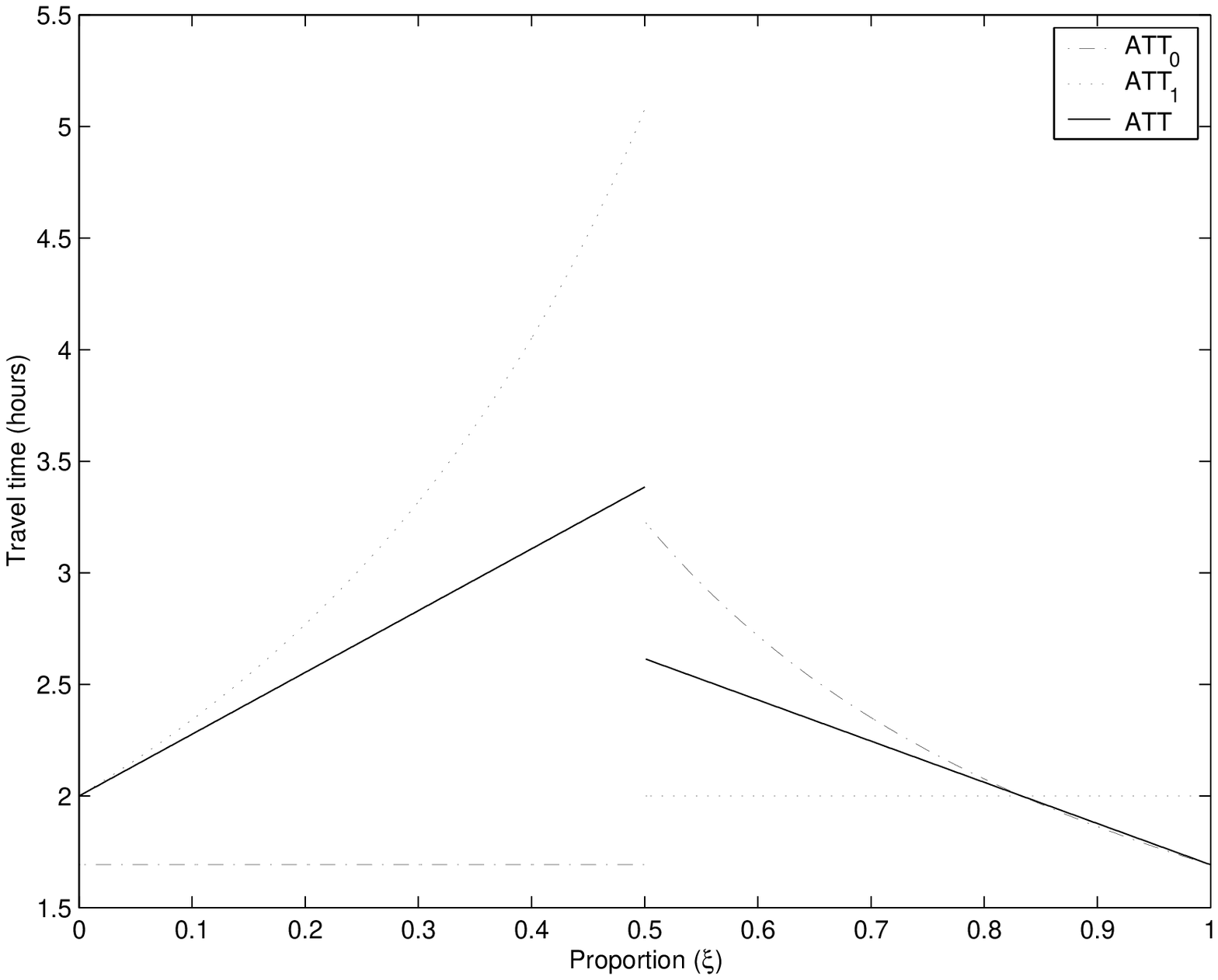}\ec\caption{Travel times at equilibrium states v.s. proportion $\xi$}\label{F_estt}\efg

From the table and figure, we can see the following properties of travel times. First, the average travel time of all commodities, which can be considered the performance function of the whole road network, consists of two pieces of lines and jumps when $\xi=0.5$. This discontinuity asserts the fact that the equilibrium state when $\xi=0.5$ is not stable. Third, the minimum average travel time $ATT$ attains its minimum when $\xi=1$; i.e., when all vehicles take the shorter path. Third, the travel times of commodities 0 and 1 are equal when $\xi=5/6$. However, $ATT$ is not at its minimum when $\xi=5/6$.

\subsection{Discussions}
In the previous subsections, we examine equilibrium states for the road network in \reffig{F_assignmentNetwork} and its performance at these states against the proportion $\xi$. The determination of $\xi$ is really part of a traffic assignment process. Therefore, this study provides another angle for checking the user-equilibrium property \citep{wardrop1952ue}. For this network, it is shown to have two user equilibrium states, one at $\xi=5/6$ and the other at $\xi=1$, with the latter has smaller path travel times. This is a clear evidence that the inclusion of physical queues destroys the uniqueness property of traffic assignment with other commonly used models, in which queues do not take up space \citep{wardrop1952ue}. Note that, in our discussions, the performance function is the average path travel time, which is different from the link performance function used in static traffic assignment studies \citep{sheffi1984networks}.

Although the connection between equilibrium states from a dynamic point of view and the user-equilibrium is still premature, these studies shed light on solving dynamic traffic assignment problems, in which the choice of a proper link performance function is still in debate \citep{daganzo1995tt}. From the MCKW simulation, we can obtain average origin/destination (O/D) travel times under an assignment strategy and use these travel times to evaluate the strategy. To make the performance function of an O/D flow more accurate, one can include the loading time at the origin, which is defined as the average waiting time before a vehicle enters the road network. If departure (or loading) flow-rate at the origin is $f(t)$ ($t\in[0,T]$), then the average loading time is
\bqn
LT&=&\frac{\int_{t=0}^T f(t) \: t \m{ dt}}{\int_{t=0}^T f(t) \m{ dt}}.
\eqn
 In the examples above, since the loading flow-rate is the same for all $\xi$, there is no difference between the loading times.

%Wenlong Jin: 2003
\section{The formation and structure of periodic oscillations in the kinematic wave model of road networks}\label{sec:app3}
In \citep{jin2003inhLWR} it is shown that there are three types of basic kinematic waves on a road link: shock (decelerating) waves, rarefaction (accelerating) waves, and transition (standing) waves. Further, for a single merge or diverge, we still have these three types of waves on each branch \citep{jin2003merge, jin2002continuousmerge,jin2002diverge, jin2001diverge2}. In this section, we show with MCKW simulation an interesting type of solutions, periodic oscillations, which can be observed in real traffic. We show that these solutions can exist in a small road network with a diverge and a merge and check their formation and structure against network characteristics.

\subsection{Network for studying periodic solutions}
In this section, we study periodic solutions on a road network shown in \reffig{F_periodic}. In this network, lengths of links 2, 3, 4, and 5 are $L_2=10$ miles, $L_3=1$ mile, $L_4=2$ miles, and $L_5=1$ mile, respectively; the number of lanes of these links are 3, 1, 2, and 2, correspondingly; all links have the same fundamental diagram for each lane:
\bqs
Q(a,\r)&=& \cas{{ll} v_{f}\r,& 0\leq\r\leq a \r_{c};\\\frac{\r_{c}}{\r_{j}-\r_{c}} v_{f} (a\r_{j}-\r),& a\r_{c}<\r\leq a\r_{j};}
\eqs
where $\r$ is the total density of all lanes, $a$ the number of lanes, the jam density $\r_j$=180 vpmpl, the critical density $\r_c$=36 vpmpl, the free flow speed $v_f$=65 mph, and the capacity of each lane $q_c=\r_c v_f$=2340 vphpl.

\bfg\bc\includegraphics[height=3cm]{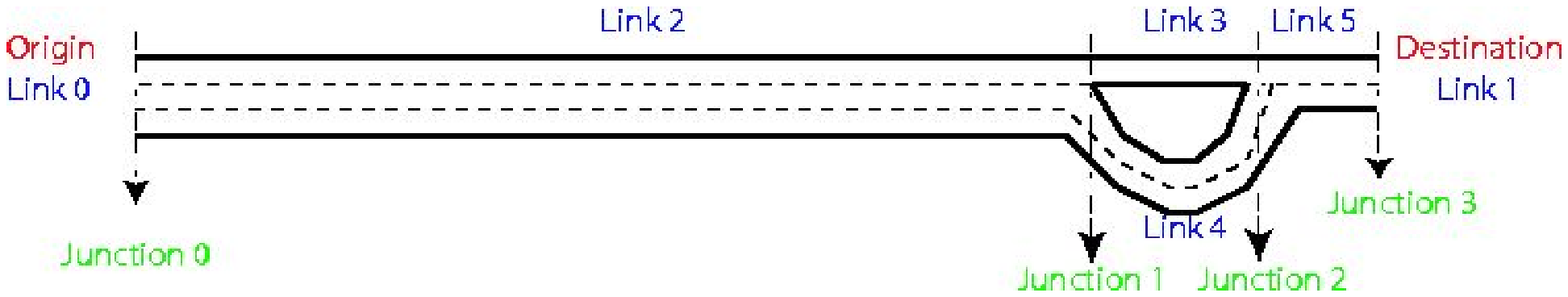}\ec\caption{Network for studying periodic oscillations}\label{F_periodic}\efg

Initially, the road network is empty. Traffic supply at the destination is always capacity of two lanes. Here, traffic demand at the origin is always the capacity of link 2, $3q_c$.

Links 2, 3, 4, and 5 are partitioned into 800, 80, 160, and 80 cells respectively, with the length of all cells as 0.0125 miles. The total simulation time of 1.4 hours is divided into $8000$ time steps, with the length of a time step $\dt=1.75\times 10^{-4}$ hours=0.63 seconds. Thus the CFL \citep{courant1928CFL} number is no bigger than $v_f \dt/dx=0.91$, which is valid for Godunov method \citep{godunov1959}. As we show in subsection \ref{convergence_MCKW}, the simulation is numerically convergent. Therefore, with this very fine partition of links, we are able to obtain results closer to theoretical solutions.

\subsection{Periodic oscillations}
As in the preceding section, we also consider the equilibrium states for a proportion of commodity 0, $\xi$. In equilibrium states, at junction 1, we have
\bqn
\ba{lcl}
f_{2,out}&=&\min\{3q_c, S_3/\xi, S_4/(1-\xi)\},\\
f_{3,in}&=&\xi f_{2,out},\\
f_{4,in}&=&(1-\xi) f_{2,out},
\ea\label{junc1}
\eqn
and at junction 2
\bqn
\ba{lcl}
f_{5,in}&=&\min\{D_3+D_4, 2q_c\},\\
f_{3,out}&=&\frac{D_3}{D_3+D_4} f_{5,in},\\
f_{4,out}&=&\frac{D_4}{D_3+D_4} f_{5,in},
\ea\label{junc2}
\eqn
where $D_i$ and $S_i$ are traffic demand and supply of link $i$ ($i=2,3,4,5$) respectively, $D_2=3q_c$, and $S_5=2q_c$. Moreover, we have
\bqn
\ba{lcl}
f_{3,in}&=&f_{3,out},\\
f_{4,in}&=&f_{4,out},
\ea\label{junc12}
\eqn
and constraints defined in \reft{es_cons}, where $q_3$ and $q_4$ are flow-rates in the upstream cells of links 3 and 4.

\btb
\bc\begin{tabular}{|c|c|c||c|c|c|}\hline
Link 3& under-critical&over-critical&Link 4&under-critical&over-critical\\\hline
$S_3$ & $q_c$& $q_3$ & $S_4$ & $2 q_c$& $q_4$ \\\hline
$D_3$ & $\leq q_c$ &$q_c$ & $D_4$ & $\leq 2 q_c$& $2 q_c$\\\hline
\end{tabular}\caption{Constraints on equilibrium states}\label{es_cons}
\ec
\etb

Since links 3 and 4 are initially empty, $S_3=q_c$ and $S_4=2q_c$. Therefore, when $\xi=0$ or $0.5\leq \xi\leq 1$, $f_{2,out}\leq 2q_c$ according to \refe{junc1}. Note that the capacity of link 5 is $2q_c$. Under this circumstance, solutions are trivial. When $0<\xi<0.5$, from \refe{junc1}, we have $f_{2,out}>2q_c$. This means that at least one of links 3 and 4 has to be over-critical. When both links are over-critical, we have from \refe{junc1}-\refe{junc12} and \reft{es_cons} that $\xi=\frac 13$. When $0<\xi<\frac 13$, link 4 is congested in equilibrium state, and the solution pattern is similar to \reffig{network64} and \reffig{network46}. Refer to subsection \ref{ana_es} for an analysis of the equilibrium states.

For the network in \reffig{F_periodic}, interesting solutions occur for $\frac 13<\xi<\frac 12$, when link 3 is over-critical and link 4 under-critical. When $\xi=0.45$, contour plots of traffic densities on the four links are shown in \reffig{periodicContour}. This figure shows the formation of periodic  solutions as follows: (i) During $t_0=0$ and $t_1=L_3/v_f$, capacity flow travels forward on link 2. When the first vehicle reaches junction 1 at $t_1$, we obtain from \refe{junc1} that $f_{2,out}=q_c/0.45$, $f_{3,in}=q_c$, and $f_{4,in}=\frac{11}9 q_c$. (ii) Thus, after $t_1$, a back-traveling shock wave forms on link 2, capacity flow travels on link 3, and free flow travels on link 4. At $t_2=t_1+L_3/v_f$, the first vehicle on link 3 reaches junction 2. Since the first commodity-1 vehicle is still half way on link 4, from \refe{junc2}, we obtain $f_{5,in}=q_c$, $f_{3,out}=q_c$, and $f_{4,out}=0$. Therefore, there is no change in traffic patterns on links 2 or 3 at $t_2$. (iii) AT $t_3=t_1+L_4/v_f$, vehicles on link 4 reach junction 2. From \refe{junc2}, we have $f_{5,in}=2q_c$, $f_{3,out}=0.9 q_c$, and $f_{4,out}=1.1 q_c$. (iv) After $t_3$, a back-traveling shock forms on link 3 and the over-critical traffic state, whose flow-rate is $0.9 q_c$, propagates upstream. On link 4, there exists a flimsy intermediate area, which is still under-critical. Thus there is no backward wave on link 4. The shock wave on link 3 travels at the speed $\frac 14 v_f$ and hits junction 1 at $t_4=t_3+4L_3/v_f$. (v) At $t_4$, we have $S_3=0.9 q_c$, which yields from \refe{junc1} that $f_{2,out}=2q_c$, $f_{3,in}=0.9 q_c$, and $f_{4,in}=1.1 q_c$. Therefore, a forward shock wave forms on link 4, whose upstream flow-rate is $1.1 q_c$ and downstream flow-rate $\frac {11} 9 q_c$. (vi) At $t_5=t_4+L_4/v_f$, the new traffic state reaches junction 2. Consequently, $f_{3,out}=2/2.1 q_c>0.9 q_c$, and $f_{4,out}=2.2/2.1 q_c<1.1 q_c$. Thus the downstream of link 3 becomes less congested, and a rarefaction wave travels upstream. Note that there is no back-traveling wave on link 4. (vii) At $t_6=t_5+4L_3/v_f$, the rarefaction wave reaches junction 1 and increases out-flux of link 2, $f_{2,out}=2.1164 q_c$, and $f_{4,in}=1.164 q_c>1.1 q_c$. (viii) When the new flow on link 4 reaches junction 2 at $t_7=t_6+L_4/v_f$, link 3 becomes more congested again. After another time period, $(4L_3+L_4)/v_f$, link 4 will go back to lower congestion. This process repeats in a periodic manner with the period $T=2(4L_3+L_4)/v_f$=0.1846 hours. After a number of periods, traffic conditions on links 2, 3, and 4 reach a stable pattern, a period of which is explained as follows.

\bfg\bc\includegraphics[height=12cm]{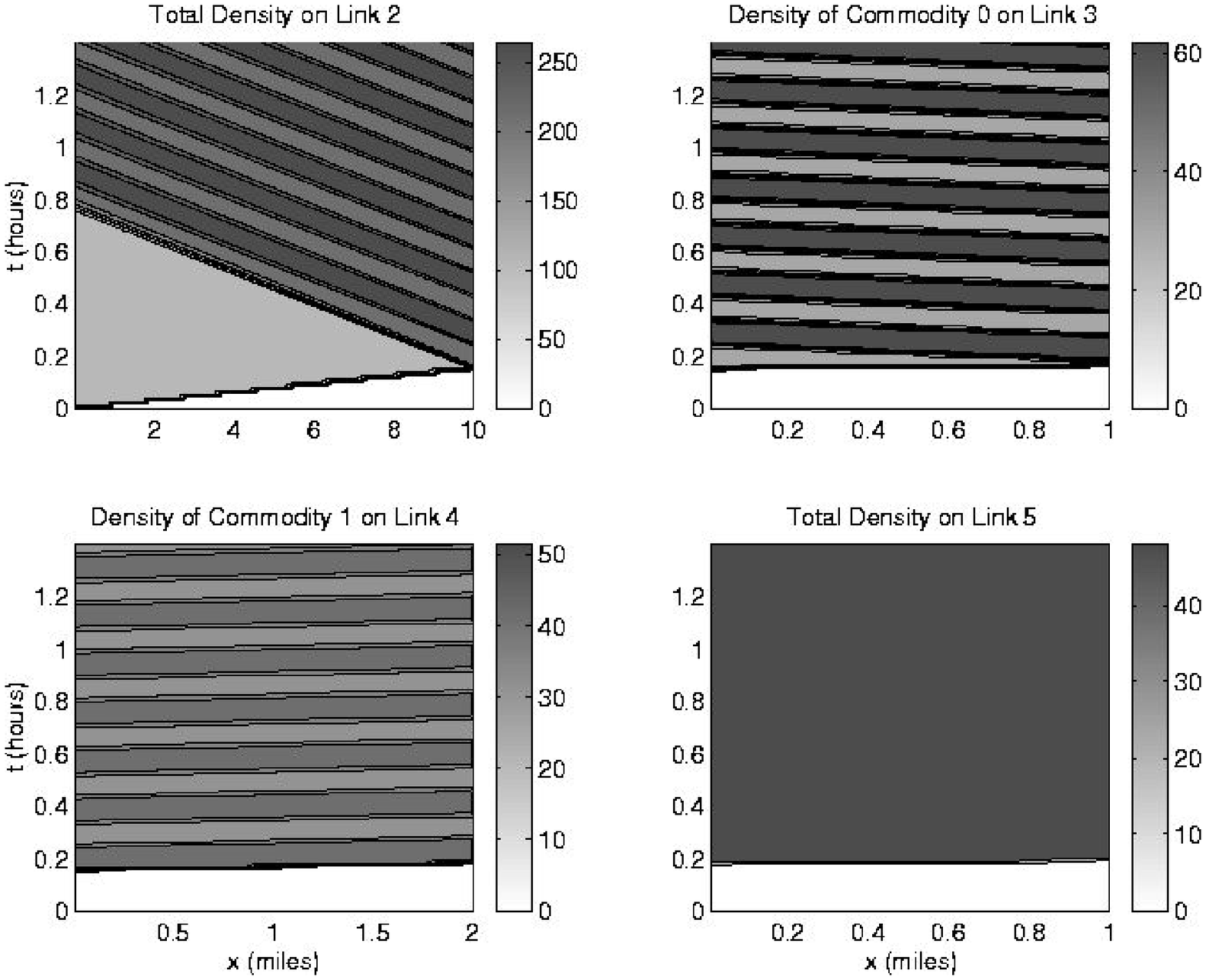}\ec\caption{Contour plots of periodic oscillations}\label{periodicContour}\efg

The formation of periodic solutions can also be observed at certain locations, such as the end of link 2, start of link 3, and start of link 4, as shown in \refe{timecurves}. From \reffig{periodicContour}, we can see that, most of the time, wave speeds on links 2, 3, and 4 are $-\frac  14 v_f$, $-\frac 14 v_f$, and $v_f$, respectively. Thus, from curves in \reffig{timecurves}, we can derive traffic pattern at other locations.

\bfg\bc\includegraphics[height=12cm]{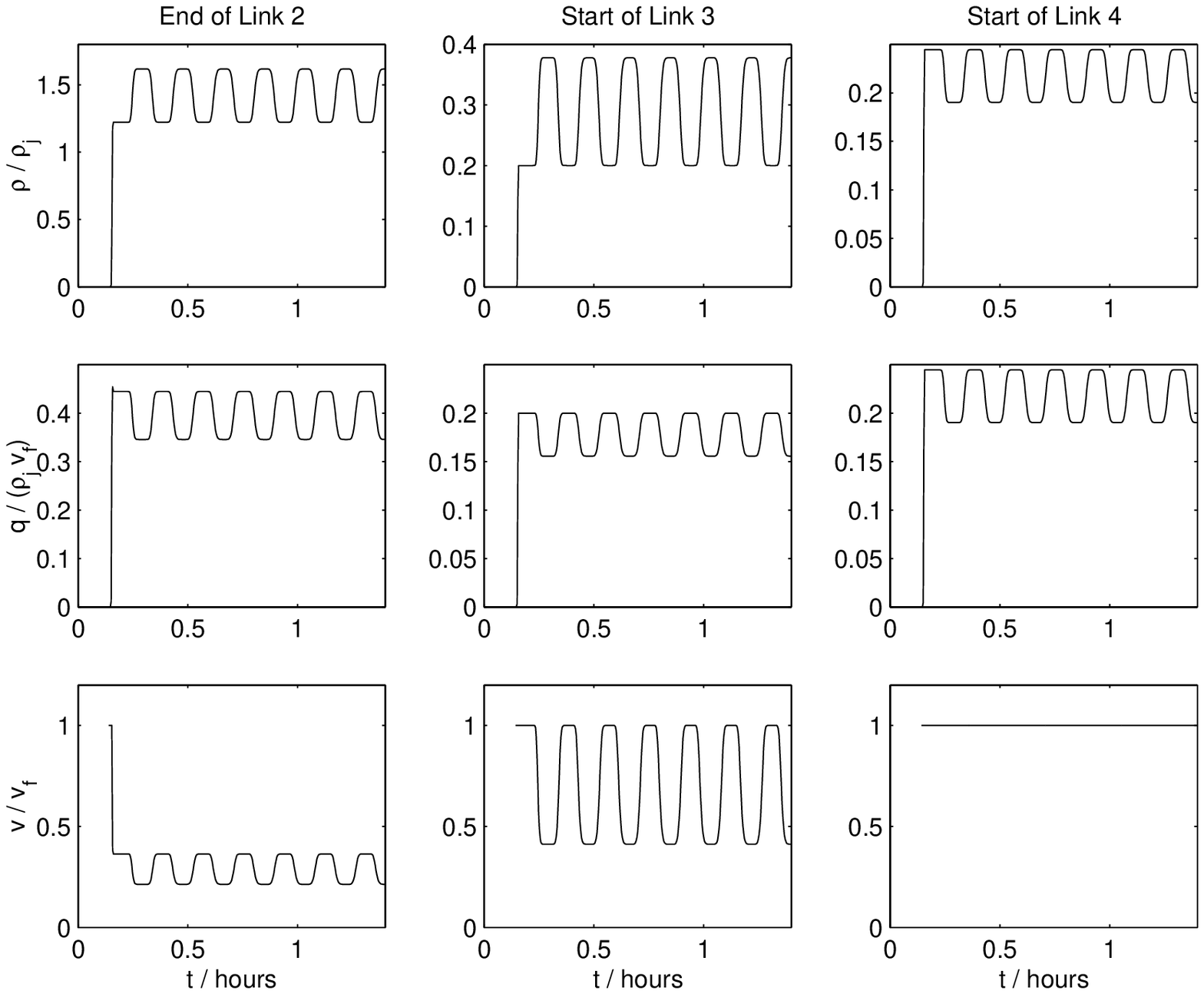}\ec\caption{Density, flow-rate, and speed at the end of link 2, start of link 3, and start of link 4}\label{timecurves}\efg

In a period, traffic conditions evolve in the following five stages. (i) Assuming the period starts at $t_0=0$ when the downstream of link 3 becomes less congested at $\bar \r_3$, which is always over-critical, at $t_0$, link 4 has a smaller density $\bar \r_4$. (ii) At $t_1=4L_3/v_f$, traffic density on link 4 starts to increase to $\hat \r_4$. (iii) After $t_2=t_1+L_4/v_f$, traffic density on link 3 increased to $\hat \r_3$. This finishes the half period. (iv) At $t_3=t_2+4L_3/v_f$, traffic density on link 3 decreases to $\bar \r_3$. (v) At $t_4=T$, traffic density on link 3 gets back to $\hat \r_3$.
\subsection{The structure of periodic solutions}
 In half a period,  solutions at the end of link 2,  start of link 3, and start of link 4 in the $(\r,q)$-plane are shown in \reffig{rqplane}, in which lines are for fundamental diagrams and dots for solution data. From the figure, we can see that solutions in each period is highly symmetric, with transition layer evenly distributed in the $(\r,q)$-plane. In half a period, $\r$ decreases and approximately satisfies the following equation:
\bqn
\frac{\m{d }\r}{\m{d t}}&=&\alpha (\r-\hat \r) (\r-\bar \r),
\eqn
where $\hat \r$ and $\bar \r$ are the maximum and minimum of $\r$ respectively, and $\alpha$ is a constant to be determined, and the other half satisfies
\bqn
\frac{\m{d }\r}{\m{d t}}&=&\alpha (\hat \r-\r) (\r-\bar \r),
\eqn
which is increasing.

\bfg\bc\includegraphics[height=12cm]{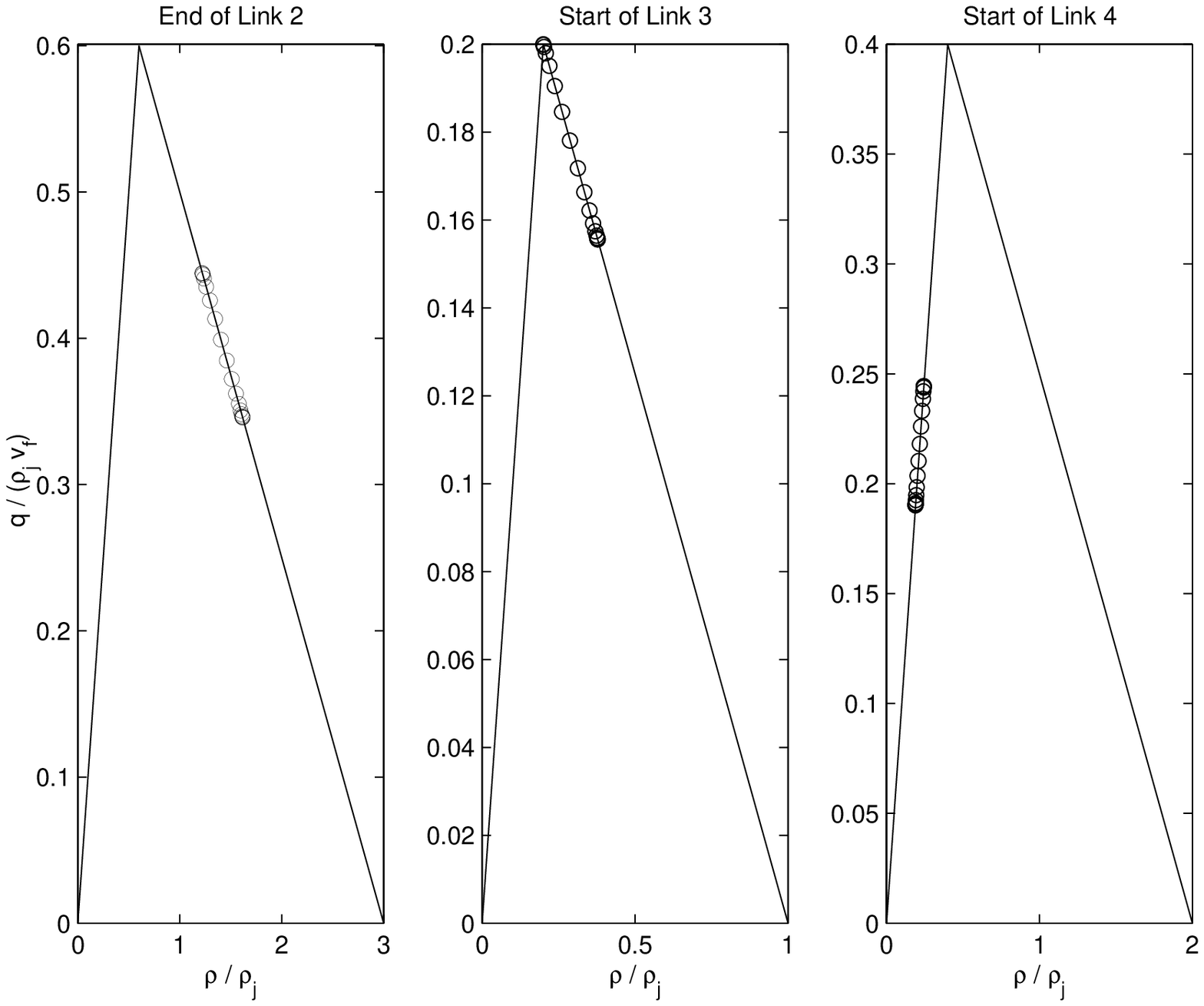}\ec\caption{Periodic solutions at the end of link 2, start of link 3, and start of link 4 in the $(\rho, q)$-plane}\label{rqplane}\efg

In a period, the solutions can be approximated by
\bqn
\r(t)&=&\cas{{ll} \hat \r-(\hat \r-\bar \r)\Big/\left\{1+\exp\left[-\alpha(\hat \r-\bar \r)(t-T/4)\right]\right\}&\m{when }t\in[0,\frac T2),\\
\hat \r-(\hat \r-\bar \r)\Big/\left\{1+\exp\left[\alpha(\hat \r-\bar \r)(t-3T/4)\right]\right\}&\m{when }t\in[\frac T2,T).} \label{periodic:rho}
\eqn
The corresponding $v$ and $q$ are solved by
\bqn
v(t)&=&\cas{{ll} \hat v-(\hat v-\bar v)\Big/\left\{1+\frac {\hat \r}{\bar \r} \exp\left[\alpha(\hat \r-\bar \r)(t-T/4)\right]\right\}&\m{when }t\in[0,\frac T2),\\
\hat v-(\hat v-\bar v)\Big/\left\{1+\frac {\hat \r}{\bar \r} \exp\left[-\alpha(\hat \r-\bar \r)(t-3T/4)\right]\right\}&\m{when }t\in[\frac T2,T),} \label{periodic:v}
\\
q(t)&=&\cas{{ll} \bar q +(\hat q-\bar q)\Big/\left\{1+\exp\left[-\alpha(\hat \r-\bar \r)(t-T/4)\right]\right\}&\m{when }t\in[0,\frac T2),\\
\bar q +(\hat q-\bar q)\Big/\left\{1+\exp\left[\alpha(\hat \r-\bar \r)(t-3T/4)\right]\right\}&\m{when }t\in[\frac T2,T).} \label{periodic:q}
\eqn
As shown in \reffig{structure}, density, flow-rate, and travel speed at the end of link 2 can be well approximated with $\alpha=500$, and the average of $q_2$ in a period is equal to $0.4 \r_j v_f\approx 2 q_c$. Therefore, this pattern is stable.

\bfg\bc\includegraphics[height=12cm]{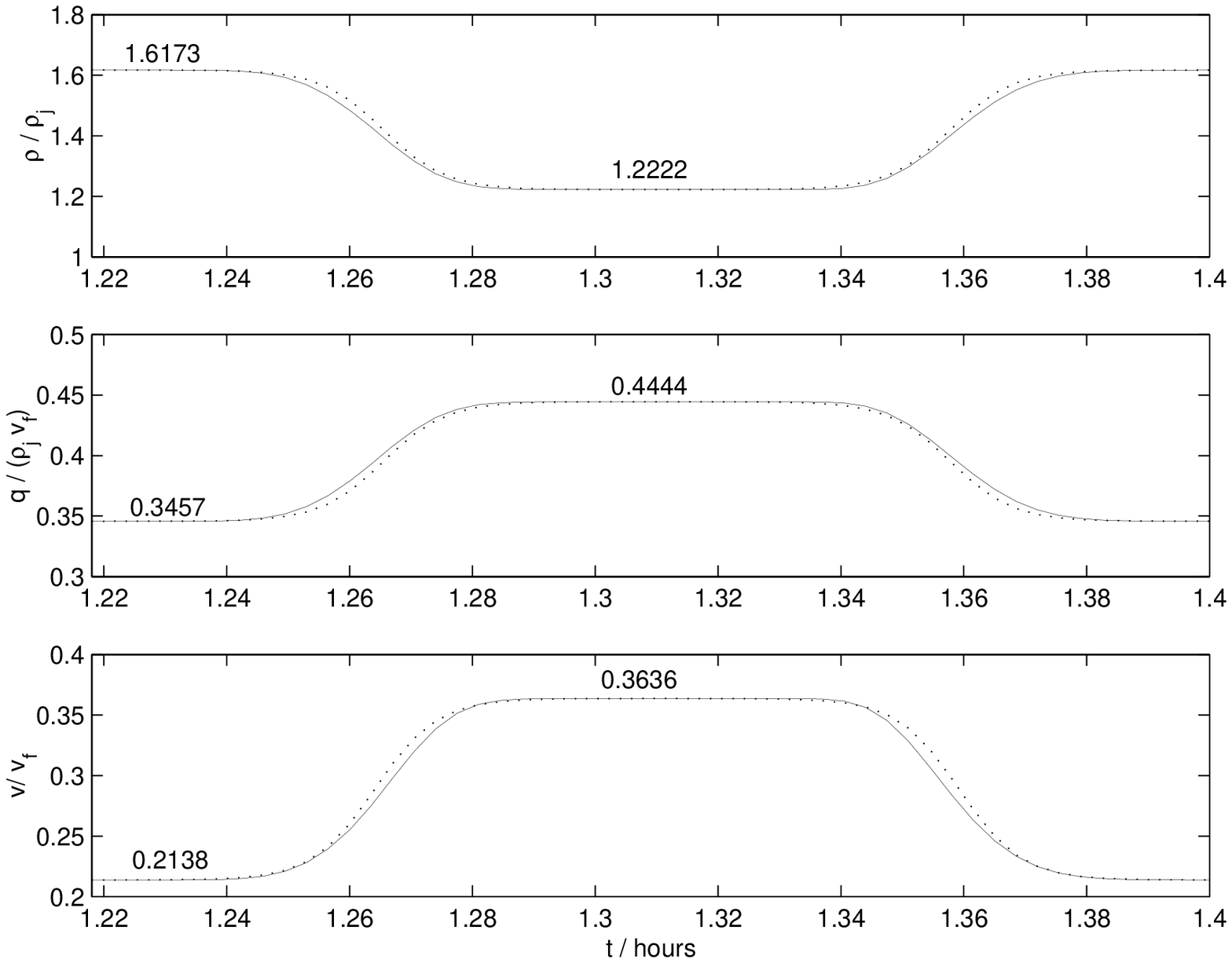}\ec\caption{Structure of periodic solutions on link 2}\label{structure}\efg

Solutions at other locations on link 2, or locations on links 3 and 4 share similar structure as in \refe{periodic:rho}-\refe{periodic:q}, but the phase or amplitude may be different.
\subsection{Discussions}
From \reffig{periodicContour}, we can see that we can observe several periods at the same time only on link 2, although traffic on links 3 and 4 are also periodic. Therefore, these periodic waves can be considered to exist on the link upstream to a diverge and a merge when the proportion satisfies a certain condition. The special road structure, i.e., a diverge and a merge, can be considered as the generator of periodic solutions. We can also see that the period is determined by the length of links between the diverge and merge. In addition, possible approximation of the periodic solutions are provided.

This study shows the effect of network structure on traffic dynamics. We expect more un-revealed solutions patterns exist in more complicated road networks. In the future, it will also be interesting to study the propagation of periodic solutions in road networks. It is not hard to check that when traffic demand decreases, the periodic solutions will disappear.

Since periodic solutions have important negative impact on driver's behavior and vehicles' emission, this study will help to resolve such kind of solutions. 
%\input{applications/4th_app}
%Wenlong Jin: 2003
\section{Conclusions}
In this chapter, with the MCKW simulation model, we studied equilibrium states in a demonstration network and oscillations caused by a special network structure.

In Section \ref{sec:app2},  we showed that traffic dyanmics on a network approach equilibrium states  with constant boundary conditions. Against different combination of traffic commodities, solutions of equilibrium states and the performance of the road network were carefully examined. We further discussed possible relationship between these equilibrium states and the user-equilibrium. Traffic equilibrium states in a road network can be used pursuing important concepts in many network traffic studies, such as traffic assignment, origin demand estimation, etc.

In Section \ref{sec:app3}, we revealed the existence of periodic solutions in a road network. We showed the structure of periodic solutions, whose period is solely dependent on the network structure. This study shows that, in a road network, special network structure may cause stop and go oscillatory waves. \footnote{Another type of vehicle clusters are shown to exist in higher-order models \citep{kerner1994cluster,jin2003cluster}.}  The periodic oscillations are in fact a new type of kinematic waves initiated from simple, jump initial conditions. From this study, we expect more complicated kinematic waves when more complicated networks are considered.

Studies in this chapter show that, with the MCKW simulation, we are able to study more complicated kinematic waves in a road network. This is important for understanding traffic dynamics as demonstrated that complicated traffic phenomena can be caused by network topology and O/D demand pattern. However, other factors may also contribute to stop and go waves in a realistic network.  

\newpage
\pagestyle{myheadings}
\markright{  \rm \normalsize CHAPTER 8. \hspace{0.5cm}
 CONCLUSIONS}
\chapter{Conclusions}
%\thispagestyle{myheadings}
%Wenlong Jin: 2003
%Wenlong Jin: 2003
\section{Summary}
The original kinematic wave theory, also known as the Lighthill-Whitham-Richards theory \citep{lighthill1955lwr,richards1956lwr}, explains traffic dynamics of single origin/destination traffic on a linear road with kinematic waves, including decelerating (shock) waves and accelerating (expansion or rarefaction) waves. In this dissertation, we study traffic phenomena in a road network within the framework of kinematic waves. We have theoretically and numerically investigate traffic dynamics for link inhomogeneities, junctions, and mixed-type vehicles. Furthermore, we have developed a simulation platform of multi-commodity network traffic and studied some network traffic phenomena.

The major results of this dissertation are as follows.
\ben
\item In Chapter 2, we reformulate  the Lighthill-Whitham-Richards model into a nonlinear resonant system for link bottlenecks such as lane-drops. We then show that there is an additional type of kinematic waves, namely, standing (transition) waves. A standing wave always stays at the inhomogeneity spot, and traffic on its both sides has the same flow-rate but different velocities. We further solve the Riemann problem by ten combinations of shock, rarefaction, and transition waves. The wave solutions are summarized and shown to be consistent with solutions by the supply-demand method. Using numerical simulations, we have also shown how traffic queues up in the region upstream to a link bottleneck.
\item In Chapter 3, we take a closer look at existing kinematic wave models for merging traffic within the supply-demand framework. In particular, we study the distribution scheme, which helps to uniquely determine traffic flow from each upstream branch to the downstream link. Further, we propose a ``fairness" condition, under which the out-flow from each upstream branch is proportional to its demand. This condition yields the simplest distribution scheme, which is the only valid distribution independent of downstream conditions. We also show that it can capture the key characteristics of a merge and leads to a merge model that is computationally efficient and easy to calibrate. Furthermore, we show that the new distribution scheme is well-defined and qualitatively sound. With numerical simulations, we demonstrate that the distribution scheme produce convergent solutions.
\item In Chapter 4, we propose a new kinematic wave model for highway diverges. Based on the assumption that, in a short time interval, diverging traffic flows are independent of each other, dynamics of traffic to a downstream link can be described by a system of non-strict conservation laws, or a nonlinear resonant system, whose Riemann problem can be solved by seven types of wave combinations. These waves are called as instantaneous kinematic waves. Further, we show that this model is equivalent to a supply-demand method with modified definitions of  traffic demands of diverging traffic. With numerical simulations, we show that this model is consistent with some existing rules on diverging traffic.
\item In Chapter 5, we study traffic dynamics with mixed-type vehicles, for which we find an additional family of kinematic waves, contact waves. We carefully study the wave solutions of the Riemann problem and develop a Godunov method for solving the model. Using simulations, we demonstrate how mixed traffic evolves on a road link and show that the First-In-First-Out (FIFO) principle is always observed.
\item In Chapter 6, we propose a multi-commodity kinematic wave (MCKW) model of network traffic flow. In this model, we apply the models studied in the previous chapters for different network components and classify vehicles into a number of commodities according to their paths. The proportions of commodities in a road cell are updated based on the fact that traffic is anisotropic. We then show that this model yields solutions, whose departure from the FIFO principle is in the order of a time interval. This model is macroscopic, but is able to provide individual vehicles' trajectories and travel times by using cumulative curve. We also propose an implementation of the MCKW simulation and carefully design the data structure for network topology, traffic characteristics, and simulation algorithms. Our numerical simulations demonstrate that the numerical results converge to FIFO solutions.
\item In Chapter 7, we first study equilibrium states in a road network with a single origin-destination (O/D) pair and two routes. We demonstrate the formation of an equilibrium state. After examining the relationship between equilibrium states and the distribution of  vehicles to different routes, we find that multiple equilibrium status may exist for the same O/D flow but different route distributions. We then show the formation of periodic oscillations for certain network structure and route distributions. The periodic oscillations are in fact a new type of kinematic wave in a road network. We then discuss their structure and properties. \een

The kinematic wave models of network vehicular traffic studied in this dissertation are theoretically rigor, numerically reliable, and computational efficient. These studies help better understand traffic dynamics in a road network, in particular the formation and propagation of traffic congestion, and establish a solid foundation for applications in traffic control and management.
\section{Future research directions}
There are three directions in which we could extend this dissertation research: further investigation of the kinematic wave theories, enhancement of the MCKW simulation model, and applications of the MCKW simulation model to solve real-world problems.

\subsection{Further investigations of the kinematic wave theories}
In studies of inhomogeneous links (Chapter 2) and diverges (Chapter 4), the kinematic wave theories are shown to yield consistent results with the supply-demand method given proper definitions of traffic supply and demand. That is, the supply-demand method can be considered as the discrete form of the corresponding continuous kinematic wave theory. It is also interesting to directly prove that the discrete model in the supply-demand method converge to the corresponding kinematic wave theory. For the inhomogeneous LWR model under jump initial conditions, for example, we can find solutions by the supply-demand method and obtain a time series for each cell. Then will these series get closer to the wave solutions predicted in Chapter 2 as the length of a time interval keeps decreasing? The answer to this problem will help to understand the continuous kinematic wave model of merging traffic, whose discrete model is simple and clear (Chapter 3). In addition, the diverging model (Chapter 4) is an instantaneous approximation, also in discrete form, whose continuous counterpart is another interesting topic. One probable alternative approach is to introduce appropriate fundamental diagram as in \citep{daganzo1997special}.

The kinematic wave models studied in this dissertation are deterministic. However, these models can be extended to include randomness in a road network. For example, non-recurrent traffic congestion is generally caused by accidents or work zones, which are reflected by a change in the number of lanes. Since the number of lanes has been included in the models we studied, these stochastic effects could be considered in the kinematic waves. As another example, weaving effects on traffic dynamics could also be studied in a similar fashion.

Finally, it is of interest to study kinematic waves that arise from the traffic dynamics of other traffic systems, such as those with special lanes, special vehicles, and so on.
\subsection{Calibration, validation, and enhancement of the MCKW simulation model}
The enhancement of MCKW simulation model will be carried out in three aspects. First, we can incorporate new kinematic wave theories so that it is capable of simulating traffic systems with more components, such as different types of vehicles and special lanes. To achieve this, traffic has to be categorized into more commodities, and more complicated data and program structures are expected. Second, we can improve the computational efficiency of the MCKW simulation for large-scale road networks. In this dissertation, the solution methods, either the supply-demand methods or the Godunov methods, are  of first order. One improvement can be made by applying higher-order methods in the  MCKW model \citep[e.g., ][]{daganzo1999lagged,collela2000numerical}, which are more efficient. Another approach is to incorporating parallel algorithms and more efficient memory management methods. Finally, we can design the programming of the simulation model in order to satisfy different requirements of different applications. For example, we need very detailed information when studying traffic dynamics and phenomena, but need less information in traffic assignment. Thus it will be better to use two different programming structure for these two applications.

The MCKW simulation model, as other simulation models, is also subject to calibration and validation with observed traffic dynamics. A conceptual approach is to calibrate and validate the simulation model first for inhomogeneous links, merges, and diverges and then for a whole road network. We can check existing databases or collect data by ourselves in the calibration and validation process. In many data sets collected by loop detectors, however, volumes (the number of vehicles passing a detector) are not conserved over a link. This inconsistency will be fatal for calibrating the kinematic wave models, since the fundamental assumption in these models is traffic conservation. Thus, it may be necessary to seek advanced measures, such as  vedio-taping or by satellite sensing techniques, to provide the ideal accuracy in data in these studies.
\subsection{Applications of the MCKW simulation model}
First, using the MCKW simulation model, we can better understand the formation and characteristics of traffic congestion. This will form a foundation to alleviate it by either expansion of infrastructure or traffic control and management.

Second, the MCKW simulation model can be applied in evaluating strategies in infrastructure expansion and traffic operations and management. For example, we could develop and evaluate on-ramp metering and arterial signal control algorithms by using the MCKW model, in which signals have been incorporated in the computation of traffic demand (see Chapter 3). In addition, in simulation-based dynamic traffic assignment, the MCKW model can be used as a loading model.

Finally, the MCKW simulation model can be used to estimate the travel demands in a road network with observed traffic conditions at a number of locations and can be further applied in regional and transportation plan.

\newpage
\phantomsection
\addcontentsline{toc}{chapter}{Bibliography}
\pagestyle{myheadings}
\markright{  \rm \normalsize BIBLIOGRAPHY}
\bibliographystyle{abbrvnat}
\bibliography{references}
\end {document}